%% file: main_-_Apr_10.tex
\documentclass[english]{smfart}

\usepackage[french,english]{babel}
\usepackage{smfthm}
\usepackage{amsmath,amssymb}
\usepackage{mathtools}
\usepackage{tikz-cd}
\usepackage{enumerate}
\usepackage[T1]{fontenc}     
\usepackage[utf8]{inputenc}
\usepackage{graphicx}
\usepackage{subcaption}
\usepackage{xcolor}
\usepackage{epsfig,color,soul}
\usepackage{wrapfig}
\usepackage{url}

\usepackage{fourier}

\theoremstyle{plain}
\newtheorem{thm}{Theorem}[section]
\newtheorem{lem}[thm]{Lemma}
\newtheorem{cor}[thm]{Corollary}
\newtheorem{proposition}[thm]{Proposition}
\newtheorem{question}[thm]{Question}

\newtheorem*{thm*}{Theorem}
\newtheorem*{prop*}{Proposition}
\newtheorem*{lem*}{Lemma}
\newtheorem*{cor*}{Corollary}

\newtheorem*{thmA}{Theorem A}
\newtheorem*{thmB}{Theorem B}

\theoremstyle{remark}
\newtheorem{ex}[thm]{\textbf{Example}}
\newtheorem{rem}[thm]{\textbf{Remark}}

\newtheorem*{rem*}{\textbf{Remark}}
\newtheorem*{ex*}{\textbf{Example}}

\newtheorem*{cl}{\textbf{Claim}}
\newtheorem*{fact}{\textbf{Fact}}

\newtheorem*{interpolating nbhd}{The interpolating neighborhoods}
\newtheorem*{construction of H_W}{The construction of $H_W$}
\newtheorem*{Scheme of the proof}{Scheme of the proof}
\newtheorem*{The Construction of $H_V$ and $H_V'$}{The construction of $H_V$ and $H_V'$}
\newtheorem*{Proof of the steps}{Proof of the steps}

\theoremstyle{definition}
\newtheorem{defn}[thm]{Definition}

\newtheorem*{defn*}{Definition}

\newtheorem*{question*}{Question}

\DeclareMathOperator{\interior}{\mathrm{int}}

\DeclareMathOperator{\spa}{\mathrm{span}}
\DeclareMathOperator{\dist}{\mathrm{dist}}

\newcommand{\B}{\mathring{B}}
\newcommand{\U}{\mathring{U}}

\newcommand{\R}{\mathbb{R}}
\newcommand{\Z}{\mathbb{Z}}
\newcommand{\D}{\mathbb{D}}
\newcommand{\s}{\mathbb{R}/\mathbb{Z}}
\newcommand{\F}{\mathcal{F}}
\newcommand{\Fs}{\mathcal{F}^\mathrm{s}}
\newcommand{\Fu}{\mathcal{F}^\mathrm{u}}
\newcommand{\Es}{E^\mathrm{s}}
\newcommand{\Cs}{C^\mathrm{s}}
\newcommand{\Cu}{C^\mathrm{u}}
\newcommand{\FFs}{F^\mathrm{s}}
\newcommand{\FFu}{F^\mathrm{u}}
\newcommand{\Ec}{E^\mathrm{c}}
\newcommand{\Eu}{E^\mathrm{u}}
\newcommand{\Ecu}{E^\mathrm{cu}}
\newcommand{\Ecs}{E^\mathrm{cs}}

\newcommand{\Ccs}{C^\mathrm{cs}}
\newcommand{\Esu}{E^\mathrm{su}}
\newcommand{\Ws}{W^\mathrm{s}}

\newcommand{\Wu}{W^\mathrm{u}}
\newcommand{\Wsloc}{W^s_\mathrm{loc}}
\newcommand{\Wuloc}{W^u_\mathrm{loc}}

\author{Mario Shannon}
\address{
Department of Mathematics\\
The Pennsylvania State University\\
State College, PA 16801, USA
}
\email{mjs8822@psu.edu}
\title{Hyperbolic models of transitive topologically Anosov flows in dimension three}
\alttitle{Modèles hyperboliques de flots topologiquement d'Anosov transitifs}

\begin{document}
\frontmatter

\begin{abstract}
We prove that every transitive topologically Anosov flow on a closed $3$-manifold is orbitally equivalent to a smooth Anosov flow, preserving an ergodic smooth volume form. 
\end{abstract}

\begin{altabstract}
Nous montrons que tout flot topologiquement d'Anosov et transitif sur une 3-variété compacte sans bord est orbitalement équivalent à un flot d'Anosov lisse, qui préserve une forme de volume ergodique. 
\end{altabstract}

\altkeywords{}

\maketitle
\setcounter{tocdepth}{1}
\tableofcontents

\mainmatter

\section{Introduction.}
A \emph{topologically Anosov flow} on a closed $3$-manifold is a continuous $\R$-action that is non-singular (i.e., has no global fixed point) and satisfies the following condition: the two partitions of the manifold into stable and unstable sets of the trajectories of the flow induce a pair of transverse foliations that intersect along the respective trajectories. 

The first family of examples in this category are the \emph{smooth Anosov flows}. The flow generated by a non-singular smooth vector field on a closed manifold is \emph{Anosov}, if its natural action on the tangent bundle of the manifold preserves a \emph{hyperbolic splitting}. In this case, the \emph{stable manifold theorem} of hyperbolic dynamics implies that the partitions of the manifold into stable and unstable sets of the trajectories form a pair of transverse, invariant, $2$-dimensional foliations, that intersect along the trajectories of the flow. Therefore, every smooth flow that is globally hyperbolic (Anosov) is, in particular, a topologically Anosov flow.

A general topologically Anosov flow has the same dynamical behavior as a smooth Anosov flow from the point of view of topological dynamics. For instance, it is \emph{orbitally expansive} and satisfy the \emph{global pseudo-orbits tracing property}, and it turns out that these two properties are enough to characterize its dynamics (cf. Section \ref{subsection_top-Anosov_defn}). The main difference between topologically Anosov and smooth Anosov is not a matter of regularity, it is about the lack of hyperbolic splitting. A smooth flow can be topologically Anosov and non-hyperbolic at the same time, hence not Anosov. Many of the properties of smooth Anosov flows derived from the existence of a hyperbolic splitting (e.g. $C^1$-structural stability, exponential contraction/expansion rates along invariant manifolds, or existence of special ergodic measures) are no longer available in the topologically Anosov framework. 

The main problem motivating this work is the following question:

\begin{question}\label{question:top-Anosov=Anosov}
Is every topologically Anosov flow on a closed $3$-manifold orbitally equivalent to a smooth Anosov flow?
\end{question} 

When a topologically Anosov flow is orbitally equivalent to a smooth Anosov flow, we say that the latter is a \emph{hyperbolic model} of the first. 

The question above has at least two traceable origins. One is related with the construction of different examples of Anosov flows on closed $3$-manifolds, where topologically Anosov flows appear as intermediate steps in a battery of \emph{3-manifold cut-and-paste} techniques called \emph{surgeries}. The most paradigmatic example is the so called \emph{Fried surgery} introduced in \cite{fried_surgery}, a kind of Dehn surgery but adapted to the pair (flow, $3$-manifold), which produces topologically Anosov flows that are not hyperbolic in a natural way (cf. Question \ref{question_almost-Anosov} below).  

The other origin is related with the study of \emph{smooth orbital equivalence classes} of Anosov flow. The existence of a hyperbolic model for a given topologically Anosov flow $(\phi,M)$ means that it is possible to endow the manifold $M$ with a smooth atlas, such that the foliation by flow-orbits is tangent to an Anosov vector field. In this sense, topologically Anosov flows can be seen as the toy models from which smooth Anosov flows can be obtained, by adding suitable smooth structures on the manifold. It is thus relevant to understand if all of these topologically Anosov toy models effectively correspond to a hyperbolic dynamic. The general problem of the existence of smooth models for a given expansive dynamic has been treated in other contexts, from where it is relevant to mention the case of \emph{pseudo-Anosov} homeomorphisms on closed surfaces treated in \cite{gerber-katok}. (See also \cite{lewowicz_anal} and \cite{farrell-jones_Anosov_exotic}).

The notion of topologically Anosov turns out to also be involved in other results about classification of \emph{partially hyperbolic} diffeomorphisms on $3$-manifolds. An example is \cite{bonatti-wilkinson}, where an affirmative answer to Question \ref{question:top-Anosov=Anosov} is conjectured.

Recall that a flow is called \emph{transitive} if there exists a dense orbit. The purpose of this work is to prove the following statement:

\begin{thmA}
Every transitive topologically Anosov flow on a closed $3$-manifold admits a hyperbolic model, which in addition preserves an ergodic smooth volume form.
\end{thmA}   

In some literature on the subject (e.g. \cite{brunella_thesis}, \cite{fried_surgery} and \cite{mosher_monograph}) we can find arguments trying to support the existence of hyperbolic models of topologically Anosov flows, specially in the transitive case. Nevertheless, as we can confirm from more recent works (e.g. \cite{bonatti-wilkinson}) there is no agreement about a complete proof.

Some special cases of this result can be derived from the works \cite{plante_sol_mafld_classification} and \cite{ghys_Classification_Anosov_circle_bundles}, if we assume that the ambient manifold is a torus bundle or a Seifert manifold. In those works, by studying the topological properties of the invariant foliations associated with an Anosov flow, it is shown that, in each case, there is essentially one possible model: \emph{suspension of a hyperbolic automorphism of the torus} or \emph{geodesic flow of a hyperbolic surface}. These results extend into the context of topologically Anosov (or even expansive, see \cite{brunella_expansivos_seifert} and \cite{paternain_geodesic_flows}) and provide the existence of hyperbolic models. 

We do not know a proof of Theorem \textbf{A} without the transitivity assumption. The proof that we present here makes use of an \emph{open book decomposition} of the supporting 3-manifold associated to the topologically Anosov flow, that is not available in the non-transitive case. 

\begin{rem}
It is worth to note that every transitive smooth Anosov flow in dimension three (or in any dimension, provided that the dimension of the strong stable manifold is one) admits a hyperbolic model which preserves a smooth volume form, as it is proved in \cite{asaoka_invariant_volumes}. Nevertheless, we do not make use of this result in the present construction, but instead we directly construct the invariant measure, since it appears naturally in the hyperbolic models that we provide at Section \ref{section_hyperbolic_models}. 
\end{rem}

Let us outline the construction of the hyperbolic models in Theorem \textbf{A}.

\subsubsection*{Transitiveness and open book decomposition}

Given a non-singular flow on a closed $3$-manifold, a \emph{Birkhoff section} is a compact surface, usually with non-empty boundary, immersed in the phase space in such a way that: (1) The interior of the surface is embedded and transverse to the flow lines; (2) the boundary components are periodic orbits of the flow; and (3) every orbit intersects the surface in an uniformly bounded time. These sections come equipped with a first return map defined on the interior of the surface and hence the flow is a suspension on the complement of a finite set of periodic orbits. 

The joint information of a flow and an associated Birkhoff section provides what is called an \emph{open book decomposition} of the $3$-manifold: in the complement of a finite union of closed curves (called the \emph{binding}) the manifold is a surface bundle over the circle, with fibers homeomorphic to the interior of the Birkhoff section (called \emph{pages}), and monodromy map corresponding to the first return. Apart from the closed curves in the $3$-manifold that form the binding and the monodromy onto the pages, the isotopy class of the embedding of the pages in a neighborhood of the binding is an important information associated with the open book decomposition, that may be encoded by a set of integer parameters.

\begin{figure}[t]
\begin{center}
\includegraphics[width=\textwidth]{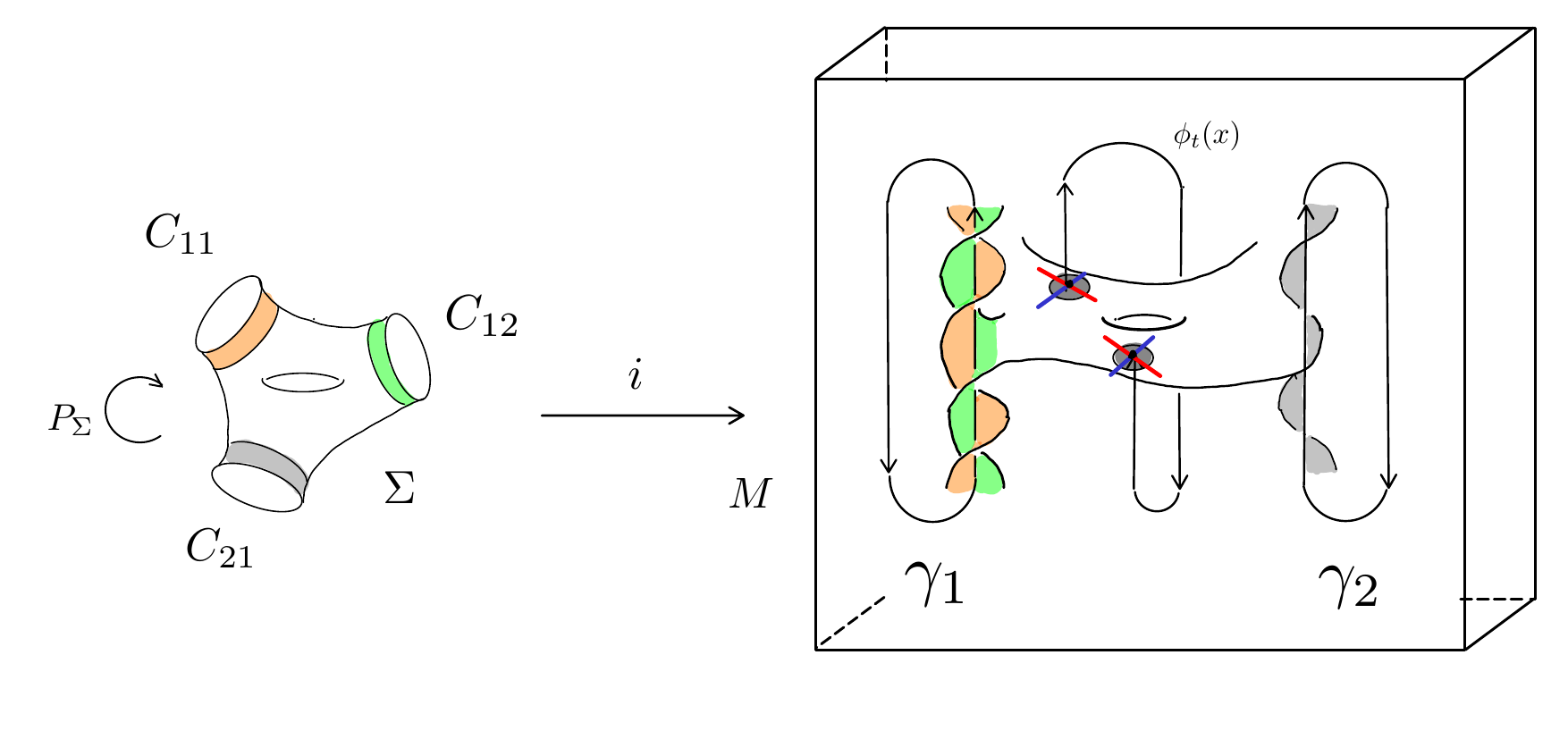}
\caption{Birkhoff section.}
\end{center}
\end{figure}

In \cite{fried_surgery}, Fried proved that every transitive Anosov flow admits a Birkhoff section whose first return map is pseudo-Anosov (in a non-closed surface). This extraordinary fact, later generalized by Brunella to any transitive expansive flow in \cite{brunella_thesis}, opens the possibility of reducing part of the analysis of transitive Anosov flows to the theory of pseudo-Anosov homeomorphisms on surfaces and open book decompositions. Observe that transitivity is a necessary condition for the existence of a Birkhoff section, since the first return is pseudo-Anosov and hence transitive.  

From a Birkhoff section associated with a topologically Anosov flow it is possible to derive many interesting informations about the flow and the topology of the ambient space. For instance, it is possible to construct \emph{Markov partitions} (\cite{brunella_thesis}) or to derive properties of the topology of the \emph{orbit space} associated to the flow (\cite{fenley_Anosov},\cite{asaoka2022oriented}). In addition, Birkhoff sections give a good framework for working with Fried surgeries. 

We use Birkhoff sections to construct the hyperbolic models of Theorem \textbf{A}. Along this construction, we need to show that the orbital equivalence classes of topologically Anosov flows are determined by the combinatorial data associated to a Birkhoff section. There is a well-known property in surface dynamics, which states that two isotopic pseudo-Anosov homeomorphisms on a closed surface are conjugated by a homeomorphism isotopic to the identity. In turn, this property means that the conjugacy class of a pseudo-Anosov on a closed surface is determined by the action of the homeomorphism on the fundamental group of the space (see Section \ref{subsection_psudo-Anosov-homeos}). We show the following:

\begin{thmB}[Theorem \ref{thmA-statement}]
The orbital equivalence class of a transitive topologically Anosov flow on a closed $3$-manifold is completely determined by the following data associated with a Birkhoff section: 
\begin{enumerate}[i.]
\item The action of the first return map on the fundamental group of the Birkhoff section,
\item The combinatorial data associated with the embedding near the boundary components.
\end{enumerate}
\end{thmB}
This theorem can be interpreted as an extension of the former property of pseudo-Anosov homeomorphisms into the context of expansive flows. (Compare with Lemma 7 of \cite{brunella_top_equivalence_criterion}.)

The main difficulty in proving this theorem can be explained as follows: Given two non-singular flows equipped with Birkhoff sections, if the first return maps (which are defined only on the interior of a compact surface) are conjugated, then we can use this conjugation to produce an orbital equivalence between the two flows, only defined in the complement of finite sets of periodic orbits. In order to extend the equivalence onto the periodic orbits on the boundary, there is an essential obstruction coming from the fact that the homeomorphism which conjugates the first return maps, in general, does not extend continuously onto the boundary of the sections. In our case, the corresponding first return applications are pseudo-Anosov but, for instance, \emph{it is not true that the action on the fundamental group determines the conjugacy class of a pseudo-Anosov homeomorphism in a surface with boundary} (cf. Section \ref{subsection_psudo-Anosov-homeos}). 

To deal with this obstruction we use techniques coming from \cite{bonatti-guelman}. In that work, the authors provide many tools for the analysis of the germ of a hyperbolic flow in the neighborhood of a boundary component of an immersed Birkhoff section. When the Birkhoff section is nicely embedded in the $3$-manifold, a property that they call \emph{tame}, it is possible to reconstruct the germ of the periodic orbit as a function of the first return map on the interior of the surface.

\begin{rem}
It has been subsequently observed that it is possible to write an alternative proof of Theorem \textbf{B}, based on completely different methods coming from \cite{barbot_charact_flots_Anosov} and \cite{fenley_Anosov}. See \cite{iakovoglou_Thesis} for a sketch of this argument. (We thank the referee for pointing out this alternative proof.)
\end{rem}

\subsubsection*{Construction of the hyperbolic models}

One strategy for the problem of finding a smooth representative of a given expansive dynamical system is as follows: First, by some modification process we construct a smooth model expected to be equivalent to the original one; and second, we prove that the smooth model is actually equivalent to the original dynamical system by some stability argument. We can put this in practice in the case of transitive topologically Anosov flows due to the existence of Birkhoff sections. 

Since in the complement of a finite set of periodic orbits the flow is orbitally equivalent to the suspension of a pseudo-Anosov map on a non-closed surface, we can endow this open and dense set with a smooth transversally affine atlas and a Riemannian metric, such that the action of the flow preserves a uniformly hyperbolic splitting in an open manifold. These smooth structures are called \emph{almost Anosov structures}. 

Each Birkhoff section produces an almost Anosov atlas only defined in the complement of the boundary orbits. It is not clear whether such a smooth atlas can be extended along the boundary of the section, in such a way that the original topologically Anosov flow becomes smooth and preserves a hyperbolic splitting.

\begin{question}\label{question_almost-Anosov}
Let $(\phi,M)$ be a transitive topologically Anosov flow equipped with an almost Anosov structure in the complement of a finite set of periodic orbits $\Gamma$ (cf. Section \ref{subsection_psudo-Anosov-homeos}). Does it exist a $C^r$-smooth Anosov flow $(\psi,N)$, where $r\geq 1$, together with an orbital equivalence $H:(\phi,M)\to (\psi,N)$, such that the restriction of $H$ onto $M\backslash\Gamma$ is a diffeomorphism onto its image in $N$ ? 
\end{question}

\begin{rem}\label{non-extension_Bonatti}
For instance, C. Bonatti has orally indicated to the author an argument that shows the answer is no for $r\geq 1+\alpha$. That is, \emph{there is no such flow $\{\psi_t:N\to N\}_{t\in\R}$ generated by a vector field $X$ of class $C^1$ with $\alpha$-H\"older differential $DX$, $0<\alpha<1$.} We plan to sketch the argument in a future work about ``Transitive Anosov flows and transverse Lorentzian structures with singularities'' See also the comment on \emph{smooth models} of pseudo-Anosov homeomorphisms at Section \ref{subsection_psudo-Anosov-homeos}.
\end{rem}

Beyond the possibility of extending or not one of these almost Anosov structures onto the whole manifold, we construct a hyperbolic model of a transitive topologically Anosov flow $(\phi,M)$ completing the following steps: 

\vspace{2pt}
\noindent
\textbf{Sketch of proof of Theorem A:}
\begin{enumerate}[(i)]
\item By making a \emph{Goodman-like} surgery (cf.\ \cite{goodman_surgery}) of an almost Anosov structure in a neighborhood of the singular orbits, we construct a volume-preserving smooth flow $\Psi$ on a smooth manifold $N$ homeomorphic to $M$. Then, using the so-called \emph{cone field criterion} (cf. \cite{katok-hasselblatt}), we show that this flow preserves a uniformly hyperbolic splitting. 

\item We show that both flows $(\phi,N)$ and $(\psi,N)$ are equipped with Birkhoff sections satisfying the criterion for orbital equivalence stated in Theorem \textbf{B}.
\end{enumerate}

\begin{rem}
We encounter the same technique in the construction of \emph{smooth models of pseudo-Anosov homeomorphisms} given in \cite{gerber-katok}. Observe that, given a pseudo-Anosov in a closed surface, the transverse invariant measures associated with the stable/unstable foliations define a smooth (translation) atlas in the complement of a finite singular set, that makes the map smooth out of these singularities. The construction of the smooth models in \cite{gerber-katok} is done by modifying the given pseudo-Anosov map in a neighborhood of its singularities to obtain an everywhere smooth map, and then showing equivalence between the original pseudo-Anosov and the new diffeomorphism. But to show equivalence between the two is not at all direct, since pseudo-Anosov homeomorphisms are not \emph{structurally stable}. Even if the smooth model is obtained by modifying the pseudo-Anosov map in a small neighborhood of the singularities, the conjugation between them is a homeomorphism that, in general, is different from the identity almost everywhere. 
\end{rem}

A phenomenon similar to the one described in the previous remark occurs in the construction of Theorem \textbf{A}. The surgery that we apply is contained in a small neighborhood of a finite set of periodic orbits. However, the identity map in the complement of that neighborhood does not extend to a global orbital equivalence. We refer Proposition 3.16 of \cite{shannon_thesis} and comments therein for an explanation. 

The idea for constructing transitive smooth Anosov flows by surgery methods goes back to the seminal work \cite{handel-thurston} of Handel and Thurston, and its general strategy has been applied several times in constructions involving Anosov flows. The construction presented in Theorem \textbf{A} is based on the techniques introduced by Goodman in \cite{goodman_surgery}, which adapt Handel-Thurston's method to remove a tubular neighborhood of a periodic orbit, and glue it back in order to obtain another Anosov flow. 

The technique used in the present work can be used to produce special \emph{models} in other contexts involving transitive Anosov (or pseudo-Anosov) flows, due to the existence of Birkhoff sections and Theorem \textbf{B}. For instance, in \cite{shannon_thesis} it has been proved that the Goodman surgery referred above only depends on its homological parameters (independent of the particular position of the surgery locus), producing an Anosov flow uniquely defined up to orbital equivalence. Other example is \cite{agol-tsang_Veering-triangulations_AGT}, where some particular models of pseudo-Anosov flows are constructed.

\subsection*{Description of the content}
This paper is based on results from the author's doctoral thesis \cite{shannon_thesis}. At some places where the proofs are standard arguments, we  refer to \cite{shannon_thesis} for more details. The paper is organized in two main parts: The proof of Theorem \textbf{B} contained in Section \ref{section_Birk_sect_and_equiv} and the proof of Theorem \textbf{A} contained in Section \ref{section_hyperbolic_models}. 

In Section \ref{section_preliminaries} we summarize several definitions and elementary properties about orbit equivalence, transverse sections, topologically Anosov flows and pseudo-Anosov homeomorphisms. This section does not contain proofs and is intended only for completeness of the material. Nevertheless, we explain some prior results in the form that are used in subsequent sections.

The proof of Theorem \textbf{B} is contained in Section \ref{section_Birk_sect_and_equiv}. For a quick reading on Theorem \textbf{B} we recommend the reader to go through: Section \ref{subsection_psudo-Anosov-homeos}, Proposition \ref{prop_combinatorial_relation_projections_first_return}, Theorem \ref{thm_fund_local} on Section \ref{subsection-Local_version_thm-A}-\ref{subsection-Local_version_thm-A}.1 and Section \ref{subsection-Proof_thm-A}. 

In Section \ref{section_almost-Anosov} we explain what is an almost Anosov atlas and we give a description of the atlas in a neighborhood of a singular orbit.  

The proof of Theorem \textbf{A} is contained in Section \ref{section_hyperbolic_models}. For a quick reading on Theorem \textbf{A} we recommend the reader to go through: Section \ref{subsection_top-Anosov_defn}, statement of Theorem \textbf{B} at the introduction of Section \ref{section_Birk_sect_and_equiv}, Section \ref{section_almost-Anosov} and \ref{subsection_proof_thm_A}, \ref{subsection_cross-shaped_neighbourhoods} and \ref{subsection_smooth_model_is_equivalent}. 

\subsection*{Acknowledgements}
The results presented in this work are a continuation of previous ones obtained by M. Brunella, C. Bonatti, F. Béguin and Bin Yu, and were developed in the course of the PhD thesis of the author. We are infinitely grateful to C. Bonatti, who has introduced us into these beautiful questions and techniques. We also want to thank T. Barbot, F. Béguin, P. Dehornoy, S. Fenley, S. Hozoori, R. Potrie, F. R. Hertz and B. Yu for many helpful conversations during the preparation of this work.

\section{Preliminaries.}\label{section_preliminaries}
This section contains the preliminary results and definitions about flows and orbital equivalence (\ref{subsection_flows_and_equivalence}), transverse sections and first return maps (\ref{subsection_transverse_sect_and_1st-return}), expansive flows (\ref{subsection_top-Anosov_defn}), and pseudo-Anosov homeomorphisms on surfaces (\ref{subsection_psudo-Anosov-homeos}). 

\subsection{Flows and orbital equivalence.}\label{subsection_flows_and_equivalence}
A \emph{continuous flow} $\phi$ on a manifold $M$ is a continuous map $\phi:\R\times M\to M$ satisfying the \emph{group property} $\phi(s,\phi(t,x))=\phi(s+t,x)$, for every $s,t\in\R$ and $x\in M$. The point $\phi(t,x)$ is called the \emph{action at time $t$ over $x$}, and is denoted by $\phi_t(x)$. For every $t\in\R$ the map $\phi_t:M\to M$ is a homeomorphism, and the group property implies that the 1-parameter family $\{\phi_t:M\to M\}_{t\in\R}$ is a subgroup of the homeomorphism group of $M$. A flow $\phi$ on a manifold $M$ can be equivalently defined as a continuous 1-parameter family of homeomorphisms $\phi=\{\phi_t\}_{t\in\mathbb{R}}$ of $M$ such that $\phi_0$ equals the identity and $\phi_{t}\circ\phi_{s}=\phi_{s+t}$, for every $s,t\in\R$. 

The \emph{orbit} of $x$ under the action of $\phi$ is the set $\mathcal{O}(x)=\{\phi_t(x):t\in\R\}$. The set of all $\phi$-orbits $\mathcal{O}=\{\mathcal{O}(x)\}_{x\in M}$ constitutes a partition of $M$. The \emph{positive and negative semi-orbits $\mathcal{O}^{+}(x)$ and $\mathcal{O}^{-}(x)$} are the subsets of $\mathcal{O}(x)$ obtained by restricting the time parameter to either $t\geq 0$ or $t\leq 0$, respectively. Given a point $x\in M$ and two real numbers $t_1<t_2$, we denote by $[\phi_{t_1}(x),\phi_{t_2}(x)]$ the orbit segment $\{\phi_s(x):t_1\leq s\leq t_2\}$. Given a non-empty open set $U\subset M$ and a point $x\in U$ we denote by $\mathcal{O}_U(x)$ the connected component of $\mathcal{O}(x)\cap U$ that contains $x$. We denote by $\mathcal{O}_U=\{\mathcal{O}_U(x)\}_{x\in U}$ the partition of $U$ induced $\phi$-orbits.  

The flow $\phi$ is \emph{non-singular} if no orbit $\mathcal{O}(x)$ is a singleton. The flow is \emph{regular} if every $x\in M$ has a neighborhood $U$ where the partition $\mathcal{O}_U$ is topologically equivalent to the partition of $\R^2\times\R$ by vertical lines $\{\{p\}\times\R\}_{p\in\R^2}$. The orbit partition $\mathcal{O}$ associated to a continuous, regular, non-singular flow $\phi$, is a foliation of $M$ by immersed $1$-manifolds. Observe that the orbits of a non-singular flow are naturally oriented by the direction of the flow, and hence $\mathcal{O}$ is an oriented foliation. 

If $X$ is a non-singular vector field of class $C^k$ in $M$, where $k\geq 1$, then there is a regular, non-singular flow $\{\phi_t:M\to M\}_{t\in\R}$ associated to the system of ordinary differential equations $\dot{x}=X(x)$, that is of class $C^k$. In this case, the foliation by flow orbits is of class $C^k$ and its leaves are immersed $1$-dimensional manifolds tangent to the vector field $X$. 

For $i=1,2$ let $\phi^i=\{\phi_t^i\}_{t\in\R}$ be a regular flow of class $C^k$ on a manifold $M_i$, where $k\geq 0$.

\begin{defn}[\textbf{Orbital equivalence}]
The flows $(\phi^1,M_1)$ and $(\phi^2,M_2)$ are $C^r$-\textit{orbitally equivalent}, where $r\geq 0$, if there exists a $C^r$-diffeomorphism $H:M_1\to M_2$ such that, for every $x\in M_1$, it sends the orbit $\mathcal{O}^1(x)$ homeomorphically onto the orbit $\mathcal{O}^2(H(x))$, preserving the orientation of these orbits. We denote it by $H:(\phi^1,M_1)\to(\phi^2,M_2).$
\end{defn}

\begin{rem}
When we use the word \emph{orbital equivalence} without making any reference to the regularity degree $r\geq 0$, it must be understood that $r=0$. In other words, we just care about \emph{homeomorphisms} preserving the oriented foliations by orbit segments, no matter the degree of regularity of the flows. 
\end{rem}

For technical reasons, we are also interested in a weaker notion that we explain here: Consider a non-empty open subset $U\subset M$. The foliation by $\phi$-orbits on $M$ induces a foliation on $U$ but, in general, the action $\R\times M\to M$ given by $\phi$ does not restrict onto an $\R$-action on the set $U$. Instead, what we obtain is a \emph{pseudo-flow} on $U$. That is, a map $(t,x)\mapsto\phi_t(x)$ defined for some couples $(t,x)\in\R\times U$. This pseudo-flow generates a partition of $U$ into orbit segments, which coincides with the foliation $\mathcal{O}_U$ previously defined. This restriction pseudo-flow is be simply denoted by $(\phi,U)$.

\begin{defn}[\textbf{Local orbital equivalence}]\label{defn_local_orbital_equivalence}
For $i=1,2$ let $U_i\subset M_i$ be a non-empty open subset. The pseudo-flows $(\phi^1,U_1)$ and $(\phi^2,U_2)$ are $C^r$-\emph{locally orbitally equivalent}, where $r\geq 0$, if there exists a $C^r$-diffeomorphism $H:U_1\to U_2$ such that, for every $x\in U_1$, it sends the orbit segment $\mathcal{O}_{U_1}^1(x)$ homeomorphically onto the orbit segment $\mathcal{O}_{U_2}^2(H(x))$, preserving the orientation of these orbit segments. We denote it by
$H:(\phi^1,U_1)\to(\phi^2,U_2).$
\end{defn}

If $\gamma$ is a periodic orbit of $(\phi,M)$ and $W,\ W'$ are two neighborhoods of $\gamma$, the pseudo-flows obtained by restriction onto these sets are the same, in the sense that they coincide over their common domain of definition $W\cap W'$. The \emph{germ of $\phi$ at $\gamma$} is the equivalence class of the pseudo-flows $\{(\phi,W):W\text{ neighborhood of }\gamma\}$ under this relation, and we denote it by $(\phi,M)_\gamma$ or $(\phi,W)_\gamma$. By \emph{abuse of terminology}, we use the word \emph{germ} for referring to a given pseudo-flow $(\phi,W)$, instead of its equivalence class. Through applications, it will be understood that the neighborhood $W$ can be freely replaced by a smaller one if needed.

\begin{defn}[\textbf{Local orbital equivalence of germs}]
For $i=1,2$ let $\gamma_i$ be a periodic orbit of a flow $(\phi^i,M_i)$. The germs $(\phi^1,M_1)_{\gamma_1}$ and $(\phi^2,M_2)_{\gamma_2}$ are $C^r$-\emph{locally orbitally equivalent}, where $r\geq 0$, if there exists a $C^r$-local orbital equivalence $H:(\phi^1,W_1)\to(\phi^2,W_2)$ defined between some neighborhoods $W_i$ of $\gamma_i$. We denote it by
$H:(\phi^1,W_1)_{\gamma_1}\to(\phi^2,W_2)_{\gamma_2}.$   
\end{defn}

\begin{defn}[\textbf{Conjugation of flows}]
If an orbital equivalence $H:(\phi^1,M_1)\to(\phi^2,M_2)$ satisfies in addition that  $\phi_t^2(H(x))=H(\phi_t^1(x))$, for every $x\in M_1$ and $t\in\R$, we say that $H$ is a \emph{conjugation} of the flow actions. 
\end{defn}

In the case where the flows $(\phi^i,M_i)$ are generated by a vector field $X_i$ of class $C^1$ and the conjugation is at least $C^1$, then conjugation between flows is equivalent to the condition
$$X_2(z)=H_*(X_1)(z)=DH(H^{-1}(z))\cdot X_1(H^{-1}(z)),\ \text{ for all }z\in U_2.$$
The same considerations apply in the case of local orbital equivalence and pseudo-flows. 

\subsection{Transverse sections and first return map.}\label{subsection_transverse_sect_and_1st-return}

Let $\phi=\{\phi_t\}_{t\in\R}$ be a continuous, regular, non-singular flow, acting on a 3-manifold $M$.

\begin{defn}[\textbf{Transverse section}]\label{transverse_sections}
A \emph{transverse section} for $(\phi,M)$ is a boundaryless embedded surface $\Sigma\subset M$, satisfying:
\begin{enumerate}[(i)]
\item The surface $\Sigma$ is \emph{topologically transverse to the flow lines}. That is, $\forall\ x\in\Sigma$ there exists a neighborhood $W$ of $x$ inside $M$ and some $\delta>0$ such that:
\begin{itemize}
\item $\Sigma\cap W$ is connected and $W\backslash\Sigma$ has two connected components;
\item  $\Sigma\cap[\phi_{-\delta}(x),\phi_{\delta}(x)]=\{x\}$ and it is verified that $\interior\left([\phi_{-\delta}(x),x]\right)$ is contained in one component of $W\backslash\Sigma$ and $\interior\left([x,\phi_{\delta}(x)]\right)$ is contained in the other.
\end{itemize}
\item The surface is \emph{proper with respect to the flow lines}. That is, every compact orbit segment $[\phi_{t_1}(x),\phi_{t_2}(x)]$ intersects $\Sigma$ in a compact set, where $t_1,t_2\in\R$ and $x\in M$.
\end{enumerate}
\end{defn}

Observe that conditions (i) and (ii) in the definition actually imply that every compact orbit segment cuts the transverse section in a finite set. 
Let $\Sigma\subset M$ be a transverse section for a flow $\phi$ acting on $M$, and assume that there exists a non-empty open subset $U\subset\Sigma$ where there is a well-defined \textit{first return map} $P:U\to\Sigma$. That is, for every $x\in U$ there exists $\tau(x)=\min\{t>0:\phi_t(x)\in\Sigma\}$ and 
\begin{itemize}
\item The function $\tau:U\to(0,+\infty)$ is continuous,
\item The first return map is given by $P(x)=\phi_{\tau(x)}(x)$.
\end{itemize}

If $x$ is a point in $\Sigma$ satisfying that there exists some $t>0$ such that $\phi_t(x)\in\Sigma$, then from the continuity of the flow it follows that there exists a first return map defined in a neighborhood of $x$. In fact, a first return map as above is always a homeomorphism from its domain onto its image. Moreover, as it follows from the \emph{implicit function theorem}, it is a $C^l$-diffeomorphism, where $l$ is the minimum between the regularities of $\phi$ and $\Sigma$. We use the following definition:

\begin{figure}[t]
\begin{center}
\scalebox{0.4}{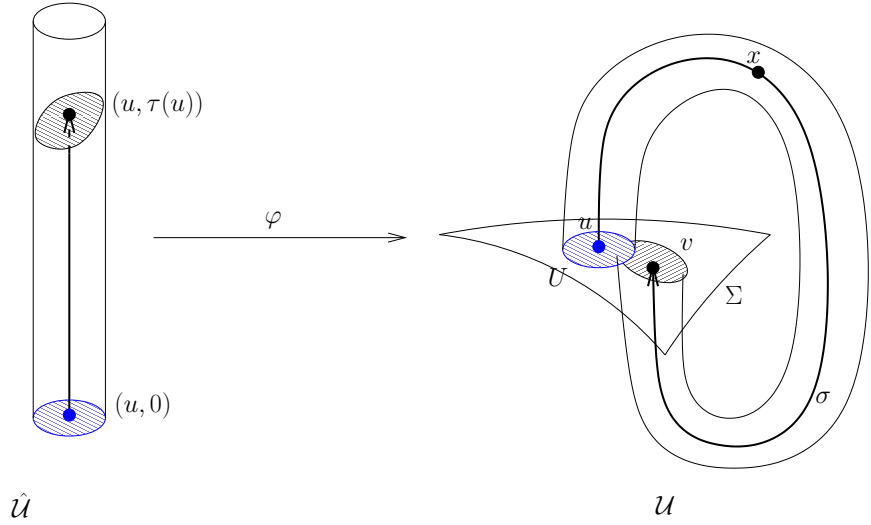}
\caption{First return saturation $\mathcal{U}$ of the open set $U\subset\Sigma$.}
\label{fig_1_2_1}
\end{center}
\end{figure}

\begin{defn}\label{saturation_by_orbit_segments}
The \emph{first return saturation} of $U$ is the set 
\begin{equation*}
\mathcal{U}=\{x\in M:\exists\ t\leq 0\text{ such that }\phi_t(x)\in U\text{ and }\interior\left([\phi_t(x),x]\right)\cap\Sigma=\emptyset\}.
\end{equation*}
\end{defn}
This set is the union of all compact orbits segments joining each point in $U$ with its first return to $\Sigma$, as we see in Figure \ref{fig_1_2_1}. The following lemma is elementary, see \cite{shannon_thesis} for a proof.

\begin{lem}[\cite{shannon_thesis}, Lemma 1.9]\label{lemma_conjugacion_vs_equiv_orbital}
Let $\{\phi^i_t:M_i\to M_i\}_{t\in\R}$, $i=1,2$ be a pair of regular, non-singular flows. Let $\Sigma_i\subset M_i$ be a transverse section for each flow and assume there exists a first return map $P_i:U_i\to\Sigma_i$ defined in an open set $U_i\subset\Sigma_i$. If there exists a homeomorphism $h:\Sigma_1\to\Sigma_2$ such that $h(U_1)=U_2$ and $h\circ P_1(x)=P_2\circ h(x)$, $\forall$ $x\in U_1$, then there exists a homeomorphism $H:\mathcal{U}_1\to\mathcal{U}_2$ such that:
\begin{enumerate}
\item For every point $x\in\mathcal{U}_1$ the map $H$ takes the orbit segment $\mathcal{O}_{\mathcal{U}_1}(x)$ onto $\mathcal{O}_{\mathcal{U}_2}(H(x))$,

\item The restriction map $H|_{U_1}$ coincides with $h|_{U_1}$.
\end{enumerate} 
\end{lem}

A \emph{global transverse section} $\Sigma$ for $(\phi,M)$ is a transverse section that is properly embedded in $M$ and for which there exists $T>0$ such that $[x,\phi_T(x)]\cap\Sigma\neq\emptyset$, for every $x\in M$. In this case, the first return map is a homeomorphism $P_\Sigma:\Sigma\to\Sigma$. Given a periodic orbit $\gamma$ of the flow, a \textit{local transverse section at $\gamma$} is a transverse section $D$, homeomorphic to a disk, such that $\{x_0\}=\gamma\cap D$ contains exactly one point. Given a local transverse section $D$ there always exists a neighborhood $U\subset D$ of $x_0$ and a first return map $P_D:U\to D$ that fixes $x_0$. The following statement is a direct corollary of Lemma \ref{lemma_conjugacion_vs_equiv_orbital} above.

\begin{proposition}\label{prop_conjugacion_vs_equiv_orbital}
Let $\{\phi^i_t:M_i\to M_i\}_{t\in\R}$, $i=1,2$, be two regular, non-singular flows.
\begin{enumerate}
\item Assume there exists a global transverse section $\Sigma_i$ for each flow $\phi^i$ and let $P_i:\Sigma_i\to\Sigma_i$ be the first return map. If there exists a homeomorphism $h:\Sigma_1\to\Sigma_2$ such that $h\circ P_1=P_2\circ h$, then there exists a homeomorphism $H:M_1\to M_2$ such that:
\begin{enumerate}
\item $H$ is an orbital equivalence between the flows,
\item $H|_{\Sigma}=h$.
\end{enumerate}
\item Let $\gamma_i$ be a periodic orbit of each flow $\phi^i$, $D_i$ a local transverse section and $P_{D_i}:U_i\to D_i$ a first return map defined in a neighborhood $U_i\subset D_i$ of the intersection point $x_i=\gamma_i\cap D_i$. If there is a homeomorphism $h:D_1\to D_2$ such that $h(U_1)\subset U_2$ and $h\circ P_{D_1}(x)=P_{D_2}\circ h(x)$, $\forall$ $x\in U_1$, then there exists a tubular neighborhood $W_i$ of each $\gamma_i$ and a homeomorphism $H:W_1\to W_2$ such that:
\begin{enumerate}
\item $H$ is a local orbital equivalence between the respective germs at each $\gamma_i$,
\item $H|_{D_1\cap W_1}=h|_{D_1\cap W_1}$.
\end{enumerate}
\end{enumerate}
\end{proposition}

\subsection{Topologically Anosov flows.}\label{subsection_top-Anosov_defn}

\begin{defn}[\textbf{Hyperbolic splitting}]\label{defn_hyperbolic_splitting}
Let $M$ be a smooth $3$-manifold equipped with a Riemannian metric $\Vert\cdot\Vert$ and let $\{\phi_t:M\to M\}_{t\in\R}$ be a flow generated by a non-singular $C^k$-vector field $X$, where $k\geq 1$. A \emph{hyperbolic splitting} of the tangent bundle $TM$ is a decomposition as Whitney sum of three line bundles $TM=\Es\oplus\Ec\oplus\Eu$, where $\Ec=\spa\{X\}$, satisfying that:
\begin{enumerate}
\item Each line bundle is invariant under the action of $D\phi_t:TM\to TM$, for every $t\in\R$;
\item There exist constants $C>0$ and $0<\lambda<1$ such that 
\begin{align}\label{contraction/expansion_rates}
& \Vert D\phi_t(x)\cdot v\Vert\leq C\lambda^t\Vert v\Vert,\ \forall\ v\in\Es(x),\ t\geq 0,\ x\in \Lambda;\\
& \Vert D\phi_t(x)\cdot v\Vert\leq C\lambda^{-t}\Vert v\Vert,\ \forall\ v\in\Eu(x),\ t\leq 0,\ x\in \Lambda.\nonumber
\end{align}
\end{enumerate}
\end{defn}

The bundle $\Es$ is called the \textit{stable bundle}, $\Eu$ is called the \textit{unstable bundle} and $\Ec$, the one who is tangent to the flow lines, is called the \textit{central bundle}. The two dimensional bundles $\Ecs=\Es\oplus\Ec$, $\Ecu=\Ec\oplus\Eu$ and $\Esu=\Es\oplus\Eu$ are respectively called the \textit{center-stable bundle} (or \textit{cs}-bundle), the \textit{center-unstable bundle} (or \textit{cu}-bundle) and the \textit{stable-unstable bundle} (or \textit{su}-bundle). 

\begin{defn}[\textbf{Anosov flow}]\label{defn_Anosov_flow}
Let $M$ be a \emph{closed} smooth $3$-manifold and consider a flow $\{\phi_t:M\to M\}_{t\in M}$ generated by a non-singular $C^k$-vector field $X$, where $k\geq 1$. The flow is \emph{Anosov} if its derivative action on $TM$ preserves a hyperbolic splitting, for some given Riemannian metric on $M$. 
\end{defn}

\begin{rem}
Observe that from the compactness of the ambient space, in the case of an Anosov flow it follows that the vectors in $\Es$ or $\Eu$ satisfy condition \eqref{contraction/expansion_rates} above for any chosen Riemannian metric on $M$ and for every reparametrization of the flow, up to modifying the constants $C>0$ and $0<\lambda<1$ if necessary (cf. \cite{katok-hasselblatt}). Thus, the definition of Anosov flow makes an auxiliary use of a Riemannian metric, but it only depends on the $C^1$-equivalence class of the flow. As well, it is not difficult to check that the decomposition of $TM$ must be unique and continuous. 

This observation is no longer true if $M$ is not compact. In this case, condition \eqref{contraction/expansion_rates} in Definition \ref{defn_hyperbolic_splitting} above strongly depends on the choice of Riemannian metric. 
\end{rem}

One of the fundamental properties of Anosov flows is the integrability of its (center-) stable and (center-)unstable bundles into foliations which are preserved by the flow action, and can be merely defined by dynamical properties. This well-known property is called \emph{stable manifold theorem}, see for example \cite{katok-hasselblatt}.

\begin{thm*}[\textbf{Stable manifold theorem}]
Let $\{\phi_t:M\to M\}_{\in\R}$ be an Anosov flow. Then each 1-dimensional bundle $\Es$ and $\Eu$ is uniquely integrable, and the partition of $M$ by integral curves respectively determines a pair of 1-dimensional foliations $\Fs$ and $\Fu$, invariant by the action of the flow. Moreover, for every $x\in M$ it is satisfied that
\begin{align*}
&\Fs(x)=\Ws(x)=\{y\in M:\dist(\phi_t(y),\phi_t(x))\to 0,\ t\to+\infty\},\\
&\Fu(x)=\Wu(x)=\{y\in M:\dist(\phi_t(y),\phi_t(x))\to 0,\ t\to-\infty\}.
\end{align*}
In addition, each of the bundles $\Ecs$ and $\Ecu$ is uniquely integrable into a $2$-dimensional foliation $\mathcal{F}^{cs}$ and $\mathcal{F}^{cu}$ respectively, and for every $x\in M$ it is satisfied that:
\begin{align*}
&\mathcal{F}^{cs}(x)=\{y\in M:\exists\ s\in\R,\text{ s.t. }\dist(\phi_t(y),\phi_{t+s}(x))\to 0,\ t\to+\infty\},\\
&\mathcal{F}^{cu}(x)=\{y\in M:\exists\ s\in\R,\text{ s.t. }\dist(\phi_t(y),\phi_{t+s}(x))\to 0,\ t\to-\infty\}.
\end{align*}
\end{thm*}

\begin{defn}[\textbf{Expansivity}]
Let $\{\phi_t:M\to M\}_{t\in\R}$ be a non-singular flow on a closed $3$-manifold. The flow is said to be \emph{orbitally expansive} if for some fixed metric on $M$ and for every $\varepsilon>0$, there exists $\alpha=\alpha(\varepsilon)>0$ such that: If two points $x$, $y$ satisfy that 
$\dist(\phi_{h(t)}(y),\phi_t(x))\leq\alpha,\ \forall\ t\in\R,$
where $h:(\R,0)\to(\R,0)$ is some increasing homeomorphism, then $y=\phi_s(x)$ for some $|s|\leq\varepsilon$.  
\end{defn}

As before, this definition is independent of the chosen metric by compactness of the ambient space. General expansive flows on closed $3$-manifolds have been studied by Paternain, Inaba and Matsumoto. In \cite{paternain_expansive_flows} and \cite{inaba-matsumoto_expansive_flows} it is shown that every orbitally expansive flow is orbitally equivalent to a \emph{pseudo-Anosov flow}. This is obtained by showing that the partitions of $M$ into stable and unstable sets of the orbits, that is
\begin{align*}
&\mathcal{W}^{\mathrm{s}}(\mathcal{O}(x))=\{y\in M:\exists s\in\R\ \text{and }h:(\R,0)\to(\R,0)\text{ s.t. }\dist(\phi_{h(t)}(y),\phi_{t+s}(x))\to 0,\ t\to+\infty\},\\
&\mathcal{W}^{\mathrm{u}}(\mathcal{O}(x))=\{y\in M:\exists s\in\R\ \text{and }h:(\R,0)\to(\R,0)\text{ s.t. }\dist(\phi_{h(t)}(y),\phi_{t+s}(x))\to 0,\ t\to-\infty\},
\end{align*}
constitute a pair of transverse foliations in $M$, possibly with singularities, which are invariant by the flow action and intersect along the flow orbits. We denoted them by $\mathcal{F}^{cs}$ and $\mathcal{F}^{cu}$ respectively. The singularities are of a special type called \emph{circle prongs}, that are obtained by suspending $k$-prong multi-saddle fixed points in a disk. See \cite{brunella_thesis} or the referred works for precise statement and definitions. 

\begin{defn}[\textbf{Topologically Anosov flow}] 
A non-singular flow on a closed $3$-manifold is \emph{topologically Anosov} if it is orbitally expansive and its invariant foliations have no singularities.
\end{defn}

As a final remark, it is possible to see that a flow is topologically Anosov if and only if it is orbitally expansive and satisfies the \emph{shadowing property of Bowen}, since in the presence of invariant foliations then multi-saddle orbits are the only obstruction to shadowing. In general, a \emph{topologically Anosov dynamical system} (discrete- or continuous-time) is one that is expansive and satisfies (globally) the shadowing property. See \cite{aoki-hiraide} for definitions and general properties of topologically Anosov dynamics.  

\subsection{Pseudo-Anosov homeomorphisms.}\label{subsection_psudo-Anosov-homeos}
Let $\Sigma$ be a closed orientable surface. 

\begin{defn}\label{defn_pseudo_Anosov}
A homeomorphism $f:\Sigma\to\Sigma$ is \emph{pseudo-Anosov} if there exists a pair of transverse, $f$-invariant, foliations $\Fs$ and $\Fu$ on $\Sigma$, respectively equipped with \emph{transverse measures} $\mu_\mathrm{s}$ and $\mu_\mathrm{u}$ and a constant $0<\lambda<1$ such that $f_*(\mu_\mathrm{s})=\lambda^{-1}\cdot\mu_\mathrm{s}$ and $f_*(\mu_\mathrm{u})=\lambda^{-1}\cdot\mu_\mathrm{u}$ The transverse measures are required to be non-atomic and with full support.
\end{defn}

\begin{figure}[t]
     \centering
	 \fbox{	     
     \begin{subfigure}[b]{0.45\textwidth}
         \centering
         \includegraphics[width=0.7\textwidth]{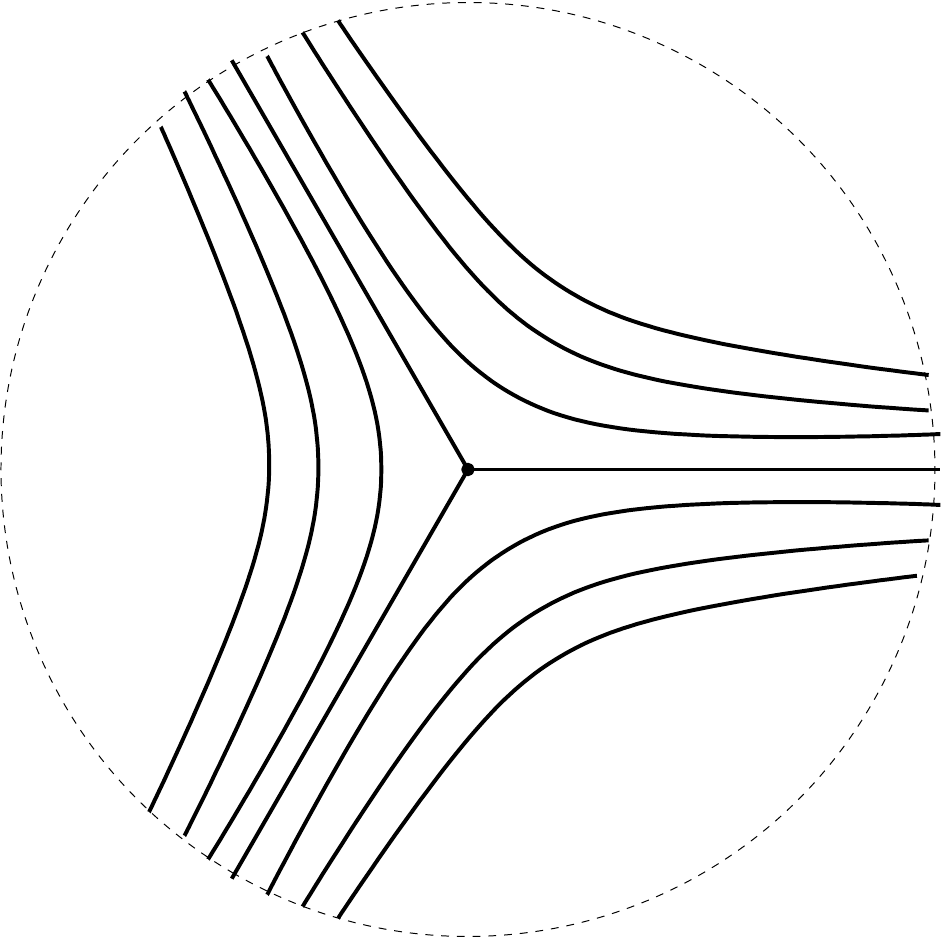}
         \caption{$k$-prong singularity, $k=3$.}
         \label{fig_k-prong}
     \end{subfigure}
	}
	\fbox{    
     \begin{subfigure}[b]{0.45\textwidth}
         \centering
         \includegraphics[width=0.7\textwidth]{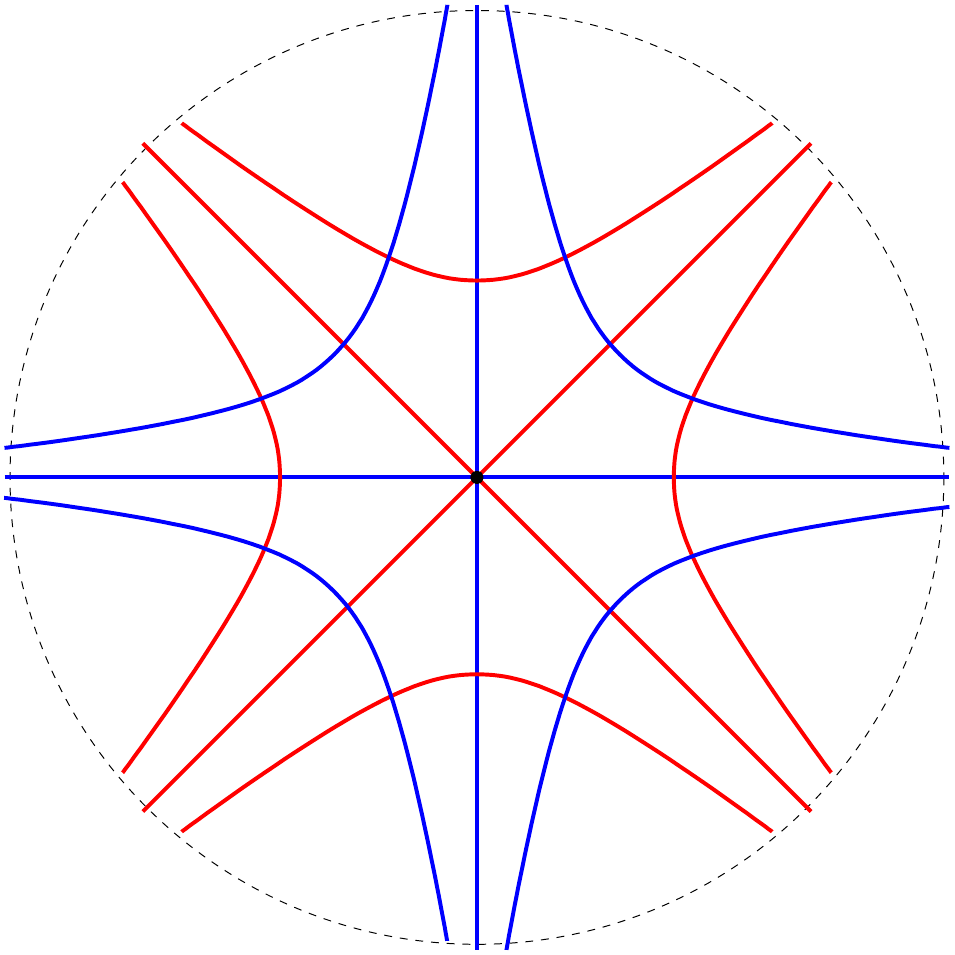}
         \caption{Transverse foliations at $k$-prong, $k=4$.}
         \label{fig_k-prong_2}
     \end{subfigure}
     }
     \caption{Local model of $k$-prong singularities}
\end{figure}

If the genus of $\Sigma$ is greater than one, then the two foliations necessarily have singularities. If this is the case, in the previous definition we just allow singularities whose local model is a \emph{$k$-prong singularity} (see Figure \ref{fig_k-prong}) and with $k\geq 3$. Since the foliations must be transverse between them, then $\Fs$ and $\Fu$ share the same (finite) set of singularities, and in a small neighborhood of each singularity the two foliations intersect as in the local model \ref{fig_k-prong_2}. 

By the Euler-Poincar\'e formula (see \cite{FLP}), the sphere is excluded from having a pseudo-Anosov homeomorphism and in the torus the foliations are non-singular. If $\textrm{gen}(\Sigma)\geq 2$, observe that the set of singularities $\Delta$ is necessarily included in the set of periodic points of $f$. 

For higher genus surfaces, pseudo-Anosov homeomorphisms can be seen as a counter-part of linear hyperbolic automorphisms of the torus. In view of the singularities, they can never be \emph{hyperbolic diffeomorphisms}, but they share some properties the latter ones. In particular, every pseudo-Anosov homeomorphism $f:\Sigma\to\Sigma$ is transitive, expansive, the set of periodic points is dense in $\Sigma$, and the topological entropy of $f$ is positive. More generally, the dynamic of a pseudo-Anosov map can be encoded using a Markovian partition constructed from its invariant foliations, as is shown in \cite{FLP}. From the symbolic point of view these maps are equivalent to subshifts of finite type. 

In \cite{lewowicz_exp} and \cite{hiraide-Expansivos}, Lewowicz and Hiraide have shown that the pseudo-Anosov maps are the only expansive homeomorphisms in a closed orientable surface. More precisely, if $f:\Sigma\to\Sigma$ is an expansive homeomorphism then $\textrm{gen}(\Sigma)\geq 1$ and if $\textrm{gen}(\Sigma)=1$ then $f$ is $C^0$-conjugated to a linear Anosov map, and in higher genus $f$ is $C^0$-conjugated to some pseudo-Anosov homeomorphism. 

\subsubsection{Smooth models.} In \cite{gerber-katok} Gerber and Katok proved that \emph{every pseudo-Anosov homeomorphism $f$ on a closed surface $\Sigma$ is $C^0$-conjugated to some diffeomorphism $g:\Sigma\to\Sigma$, which in addition preserves an ergodic smooth measure.} See also \cite{lewowicz_anal} for an analytic version. Let us point out some interesting facts about this result. 

If $f:\Sigma\to\Sigma$ is a pseudo-Anosov homeomorphism, then the system of transverse foliations equipped with transverse measures defines a \emph{translation atlas} $\mathcal{A}_\Delta$ in the complement of the set of singularities $\Delta$ (see e.g. \cite{FLP}). Choose some smooth atlas $\mathcal{A}$ on $\Sigma$ such that the inclusion map $(\Sigma\backslash\Delta,\mathcal{A}_\Delta)\hookrightarrow(\Sigma,\mathcal{A})$ induces a diffeomorphism onto its image. From now on we consider $\mathcal{A}$ as a fixed smooth structure on $\Sigma$. 

Directly from the construction we can see that the restriction $f:\Sigma\backslash\Delta\to\Sigma\backslash\Delta$ is a smooth diffeomorphism, and that it is not differentiable on the singular points $x\in\Delta$. In \cite{gerber-katok}, Proposition at page 177, it is showed that \emph{$f$ cannot be $C^0$-conjugated to a diffeomorphism $g:\Sigma\to\Sigma$, via a homemorphism that is a $C^1$-diffeomorphism except at the singularities of $f$}. This result says that singularities of the smooth atlas $\mathcal{A}_\Delta$ are essential for the dynamics of $f$, in the sense that it cannot be extended onto a globally defined $C^1$-atlas, in such a way that $f$ turns into a diffeomorphism. The principal result in \cite{gerber-katok} gives a smooth atlas $\mathcal{B}$ on $\Sigma$ such that, in the coordinates of this atlas, $f$ is a smooth diffeomorphism. Nevertheless, the atlas $\mathcal{B}$ is nowhere compatible with $\mathcal{A}_\Delta$.

\begin{rem}
Question \ref{question:top-Anosov=Anosov} is similar to the question treated by Gerber and Katok in the referred work, but translated into the context of topologically Anosov flows. It should be noted that our procedure for constructing hyperbolic models of topologically Anosov flows follows the same general strategy of that in \cite{gerber-katok} as we explain in Section \ref{section_hyperbolic_models}. In our context, the role of the translation atlas with singularities is played by what we call an \emph{almost Anosov atlas} in Section \ref{section_almost-Anosov}. However, we have not been able to obtain an analogous $C^1$ non-extension result as that in \cite{gerber-katok}, see Remark \ref{non-extension_Bonatti}. 
\end{rem}

\subsubsection{Conjugacy classes of pseudo-Anosov homeomorphisms.}
A remarkable property about pseudo-Anosov homeomorphisms is that if two of these maps are isotopic, then they are conjugated by a homeomorphism isotopic to the identity.

\begin{thm}[\cite{FLP}, Expos\'e XII, Theorem 12.5]\label{isotopic_pseudo_Anosov}
Let $\Sigma$ be a closed orientable surface and let $f$ and $g$ be two pseudo-Anosov homeomorphisms. If $g$ is isotopic to $f$ then there exists a homeomorphism $h:\Sigma\to\Sigma$, isotopic to the identity, such that $f\circ h=h\circ g$.
\end{thm}

This theorem, in combination with the \emph{Dehn-Nielsen-Baer theorem} about mapping class groups, allows to decide if two pseudo-Anosov homeomorphisms are conjugated by looking at their actions on fundamental groups. In the following we explain this fact in the way that it will be used later.

\paragraph*{The action on the fundamental group.}
Let $x_0$ be a point in $\Sigma$. Every homeomorphism $f\in\textrm{Homeo}(\Sigma)$ induces an automorphism of $\pi_1(\Sigma,x_0)$ in the following way: Let $\beta:[0,1]\to\Sigma$ be an arc that connects $x_0=\beta(0)$ with $f(x_0)=\beta(1)$. Given a class $[\gamma]\in\pi_1(\Sigma,x_0)$ represented by a curve $\gamma:[0,1]\to\Sigma$ we define 
$$f_*^\beta:[\gamma]\mapsto\left[\bar{\beta}\cdot f(\gamma)\cdot\beta\right],\text{ where }\bar{\beta}\text{ is }\beta\text{ parametrized with inverse sense}.$$
The map $f_*^\beta$ is a well-defined automorphism of $\pi_1(\Sigma,x_0)$ which depends on the particular election of the arc $\beta$. If we choose another arc $\beta'$ connecting $x_0=\beta'(0)$ with $f(x_0)=\beta'(1)$, then $f_*^{\beta'}=[\alpha]^{-1}\cdot f_*^\beta \cdot[\alpha]$ where $[\alpha]=[\bar{\beta}\cdot\beta']\in\pi_1(\Sigma,x_0)$. Thus, changing the arc $\beta$ has the effect of conjugating $f_*^\beta$ by an inner automorphism of the fundamental group $\pi_1(\Sigma,x_0)$. 

\begin{defn}
Given a pair of homeomorphisms $f_i:\Sigma_i\to\Sigma_i$, $i=1,2$, where $\Sigma_1$ and $\Sigma_2$ are two homeomorphic closed orientable surfaces, we say that $f_1$ and $f_2$ are \emph{$\pi_1$-conjugated} if there exist points $x_i\in\Sigma_i$, induced actions $(f_i)_*^{\beta_i}:\pi_1(\Sigma_i,x_i)\to\pi_1(\Sigma_i,x_i)$ and an isomorphism $\phi:\pi_1(\Sigma_1,x_1)\to\pi_1(\Sigma_2,x_2)$ such that $(f_2)_*^{\beta_2}\circ\phi=\phi\circ (f_1)_*^{\beta_1}$.
\end{defn}

Observe that if $f_1$ and $f_2$ are $\pi_1$-conjugated then, for every pair of points $x_i$, every pair of induced actions $(f_i)_*^{\beta_i}$ on $\pi_1(\Sigma_i,x_i)$ are conjugated by an isomorphism $\phi:\pi_1(\Sigma_1,x_1)\to\pi_1(\Sigma_2,x_2)$.

\begin{proposition}[\cite{shannon_thesis}, Proposition 1.28]\label{pi_1-conjugated_pseudo_Anosov}
For $i=1,2$ consider a pseudo-Anosov homeomorphism $f_i:\Sigma_i\to\Sigma_i$ defined in a closed orientable surface $\Sigma_i$. If $f_1$ and $f_2$ are $\pi_1$-conjugated, then there exists a homeomorphism $h:\Sigma_1\to\Sigma_2$ such that $f_2\circ h=h\circ f_1$.
\end{proposition}

\subsubsection{The case of punctured surfaces.}\label{section_pseudo-Anosov_punctured_surfaces}
Consider two pseudo-Anosov homeomorphisms $f_i:\Sigma_i\to\Sigma_i$, where $i=1,2$, defined in two closed, orientable, homeomorphic surfaces. For each $f_i$ consider a finite collection of periodic orbits $\mathcal{O}^i_1,\dots,\mathcal{O}^i_N$. That is, for each $k=1,\dots,N$ let
$\mathcal{O}^i_k=\{x^i_{k1},\dots,x^i_{kp_k}\}$
where $x^i_{kn}=f_i^{n-1}(x^i_{k1})$ and $p_k\geq 1$ is the period of the orbit. We are interested in knowing when $f_1$ is conjugated to $f_2$ by a homeomorphism $h:\Sigma_1\to\Sigma_2$, with the additional property that $h$ sends each orbit $\mathcal{O}^1_k$ to the orbit $\mathcal{O}^2_k$. 

A \emph{finite type punctured surface} is the data of a compact surface $\Sigma$ together with a finite subset $\mathcal{O}\subset\Sigma$. In this text we just consider the case where the surface is closed. The \emph{mapping class group of the punctured surface $(\Sigma;\mathcal{O})$} is defined to be the set homeomorphisms $f:\Sigma\to\Sigma$ preserving the set $\mathcal{O}$, modulo isotopies fixing $\mathcal{O}$. 

Denote by $\pi_1(\Sigma,x_0;\mathcal{O})$ the fundamental group of $\Sigma\backslash\mathcal{O}$ based in a point $x_0$ not in $\mathcal{O}$. If $f$ is a homeomorphism of $\Sigma$ preserving $\mathcal{O}$, then it induces a permutation of this finite set as well as an action on $\pi_1(\Sigma,x_0;\mathcal{O})$, uniquely defined up to conjugation by inner automorphisms. 

Let $\mathcal{O}=\{x_1,\dots,x_R\}$. For each of the points $x_l$ consider a closed curve homotopic to the puncture $x_l$ and joined to $x_0$ with an arbitrary path. This curve determines an element of $\pi_1(\Sigma,x_0;\mathcal{O})$ that we denote by $c_l$. We denote by $\Gamma(\mathcal{O})$ the set of conjugacy classes of the elements $c_l$. That is,
$$\Gamma(\mathcal{O})=\{\bar{\gamma}\cdot c_l\cdot\gamma:\text{ where }l=1,\dots,R\text{ and }\gamma\in\pi_1(\Sigma,x_0;\mathcal{O})\}.$$
The action $f_*$ leaves invariant the set $\Gamma(\mathcal{O})$, since it comes from a homeomorphism of the surface. 

\begin{proposition}[\cite{shannon_thesis}, Proposition 1.29]\label{pi_1-conjugated_pseudo_Anosov_punctures}
For $i=1,2$ consider a pseudo-Anosov homeomorphism $f_i:\Sigma_i\to\Sigma_i$ defined in a closed orientable surface $\Sigma_i$ and a finite collection of periodic orbits $\mathcal{O}^i_1,\dots,\mathcal{O}^i_N$ of periods $p_1,\dots,p_N$, respectively. 
Assume there exists an isomorphism $$\phi:\pi_1(\Sigma_1;\mathcal{O}^1_1\cup\dots\cup\mathcal{O}^1_N)\to\pi_1(\Sigma_2;\mathcal{O}^2_1\cup\dots\cup\mathcal{O}^2_N)$$ such that:
\begin{enumerate}
\item $\phi$ conjugates the actions $(f_i)_*$ induced on fundamental groups of $\Sigma_i\backslash\mathcal{O}^i_1\cup\cdots\cup\mathcal{O}^i_N$,
\item $\phi(\Gamma(\mathcal{O}^1_k))=\Gamma(\mathcal{O}^2_k)$ for every $k=1,\dots,N$.
\end{enumerate}

Then, there exists a homeomorphism $h:\Sigma_1\to\Sigma_2$ such that $f_2\circ h=h\circ f_1$ which in addition satisfies $h(\mathcal{O}^1_k)=\mathcal{O}^2_k$, $\forall\ k=1,\dots,N$.
\end{proposition}

\subsubsection{Pseudo-Anosov on non-closed surfaces.}\label{section_pseudo-Anosov_non-closed_surfaces}
We will use Proposition \ref{pi_1-conjugated_pseudo_Anosov_punctures} in the course of the proofs of theorems \textbf{A} and \textbf{B}. By the way, since we work with non-closed surfaces, we make some remarks in the sequel. Let $\Sigma$ be a compact orientable surface.

\begin{defn}[\textbf{Pseudo-Anosov on non-closed surfaces}]
Let $f:\Sigma\to\Sigma$ be a homeomorphism. Let $\widehat{\Sigma}$ the surface obtained by collapsing each boundary component into a point and let $\hat{f}$ be the corresponding induced map on the quotient. We say that $f$ is pseudo-Anosov if $\hat{f}$ is pseudo-Anosov according to Definition \ref{defn_pseudo_Anosov}.
\end{defn}

\begin{figure}[t]
     \centering
	 \fbox{	     
     \begin{subfigure}[b]{0.4\textwidth}
         \centering
         \includegraphics[width=0.7\textwidth]{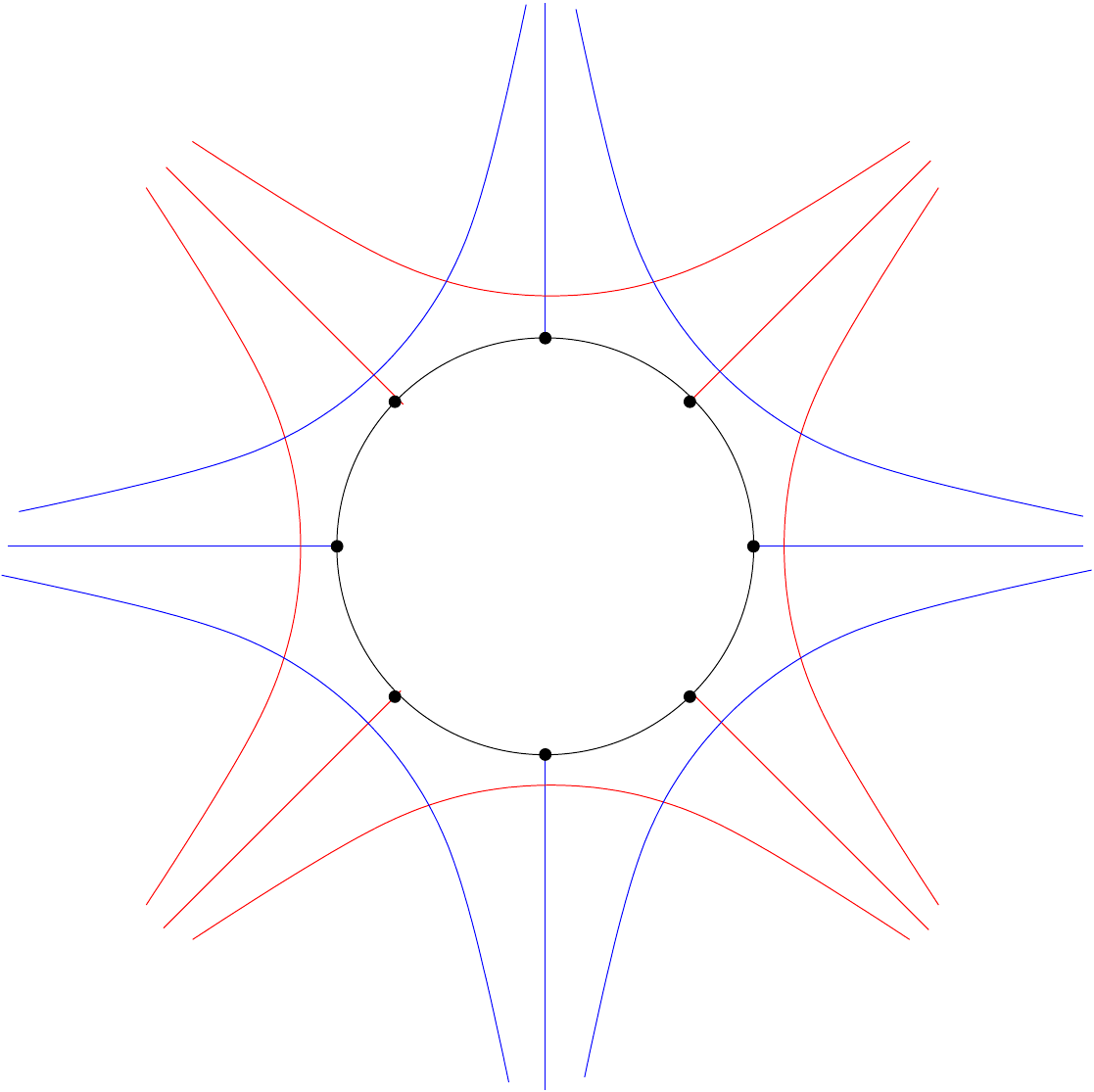}
         \caption{Tame}
         \label{fig_tame_boundary}
     \end{subfigure}
	}
	\fbox{    
     \begin{subfigure}[b]{0.4\textwidth}
         \centering
         \includegraphics[width=0.7\textwidth]{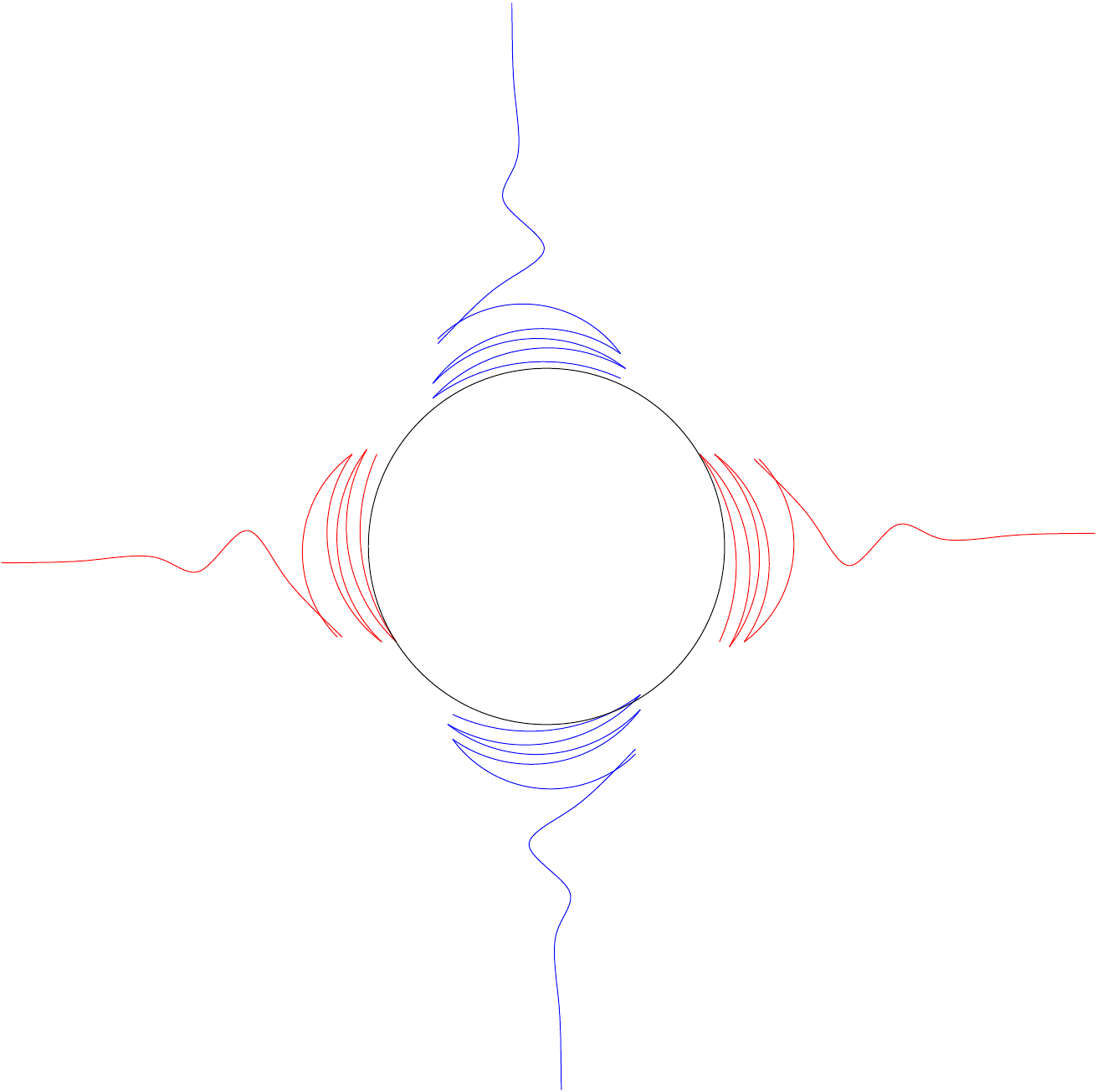}
         \caption{Wild}
         \label{fig_wild_boundary}
     \end{subfigure}
     }
     \caption{Local model of $k$-prong singularities}
\end{figure}

Let $C$ be a boundary component of $\Sigma$ and $p\in\widehat{\Sigma}$ the point obtained after collapsing $C$. Each invariant foliation of $\hat{f}$ has a finite number of leaves which accumulate on $p$, which are usually called \emph{branches}. When lifted to $\Sigma$, these branches do not necessarily converge to a point in $C$, as is depicted in Figure \ref{fig_wild_boundary}. We say that the foliations are \emph{tame} if the local model in a neighborhood of the boundary component $C$ is as in Figure \ref{fig_tame_boundary}. In this case, each branch converges to a point in $C$, which is necessarily periodic for~$f$. 

Collapsing each boundary component into a point provides a semi-conjugation from $f$ to $\hat{f}$, which is actually a conjugation on the interior of the surface. So, most of the dynamical properties of $\hat{f}$ are also available for~$f$. However, we want to point out the following:

\begin{rem} 
In the case of non-closed surfaces, Theorem \ref{isotopic_pseudo_Anosov} as well as Proposition \ref{pi_1-conjugated_pseudo_Anosov} are no longer available.
\end{rem}

We want to point out an obstruction to conjugacy taken from \cite{Palis_saddle-connection-stability}, which depends on the behavior of the map in a neighborhood of the boundary. 

\begin{figure}[t]
\begin{center}
\includegraphics[width=0.4\textwidth,keepaspectratio,angle=0]{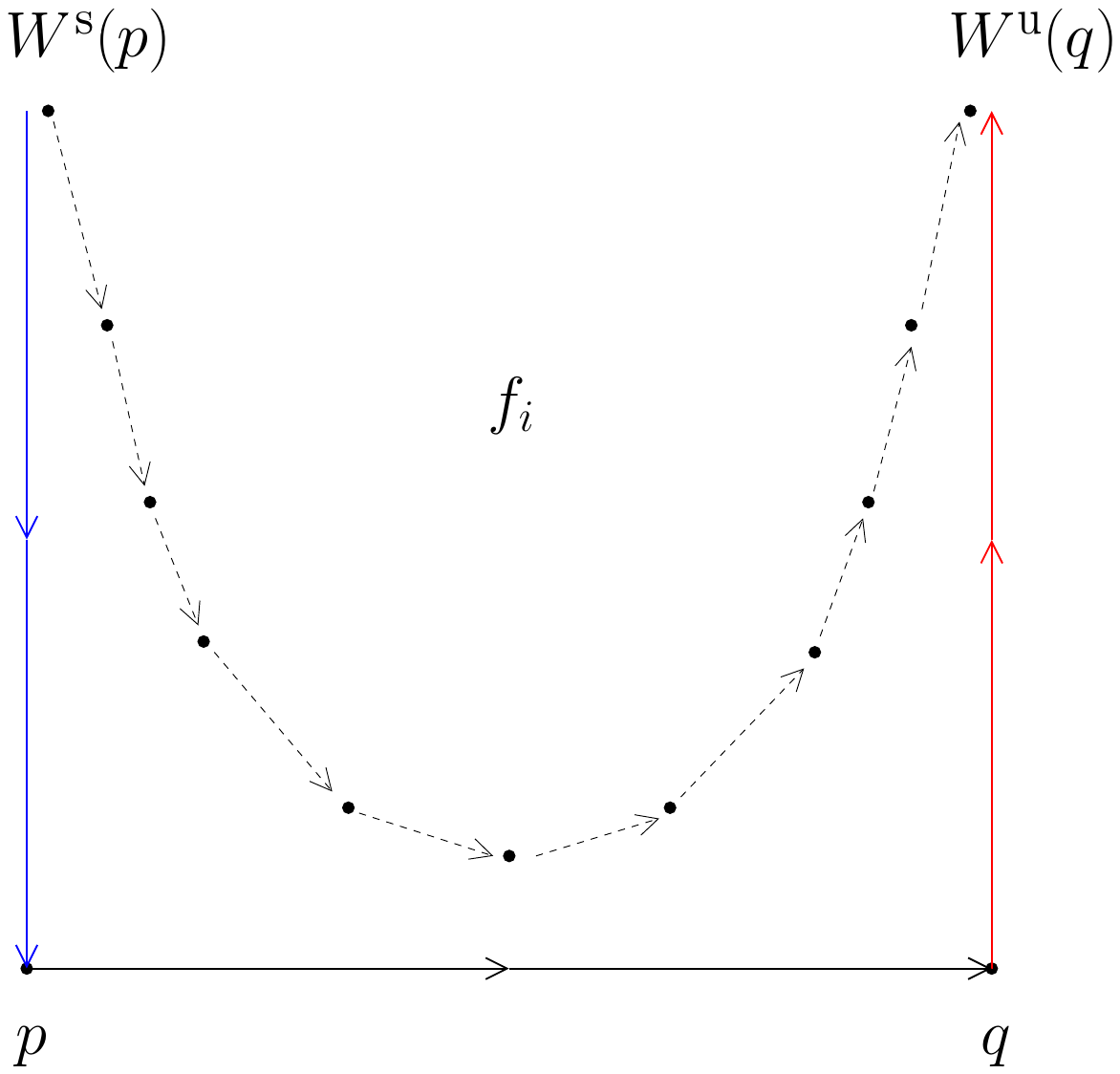}
\caption{The homeomorphisms in Example \ref{example_non_conjugacy}.}
\label{non-conjugated_boundary}
\end{center}
\end{figure}

\begin{ex}\label{example_non_conjugacy}
Consider two $C^1$-diffeomorphisms $f$ and $g$ defined in the band $[-1,1]\times[0,+\infty)$, whose phase portrait is as follows:
\begin{itemize}
\item The non-wandering set consists in the corner points $p=(-1,0)$ and $q=(1,0)$, which are saddle type hyperbolic fixed points. 
\item $\Ws(p)=\{-1\}\times[0,+\infty)$,
\item $\Wu(q)=\{1\}\times[0,+\infty)$,
\item The segment $[-1,1]\times\{0\}$ is a saddle connection between $p$ and $q$.
\end{itemize}

This is illustrated in Figure \ref{non-conjugated_boundary}. For each fixed point $x=p,q$ we have  
$$Df(x)=
\begin{pmatrix}
\lambda_x(f) & 0\\
0 & \mu_x(f)
\end{pmatrix}
\text{ and }
Dg(x)=
\begin{pmatrix}
\lambda_x(g) & 0\\
0 & \mu_x(g)
\end{pmatrix}.
$$
\end{ex}

\begin{prop*}[Palis \cite{Palis_saddle-connection-stability}]
If there exists a homeomorphism $h:[-1,1]\times[0,+\infty)\to[-1,1]\times[0,+\infty)$ such that $g\circ h=h\circ f$, then
\begin{equation*}
\frac{\log(\mu_q(f))}{\log(\mu_p(f))}=\frac{\log(\mu_q(g))}{\log(\mu_p(g))}. 
\end{equation*} 
\end{prop*}

In particular, general homeomorphisms (even $C^1$-diffeomorphisms) whose phase portrait is as in Example \ref{example_non_conjugacy} are not $C^0$-conjugated. Observe that there always exists a conjugation between these dynamics in the complement of the segment $[-1,1]\times\{0\}$. This can be seen by dividing the band into adequate fundamental domains for the action of each map. The obstruction appears when we try to extend the conjugation to the segment that connects the two saddles. 

This dynamical behavior is what we encounter in the neighborhood of a boundary component of a surface $\Sigma$, when we look at the action of a pseudo-Anosov map. If we choose a power of $f$ that fixes the boundary component, then we can decompose a neighborhood of this component into a finite number of bands homeomorphic to $[-1,1]\times[0,+\infty)$ where the dynamic looks like in the Example \ref{example_non_conjugacy}. As a final remark, observe that if $\hat{f}:\widehat{\Sigma}\to\widehat{\Sigma}$ is a pseudo-Anosov in a closed surface, we can construct a pseudo-Anosov $f:\Sigma\to\Sigma$ in a non-closed surface by blowing up $\hat{f}$ along a periodic orbit. But, in view of \ref{example_non_conjugacy}, different ways of blowing up could lead to non-conjugated maps, even if all the actions on $\pi_1(\Sigma)$ are the same.

\section{Birkhoff sections and orbital equivalence classes.}\label{section_Birk_sect_and_equiv}
Let $\{\phi_t:M\to M\}_{t\in\R}$ be a regular, non-singular, continuous flow on a closed $3$-manifold.

\begin{defn}\label{defn_Birkhoff_section}
A \textit{Birkhoff section} for $(\phi,M)$ is an immersion $\iota:(\Sigma,\partial\Sigma)\to(M,\Gamma)$, where $\Sigma$ is a compact surface, sending $\partial\Sigma$ onto a finite set $\Gamma$ of periodic orbits of $\phi$, such that:
\begin{enumerate}
\item The restriction of $\iota$ to each component of $\partial\Sigma$ is a covering map onto a closed curve in $\Gamma$.

\item The restriction of $\iota$ to the interior $\mathring{\Sigma}$ is an embedding inside $M\backslash\Gamma$, and the submanifold $\iota(\mathring{\Sigma})$ is transverse to the foliation by trajectories of the flow.

\item There exists a real number $T>0$ such that $[x,\phi_T(x)]\cap\iota\left(\Sigma\right)\neq\emptyset$, $\forall$ $x\in M$.

\end{enumerate}  
\end{defn}

We indistinctly use $\mathring{\Sigma}$ to denote either $\Sigma\backslash\partial\Sigma$ or its inclusion $\iota(\mathring{\Sigma})$ inside the open manifold $M_\Gamma\coloneqq M\backslash\Gamma$. Conditions $2.$ and $3.$ above imply that there is a well-defined \emph{first return map} $P:\mathring{\Sigma}\to\mathring{\Sigma}$, which is a homeomorphism of the form $x\mapsto P(x)=\phi(\tau(x),x)$, for some positive and bounded continuous function $\tau:\mathring{\Sigma}\to\R^+$. In particular, the surface $\mathring{\Sigma}$ is a global transverse section for the flow $\phi$ restricted to $M_\Gamma$. Hence, $(\phi,M_\Gamma)$ is topologically equivalent to the suspension flow generated by $P:\mathring{\Sigma}\to\mathring{\Sigma}$. 

The main connection between transitive expansive flows and Birkhoff sections is given by the following theorem:
\begin{thm}[Fried \cite{fried_surgery}, Brunella \cite{brunella_thesis}]\label{thm_Fried-Brunella}
Every transitive and orbitally expansive non-singular continuous flow on a compact $3$-manifold admits a Birkhoff section. Moreover, the first return map is obtained from a pseudo-Anosov homeomorphism on a closed surface, by removing a finite set of periodic orbits.  
\end{thm}

It turns out that the orbital equivalence class of a transitive topologically Anosov flow is completely determined by the action of the first return map on a given Birkhoff section, and the homotopy class of the embedding of the surface in the $3$-manifold in a neighborhood of the boundary. This is the content of Theorem \ref{thmA-statement} below, and all this section is devoted to prove it. Before, we state the main definitions related to Birkhoff sections.

In the following, we set $(\phi,M)$ to be a transitive topologically Anosov flow. 
 
\subsubsection{Homological coordinates near the boundary.} 
Denote by $\gamma_1,\dots,\gamma_k$ the periodic orbits in $\Gamma$. Since the flow is topologically Anosov, the germ of the flow in a neighborhood of any curve $\gamma_i$ is equivalent to the suspension of a linear transformation of the plane, given by a diagonal matrix with eigenvalues $\lambda,\mu$ satisfying $0<|\lambda|<1<|\mu|$. The local invariant manifolds $\Wsloc(\gamma_i)$ and $\Wuloc(\gamma_i)$ are then homeomorphic to \emph{cylinders} or \emph{Möbius bands}, according to the signature of the eigenvalues. We say that $\gamma_i$ is a \emph{saddle type} periodic orbit. 

\begin{rem}\label{remark_existence-orientable_loc_inv_mflds}
As it is shown along the proof of theorem 2.1 in \cite{brunella_thesis}, it is always possible to find Birkhoff sections satisfying the following additional hypothesis:
\begin{itemize}
\item[4.] \emph{Each $\gamma_i\in\Gamma$ is a saddle type periodic orbit whose local invariant manifolds are orientable.} 
\end{itemize}   
In this statement it is implicit that a topologically Anosov flow always has periodic orbits with orientable invariant manifolds, and the statement holds for general transitive expansive flows as well. 
\end{rem}

From now on, we will always work with Birkhoff sections under the previous condition. The preimage by $\iota$ of each $\gamma_i$ may consists in many boundary components of $\Sigma$. Let us denote the components of $\iota^{-1}(\gamma_i)$ by $C^i_1,\dots,C^i_{p_i}$, in order to have
\begin{equation}
\partial\Sigma= C_1^1 \cup\cdots\cup C_{p_1}^1 \cup\cdots\cup C_1^n \cup\cdots\cup C_{p_n}^n.
\end{equation}
\begin{defn}
We say that $p_i=p_i(\gamma_i,\Sigma)$ is the \emph{number of connected components at $\gamma_i$}, $i=1,\dots,k$.
\end{defn}

For each $\gamma_i$ consider a (small) tubular neighborhood $W_i$ of $\gamma_i$. Since we assume that $M$ is an oriented manifold, each $W_i$ is an oriented open set. We regard $\gamma_i$ as an oriented closed curve, the orientation being provided by the forward-time action of the flow. The local invariant manifolds give a frame along the oriented closed curve $\gamma_i$, that may be used to define a \emph{meridian/longitude} basis of the homology of the punctured neighborhood $W_i\backslash\gamma_i$. More precisely,
\begin{enumerate}
\item \textbf{Longitude:} Let $\beta$ be an oriented simple closed curve contained in $\Wuloc(\gamma_i)\cap(W_i\backslash\gamma_i)$, that is orientation preserving isotopic to $\gamma_i$ inside the cylinder $\Wuloc(\gamma_i)$. We call \textit{longitude} to the homology class 
$$b=[\beta]\in H_1(W_i\backslash\gamma_i).$$

\item \textbf{Meridian:} Let $\alpha\subset W_i\backslash\gamma_i$ be a simple closed curve that is the boundary of an embedded closed disk $D\subset W_i$, which is transverse to $\gamma_i$ and intersects it in exactly one point. Here, the orientation in $D$ is induced by the co-orientation associated with $\gamma_i$, and $\alpha=\partial D$ is endowed with the boundary orientation. We call \textit{meridian} to the homology class
$$a=[\alpha]\in H_1(W_i\backslash\gamma_i).$$  
\end{enumerate}

\begin{figure}[t]
\begin{center}
\includegraphics[width=0.8\textwidth]{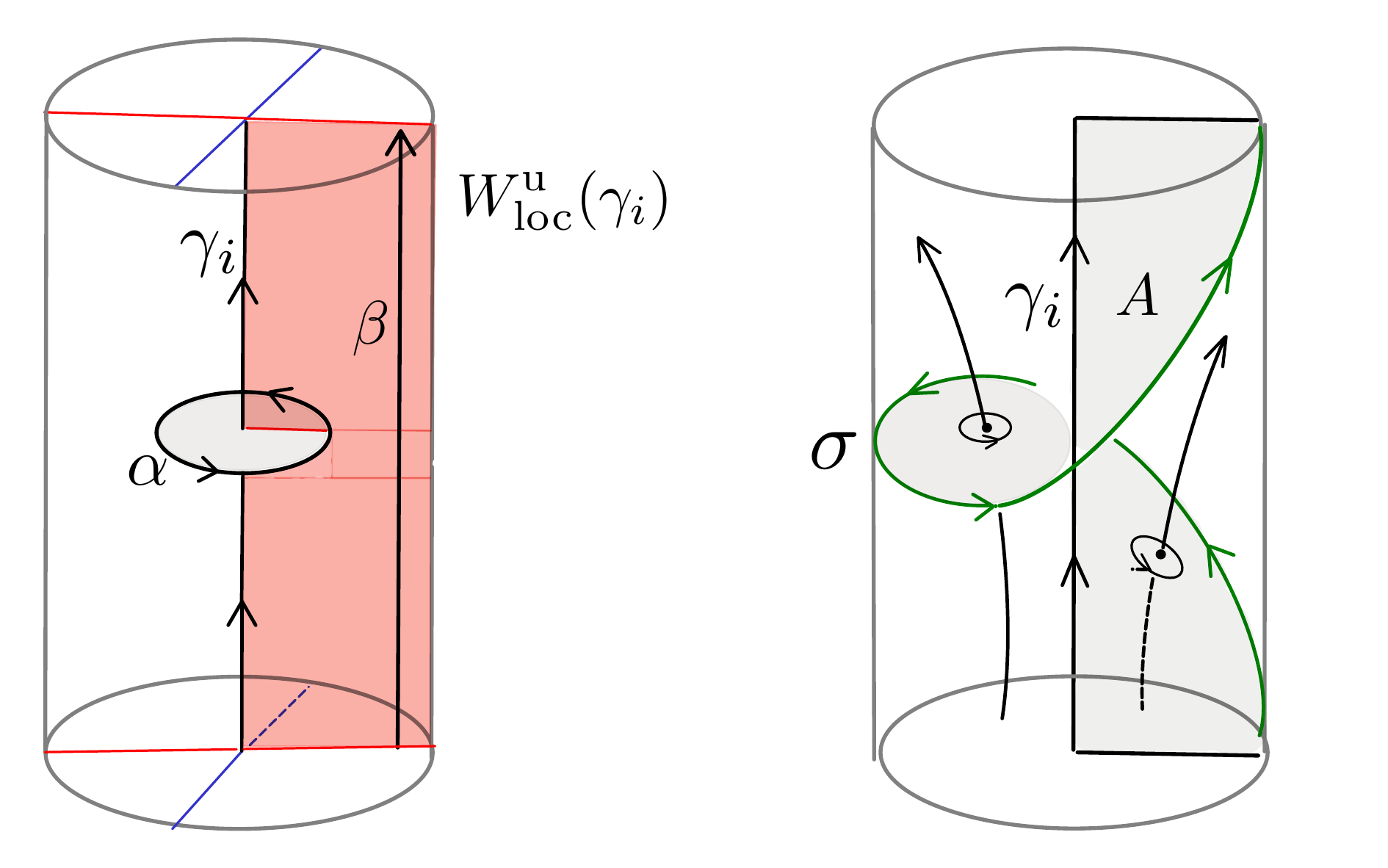}
\caption{Meridian/longitude basis of a punctured tubular neighborhood $W_i\backslash\gamma_i$.}
\end{center}
\end{figure}

If the neighborhood $W_i$ is small enough, then $(W_i\backslash\gamma_i)\cap\Sigma$ splits as $p_i$ different annuli $B_1^i,\dots,B_{p_i}^i$, properly embedded in the punctured solid torus $W_i\backslash\gamma_i$ and pairwise isotopic. Let $\sigma$ be a simple closed curve in 
$B_j^i\backslash\gamma_i$ that generates the fundamental group, oriented in the following way: The coordinates of $[\sigma]$ in the meridian-longitude basis $\{a,b\}$ of $H_1(W_i\backslash\gamma_i)$ are two integers $n=n(\gamma_i,\Sigma)$ and $m=m(\gamma_i,\Sigma)$ satisfying that 
$$[\sigma]=n(\gamma_i,\Sigma)\cdot a+m(\gamma_i,\Sigma)\cdot b.$$
We choose an orientation of $\sigma$ such that $n$ is non-negative. These integers are independent of the particular annulus $B_j^i$ that we have chosen for $j=1,\dots,p_i$.
\begin{defn}\label{defn_linking_number_multiplicity}
The integers $n(\gamma_i,\Sigma)$ and $m(\gamma_i,\sigma)$ are called the \emph{linking number} and the \emph{multiplicity} of $\Sigma$ at $\gamma_i$, respectively. 
\end{defn}

\begin{proposition}\label{prop_linking-multiplicity} With the conventions stated above, it is satisfied that:
\begin{itemize}
\item $n(\gamma_i,\Sigma)\geq 1$;
\item $m(\gamma_i,\Sigma)\neq 0$ and $m(\gamma_i,\Sigma)=\pm 1$ if and only if $\Sigma$ is embedded in a neighborhood of $\gamma_i$;
\item $\gcd(n(\gamma_i,\Sigma),m(\gamma_i,\Sigma))=1$.
\end{itemize}
\end{proposition}

\begin{proof}
We postpone the verification of the first item for Section \ref{subsubsection_projection_along_flow_lines}, see the proof of Proposition \ref{prop_quadrants_of_Birkhoff_section} there. For the second item observe that, from the one side, the map $H_1(W_i\backslash\gamma_i)\to H_1(W_i)$ induced by inclusion sends $b\mapsto[\gamma_i]$ and $a\mapsto 0$, so $[\sigma]\mapsto m\cdot[\gamma_i]$. From the other side, the curve $\sigma$ is isotopic to a boundary component $C_j^i$ of the Birkhoff section where $\iota:C_j^i\to\gamma_i$ is a covering map. Therefore $m$ is the degree of the covering $\iota:C_j^i\to\gamma_i$, so it must be non-zero and $m=1$ if and only if the section is embedded at $\gamma_i$. For the last item, observe that $\sigma$ is isotopic to $B\cap\partial W_i$. Up to shrinking the tubular neighborhood $W_i$ if necessary, the last one is a simple closed curve in $\partial W_i$. Since it is simple (no self-intersections) then its coordinates $n\cdot a+m\cdot b$ satisfy $\gcd(n,m)=1$. 
\end{proof}

\subsubsection{Blow-down operation.}
Associated with the first return map $P:\mathring{\Sigma}\to\mathring{\Sigma}$ there is a construction called \emph{blow-down} that consists in the following: Let $\widehat{\Sigma}$ be the surface obtained by collapsing each boundary component of $\Sigma$ into a point, and denote by $x_j^i$ the point obtained when collapsing $C^i_j$. If $W_i$ is a small tubular neighborhood of $\gamma_i$ then the first return map induces a cyclic permutation of the components of $\Sigma\cap(W_i\backslash\gamma_i)$ and hence of the closed curves $C_1^i,\dots,C_{p_i}^i$. Therefore, there is an associated homeomorphism $\widehat{P}:\widehat{\Sigma}\to\widehat{\Sigma}$, and each set $\{x_1^i,\dots,x_{p_i}^i\}$ constitutes a periodic orbit of $\widehat{P}$ of period $p_i=p(\gamma_i,\Sigma)$, $i=1,\dots,k$. Denote by 
$\Delta=\{x_j^i\ :\  i=1,\dots,k,\ j=1,\dots,p_i\}.$ We say that $\widehat{P}$ is a \emph{marked homeomorphism} on the surface $\widehat{\Sigma}$, and the finite set $\Delta$ is called the set of \emph{marked periodic points}.
\begin{defn}
The marked homeomorphism $(\widehat{P},\widehat{\Sigma},\Delta)$ is called the \emph{blow-down} associated with the Birkhoff section $\iota:(\Sigma,\partial\Sigma)\to(M,\Gamma)$. 
\end{defn}
This homeomorphism is expansive and thus it is conjugated with a pseudo-Anosov homeomorphism on $\widehat{\Sigma}$. The intersection of the invariant foliations of $\phi$ with $\mathring{\Sigma}$ determine the stable and unstable foliations of $\widehat{P}$, which have no singularities on $\widehat{\Sigma}\backslash\Delta$. Under the hypothesis that every $\gamma_i\in\Gamma$ is saddle type with orientable local invariant manifolds, the following condition holds (cf. \cite{fried_surgery}): 

\subsubsection{2-cycle condition.}
\emph{On each $x_j^i$ the stable (resp. unstable) foliation has $k_i=2\cdot n(\gamma_i,\Sigma)\geq 2$ prongs and $\widehat{P}^{p_i}$ permutes the set of prongs of $W^\mathrm{s}_\varepsilon(x_j^i)$ (resp. unstable) generating a partition in two subsets.} 

The \emph{multi-saddle} periodic points of $\widehat{P}$ with at least $k\geq 3$ prongs (i.e. the \emph{singularities} of the invariant foliations) are contained in $\Delta$, although points in $\Delta$ may be regular periodic points. 

\subsubsection{Orbital equivalence class and Birkhoff sections.}
Consider two transitive topologically Anosov flows $(\phi^i,M_i)$, $i=1,2$, each one equipped with a Birkhoff section $\iota_i:(\Sigma_i,\partial\Sigma_i)\to(M_i,\Gamma_i)$ and a corresponding first return map $P_i:\mathring{\Sigma}_i\to\mathring{\Sigma}_i$. Assume that the first return maps are conjugated via a homeomorphism $h:\mathring{\Sigma}_1\to\mathring{\Sigma}_2$. This condition automatically implies that the corresponding blow-down $(\widehat{P}_i,\widehat{\Sigma}_i,\Delta_i)$, $i=1,2$ are conjugated and thus, $h$ induces a correspondence $\hat{h}:\Gamma_1\to\Gamma_2$ such that
\begin{align*}
& p(\gamma,\Sigma_1)=p\left(\hat{h}(\gamma),\Sigma_2\right),\ \forall\ \gamma\in\Gamma_1,\\
& n(\gamma,\Sigma_1)=n\left(\hat{h}(\gamma),\Sigma_2\right),\ \forall\ \gamma\in\Gamma_1.
\end{align*}
Assume in addition that
\begin{align*}
& m(\gamma,\Sigma_1)=m\left(\hat{h}(\gamma),\Sigma_2\right),\ \forall\ \gamma\in\Gamma_1.\\
\end{align*}

We want to address the question of whether or not these conditions are sufficient to conclude orbital equivalence between $(\phi^1,M_1)$ and $(\phi^2,M_2)$. Observe that, by \emph{suspending} the homeomorphism $h$, we obtain an orbital equivalence $H_0:(\phi^1,M_1\backslash\Gamma_1)\to(\phi^2,M_2\backslash\Gamma_2)$ in the complement of finite sets of periodic orbits, what is also called an \emph{almost orbital equivalence}. In addition, the previous assumptions on the combinatorial parameters associated with the embedding on the boundary imply that $M_1$ is homeomorphic to $M_2$. However, we remark the following:

\begin{rem}\label{remark_non-extension_of_almost-equiv}\ 
The homeomorphism $H_0:M_1\backslash\Gamma_1\to M_2\backslash\Gamma_2$ rarely extends onto a global homeomorphism $M_1\to M_2$. If this happens, observe that then the map $h:\mathring{\Sigma}_1\to\mathring{\Sigma}_2$ conjugating the two (pseudo-Anosov) first return maps must extend to the boundary of the surfaces. The latter rarely happens, as we explained in Section \ref{section_pseudo-Anosov_non-closed_surfaces}.  
\end{rem}

Despite the previous remark, under the assumption that the flows are topologically Anosov or expansive, the conditions stated above are enough to guarantee orbital equivalence, according to Theorem \ref{thmA-statement} below. We state this property in the form that is used along this work.

\begin{thm}[\textbf{Theorem B}]\label{thmA-statement}
Consider two topologically Anosov flows $(\phi^i,M_i),\ i=1,2$ equipped with homeomorphic Birkhoff sections $\iota_i:(\Sigma_i,\partial\Sigma_i)\to(M_i,\Gamma_i)$ satisfying that each curve $\gamma\in\Gamma_i$ has orientable local invariant manifolds. Assume that there exists a homeomorphism $\Psi:\mathring{\Sigma}_1\to\mathring{\Sigma}_2$ such that:
\begin{enumerate}
\item The induced homomorphism $[\Psi]:\pi_1(\mathring{\Sigma}_1)\to\pi_1(\mathring{\Sigma}_2)$ conjugates the actions on fundamental groups $[P_i]:\pi_1(\mathring{\Sigma}_i)\to\pi_1(\mathring{\Sigma}_i)$.
\item For each $\gamma\in\Gamma_1$ it is verified that $m(\gamma,\Sigma_1)=m(\Psi(\gamma),\Sigma_2)$.
\end{enumerate}
Then, there exists an orbital equivalence $H:(\phi^1,M_1)\to(\phi^2,M_2)$, which sends each curve $\gamma\in\Gamma_1$ homeomorphically onto $\Psi(\gamma)$.
\end{thm} 

First of all, we remark that Theorem \ref{thmA-statement} is equally valid in the more general case of pseudo-Anosov flows. We restrict ourselves to the case of topologically Anosov (i.e. non-singular invariant foliations), satisfying the additional condition of orientability of invariant manifolds on the boundary of the sections, which simplifies much of the combinatorial analysis that must be done. We assume as well, for the sake of simplicity, that the closed $3$-manifold $M$ is always orientable. 

Using the fact that \emph{conjugacy classes of pseudo-Anosov homeomorphisms on closed surfaces} are determined by the action on fundamental group (cf. Theorem \ref{isotopic_pseudo_Anosov} on Section \ref{subsection_psudo-Anosov-homeos}) we can directly obtain the existence of a homeomorphism $h:\mathring{\Sigma}_1\to\mathring{\Sigma}_2$ satisfying $P_2\circ h=h\circ P_1$. By suspending $h$, we obtain an orbital equivalence $H_0:(\phi^1,M_1\backslash\Gamma_1)\to(\phi^2,M_2\backslash\Gamma_2)$. The difficulty resides in the fact that many of the homeomorphisms obtained in this way do not extend continuously onto the boundary of the sections as we pointed in Remark \ref{remark_non-extension_of_almost-equiv} above.

Nevertheless, we show that if the Birkhoff sections are \emph{well-positioned}, then $H_0$ can be modified in an arbitrarily small neighborhood of $\Gamma_1$ by pushing along the flow-lines, in order to obtain a homeomorphism $H:M_1\backslash\Gamma_1\to M_2\backslash\Gamma_2$ that extends onto the whole $M_1\to M_2$. This notion of well-positionedness is called called \emph{tameness} and was taken from \cite{bonatti-guelman}. Under this assumption, the germ of the flow in a neighborhood of each $\gamma\in\Gamma_i$ can be described using the information associated with the section (first return and homological coordinates), and Theorem \ref{thmA-statement} can be reduced to a local problem near the sets $\Gamma_i$. For going from tame sections to arbitrary ones, we make use of Proposition \ref{pi_1-conjugated_pseudo_Anosov_punctures}.

The rest of the section is organized as follows: In Section \ref{subsection-local_Birkhoff_sectiions} we summarize all the information from \cite{bonatti-guelman} needed for proving Theorem \ref{thmA-statement}. In Section \ref{subsection-Local_version_thm-A} we state and prove a local version of \ref{thmA-statement}. We give the proof of \ref{thmA-statement} in Section \ref{subsection-Proof_thm-A}. In a first reading, it is possible to skip Sections \ref{subsection-local_Birkhoff_sectiions} and \ref{subsection-Local_version_thm-A} and go directly to the proof of \ref{thmA-statement} in Section \ref{subsection-Proof_thm-A}.

\subsection{Local Birkhoff sections.}\label{subsection-local_Birkhoff_sectiions}
Let $\{\phi_t:M\to M\}_{t\in\R}$ be a regular flow in an oriented $3$-manifold, having a saddle type periodic orbit $\gamma\subset M$ whose local invariant manifolds $\Wsloc(\gamma)$ and $\Wuloc(\gamma)$ are cylinders.  

\begin{defn}\label{defn_local_Birkhoff_section}
A \textit{local Birkhoff section at $\gamma$} is an immersion $\iota:[0,1)\times\R/\Z\to M$ such that:
\begin{enumerate}
\item $\gamma=\iota\left(\{0\}\times\R/\Z\right)$;
\item The restriction of the map $\iota$ to $(0,1)\times\R/\Z$ is an embedding and the submanifold $\iota\left((0,1)\times\R/\Z\right)$ is transverse to the flow lines;
\item The exists a real number $T>0$ and a neighborhood $W$ of $\gamma$ such that $[x,\phi_T(x)]\cap\iota\left([0,1)\times\R/\Z\right)\neq\emptyset$, $\forall$ $x\in W$.
\end{enumerate}
\end{defn}

Let $B\subset M$ be the image of $\iota:[0,1)\times\R/\Z\to M$. When confusion is not possible, we just use $B$ for denoting the local Birkhoff section. We denote by $\B$ the set $B\backslash\gamma$. Given a local Birkhoff section $B$ at $\gamma$, the last property in the previous definition implies that there exists a collar neighborhood $U\subset B$ of $\gamma$ and a first return map $P_B:\U\to\B$, where $\U=U\backslash\gamma$. In general we will be concerned with the germ of the first return map to a local Birkhoff section, that is, we consider the map $P_B$ up to changing $U$ for a smaller collar neighborhood if needed.

Let $W\subset M$ be a regular tubular neighborhood of $\gamma$ satisfying condition 3 in the definition above. We are interested in describing the germ $(\phi,W)_\gamma$ of the foliation by flow orbits restricted to the neighborhood $W$ (or a smaller one), using the information associated with the Birkhoff section. Since the periodic orbit $\gamma$ is saddle type, there is an associated meridian-longitude basis and thus the homotopy class of the embedding $\B\hookrightarrow W\backslash\gamma$ has an associated linking number and multiplicity (cf. definition  \ref{defn_linking_number_multiplicity}) satisfying Proposition \ref{prop_linking-multiplicity}.

\subsubsection{Tameness.}
Along this part we make the following assumption about the embedding $B\hookrightarrow W$. 

\begin{defn}[\textbf{Tameness}]
The local Birkhoff section $B$ is \emph{tame} if there exists a collar neighborhood $U\subset B$ of $\gamma$ such that the sets $U\cap\Ws_{loc}(\gamma)$ and $U\cap\Wu_{loc}(\gamma)$ consists of the union of $\gamma$ with finitely many compact segments, each of them intersecting $\gamma$ exactly at one of its extremities. (cf. Figure \ref{fig_partition into_quadrants}.) 
\end{defn}

\subsubsection{Partition into quadrants.}
\begin{figure}[t]
\begin{center}
\includegraphics[height=0.45\textheight, keepaspectratio]{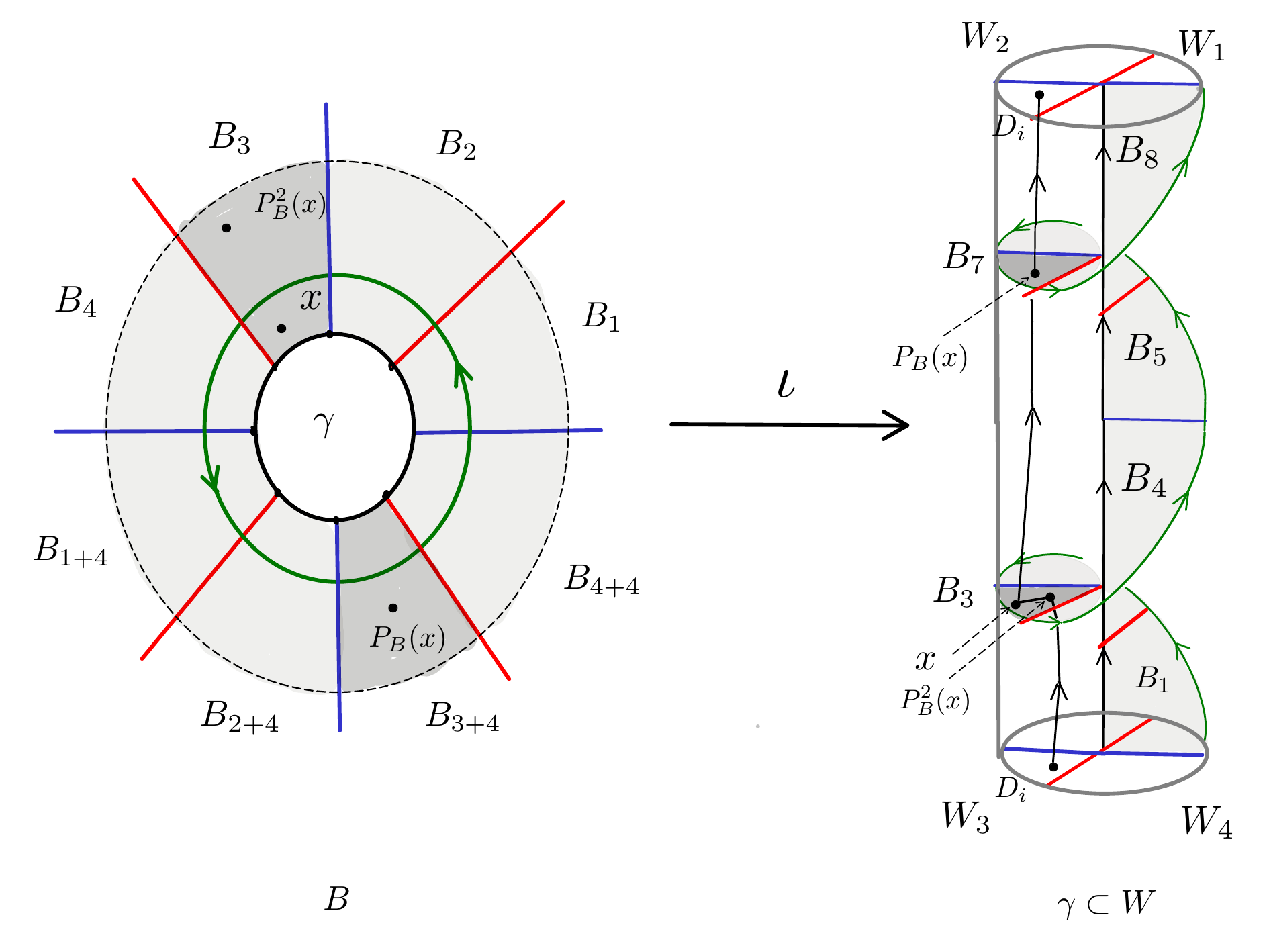}
\caption{Partition into quadrants of a tame Birkhoff section with $n(\gamma,B)=2$ and $m(\gamma,B)=1$.}
\label{fig_partition into_quadrants}
\end{center}
\end{figure}

Let us assume that $W$ is small enough such that the union of the local stable and unstable manifolds separate $W$ in four quadrants. Observe that the orientation of the meridian $a\in H_1(W\backslash\gamma)$ induces a cyclic order on this set of connected components. The four quadrants are denoted by $W_i$, $i=1,\dots,4$, where the indices are chosen to respect the cyclic order of the quadrants, as in Figure \ref{fig_partition into_quadrants}. 

Let $D\subset W$ be a local transverse section and $x_D$ its intersection point with $\gamma$. Then, the quadrants $W_1,\dots,W_4$ determine four quadrants $D_i=D\cap W_i$ in $D$. Observe that the first return map $P_D:V\to D$ defined in a neighborhood $V$ of $x_D$ preserves the quadrants, i.e. $P_D(V\cap D_i)\subset D_i$ for every $i=1,\dots,4$.

Let $B\subset W$ be a tame local Birkhoff section at $\gamma$ with linking number $n$ and multiplicity $m$. Then, because of the tameness hypothesis the four quadrants $W_i$ also determine a partition of the annulus $B$ into quadrants. Each of these quadrants can be though of as a rectangle, whose boundary contains a segment of $\gamma$ and two segments which are connected components of $\B\cap\left(\Wsloc(\gamma)\cup\Wuloc(\gamma)\right)$. Observe that the first return map $P_B:\U\to\B$ defined in a collar neighborhood $U$ of $\gamma$ sends quadrants of $B$ into quadrant of $B$.

\begin{proposition}\label{prop_quadrants_of_Birkhoff_section}
There are exactly $4n$ quadrants of $B$, and each $W_i$ contains $n$ of them.
\end{proposition}

If we choose a quadrant of $B$ which lies in $W_1$ and we call it $B_1$, then we can inductively label the quadrants of $B$ as $B_1,\dots,B_{4n}$ by declaring that $\forall j=1,\dots,4n$, if $W_i$ contains $B_j$ then $B_{j+1}$ is the quadrant adjacent to $B_j$ which lies in $W_{i+1}$, $i=1,\dots,4$. We always use this labeling for the quadrants of a Birkhoff section. 

We postpone the proof of Proposition \ref{prop_quadrants_of_Birkhoff_section} for the next subsection. 

\subsubsection{Projections along flow lines.}\label{subsubsection_projection_along_flow_lines}
The foliation $\mathcal{O}$ by $\phi$-orbits of $M$ induces a foliation $\mathcal{O}_W$ by orbit segments on the tubular neighborhood $W$, where for every $x\in W$, $\mathcal{O}_W(x)$ is the connected component of $\mathcal{O}(x)\cap W$ that contains $x$. Each segment $\mathcal{O}_W(x)$ is parametrized by the action of $\phi$, that is, 
$$\mathcal{O}_W(x)=\{\phi_t(x):-\infty\leq a_x<t<b_x\leq+\infty\}$$ for some $a_x,b_x$. The periodic orbit $\gamma$ is the unique $\phi$-orbit that is entirely contained in $W$.

The (germ of the) $\phi$-foliation restricted on the tubular neighborhood $W$ can be obtained by suspending the (germ of the) first return map $P_D:V\to D$ onto the transverse disk $D$. On the punctured neighborhood $W_\Gamma\coloneqq W\backslash\gamma$, the foliation by flow orbits can be written as both:
\begin{itemize}
\item The suspension of $P_D$ on the punctured disk $D^*=D\backslash\{x_D\}$,
\item The suspension of $P_B$ on the interior of the Birkhoff section $\B=B\backslash\gamma$.
\end{itemize}
Our aim is to relate these two representations of the foliation induced by $\phi$-orbits on $W_\Gamma$. This can be done by projecting points in $\B$, along the flow orbits, onto points in $D^*$. Although the surfaces $D^*$ and $\B$ are not homotopic in $W_\Gamma$, these maps can be well-defined over simply connected regions of $\B$ (in particular, over the quadrants of the Birkhoff section).

\begin{figure}[t]
\begin{center}
\includegraphics[height=0.45\textheight, keepaspectratio]{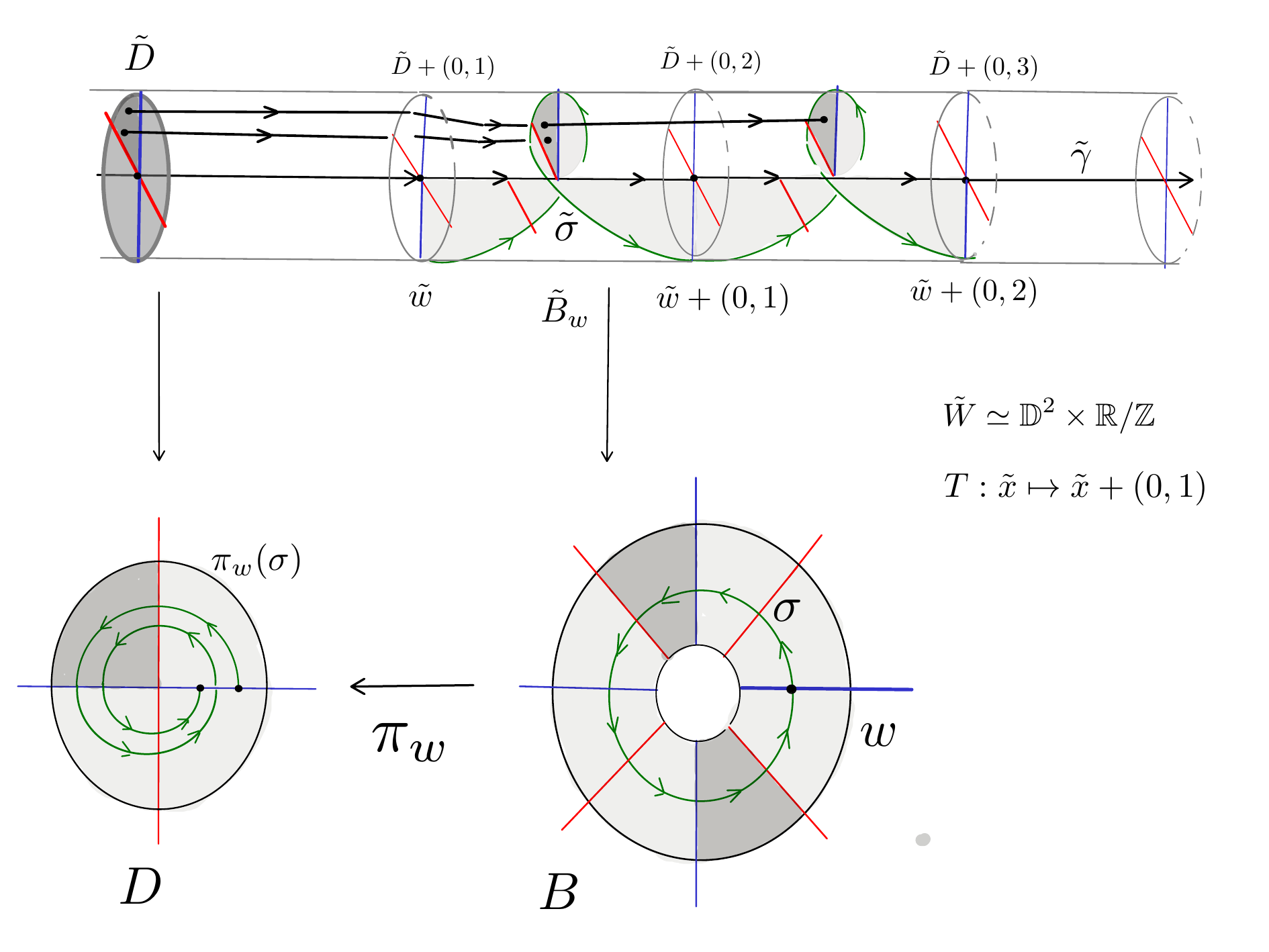}
\caption{Projection along the flow lines from a local Birkhoff section onto a local transverse section.}
\end{center}
\end{figure}

\paragraph*{Construction.}
Let $w$ be a segment which is the closure of a connected component of $\B\cap\left(\Wsloc(\gamma)\cup\Wuloc(\gamma)\right).$ Since the Birkhoff section is an immersed annulus, we can cut along $w$ and obtain an immersed compact strip inside $W$ (embedded on the interior), with two opposite sides that are naturally identified with $w$. We denote this strip by $B_w$. 

Consider a universal cover $\tilde{W}$ of $W$, which is homeomorphic to $\D^2\times\R$ since $W\simeq\D^2\times\R/\Z$. By lifting the foliation $\mathcal{O}_W$ on $W$, we obtain a foliation $\tilde{\mathcal{O}}_{W}$ on $\tilde{W}$ by parametrized segments $\tilde{\phi}_{t}(\tilde{x})$ with $\tilde{x}\in\tilde{W}$ and 
$-\infty\leq a_{\tilde{x}}<t<b_{\tilde{x}}\leq +\infty$. Let $\tilde{B}_w$, $\tilde{D}$ and $\tilde{\gamma}$ be lifts to the universal of $B_w$, $D$ and $\gamma$ respectively. By the continuity of the flow there exists some neighborhood $O$ of the compact segment $\tilde{\gamma}\cap\tilde{B}_w$ and some $T>0$ such that, for every $\tilde{x}\in O$, the $\tilde{\phi}$-segment $[\tilde{\phi}_{-T}(\tilde{x}),\tilde{\phi}_{T}(\tilde{x})]$ intersects the disk $\tilde{D}$ and exactly in one point. This allows to consider a collar neighborhood $\tilde{U}_w$ of $\tilde{\gamma}\cap\tilde{B}_w$ inside $O\cap\tilde{B}_w$ and to define a map $\tilde{\pi}_w:\tilde{U}_w\to\tilde{D}$ of the form $\tilde{\pi}_w(x)=\tilde{\phi}(s(x),x)$, where $s:\tilde{U}_w\to\R$ is continuous, and is bounded as a consequence of the tame hypothesis (cf. \cite{bonatti-guelman}). Observe that all the points in $\tilde{\gamma}\cap\tilde{U}_w$ are mapped over the intersection point $\tilde{x}_D$ of $\tilde{D}$ with $\tilde{\gamma}$, but we are not interested in these points. So, we just consider the restriction $\tilde{\pi}_w:(\tilde{U}_w\backslash\tilde{\gamma})\to(\tilde{D}\backslash\{\tilde{x}_D\})$. 

Since the universal covering map provides identifications $\tilde{D}\to D$ and $\tilde{B}_w\to B_w$, we can think of $\tilde{\pi}_w$ as a map $\pi_w:\B_w\to D^*$ of the form $\pi_w(x)=\phi(s(x),x)$, defined for points $x$ in the strip $B_w$ which are sufficiently close to $\gamma$. 
\begin{defn}\label{defn_projection_along_flow_lines}
Let $B$ be a tame local Birkhoff section at $\gamma$, $D$ a local transverse section which intersects $\gamma$ at the point $x_D$, and let $w$ be a connected component of $\B\cap\left(\Wsloc(\gamma)\cup\Wuloc(\gamma)\right)$. A \textit{local projection along the flow from $B$ onto $D$} is a map $\pi_w:\U_w\to D^*$ of the form $\pi_w(x)=\phi(s(x),x)$, where 
\begin{itemize}
\item $U\subset B$ is a collar neighborhood of $\gamma$ and $\U=U\backslash\gamma$,
\item $\U_w$ is the strip obtained by cutting $\U$ along $w$,
\item $s:\U_w\to\R$ is continuous and bounded.  
\end{itemize}
\end{defn}

If we choose a quadrant $B_j$ of $B$ that lies in $W_i$ and we consider the restriction $\pi_j=\pi_w|_{B_j}$, this map is a homeomorphisms between $\U\cap B_j$ and a open set $V_i$ in $D_i$ that accumulates in $x_D$. The main interest about this map is that it conjugates the first return maps on the quadrant surfaces $B_j$ and $D_i$. Observe that a projection along the flow depends on the particular choice of the segment $w$, as well as the particular choices of the the lifts of $D$ and $B_w$. 

The following proposition summarizes different properties that we will need later. 

\begin{proposition}\label{prop_combinatorial_relation_projections_first_return}
Let $B$ be a tame local Birkhoff section at $\gamma$ with linking number $n\geq 1$ and multiplicity $0\neq m\in\Z$. Let $\pi=\pi_w:\U_w\to D$ be a projection along the flow onto a local transverse section $D$ as defined above. Here, $w$ is a segment of the intersection of $\B$ with $\Wsloc(\gamma)\cup\Wuloc(\gamma)$.  We enumerate the quadrants of $B$ as $B_1,\dots,B_{4n}$ in such a way that $B_1$ and $B_{4n}$ intersect along $w$. Observe that for all the quadrants $B_j$ contained in $W_i$ it is satisfied that $j=i+4r$, where $0\leq r\leq n-1$ and $1\leq i\leq 4$. 
Then, we have:

\begin{enumerate}
\item The map $P_B$ permutes cyclically all the quadrants of $B$ that are contained in the same component $W_i$, and in particular we have that $P_B^{n}$ preserves each quadrant $B_j$. Moreover, in the case that $n>1$, let $1\leq k,l\leq n-1$ be such that $k\equiv m$ ($\bmod n$) and $l\equiv m^{-1}$ ($\bmod n$). Then, the first return map to $B$ permutes the quadrants in the following fashion: 

\begin{enumerate}
\item 
$P_B$ takes the quadrant $B_j$ into 
$\left\{
\begin{array}{ll}
B_{j+4l} & \text{if }  m>0, \\
B_{j-4l} & \text{if }  m<0. 
\end{array}
\right.$

\item 
$P_B^{k}$\text{ takes the quadrant }$B_j$\text{ into }$\left\{
\begin{array}{ll}
B_{j+4} & \text{if }  m>0, \\
B_{j-4} & \text{if }  m<0. 
\end{array}
\right.$
\end{enumerate}

\item Let $B_j$ be a quadrant of $B$ where $j=i+4r$, $1\leq i\leq 4$ and $0\leq r\leq n-1$. Let us denote by $\U_j=\U\cap B_j$ and $\pi_j$ the restriction $\pi|_{\U_j}$. Then, the $\pi_j$ takes $\U_j$ homeomorphically onto its image in $D_i-\{x_D\}$ and it is satisfied that:

\begin{enumerate}

\item The homeomorphism $\pi_j$ induces a local conjugacy between $P_B^{n}$ and $P_D$. That is, 
\begin{equation*}
\pi_j\circ P_B^{n}(z)=P_D\circ\pi_j(z),\ \forall j=1,\dots,4n,\ \text{and}\ z\ \text{sufficiently close to}\ \gamma.
\end{equation*}

\item The holonomy defect over $w$ is given by 
\begin{equation*}
\pi_1\circ\pi_{4n}^{-1}(z)=P_D^m(z),\ \forall\ z\in w=B_1\cap B_{4n}\ \text{sufficiently close to}\ \gamma. 
\end{equation*}

\end{enumerate}

\item The projection along the flow depends on $w$ and on a particular choice of a lift of $B_w$ to the universal cover. For a fix segment $w$ we have that if $\pi_w$ and $\pi_w'$ are two such projections then $\exists$ $k\in\mathbb{Z}$ such that $\pi_w'=\pi_w\circ P_B^k$, in a suitable common domain of definition. 
\end{enumerate}
\end{proposition}

The proof of this proposition (as well as Proposition \ref{prop_quadrants_of_Birkhoff_section} and the fact that the linking number is not zero in Proposition \ref{prop_linking-multiplicity}) resides in the following two facts. 

First, let $B_j$ be a quadrant of $B$ contained in some component $W_i$ of $W\backslash(\Wsloc(\gamma)\cup\Wuloc(\gamma))$, and let $D_i$ be the corresponding quadrant of $D$. Let $\U_j=\U\cap B_j$ and $\pi_j:\U_j\to D_i$ be the restriction to the quadrant $B_j$ of the projection along the flow. The map $\pi_j$ is a homeomorphism onto its image. Let us denote $Q\coloneqq B_j$ and let $P_Q:\U_j\to Q$ be the first return map \emph{onto the quadrant $Q$}. Let $z\in\U_j$ be a point such that $z'=P_Q(z)\in\U_j$, and let $\eta_0:[0,1]\to\U_j$ be an arc contained in $Q$, connecting $\sigma_0(0)=z'$ with $\sigma_0(1)=z$. Define
\begin{itemize}
\item $\beta_0=[z,z']\cdot\eta_0$ is the curve obtained by concatenating the $\phi$-orbit segment $[z,z']$ with the arc $\eta_0$.
\item $x=\pi_j(z)$, $x'=\pi_j(z')$, $\eta_1=\pi_j\circ\eta_0$. (Observe that $z,\ x,\ z',\ x' $ all lie in the same $\phi$-orbit.) 
\item $\beta_1=[x,x']\cdot\eta_1$ is the curve obtained by concatenating the $\phi$-orbit segment $[x,x']$ with the arc $\eta_1$. 
\end{itemize}
\begin{fact}[1]
If $z$ is chosen sufficiently close to $\gamma$, then $\beta_0$ is homotopic to $\beta_1$ inside $W_i$.
\end{fact}

\begin{figure}[t]
\begin{center}
\includegraphics[scale=0.35]{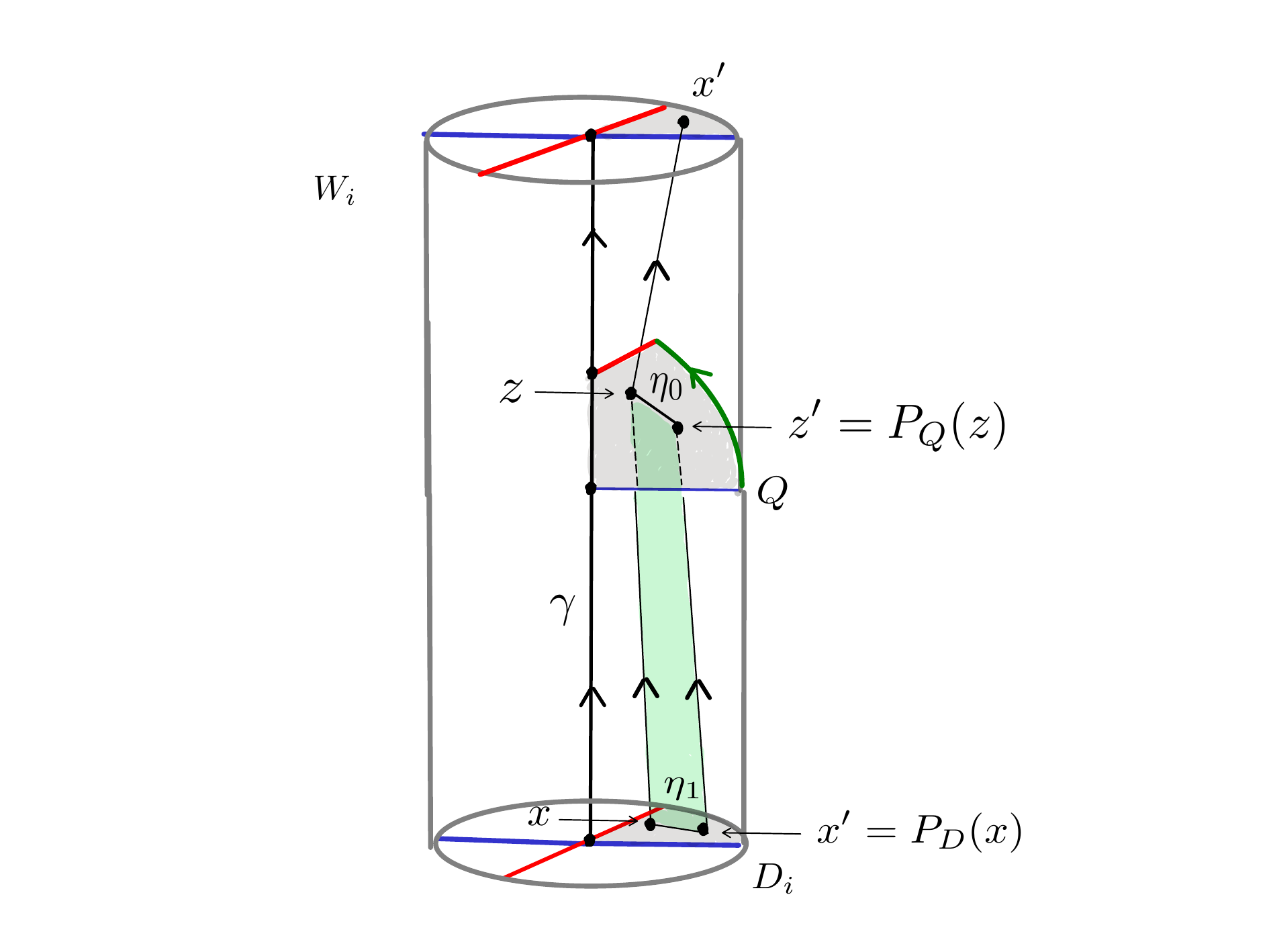}
\caption{The curves $\beta_0$ and $\beta_1$ of Fact (1).}
\label{fig_flow-conjugation}
\end{center}
\end{figure}

\begin{proof}
This follows from the fact that the formula $\pi_j^\theta:x\mapsto\phi(\theta s(x),x)$ with $0\leq\theta\leq 1$ defines a proper isotopy between $\pi_j$ and the inclusion $Q\hookrightarrow W_i$. See Figure \ref{fig_flow-conjugation}.
\end{proof}

Second, since $\B\hookrightarrow W\backslash\gamma$ is a \emph{properly embedded} surface (for a small choice of the regular tubular neighborhood $W$ of $\gamma$) there is a well-defined notion of \emph{algebraic intersection} between the homology class $[\alpha]\in H_1(W\backslash\gamma)$ of a closed curve $\alpha$ and the relative homology class $[B]\in H_2(W,\partial W\cup\gamma)$ of the annulus $\B$. This intersection can be defined by counting oriented intersections of a representative $\alpha$ that is transverse to $\B$. It is satisfied that:
\begin{fact}[2]
Let $\{a,b\}$ be the meridian/longitude basis of $H_1(W\backslash\gamma)$ defined in \ref{defn_linking_number_multiplicity} above, and let $[\alpha]=p\cdot a +q\cdot b$ be the homology class of a closed curve $\alpha$, where $p,q\in\Z$. Then
\begin{equation}\label{eq_homological_intersection}
[\alpha]\cdot[B]=-pm(\gamma,B)+qn(\gamma,B).
\end{equation}
\end{fact}

This can be checked by showing that $m=-a\cdot[B]$ and $n=b\cdot[B]$ are the multiplicity and linking number of the local Birkhoff section $B$, respectively. See \cite[Lemma 2.5]{shannon_thesis} for a proof of Fact (2). 

\begin{proof}[Proof of Proposition \ref{prop_quadrants_of_Birkhoff_section}]
Each quadrant $W_i$ is homeomorphic to a solid torus, and the map $H_1(W_i)\to H_1(W\backslash\gamma)$ induced by the corresponding inclusion sends a generator of $H_1(W_i)$ to the longitude class $b\in H_1(W\backslash\gamma)$. The homology class $[\beta_0]$ in $H_1(W_i)$ of Fact (1) above is a generator (since it has intersection number one with the properly embedded disk $Q$) and thus $[\beta_0]\cdot[B]=b\cdot[B]$. Since the intersection number $b\cdot[B]$ equals the linking number $n(B,\gamma)$ by Fact (2), and since $\beta_0$ cuts once every quadrant $B_j$ of $B$ with the same orientation, we conclude that there are $n$ quadrants of $B$ inside $W_i$. This also shows that $n\neq 0$ (Proposition \ref{prop_linking-multiplicity}).  
\end{proof}

\begin{proof}[Proof of Proposition \ref{prop_combinatorial_relation_projections_first_return}]\ 
\begin{enumerate}

\item Following the explanation of the previous proof, $P_B$ permutes cyclically the $n$ quadrants of $B$ inside $W_i$. Thus, $P_B^n$ fixes each quadrant $B_j$. In addition we have:

\begin{enumerate}
\item Consider a point $z\in\B_j$ sufficiently close to $\gamma$ such that $P_B^{n}(z)\in\B_j$. Assume first that $m>0$. Since the orbit segment $[z,P_B^{n}(z)]$ intersects once each quadrant inside $W_i$ we know that there exists some $1\leq l\leq n-1$ such that $P_B(z)\in\B_{j+4l}$ and we want to determine $l$. Let $[P_B(z),z]$ denote the curve that is obtained by reparametrizing with inverse orientation the orbit segment joining $z$ with $P_B(z)$. Consider a curve $\alpha:[0,1]\to B_j\cup B_{j+1}\cup\dots \cup B_{j+4l}$ connecting $\alpha(0)=z$ with $\alpha(1)=P_B(z)$ and let us define $\eta$ as the closed path that is obtained by concatenation of $\alpha$ with $[P_B(z),z]$.
\begin{cl}
The homology class of $\eta$ in $H_1(N\backslash\gamma)\simeq\Z^2$ is $(l,q)$ in the basis given by the meridian and the longitude, where $q$ is some integer. 
\end{cl}
Since the segment $[P_B(z),z]$ is transverse to $B$ at its endpoints and cuts it with negative orientation, and since $\alpha$ is tangent to $B$, then the homological intersection number satisfies $[B]\cdot [\eta]=-1$. Using that $[B]\cdot [\eta]=-1=-l\cdot m+q\cdot n$ (Fact (2) above) we conclude that $l\equiv m^{-1}$ (mod $n$) when $m>0$. The case $m<0$ is analogous. 

To prove the claim, observe that up to homotopy we can assume that $\alpha$ is a concatenation $\alpha=\alpha_1\cdots\alpha_{4l}$ of straight segments $\alpha_k$ contained in $B_{j+k}$ connecting the two boundaries of the quadrant. Let $S$ be a connected component of $\Wsloc(\gamma)\backslash\gamma$. Then, since the curve $\eta$ is the concatenation of $\alpha$ with a segment contained in $\interior(N_i)$, it cuts $S$ exactly $l$ times and with positive orientation, from where it follows the claim.

\item The proof of $1(b)$ is similar to the previous one. Assume that $m>0$. We now that there exists some $1\leq k\leq n-1$ such that $P_B^k(z)\in\B_{j+4}$, and we want to determine $k$. Let $[P_B^k(z),z]$ denote the curve that is obtained by reparametrizing with inverse orientation the orbit segment joining $z$ with $P_B^k(z)$. Consider a curve $\alpha:[0,1]\to B_j\cup B_{j+1}\cup B_{j+2}\cup B_{j+3}\cup B_{j+4}$ connecting $\alpha(0)=z$ with $\alpha(1)=P_B^k(z)$ and let us define $\eta$ as the closed path that is obtained as a concatenation of $\alpha$ with $[P_B^k(z),z]$. It follows that the homology class of $\eta$ in $H_1(N\backslash\gamma)$ is $(1,q)$, where $q$ is some integer. Using that $[B]\cdot [\eta]=-k=-m+q\cdot n$ we have $k\equiv m$ ($\bmod n$). The case $m<0$ is analogous.

\end{enumerate} 

\item Denote $Q=B_j$ and consider $z$,$z'\in Q$, $x,x'\in D_i$ and $\beta_0$, $\beta_1$ defined as in Fact (1) above. 

\begin{enumerate}
\item
Since $z'=P_Q(z)$ is the first return onto $Q$ then the intersection number $[\beta_0]\cdot[Q]$ equals one, and since $Q$ is properly homotopic to $D_i$ (inside $W_i$) we have also that $[\beta_1]\cdot[D_i]=1$. Since $\beta_1$ is obtained concatenating the oriented orbit segment $[x,x']$ (that always cross $D_i$ with the same orientation) and a tangent arc $\eta_1$, we conclude that $\beta_1$ has only one transverse intersection with $D_i$, so $x'=P_D(x)$. This shows that:
\begin{equation}\label{xrmptkl}
\text{For}\ z\ \text{sufficiently close to}\ \gamma\ \text{we have}\ \pi_j\circ P_Q(z)=P_D\circ\pi_j(z).
\end{equation}   
The statement of $2(a)$ follows, since $P_Q=P_B^n$ by item $1$. 

\begin{figure}[t]
\begin{center}
\includegraphics[height=0.35\textheight, keepaspectratio]{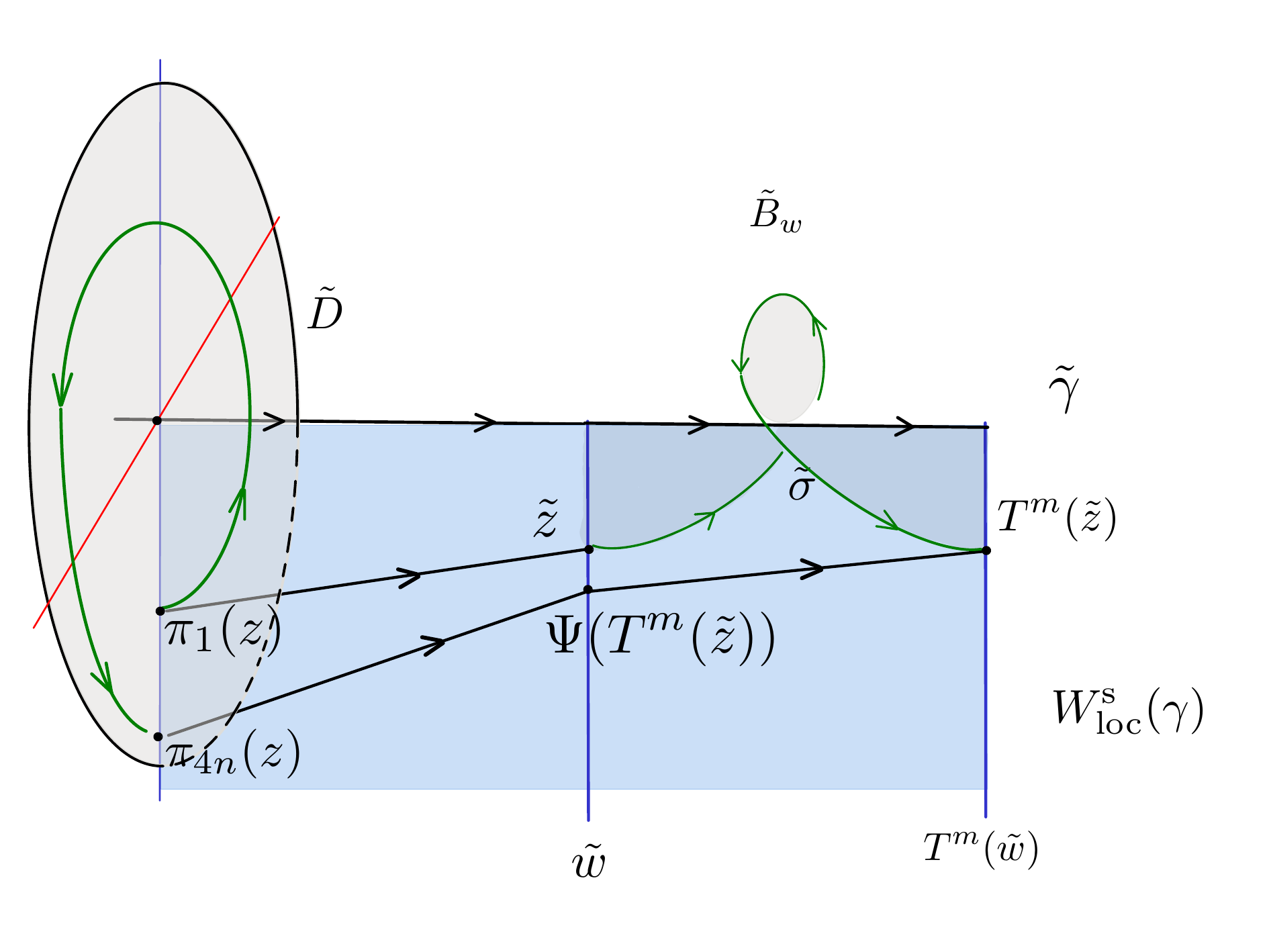}
\caption{Holonomy defect of $\pi_w$ along $w$.}
\end{center}
\end{figure}

\item
To prove $2(b)$, let $\sigma:[0,1]\to\B$ be a simple closed curve generating $H_1(B)$ such that $\sigma(0)=\sigma(1)=z\in w$ and $\sigma\cap w=\{z\}$. Let $\tilde{\sigma}:[0,1]\to\tilde{B}_w$ be a lift to the universal cover. Since $[\sigma]=n\cdot a+m\cdot b$ it follows that $\tilde{\sigma}(1)=T^m(\tilde{\sigma}(0))$, where $T:\tilde{W}\to\tilde{W}$ is the generator of the deck transformation group. Hence, the strip $\tilde{B}_w$ is delimited by two lifts of $w$, namely $\tilde{w}$ and $T^m(\tilde{w})$. Recall that we have defined first a map $\tilde{\pi}_w:\tilde{B}_w\to\tilde{D}$ (only defined for points near $\tilde {\gamma}$) and then $\pi_w:\B_w\to D^*$ by passing to the quotient $\tilde{W}\to W$. By construction we have
\begin{align}\label{eq_holonomia-borde_1}
& \pi_{4n}(z)=\tilde{\pi}_w(\tilde{\sigma}(1))\\
& \pi_1(z)=\tilde{\pi}_w(\tilde{\sigma}(0)).
\end{align}
Let $\psi:T^m(\tilde{w})\to\tilde{w}$ be the projection along the flow lines in the universal cover, so that 
$$\tilde{\pi}_w|_{T^m(\tilde{w})}=\tilde{\pi}_w|_{\tilde{w}}\circ\psi.$$ 
Recall that $w$ is a segment contained in $B\cap\left(\Wsloc(\gamma)\cup\Wuloc(\gamma)\right)$ and transverse to the flow, so there is a first return map $P_w:w\to w$ and $\pi_1\circ P_w=P_D\circ\pi_1$ as in \eqref{xrmptkl} above. Since $m>0$, the map $\psi$ from $T^m(\tilde{w})$ to $\tilde{w}$ on the universal cover induces the $-m$-th power of $P_w$ on the segment $w$. We can conclude that
$$\pi_{4n}(z)=\pi_1\circ P_w^{-m}(z)=P_D^{-m}\circ\pi_1(z).$$
\end{enumerate}

\item Fix the segment $w$ in $B\cap(\Wsloc(\gamma)\cup\Wuloc(\gamma))$ and chose flow projections $\pi_w$ and $\pi_w'$ from $B_w$ onto a transverse section $D$. Denote by $\tilde{D}$, $\tilde{B}_w$ and $\tilde{D}'$, $\tilde{B}_w'$ the corresponding lifts of $D$ and $B_w$ used to define $\pi_w$ and $\pi_w'$, respectively. By the continuity of the flow, we can define projections along the flow lines $\psi_D:\tilde{D}\to\tilde{D}'$ and $\psi_B:\tilde{B}_w'\to\tilde{B}_w$ (in suitable neighborhoods of $\tilde{\gamma}$) so that $\tilde{\pi}_w'=\psi_D\circ\tilde{\pi}_w\circ\psi_B$. These projections induce on the quotient return maps $P_D^r:D\to D$ and $P_B^l:\B\to\B$, where $r,l$ are integers such that $\tilde{D}'=T^r(\tilde{D})$ and $\tilde{B}_w=T^l(\tilde{B}_w')$. Hence, using $2(a)$ we conclude that
$$\pi_w'=P_D^r\circ\pi_w\circ P_B^l=\pi_w\circ P_B^{k},\ \text{where}\ k=nr+l.\qedhere$$
\end{enumerate}
\end{proof}

\subsubsection{Deformation by flow isotopies, equivalence and existence of tame Birkhoff sections.}
Given two embedded transverse sections of a non-singular flow, the natural way to relate these surfaces together with its first return maps (if defined) is via flow-isotopies, which are homeomorphisms between the surfaces obtained by pushing along flow lines. We explain this for the case of local Birkhoff sections. Let $\phi=\{\phi_t:M\to M\}_{t\in\R}$ be a non-singular regular flow having a saddle type periodic orbit $\gamma$ with orientable invariant manifolds. 

\begin{defn} 
Two local Birkhoff sections  $B$ and $B'$ at $\gamma$ are $\phi$\emph{-isotopic} if there exist collar neighborhoods $U\subset B$ and $U'\subset B'$ of $\gamma$ and a continuous and bounded function $s:\U\to\R$, such that the map $\psi$ given by $\psi(x)=\phi(s(x),x)$, $x\in \U$, defines a homeomorphism between $\U$ and $\U'$. (Recall the notation $\U=U\backslash\gamma$.)
\end{defn}

\begin{lem}
Let $B$ and $B'$ be two local Birkhoff sections at $\gamma$ which are $\phi$-isotopic, and let $\psi:\U\to\U'$ be a $\phi$-isotopy. Let us also assume that the first return map $P_B$ is defined for every point in $\U$. Then:
\begin{enumerate}
\item Up to shrinking the neighborhoods $U$ and $U'$ if necessary, it is satisfied that $\psi\circ P_{B}(x)=P_{B'}\circ\psi(x)$, for every $x\in\U$;
\item If $\psi'$ is another $\phi$-isotopy then $\exists$ $N\in\mathbb{Z}$ such that $\psi'(x)=\psi\circ P_{B}^N(x)$, for every $x\in\U$ sufficiently close to $\gamma$.
\end{enumerate} 
\end{lem}

\begin{proof}
Let $x\in\mathring{U}$ be such that $y=P_B(x)\in\mathring{U}$, and define $x'=\psi(x)$, $y'=\psi(y)$. Observe that all $x,y,x',y'$ lie in the same 
$\phi$-orbit segment, so in particular $y'=P_{B'}^k(x')$ for some $k\in\Z$ and we want to show that $k=1$. Let $\alpha$ be a path in $\mathring{U}$ connecting $y$ with $x$, and consider the oriented closed curve $\eta=[x,P_B(x)]\cdot\alpha$. Then the algebraic intersection number of $[\eta]\in H_1(W\backslash\gamma)$ with $[B]\in H_2(W,\partial W\cup\gamma)$ (where $W$ is a small tubular neighborhood of $\gamma$) is $[\eta]\cdot[B]=1$. Let $\alpha'=\psi(\alpha)$, which is a path in $B'$ connecting $x'=\psi(x)$ with $y'=\psi(y)$. Let $\eta'=[x',y']\cdot\alpha'$. Then $\eta'$ is homotopic to $\eta$ and thus $[\eta']\cdot[B']=[\eta]\cdot[B]=1$, so we conclude that $y'$ is the first intersection with $B'$ of the $\phi$-orbit segment starting at $x'$. This proves item 1. To prove item 2, observe that both $x'=\psi(x)$ and $x''=\psi'(x)$ are points in $\mathring{U}'$ lying in the same $\phi$-orbit segment, so there exists $t'\in\R$ such that $\phi_{t'}(x')=x''$ and thus $x''=P_{B'}^N(x')$ for some $N\in\Z$. The transversality of the $\phi$-orbits with $\mathring{B}'$ implies that $N$ varies continuously, so it is constant. 
\end{proof}

We need the following two propositions taken from \cite{bonatti-guelman}. The first one states that linking number and multiplicity determine the $\phi$-isotopy classes among the tame sections. Observe that a $\phi$-isotopy is defined by pushing along flow lines with a continuous and bounded function $s:\U\to\R$, so to obtain a bounded function it is essential the hypothesis of tameness.

\begin{proposition}[\cite{bonatti-guelman}, Lemma 3.6]\label{BG3.6}
Let $B$ and $B'$ be two \emph{tame} local Birkhoff sections on a saddle type periodic orbit $\gamma$. Then $B$ is $\phi$-isotopic to $B'$ if and only if $n(\gamma,B)=n(\gamma,B')$ and $m(\gamma,B)=m(\gamma,B')$.
\end{proposition}

The following proposition is used for the proof of Theorem \ref{thmA-statement}. It states that given two $\phi$-isotopic local Birkhoff sections $B$ and $B'$ at $\gamma$, it is possible to interpolate them to create a new local Birkhoff section that coincides with $B'$ near $\gamma$ and with $B$ outside a neighborhood of $\gamma$. 

\begin{proposition}[\cite{bonatti-guelman}, Lemma 3.7]\label{BG3.7}
Let $B$ and $B'$ be two local Birkhoff sections on a saddle type periodic orbit $\gamma$, with the same linking number and multiplicity. Then, there exists a neighborhood $W$ of $\gamma$ such that: For any neighborhood $O\subset W$ there exist another neighborhood $O'\subset O$ and a continuous and bounded function $s:W\backslash\gamma\to\R$ such that the map $\psi(u)=\phi(s(u),u)$ defines an $\phi$-isotopy onto its image, and
\begin{enumerate}
\item $\psi(u)\in\B'$, for all $u\in\B\cap O'$,
\item $\psi(u)=u$, for all $u\in\B \cap W\backslash O$.
\end{enumerate}  
\end{proposition}

\begin{rem}
In \cite{bonatti-guelman} the proofs of \ref{BG3.6} and \ref{BG3.7} are done just for the case when $m(\gamma,B)=m(\gamma,B')=1$. Nevertheless, the same proofs work in the general case, with the only difference that Lemma 3.5 in \cite{bonatti-guelman} must be replaced by item 2 of Proposition \ref{prop_combinatorial_relation_projections_first_return} of the present article. 
\end{rem}

Concerning the tameness condition defined at the beginning, general (local) Birkhoff sections need not to be tame. It is possible to construct smoothly embedded local Birkhoff sections in a neighborhood of a periodic orbit $\gamma$ that intersect $\Wsloc(\gamma)$ or $\Wuloc(\gamma)$ in a very narrow way, for example, along an infinite curve accumulating over an interval $I\subset\gamma$ not reduced to a singleton. Nevertheless, the previous proposition can be used to modify any given Birkhoff section in a neighborhood of the boundary components and obtain a tame one. (Observe that tameness is not demanded in the hypothesis of Proposition \ref{BG3.7}.) We get the following corollary:

\begin{cor}\label{prop_existence_tame_Birkhoff_sections}
Let $\{\phi_t:M\to M\}_{t\in\R}$ be a non-singular flow on a closed $3$-manifold and let $\iota:(\Sigma,\partial\Sigma)\to(M,\Gamma)$ be a Birkhoff section, such that every $\gamma\in\Gamma$ is a saddle type periodic orbit. Then, given a neighborhood $W$ of $\Gamma$ (that one can think of as a finite union of tubular neighborhoods around each curve in $\Gamma$) there exists a Birkhoff section $\iota':(\Sigma,\partial\Sigma)\to(M,\Gamma)$, with image $\Sigma'=\iota'(\Sigma)$, satisfying that:
\begin{enumerate}
\item $\mathring{\Sigma}\cap (M\backslash W)\equiv\mathring{\Sigma}'\cap (M\backslash W)$,
\item $\Sigma'$ is tame at each boundary curve $\gamma\in\Gamma$.
\end{enumerate}
\end{cor}

\subsection{Local version of Theorem \ref{thmA-statement}.}\label{subsection-Local_version_thm-A}
For $i=1,2$ consider a topological saddle type periodic orbit $\gamma^i$ of a flow $\{\phi_t^i:M_i\to M_i\}_{t\in\R}$ such that their local stable and unstable manifolds are orientable. Let $B_i\hookrightarrow M_i$ be a tame local Birkhoff section at $\gamma^i$, and assume that there exists a local conjugation $h:(B_1,P_{B_1})_{\gamma^1}\to(B_2,P_{B_2})_{\gamma^2}$ between the first return maps $P_{B_i}$. 

Let us recall that $h:\B_1\to \B_2$ is a homeomorphism, where $\B_i=B\backslash\gamma^i$, such that $P_{B_2}\circ h(x)=h\circ P_{B_1}(x)$, for every $x\in\B_1$ sufficiently close to $\gamma^1$. 

By Lemma \ref{lemma_conjugacion_vs_equiv_orbital}, we have that for every sufficiently small neighborhood $W_1$ of $\gamma^1$, the homeomorphism $h$ induces a homeomorphism 
$H:W_1\backslash\gamma^1\to W_2\backslash\gamma^2,$
where $W_2$ is a neighborhood of $\gamma^2$, which is a local orbital equivalence $(\phi^1,W_1\backslash\gamma_1)\to(\phi^1,W_1\backslash\gamma_1)$ in the sense of Definition \ref{defn_local_orbital_equivalence}, and whose restriction onto $\B_1\cap W_1$ coincides with $h$.  

\begin{thm}\label{thm_fund_local}
Consider a homeomorphism $H:W_1\backslash\gamma^1\to W_2\backslash\gamma^2$, where $W_1$ and $W_2$ are neighborhoods of $\gamma^1$ and $\gamma^2$ respectively, which verifies that
\begin{enumerate}[(a)]
\item $H:(\phi^1,W_1\backslash\gamma^1)\to(\phi^2,W_2\backslash\gamma^2)$ is an orbital equivalence,
\item $H(x)=h(x)$ for every $x\in\B_1\cap W_1$.
\end{enumerate} 
If it is satisfied that $m(\gamma^1,B_1)=m(\gamma^2,B_2)\text{ and } n(\gamma^1,B_1)=n(\gamma^2,B_2)$, then for every neighborhood $N\subset W_1$ there exists a homeomorphisms $H_N:W_1\to W_2$ such that:
\begin{enumerate}[(a)]
\item $H_N$ is a local orbital equivalence between $(\phi^1,W_1)_{\gamma^1}$ and $(\phi^2,W_2)_{\gamma^2}$,
\item $H_N(x)=H(x)$, for every $x\in W_1\backslash N$.
\end{enumerate}
\end{thm}

Cf. Remark \ref{remark_non-extension_of_almost-equiv}.

To prove Theorem \ref{thm_fund_local} we will show that given some neighborhood $N\subset W_1$ of $\gamma^1$, under the assumption that the sections are tame, it is possible to modify $H$ inside this neighborhood by pushing along flow lines, and obtain a new homeomorphism $H_N$ that extends as an orbital equivalence over the whole sets $W_i$. 

\subsubsection{Proof of Theorem \ref{thm_fund_local}.}
Consider an orbital equivalence
\begin{equation}\label{H}
H:(\phi^1,W_1\backslash\gamma^1)\to(\phi^2,W_2\backslash\gamma^2)
\end{equation} 
such that its restriction to the interior of the Birkhoff section $B_1$ coincides with $h:\B_1\to\B_2$, as in the hypothesis of \ref{thm_fund_local}. Observe that if we replace the neighborhoods $W_i$ by smaller neighborhoods $W_i'\subset W_i$, then it is enough to prove the theorem for these new neighborhoods. We will shrink the size of the $W_i$ several times in the course of the proof. Theorem \ref{thm_fund_local} relies on the following proposition.

\begin{proposition}\label{prop_fund}
For $i=1,2$ consider a saddle type periodic orbit $\gamma^i$ of a flow $\{\phi_t^i:M_i\to M_i\}_{t\in\R}$ such that their local stable and unstable manifolds are orientable. Consider:
\begin{itemize}
\item $B_i$ a tame local Birkhoff section at $\gamma^i$ and assume that there exists a local conjugation $h:(B_1,P_{B_1})_{\gamma^1}\to(B_2,P_{B_2})_{\gamma^2}$ between the first return maps $P_{B_i}$,
\item For each orbit $\gamma^i$ let $D_i$ be a local transverse section. Let $x^i=\gamma^i\cap D_i$ and let $\pi^i:(\B_i)_{w_i}\to D_i\backslash\{x^i\}$ be projections along the flow as defined in \ref{defn_projection_along_flow_lines} above, where $w_i$ is a fixed segment of the intersection of $B_i$ with the local invariant manifolds of $\gamma^i$.  
\end{itemize}
If it is satisfied that
$m(\gamma^1,B_1)=m(\gamma^2,B_2)\text{ and }n(\gamma^1,B_1)=n(\gamma^2,B_2)$,
then there exists a homeomorphism $h_D:D_1\to D_2$ satisfying:
\begin{enumerate}[(a)]
\item $h_D$ is a local conjugation between $(D_1,P_{D_1})_{x^1}$ and $(D_2,P_{D_2})_{x^2}$,
\item there exists a collar neighborhood $U_1\subset B_1$ of the curve $\gamma^1$ such that $h_D\circ\pi^1(x)=\pi^2\circ h(x)$, $\forall\ x\in U_1$.
\end{enumerate}
\end{proposition}

We postpone the proof of this proposition to the end. Now, with this result we can deduce Theorem \ref{thm_fund_local} in the following way: For each orbit $\gamma^i$ consider a local transverse section $D_i$  and a projection along the flow $\pi^i:(\B_i)_{w_i}\to D_i\backslash\{x_i\}$ as the previous proposition. Since we assume in the hypothesis of Theorem \ref{thm_fund_local} that linking number and multiplicities coincide for both $B$ and $B'$, we can apply Proposition \ref{prop_fund} and obtain a homeomorphism $h_D:D_1'\to D_2'$ satisfying $(a)$ and $(b)$, where $D_i'\subset D_i$ are smaller transverse sections. Using Proposition \ref{prop_conjugacion_vs_equiv_orbital} we see that there exist a tubular neighborhood $W_i'$ of each $\gamma^i$ and a local orbital equivalence 
\begin{equation}\label{H_D}
H_D:(\phi^1,W_1')_{\gamma^1}\to(\phi^2,W_2')_{\gamma^2}
\end{equation}
such that $H_D(x)=h_D(x)$, for every $x\in D_1'\cap W_1'$. Without loss of generality we can assume that $W_1=W_1'$. That is, we can assume that the two homeomorphisms $H$ and $H_D$ have the same domain.  

Theorem \ref{thm_fund_local} follows directly from the following proposition.

\begin{proposition}\label{thm_falso}
For every neighborhood $N\subset W_1$ of $\gamma^1$ there exists another neighborhood $N'\subset N$ and a local orbital equivalence 
$H_N:(\phi^1,W_1)_{\gamma^1}\to(\phi^2,W_2)_{\gamma^2}$
such that:
\begin{enumerate}[(a)]
\item $H_N(x)=H(x)$ for every $x\in W_1\backslash N$,
\item $H_N(x)=H_D(x)$, for every $x\in N'$.
\end{enumerate}
\end{proposition}

This interpolation between $H$ and $H_D$ is possible due to item (b) in Proposition \ref{prop_fund}. We dedicate the rest of this part to prove Propositions \ref{thm_falso} and \ref{prop_fund}.

\subsubsection{Proof of Proposition \ref{thm_falso}.}
We start by writing a scheme of the steps in the proof, and then we show each step.

\ 

\begin{Scheme of the proof}\ 

\noindent
Consider a neighborhood $N\subset W_1$. Let us denote $N_1=N$. 
\paragraph*{Step 1}
For every neighborhood $O_1\subset W_1$ of $\gamma^1$ we can find a smaller neighborhood $O_1'\subset O_1$ and an orbital equivalence $F:(\phi^1,W_1\backslash\gamma^1)\to(\phi^2,W_2\backslash\gamma^2)$ such that:
\begin{enumerate}[(a)]
\item $F(x)=H(x)$ for every $x\in B_1\cap(W_1\backslash O_1)$,
\item $F(x)=H_D(x)$ for every $x\in\B_1\cap O_1'$.
\end{enumerate}
\paragraph*{Step 2}
We will find a collection of neighborhoods
\begin{itemize}
\item $N_1'\subset O_1'\subset O_1\subset N_1\subset W_1,$
\item $N_2'\subset O_2'\subset O_2\subset N_2\subset W_2.$
\end{itemize}
Let us define $V_i=N_i\backslash O_i$ and $V_i'=O_i'\backslash N_i'$. If these neighborhoods are suitably chosen, we will be able: 
\begin{itemize}
\item to make an interpolation along the flow lines between $H$ and $F$ supported in the region $V_1$, and obtain an orbital equivalence $H_V:(\phi^1,V_1)\to(\phi^2,V_2)$ satisfying that:
\begin{enumerate}[(a)]
\item $H_V(x)=H(x)$ for every $x\in\partial N_1$,
\item $H_V(x)=F(x)$ for every $x\in\partial O_1$;
\end{enumerate}
\item to make an interpolation along the flow lines between $H_D$ and $F$ supported in the region $V_1'$, and obtain an orbital equivalence $H_V':(\phi^1,V_1')\to(\phi^2,V_2')$ satisfying that:
\begin{enumerate}[(a)]
\item $H_V'(x)=F(x)$ for every $x\in\partial O_1'$,
\item $H_V'(x)=H_D(x)$ for every $x\in\partial N_1'$.
\end{enumerate}
\end{itemize} 
\paragraph*{Step 3}
We will define $H_N$ in the following way:
\begin{equation}\label{H_N}
H_N(x)=\left\{\begin{array}{lcc}
              H(x)   & \text{if} & x\in W_1\backslash N_1,\\
            \\H_V(x) & \text{if} & x\in N_1\backslash O_1,\\
			\\F(x)   & \text{if} & x\in O_1\backslash O_1',\\
			\\H_V'(x)& \text{if} & x\in O_1'\backslash N_1',\\
			\\H_D(x) & \text{if} & x\in N_1'.\\
			\end{array}
     \right.
\end{equation}
Observe that $H_N$ is well-defined in all the boundaries $\partial N_1'$, $\partial O_1'$, $\partial O_1$, $\partial N_1$ and gives rise to a local orbital equivalence satisfying the conclusion of \ref{thm_falso}.
\end{Scheme of the proof}

\paragraph*{Step 1: The construction of $F$}

\begin{lem}\label{homeo_auxiliar}
Let $O_1\subset W_1$ be a tubular neighborhood of $\gamma^1$. Then, there exists another tubular neighborhood $O_1'\varsubsetneq O_1$ and a homeomorphism $F:W_1\backslash\gamma^1\to W_2\backslash\gamma^2$ satisfying that:
\begin{enumerate}[(a)]
\item $F:(\phi^1,W_1\backslash\gamma^1)\to(\phi^2,W_2\backslash\gamma^2)$ is an orbital equivalence,
\item $F(x)=H(x)$ for every $x\in B_1\backslash O_1$,
\item $F(x)=H_D(x)$ for every $x\in\B_1\cap O_1'$.
\end{enumerate}
\end{lem}

Let us define 
\begin{equation}
S\coloneqq \left\{H_D(x):x\in B_1\cap W_1\right\}.
\end{equation}
Then $S$ is a local Birkhoff section at $\gamma^2$. Observe in addition that $S$ satisfies that $m(\gamma^1,S)=m(\gamma^2,B_2)$ and $n(\gamma^1,S)=n(\gamma^2,B_2).$ Following Proposition \ref{BG3.6} we see that $B_2$ and $S$ are $\phi^2$-isotopic. Our next goal is to deform $B_2$ pushing along flow lines and obtain a new local Birkhoff section $B_2'$, that coincides with $B_2$ outside a tubular neighborhood $O_2$ of the orbit $\gamma^2$ and coincides with $S$ in a smaller neighborhood $O_2'\subset O_2$.

\begin{lem}\label{seccion_interpolada}
Let $O_2\subset W_2$ be a neighborhood of $\gamma^2$. Then, there exist a smaller neighborhood $O_2'\subset O_2$, a local Birkhoff section $B_2'$ at $\gamma^2$ and a homeomorphism $\psi:\B_2\to\B_2'$ such that
\begin{enumerate}
\item $B_2'\backslash O_2\equiv B_2\backslash O_2$ and $B_2'\cap O_2'\equiv S\cap O_2'$,
\item $\psi$ is an isotopy along the $\phi^2$-orbits and satisfies that $\psi(y)=y$, $\forall$ $y\in B_2\backslash O_2$,  
\item if $x\in\B_1$ satisfies that $h(x)\in B_2\cap O_2'$ then $\psi(h(x))=H_D(x)$.
\end{enumerate}  
\end{lem}

\begin{proof}
As usual we denote $\mathring{S}=S\backslash\gamma^2$. Let us start by constructing a $\phi^2$-isotopy from $\B_2$ to $\mathring{S}$. By Proposition \ref{BG3.6} there exists a collar neighborhood $U_2\subset B_2$ of the curve $\gamma^2$ and a continuous and bounded function $s:\U_2\to\R$ such that the map
\begin{equation}
\varphi:\U_2\to\mathring{S}\backslash\gamma^2\text{ defined by }\varphi(y)=\phi^2(y,s(y))
\end{equation} 
is a $\phi^2$-isotopy from $\U_2$ into $\mathring{S}$. 

Consider $U_1=h^{-1}(U_2)\subset B_1$. Chose the neighborhoods $U_i$ sufficiently small, such that they are contained in the domain of definition of the projections $\pi^i:(B_i)_{\gamma^i}\to(D_i)_{x^i}$. Let $P_S$ be the first return to the local Birkhoff section $S$.
\begin{lem}\label{ajustando_la_isotopia}
There exists $k\in\Z$ such that $\varphi(h(x))=P_S^k\circ H_D(x)$, for every $x\in\U_1$.
\end{lem}
\noindent
We postpone the proof of this lemma to the end. As a consequence this lemma, up to shrinking the size of $U_i$ if necessary and composing with some power of $P_S$ on the left, we can assume that $\varphi(h(x))=H_D(x)$ for every $x\in\U_1$. 

Given a neighborhood $O_2\subset W_2$ of $\gamma^2$, we can use Proposition \ref{BG3.7} and find another neighborhood $O_2'\subset O_2$ and a continuous and bounded function $s':\B_2\to\R$, such that the map 
\begin{equation}
\psi:y\mapsto\phi^2(y,s'(y)),\ y\in\B_2 
\end{equation}
satisfies the following
\begin{enumerate}
\item The image $B_2'\coloneqq\psi(B_2)$ is a local Birkhoff section and $\psi:\B_2\to\B_2'$ is a flow isotopy,
\item $\psi(y)=y$ for every $y\in B_2\backslash O_2$,
\item $\psi(y)=\varphi(y)$ for every $y\in\B_2\cap O_2'$.  
\end{enumerate}
The neighborhood $O_2'$, the section $B_2'$ and the map $\psi$ satisfy the properties claimed in Lemma \ref{ajustando_la_isotopia}. 
\end{proof}

To complete the proof it remains to prove Lemma \ref{ajustando_la_isotopia}.

\begin{proof}[Proof of Lemma \ref{ajustando_la_isotopia}]
Recall that the homeomorphism $h_D$ satisfies properties $(a)$ and $(b)$ of \ref{prop_fund} and that $H_D$ coincides with $h_D$ over the transverse section $D_1\cap W_1$. Recall also that for every point $z\in W_i$ we denote the connected component of $\mathcal{O}^i(z)\cap W_i$ that contains $z$ as $\mathcal{O}^i_{W_i}(z)$. We claim that if $x\in\U_1$ then $\mathcal{O}^2_{W_2}(H_D(x))$ coincides with $\mathcal{O}^2_{W_2}(h(x))$. Since the projections along the flow preserve orbit segments, it is satisfied that 
\begin{equation*}
\mathcal{O}^i_{W_i}(z)=\mathcal{O}^i_{W_i}(\pi^i(z)),\ \forall z\in\U_i. 
\end{equation*}
So we have that
\begin{equation*}
\mathcal{O}^2_{W_2}(H_D(x))=H_D(\mathcal{O}^1_{W_1}(x))=H_D(\mathcal{O}^1_{W_1}(\pi^1(x)))=\mathcal{O}^2_{W_2}(H_D\circ\pi^1(x)).
\end{equation*}
Since $H_D\circ\pi^1(x)=h_D\circ\pi^1(x)$ and by \ref{prop_fund}-$(b)$, the last term of the previous equality is equal to
\begin{equation*}
\mathcal{O}^2_{W_2}(h_D\circ\pi^1(x))=\mathcal{O}^2_{W_2}(\pi^2\circ h(x))=\mathcal{O}^2_{W_2}(h(x)),
\end{equation*}
so the claim follows. Now, since $\varphi(h(x))$ is a point in $\mathring{S}$ whose orbit segment inside $W_2$ equals that of $H_D(x)$ we deduce that $\varphi(h(x))=P_S^k\circ H_D(x)$ for some $k\in\Z$. By the continuity of the flow this integer must vary continuously with respect to $x$, so it is constant. This completes the proof of the lemma.
\end{proof}

The statement of Lemma \ref{homeo_auxiliar} is a consequence of Lemma \ref{seccion_interpolada}.

\begin{proof}[Proof of Lemma \ref{homeo_auxiliar}]
Let $V\subset B_1$ be a collar neighborhood of $\gamma^1$ contained in the domain of definition of the first return map $P_{B_1}$. Let $\mathcal{V}$ be the union of all the compact orbit segments connecting each point in $V$ with its first return to $B_1$. Up to shrinking the size of the neighborhoods $W_i$ if necessary, we can assume that $W_1\subset\interior(\mathcal{V})$.

Given $O_1\subset W_1$ consider $O_2=H(O_1)$. Then, Lemma \ref{seccion_interpolada} gives another neighborhood $O_2'\subset O_2$, a section $B_2'$ and a $\phi^2$-isotopy $\psi:\B_2\to\B_2'$. Let us define $O_1'=H_D^{-1}(O_2')$ and $h':\B_1\to\B_2'$ given by $h'(x)=\psi\circ h(x)$. Then, the homeomorphism $h'$ is a local conjugation between the first return maps to the Birkhoff sections $B_1$ and $B_2'$ respectively. So by Proposition \ref{prop_conjugacion_vs_equiv_orbital} it induces an orbital equivalence $F$ with domain $\mathcal{V}\backslash\gamma^1$ that coincides with $h'$ over $\B_1\cap\mathcal{V}$. Since $W_1\subset\interior(\mathcal{V})$ we can consider its restriction to $W_1\backslash\gamma^2$, that is $F:W_1\backslash\gamma^1\to W_2\backslash\gamma^2$. It is direct that $F$ coincides with $H$ over $B_1\backslash O_1$ and with $H_D$ over $\B_1\cap O_1'$. 
\end{proof}

\paragraph*{Step 2: The interpolation}\ 

We start by describing how to choose the neighborhoods $N_i'\subset O_i'\subset O_i\subset N_i$. Given a regular tubular neighborhood $O_1\subset N_1$ of $\gamma^1$ consider the annulus $A=B_1\cap N_1\backslash O_1$. We claim that if $O_1$ is sufficiently small, then there exists a compact annulus $K$ with non-empty interior and contained in $\interior(A)$, such that $P_{B_1}(K)\subset\interior(A)$. The claim follows directly by examining the first return map in $\B_1$ as in Figure \ref{fig_thm_global_1}.

A tubular neighborhood is said to be regular if its closure is a submanifold homeomorphic to a compact disk times an interval.  Given $N_1=N$ we chose a family of regular tubular neighborhoods $N_1'\subset O_1'\subset O_1\subset N_1$ in the following way:
\begin{enumerate}
\item Choose $O_1\subset N_1$ such that there exists a compact annulus $K$ with non-empty interior contained in $A_1=B_1\cap N_1\backslash O_1$, which satisfies that $P_{B_1}(K)\subset\interior(A_1)$,
\item Choose $O_1'\subset O_1$ given by Lemma \ref{homeo_auxiliar},
\item Choose $N_1'\subset O_1'$ such that there exists a compact annulus $K'$ with non-empty interior contained in $A_1'=B_1\cap O_1'\backslash N_1'$, which satisfies that $P_{B_1}(K')\subset\interior(A_1')$.
\end{enumerate}
Let $F:W_1\backslash\gamma^1\to W_2\backslash\gamma^2$ given by Lemma \ref{homeo_auxiliar} for the chosen neighborhood $O_1$. We define as well:
\begin{enumerate}
\item $N_2=H(N_1)$,
\item $O_2=F(O_1)$ and $O_2'=F(O_1')$,
\item $N_2'=H_D(N_1')$.
\end{enumerate}
Let us define $V_i=\bar{N}_i\backslash\interior(O_i)\subset M_i$  and $V_i'=\bar{O}_i'\backslash\interior(N_i')\subset M_i$. Since the neighborhoods that we consider are regular tubular neighborhoods it follows that $V_i$ and $V_i'$ are homeomorphic to a closed annulus times a circle. The topological equivalence $H_N$ will be an interpolation between $H$ and $F$ over the set $V_1$ and between $H_D$ and $F$ over the set $V_1'$. We describe first the topology of these interpolating sets and then we indicate how to make these interpolations over $V_1$ and $V_2$. 

\ 

\begin{interpolating nbhd}\ 

\noindent
Consider the compact sets 
\begin{align}
&V_i=\bar{N}_i\backslash\interior(O_i)\subset M_i\\
&V_i'=\bar{O}_i'\backslash\interior(N_i')\subset M_i.
\end{align}
From now on we concentrate just on $V_i$, $i=1,2$, since all the arguments are analogous for $V_i'$. The boundary components of each $V_i$ are $\partial N_i$ and $\partial O_i$. The compact annulus $A_i=B_i\cap V_i$ is a properly embedded surface in $V_i$. Consider the compact annuli with non-empty interior $K_1=K\subset\interior(A_1)$ and $K_2=h(K)\subset\interior(A_2)$. Each $K_i$ divides $A_i$ into three annuli as in Figure \ref{fig_thm_global_1}. We name the boundaries of $K_i$ as $\alpha_i$ and $\beta_i$ according to this figure. The map $h$ restricts to a homeomorphism $h:A_1\to A_2$ which defines a conjugacy between the return maps $P_{K_i}:K_i\to A_i$. 

\begin{figure}[t]
\begin{center}
\includegraphics[width=0.5\textwidth, height=\textwidth, angle=90, keepaspectratio]{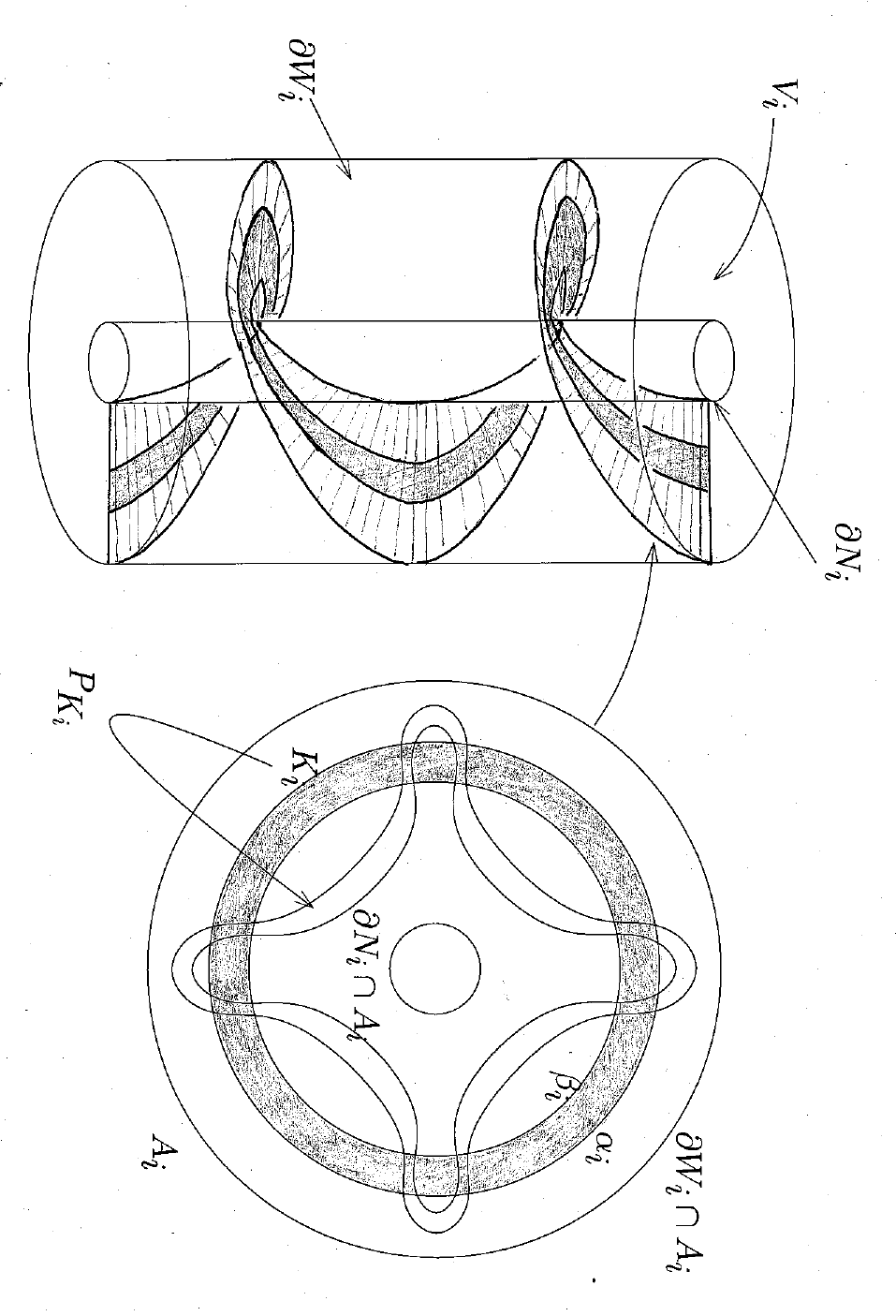}
\caption{The sets $V_i$ and the annuli $A_i$.}
\label{fig_thm_global_1}
\end{center}
\end{figure}

Consider the set
\begin{equation}
\mathcal{K}_i=\{\phi_t^i(u):u\in K_i, 0\leq t\leq\tau^i(u)\}
\end{equation}
where $\tau^i(u)$ is the time of the first return to $B_i$ of a point $u\in K_i$. This set is the union of all the compact orbit segments joining a point $u\in K_i$ with its first return $P_{K_i}(u)=\phi^i(u,\tau^i(u))$. Observe that these orbit segments are disjoint from the boundary components of $V_i$ so it is satisfied that $\mathcal{K}_i\subset\interior(V_i)$. Since we have that $H{|}_{K_1}=F{|}_{K_1}=h{|}_{K_1}$ then it is verified that $\mathcal{K}_2=H(\mathcal{K}_1)=F(\mathcal{K}_1)$.

The annulus $A_i$ is an essential surface in $V_i$ and the complement of $\mathcal{K}_i\cup A_i$ has two connected components. We call $\mathcal{C}_i\text{ and }\mathcal{D}_i$ 
to the closure of these components, where the first one is the component that contains $\partial N_i$ in its boundary and the second one is the one that contains $\partial O_i$ in its boundary. So, we have a decomposition 
\begin{equation*}
V_i=\mathcal{C}_i\cup\mathcal{K}_i\cup\mathcal{D}_i
\end{equation*}
of the neighborhood $V_i$ into three closed sets.

If we cut the set $V_i$ along $A_i$ we obtain a manifold $\widetilde{V}_i$ homeomorphic to the product of $A_i$ with a closed interval that we have depicted in Figure \ref{fig_thm_global_2}. This manifold is equipped with a map $\widetilde{V}_i\to V_i$ which corresponds to glue back the two copies of $A_i$. The three sets $\mathcal{C}_i$, $\mathcal{K}_i$ and $\mathcal{D}_i$ lift into $\widetilde{V}_i$ and gives a decomposition
\begin{equation*}
\widetilde{V}_i=\widetilde{\mathcal{C}}_i\cup\widetilde{\mathcal{K}}_i\cup\widetilde{\mathcal{D}}_i
\end{equation*} 
into three compact sets, each one homeomorphic to an annulus times an interval. The components $\widetilde{\mathcal{C}}_i$ and $\widetilde{\mathcal{D}}_i$ are disjoint, and they intersect $\widetilde{\mathcal{K}}_i$ along the annuli $L_{\alpha_i}$ and $L_{\beta_i}$ respectively, as in Figure \ref{fig_thm_global_2}. 

Observe that the foliation by orbit segments in $V_i$ lift into a foliation by segments in $\widetilde{V}_i$ which are transverse to the copies of $A_i$. Let us denote by $\widetilde{A}_i^0$ to the copy of $A_i$ where the lifted orbits point inward the manifold $\widetilde{V}_i$ and by $\widetilde{A}_i^1$ to the other one where the orbits point outward. For every point $u\in K_i$ the orbit segment connecting $u$ with its first return to $A_i$ is parametrized by $s\mapsto\phi^i(u,s)$, $s\in[0,\tau^i(u)]$, and it lifts into $\widetilde{\mathcal{K}}_i$ as a compact interval connecting the two copies of $A_i$ inside $\widetilde{V}_i$. So the set $\widetilde{\mathcal{K}}_i$ is a union of compact segments joining the two copies of $A_i$, and we can put coordinates 
\begin{equation}\label{coordenadas}
{\widetilde{\mathcal{K}}}_i\to \{(u,s)\in K_i\times[0,+\infty):0\leq s\leq\tau^i(u)\}.
\end{equation}

\begin{figure}[t]
\begin{center}
\includegraphics[width=0.5\textwidth, height=\textwidth, keepaspectratio]{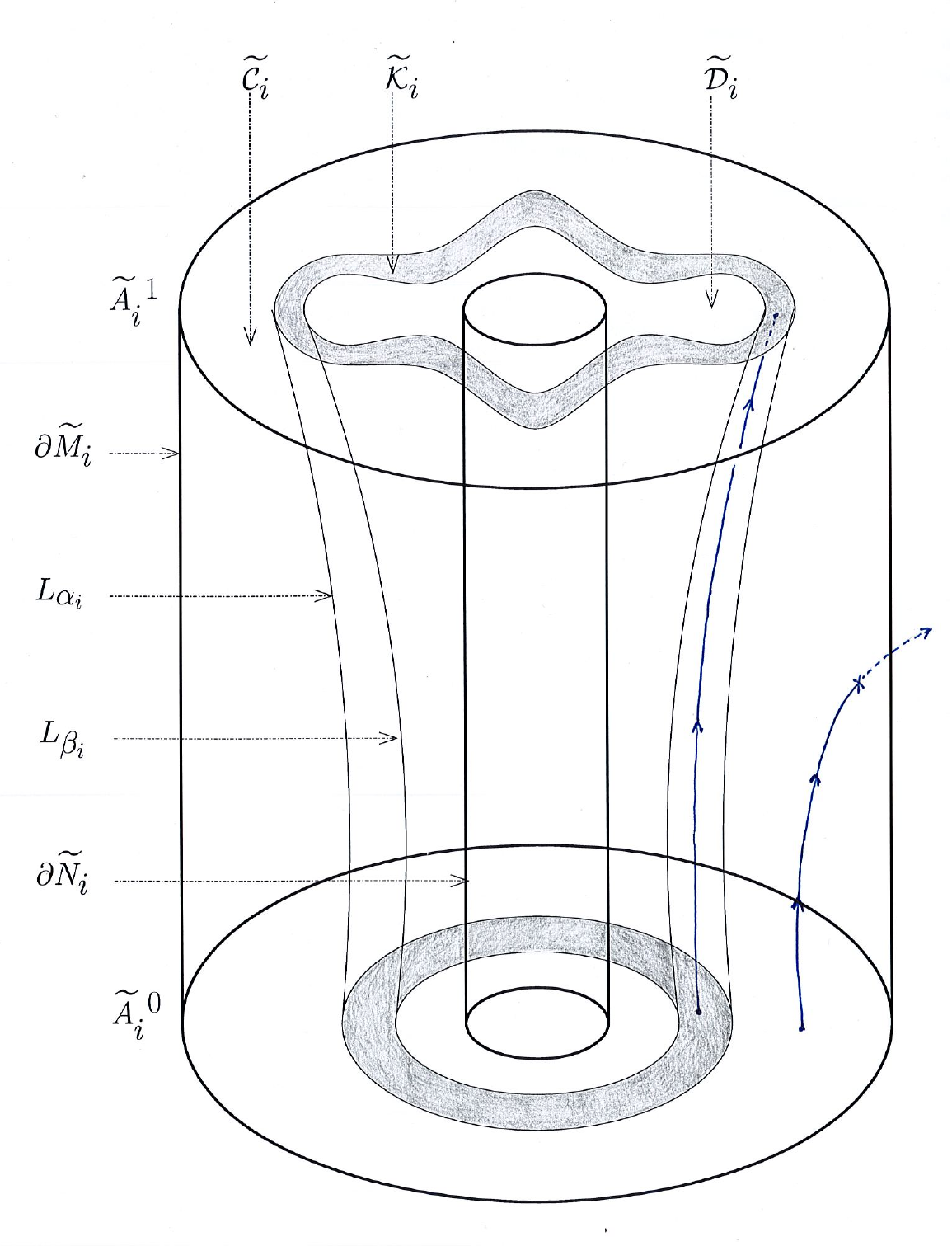}
\caption{The sets $\widetilde{V}_i$. We have depicted two orbits segments lifted into $\widetilde{V}_i$, one inside the set $\widetilde{\mathcal{K}}_i$ which connects $\widetilde{A}_i^0$ with $\widetilde{A}_i^1$ and the other one inside the components $\widetilde{C}_i$.}
\label{fig_thm_global_2}
\end{center}
\end{figure}
\end{interpolating nbhd}

\ 

\begin{The Construction of $H_V$ and $H_V'$}
\end{The Construction of $H_V$ and $H_V'$}

\begin{lem}\label{home_interpolado}
For the neighborhoods $N_i'\subset O_i'\subset O_i\subset N_i$ previously chosen, there exists homeomorphisms $H_V:V_1\to V_2$ and $H_V':V_1'\to V_2'$ satisfying that:
\begin{enumerate}[(a)]
\item $H_V:(\phi^1,V_1)\to(\phi^2,V_2)$ is an orbital equivalence,
\item $H_V(x)=H(x)$ for every $x\in\partial N_1$,
\item $H_V(x)=F(x)$ for every $x\in\partial O_1$;
\end{enumerate}
and
\begin{enumerate}[(a)]
\item $H_V':(\phi^1,V_1')\to(\phi^2,V_2')$ is an orbital equivalence,
\item $H_V'(x)=F(x)$ for every $x\in\partial O_1'$,
\item $H_V'(x)=H_D(x)$ for every $x\in\partial N_1'$.
\end{enumerate}
\end{lem}
\begin{proof}
We just do the construction of $H_V$, the other one being analogous. The key fact to prove \ref{home_interpolado} is that $H(x)=F(x)=h(x)$ for every $x\in B_1\cap V_1$. Observe that, since $\mathcal{K}_2\cup A_2=H(\mathcal{K}_1\cup A_1)=F(\mathcal{K}_1\cup A_1)$ and $H|_{A_1}\equiv F|_{A_1}\equiv h|_{A_1}$, it follows that $H(\mathcal{C}_1)=\mathcal{C}_2$ and $F(\mathcal{D}_1)=\mathcal{D}_2$. We are interested in the homeomorphisms 
\begin{align*}
& H:\mathcal{C}_1\cup\mathcal{K}_1\to\mathcal{C}_2\cup\mathcal{K}_2 \\
& F:\mathcal{K}_1\cup\mathcal{D}_1\to\mathcal{K}_2\cup\mathcal{D}_2
\end{align*}
obtained by restriction of $H$ and $F$ to the sets $\mathcal{C}_1\cup\mathcal{K}_1$ and $\mathcal{K}_1\cup\mathcal{D}_1$, respectively. We interpolate them over the closed set $\mathcal{K}_1$. Let
\begin{align*}
& \widetilde{H}:\widetilde{\mathcal{C}}_1\cup\widetilde{\mathcal{K}}_1\to\widetilde{\mathcal{C}}_2\cup\widetilde{\mathcal{K}}_2 \\ 
& \widetilde{F}:\widetilde{\mathcal{K}}_1\cup\widetilde{\mathcal{D}}_1\to\widetilde{\mathcal{K}}_2\cup\widetilde{\mathcal{D}}_2
\end{align*}
be the lifts of these maps to $\widetilde{V}_i$. Since $H$ and $F$ preserve the oriented orbit segments and coincide with $h$ over $K_1$, we can use the coordinates \ref{coordenadas} and write, for every point $(u,s)\in\widetilde{\mathcal{K}}_1$,   
\begin{align}
& \widetilde{H}(u,s)=(h(u),\theta(u,s)) \\ 
& \widetilde{F}(u,s)=(h(u),\eta(u,s)),
\end{align} 
where each function $\theta(u,\cdot)$, $\eta(u,\cdot)$ is an increasing homeomorphism between the segments $[0,\tau^1(u)]$ and $[0,\tau^2(h(u))]$, continuously parametrized over $u\in K_1$. Observe that for every $0\leq r\leq 1$, it follows that the convex combination $r\cdot\theta(u,\cdot)+(1-r)\cdot\eta(u,\cdot)$ is also an increasing homeomorphism from $[0,\tau^1(u)]$ to $[0,\tau^2(h(u))]$. 

Let $\rho:K_1\to[0,1]$ be a continuous function such that $\rho\equiv 1$ in a neighborhood of $\alpha_1$ and $\rho\equiv 0$ in a neighborhood of $\beta_1$. We define a map $\widetilde{H}_{V}:\widetilde{V}_1\to\widetilde{V}_2$ as follows: for every $\tilde{x}\in\widetilde{V}_1$,
\begin{equation}
\widetilde{H}_V(\tilde{x})=\left\{\begin{array}{lcc}
               \widetilde{H}(\tilde{x}) & \text{if}  & \tilde{x}\in\widetilde{\mathcal{C}}_1, \\
            \\ (h(u),\rho(u)\cdot\theta(u,s)+(1-\rho(u))\cdot\eta(u,s)) & \text{if}  & \tilde{x}=(u,s)\in\widetilde{\mathcal{K}}_1,\\
			\\ \widetilde{F}(\tilde{x}) & \text{if} & \tilde{x}\in\widetilde{\mathcal{D}}_1.          
             					  \end{array}
   						   \right.
\end{equation}
The map $\widetilde{H}_V$ is a well-defined homeomorphisms that preserves the foliations by (the lifts of the) orbit segments. It coincides with $\widetilde{H}$ over $\widetilde{\mathcal{C}}_1$ and with $\widetilde{F}$ over $\widetilde{\mathcal{D}}_1$. Observe also that, because of the particular election of the function $\rho$, it follows that $\widetilde{H}_V$ coincides with $\widetilde{H}$ in a neighborhood of $L_{\alpha_1}=\widetilde{\mathcal{C}}_1\cap\widetilde{\mathcal{K}}_1$ and with $\widetilde{F}$ in a neighborhood of $L_{\beta_1}=\widetilde{\mathcal{K}}_1\cap\widetilde{\mathcal{D}}_1$. By construction, over the union $\widetilde{A}_1^0\cup\widetilde{A}_1^1$ the map $\widetilde{H}_V$ coincides with the homeomorphism $\tilde{h}:\widetilde{A}_1^0\cup\widetilde{A}_1^1\to\widetilde{A}_2^0\cup\widetilde{A}_2^1$ that is obtained by lifting $h:A_1\to A_2$ to $\widetilde{V}_1$. So, if we glue back the two copies of $A_i$ inside $\widetilde{V}_i$, then $\widetilde{H}_V$ induces a homeomorphism 
\begin{equation}
H_V:V_1\to V_2
\end{equation}
which satisfies:
\begin{enumerate}[(a)]
\item for each $x\in V_1$ the map $H_V$ takes each oriented orbit segment $\mathcal{O}^1_{V_1}(x)$ homeomorphically onto the orbit segment $\mathcal{O}^2_{V_2}(H_V(x))$ preserving orientations,
\item $H_V$ coincides with $H$ on the set $\mathcal{C}_1$,
\item $H_V$ coincides with $F$ on the set $\mathcal{D}_1$.\qedhere
\end{enumerate}
\end{proof}

\paragraph*{Step 3: The construction of $H_N$}\ 

To finish the proof of \ref{thm_falso} just observe that the original neighborhood $W_1$ can be decomposed as the union of five compact manifolds which intersect along boundary tori, i.e. 
$$W_1=W_1\backslash N_1\cup V_1\cup O_1\backslash O_1'\cup V_1'\cup N_1'.$$ The homeomorphisms $H$, $H_V$, $F$, $H_V'$ and $H_D$ match well along the boundaries and give rise to the homeomorphism $H_N$ as we defined in \eqref{H_N}. It is an orbital equivalence, since it is when restricted to each piece of the decomposition of $W_1$, and clearly satisfies the properties stated in \ref{thm_falso}. This concludes the proof of \ref{thm_falso} as well as the proof of \ref{thm_fund_local}.

\subsubsection{Proof of Proposition \ref{prop_fund}.}
For each local Birkhoff section $B_i$ at the curve $\gamma^i$, $i=1,2$, consider a projection along the flow 
$\pi^i:(B_i)_{w^i}\to D_i\backslash\{x^i\}$ satisfying the hypotheses of Proposition \ref{prop_fund}, with $D_i$ a local transverse disk and $w^i$ a fixed segment of $B_i\cap\left(\Wsloc(\gamma^i)\cup\Wuloc(\gamma^i)\right)$. 

Let $n=n(\gamma^1,B_1)=n(\gamma^2,B_2)$ and $m=m(\gamma^1,B_1)=m(\gamma^2,B_2)$. Let $U_1\subset B_1$ be a collar neighborhood of $\gamma^1$ and let $U_2=h(U_1)$. We choose $U_1$ sufficiently small such that each $\U_i=U_i\backslash\gamma^i$ is contained in the domain of definition of $P_{B_i}$ and in the domain of definition of $\pi^i$. Consider the segment $v_i\subset D_i$ that equals the intersection of $D_i$ with the branch of $\Wsloc(\gamma^i)\cup\Wuloc(\gamma^i)$ that contains $w_i$. Observe that $\pi^i$ projects the points in the segment $w_i$ into the segment $v_i$.

We prove first the proposition when $n=m=1$, and then we comment how to deduce the general case using Proposition \ref{prop_combinatorial_relation_projections_first_return}.

\paragraph*{Case $n=m=1$}\ 

\begin{figure}[t]
\begin{center}
\includegraphics[height=\textwidth, angle=-90]{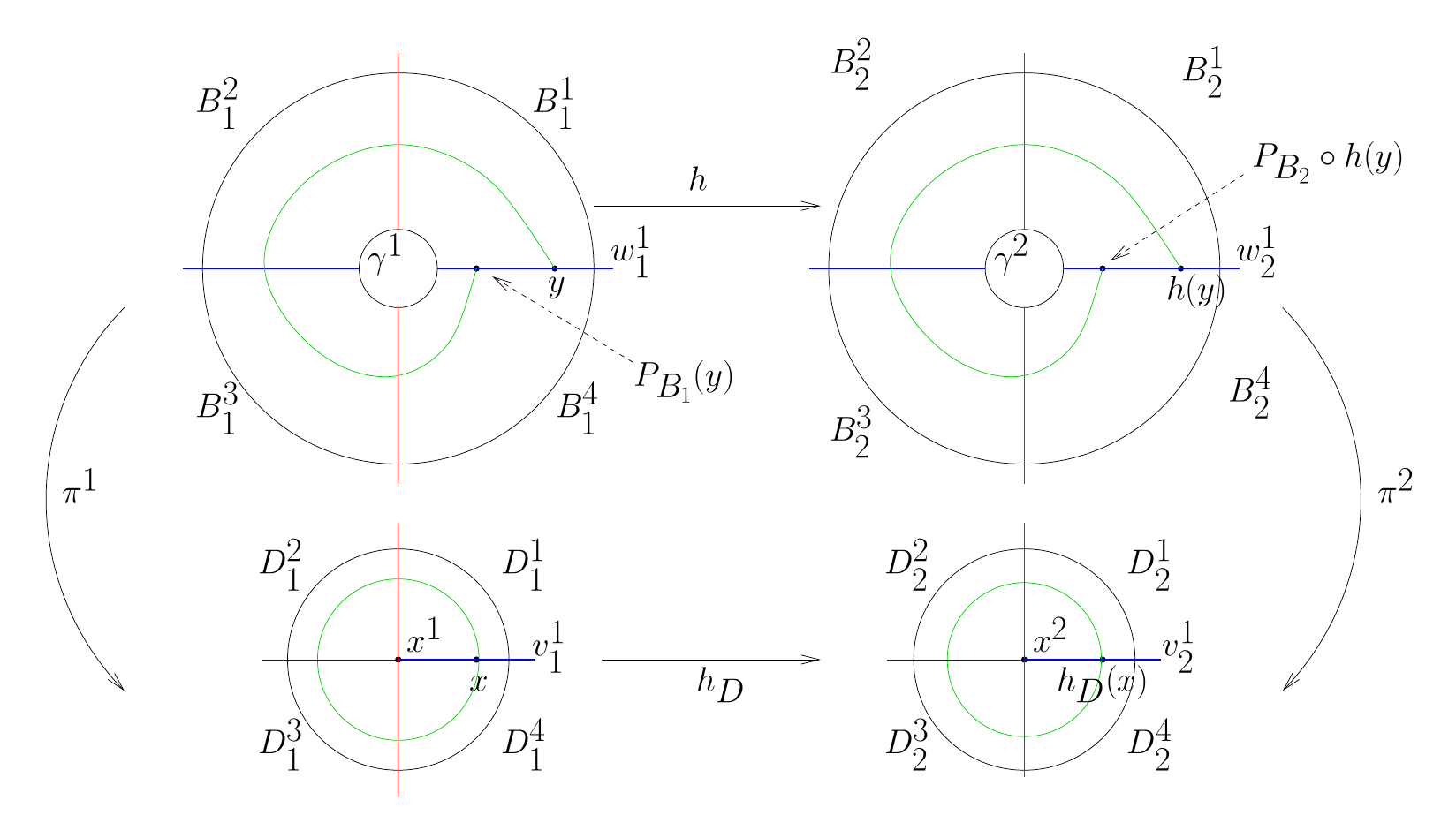}
\caption{The closed curve starting at $x\in v_1^1$ lift by the maps $\eta_s$ to a curve that connects $y=\eta_1(x)$ with $P_{B_1}(y)=\eta_4(x)$.}
\label{fig_prop_fund}
\end{center}
\end{figure}

Since $n=1$ we have that each Birkhoff section $B_i$ can be partitioned into four quadrants. Following the convention of Proposition \ref{prop_combinatorial_relation_projections_first_return}, we label the quadrants of $D_i$ and the quadrants of $B_i$ as 
\begin{align*}
&D^s_i,\ s=1,\dots,4,\\&B^s_i,\ s=1,\dots,4
\end{align*}
in such a way that $D_i^1$ and $D_i^4$ intersect along the segment $v_i$ and $B_i^1$ and $B_i^4$ intersect along $w_i$. See Figure \ref{fig_prop_fund}. We consider as well the boundaries of the quadrants
\begin{align*}
&v_i=v_i^1=D_1^4\cap D_1^1,\ v_i^2=D_1^1\cap D_1^2,\ v_i^3=D_1^2\cap D_1^2,\ v_i^4=D_1^3\cap D_1^4\\
&w_i=w_i^1=B_1^4\cap B_1^1,\ w_i^2=B_1^1\cap B_1^2,\ w_i^3=B_1^2\cap B_1^2,\ w_i^4=B_1^3\cap B_1^4.
\end{align*}
Observe that for every $s=1,\dots,4$ the map $h$ takes points in $B_1^s$, $D_1^s$, $w_1^s$ and $v_1^s$ into $B_2^s$, $D_2^s$, $w_2^s$ and $v_2^s$, respectively. The restriction of $\pi^i$ to the quadrant $B^s_i$ gives a map 
$$\pi^i_s:B^s_i\cap\U_i\to D_i^s\backslash\{x^i\}$$ 
which is a homeomorphism onto its image and takes points in $w_i^s$ into $v_i^s$. Let us denote by $\eta_s$ to the inverse map 
$$\eta_s=(\pi^1_s)^{-1}:V\cap D^s_1\backslash\{x^1\}\to B_1^s.$$ 

Let $V\subset\pi^1(\U_1)\cup\{x^1\}$ be a neighborhood of $x^1$. We start by constructing $h_D$ in each quadrant $V\cap D_1^s$. For every $s=1,\dots,4$ let 
$h_D^s:V\cap D_1^s\to D_2^s$
be defined in the following way:
\begin{equation}
h_D^s(x)=\left\{\begin{array}{lcc}
              \pi^2_s\circ h\circ\eta_s(x)  & \text{if} & x\neq x^1,\\
            \\x^2                           & \text{if} & x=x^1.\\
              \end{array}
     \right.
\end{equation}
We claim that each $h_D^s$ is a homeomorphism onto its image. Observe that $\eta_s$ takes points in $D_1^s\backslash\{x^1\}$ into the quadrant $B_1^s$, $h$ takes points in $B_1^s$ into the quadrant $B_2^s$ and then $\pi^2$ projects $B_2^s$ into the punctured quadrant $D_2^s\backslash\{x^2\}$. Since each map is a homeomorphism onto its image then $h_D^s$ takes $V\cap D_1^s\backslash\{x^1\}$ homeomorphically onto its image in $D_2^s\backslash\{x^2\}$. Since the projections along the flow send points near $\gamma^i$ in the Birkhoff section to points near $x^i$ we see that $h_D^s$ is continuous in $x^1$ and the claim follows.

We define now 
$h_D:V\to D_2\text{ such that }$ 
\begin{equation}
h_D(x)=h_D^s(x)\text{ if }x\in D_1^s,\ s=1,\dots,4.
\end{equation}
We will show that $h$ is a well-defined map, that is a homeomorphism onto its image and conjugates the first return maps $P_{D_i}$ for points close to $x^i$. To see that $h_D$ is well-defined we have to check the definition of $h_D$ over the boundaries $v_1^s$, $s=1,\dots,4$ of the quadrants $D_1^s$. If $x$ belongs to some segment $v_1^s$ for $s=2,3,4$ then we have that $\eta_{s-1}(x)=\eta_s(x)$. This is because the maps $\pi^i_{s-1}$ and $\pi^i_s$ coincides over the segment $w_i^s$ which separates the quadrants $B_i^{s-1}$ and $B_i^s$. It follows that 
\begin{equation*}
h_D^s(x)=\pi^2_s\circ h\circ\eta_s(x)=\pi^2_{s-1}\circ h\circ\eta_{s-1}(x)=h_D^{s-1}(x),\text{ for every }x\in v_1^s,
\end{equation*}
so the map $h_D$ is well-defined in the segments $v_1^2$, $v_1^3$ and $v_1^4$. 

We use now Proposition \ref{prop_combinatorial_relation_projections_first_return} to show that $h_D$ is well-defined and continuous over the segment $v_1^1$. By this proposition we have that $\pi^i_4(z)=P_{D_i}^{-1}\circ\pi^i_1(z)=\pi^i_1\circ P_{B_i}^{-1}(z)$, for every $z\in w_i^1$. So it follows that $\eta_4(x)=P_{B_i}\circ\eta_1(x)$ for every $x\in v_1^1$, as we have illustrated in Figure \ref{fig_prop_fund}. Using that $h$ conjugates the maps $P_{B_i}$ we have
\begin{align}\label{continui_v1}
h_D^4(x)&=\pi_4^2\circ h\circ\eta_4(x)=\pi_4^2\circ h\circ P_{B_1}\circ\eta_1(x)\\
		&=\pi_4^2\circ P_{B_2}\circ h\circ\eta_1(x)=\pi_1^2\circ h\circ\eta_1(x)=h_D^1(x).
\end{align}
Since $h_D^4(x)=h_D^1(x)$ for every $x\in v_1^1$, we conclude that $h_D:V\to D_2$ is well-defined, and is a homeomorphism onto its image. 

To finish the proof, consider a disk $D_1'\subset V$ that contains $x^1$ and define $D_2'=h(D_1')$. We have to show that the homeomorphism $h:D_1'\to D_2'$ is a local conjugation between $P_{D_1}$ and $P_{D_2}$. Let $x\in D_1'$ be a point which belongs to some quadrant $D_1^s$. Then
\begin{align}\label{conjugat}
h_D\circ P_{D_1}(x)&=\pi^2_s\circ h\circ\eta_s\circ P_{D_1}(x)=\pi^2_s\circ h\circ P_{B_1}\circ\eta_s(x)\\
				   &=\pi^2_s\circ P_{B_2}\circ h\circ\eta_s(x)=P_{D_2}\circ\pi^2_s\circ h\circ\eta_s(x)=P_{D_2}\circ h_D(x).
\end{align}

\paragraph*{General case}\ 

The general case follows the same proof that the previous case, the only difference is when checking the continuity along the segment $v^1_s$ and the conjugation. Using Proposition \ref{prop_combinatorial_relation_projections_first_return} we can see how to modify equations \eqref{continui_v1} and \eqref{conjugat} above and obtain the desired result. 

\subsection{Proof of Theorem B (Theorem \ref{thmA-statement}).}\label{subsection-Proof_thm-A}
Take two flows $\{\phi_t^i:M_i\to M_i\}_{t\in\R}$, $i=1,2$ equipped with Birkhoff sections $\iota_i:(\Sigma_i,\partial\Sigma_i)\to (M_i,\Gamma_i)$ as in the statement of Theorem \ref{thmA-statement}. We prove Theorem \ref{thmA-statement} assuming that $\Gamma_i$ consists in exactly one periodic orbit $\gamma^i\subset M_i$, since the general case follows by repeating the same argument in a neighborhood of each periodic orbit in $\Gamma_i$. 

By hypothesis $\Sigma_1$ and $\Sigma_2$ are homeomorphic, so each surface has the same number $p>0$ of boundary components. We call them $\partial\Sigma_i=C_1^i\cup\cdots\cup C_p^i$. By blowing down each surface $\Sigma_i$ we obtain a closed surface $\widehat{\Sigma}_i$ with a set $\Delta_i=\{x_1^1,\dots,x_p^i\}$ of marked points, where each $x^i_j$ is the point obtained by collapsing the boundary component $C^i_j$, and a pseudo-Anosov homeomorphism $\hat{P}_i$ on $\widehat{\Sigma}_i$ fixing the set $\Delta_i$. Since the first return map of a Birkhoff induces a cyclic permutation of all the boundary components arriving to the same $\gamma^i\in\Gamma_i$ (Proposition \ref{prop_combinatorial_relation_projections_first_return}) then the finite set $\Delta_i$ constitutes a single periodic orbit of period $p$. 

Since the actions induced by $P_i$ on $\pi_1(\Sigma_i)$, $i=1,2$ are conjugated by the action of a homeomorphism $\Psi:\Sigma_1\to\Sigma_2$ (thus sending conjugacy classes of curves in $\partial\Sigma_1$ onto conjugacy classes of curves in $\partial\Sigma_2$) they induce conjugated elements $[P_i]$ in the mapping class group $\textsl{MCG}(\Sigma_i,\Delta_i)$ of the punctured surface. We recall that this is the group of homeomorphisms of $\Sigma_i$ fixing $\Delta_i$, up to homotopies fixing $\Delta_i$, see Section \ref{section_pseudo-Anosov_punctured_surfaces}. Since the $\hat{P}_i$ are pseudo-Anosov then, by Theorem \ref{isotopic_pseudo_Anosov} there exists a homeomorphism $h:\widehat{\Sigma}_1\to\widehat{\Sigma}_2$ such that $\hat{P}_2\circ h=h\circ\hat{P_1}$. In particular, we deduce equality on linking numbers: $n(\gamma^1,\Sigma_1)=n(\gamma^2,\Sigma_2)=n$.

To obtain the desired orbital equivalence $(\phi^1,M_1)\to(\phi^2,M_2)$, we have to consider the following situations:

\paragraph*{Case I: Assume that each $\Sigma_i$ is tame.} 
The homeomorphism $h:\mathring{\Sigma}_1\to\mathring{\Sigma}_2$ that conjugates the first return maps $P_i$ produces an orbital equivalence $H:(\phi^1,M_1\backslash\gamma^1)\to(\phi^2,M_2\backslash\gamma^2)$, which satisfies that $H(x)=h(x)$, for every $x\in\mathring{\Sigma_1}$. 
\begin{lem}\label{thm_fund_global-tame}
Under the assumption $m(\gamma^1,\Sigma_1)=m(\gamma^2,\Sigma_2)$ then for every neighborhood $W$ of $\gamma^1$ there exists a homeomorphism $H_W:M_1\to M_2$ such that: 
\begin{enumerate}[(a)]
\item $H_W$ is a topological equivalence between $(\phi^1,M_1)$ and $(\phi^2,M_2)$;
\item $H_W(x)=h(x)$, for every $x\in\Sigma_1\backslash W$.
\end{enumerate}
\end{lem}

\begin{proof}
Let $W_1=W$, $W_2:=H(W_1\backslash\gamma^1)\cup\gamma^2$ and consider the restriction $H:W_1\backslash\gamma^1\to W_2\backslash\gamma^2$. Then, this homeomorphism is an orbital equivalence between the (open) sets $W_i\backslash\gamma^i$. More over, if we call $B_1^1$ to one of the components of $\iota_1(\Sigma_1)\cap W_1$ and $B_1^2=\tilde{h}(B_1^1)$, then $H$ coincides with $h$ over $B_1\backslash\gamma^1$ and $h$ is a local conjugation between the first return maps to the local sections $B_1^i$. So we are in the hypothesis of Theorem \ref{thm_fund_local}. Given an arbitrary neighborhood $N\subset W_1$, there exists then a local orbital equivalence $H_N:W_1\to W_2$ such that $H_N(x)=H(x)$, for every $x\in W_1\backslash N$. We define $H_W:M_1\to M_2$ such that
\begin{equation*}
H_W(x)=\left\{\begin{array}{lcc}
		&H_N(x),\ &\text{ if }x\in W_1\\
		&H(x),\ &\text{ if }x\in M_1\backslash N.	
		\end{array}
		\right.
\end{equation*} 
Then $H_W$ is a well-defined homeomorphism, and is an orbital equivalence. Observe that, since $H_W$ coincides with $H$ outside $N$, then $H_W(x)=h(x)$, for every $x\in \Sigma_1\backslash W$. This finishes the proof in the case of tame sections.
\end{proof}

\paragraph*{Case II: General case.} 
If the sections $\Sigma_i$, $i=1,2$ given in the statement of Theorem \ref{thmA-statement} are not tame, we can modify them in a neighborhood of the boundary and produce tame Birkhoff sections $\Sigma_i'$, without altering the action of the first return map on the fundamental group.

More concretely, for each $\gamma_i$ consider two tubular neighborhoods $N_i\subset W_i$ satisfying that for every $x\in\Sigma_i\backslash W_i$ then the $\phi^i$-orbit segment $[x,P_i(x)]$ connecting $x$ with its first return to $\Sigma_i$, is completely contained in $M_i\backslash N_i$. By Proposition \ref{prop_existence_tame_Birkhoff_sections} there exists a Birkhoff section $\Sigma_i'$ for $(\phi^i,M_i)$ that is tame and coincides with $\Sigma_i$ in the complement of $N_i$. Denote $\Sigma_0=\Sigma_i'\backslash W_i\equiv \Sigma_i\backslash W_i$. Denote by $P_i':\Sigma_i'\to\Sigma_i'$ the first return map. Then, it follows that $P_i'(x)=P_i(x)$, for every $x\in\Sigma_0$. Moreover, since $\Sigma_0$ is a deformation retract of both surfaces $\Sigma_i$ and $\Sigma_i'$, every homotopy class of closed curve in any of the surfaces has a representative $\alpha$ completely contained in $\Sigma_0$ and it follows that $P_i'(\alpha)=P_i(\alpha)$. Thus, $P_i$ and $P_i'$ induce equivalent actions on $\pi_1(\Sigma_i)$ and $\pi_1(\Sigma_i')$, respectively.

We obtain tame Birkhoff sections satisfying the statement of Theorem \ref{thmA-statement}, and we conclude as in the previous case.

\section{Almost Anosov structures}\label{section_almost-Anosov}
Let $(\phi,M)$ be a transitive topologically Anosov flow. In the complement of some finite collections of periodic orbits, it is possible to define a smooth atlas and a Riemannian metric such that the restriction of the flow onto this set is orbitally equivalent to a smooth flow preserving a uniformly hyperbolic splitting. This follows from the existence of Birkhoff sections, Theorem \ref{thm_Fried-Brunella}. Since $\phi$ is transitive, there exists a Birkhoff section $\iota:(\Sigma,\partial\Sigma)\rightarrow (M,\Gamma).$ In the complement of $\Gamma$ the restriction of the flow is orbitally equivalent to the suspension of the first return map $P:\mathring{\Sigma}\to\mathring{\Sigma}$. Since this map is obtained by blowing-up a pseudo-Anosov homeomorphism $\widehat{P}:\widehat{\Sigma}\to\widehat{\Sigma}$ on a closed surface, it is possible to define a convenient smooth structure and Riemannian metric on the (open) manifold $M\backslash\Gamma$, as we explain below in \ref{subsection_almost-Anosov}. 
We are interested in a description of this smooth atlas in a neighborhood of each orbit $\gamma\in\Gamma$, where the atlas is not defined. This is the content of Section \ref{subsection_normal_form}.

\subsection{Almost Anosov atlas induced by a Birkhoff section}\label{subsection_almost-Anosov}
Given an arbitrary finite set $\Gamma$ of periodic orbits of a topologically Anosov flow $(\phi,M)$, we denote $M_\Gamma=M\setminus\Gamma$ and, in general, we use $\Gamma$ as a sub-index for referring to the objects associated to the restriction of $\phi$ onto $M_\Gamma$.   

\begin{defn}
\label{defn_almost_anosov_atlas}
An \emph{almost Anosov structure} associated to $(\phi,M)$ is a smooth atlas $\mathcal{D}_\Gamma$ defined in the open 3-manifold $M_\Gamma$, where $\Gamma$ is some finite set of $\phi$-periodic orbits, satisfying:
\begin{enumerate}[(i)]
\item The orbit foliation $\mathcal{O}_\Gamma$ on $M_\Gamma$ is tangent to a smooth non-singular vector field $X_\Gamma$.
\item $D\phi_t^{X_\Gamma}:TM_{\Gamma}\to TM_{\Gamma}$ preserves a splitting $TM_\Gamma=\Es_\Gamma\oplus\Ec_\Gamma\oplus\Eu_\Gamma$, where $\phi_t^{X_\Gamma}$ denotes the flow generated by $X_\Gamma$ and the bundle $\Ec_\Gamma$ is collinear with this vector field.  
\item There exists a Riemannian metric $\vert\cdot\vert_\Gamma$ on $M_\Gamma$ for which the splitting is uniformly hyperbolic. That is, there are constants $0<\lambda<1$ and $C>0$ such that: 
\begin{align*}
& \vert D\phi_t^{X_\Gamma}(p)\cdot v\vert_\Gamma\leq C\cdot\lambda^t\cdot\vert v\vert_\Gamma,\ \forall\ v\in\Es\text{ and }t\geq 0, \\
& \vert D\phi_t^{X_\Gamma}(p)\cdot v\vert_\Gamma\leq C\cdot\lambda^{-t}\cdot\vert v\vert_\Gamma,\ \forall\ v\in\Eu\text{ and }t\leq 0,
\end{align*}
\end{enumerate} 
\end{defn}

\begin{proposition}
\label{proposition_existence_almost-Anosov_structure}
If $(\phi,M)$ is a topologically Anosov flow and $\Gamma$ is the boundary of a Birkhoff section $\iota:(\Sigma,\partial\Sigma)\rightarrow (M,\Gamma)$, then there exists an almost Anosov structure $\mathcal{D}_\Gamma$ on $M_\Gamma$.   
\end{proposition}

\begin{proof}
Let $\widehat{P}:\widehat{\Sigma}\to\widehat{\Sigma}$ be the blow down of the first return to the Birkhoff section, and consider the open surface 
$\Sigma_\Delta=\widehat{\Sigma}\backslash\Delta$, where $\Delta$ is the finite set obtained by blowing down the components of $\partial\Sigma$. Since $\Sigma_\Delta$ is canonically identified with the embedded surface $\iota(\Sigma)\setminus\Gamma$, the restricted flow $(\phi,M_\Gamma)$ is orbitally equivalent to the suspension flow 
$(\mathring{\phi}_t,\mathring{M})=\textrm{suspension}\left(P:\Sigma_\Delta\to\Sigma_\Delta\right),$
meaning that there exists a homeomorphism $H:M_\Gamma\to\mathring{M}$ preserving the corresponding orbit foliations. To prove the proposition, we show that the manifold $\mathring{M}$ can be endowed with an atlas $\mathcal{D}_\circ$ and a Riemannian metric $\vert\cdot\vert_\circ$, that make the flow 
$(\mathring{\phi}_t,\mathring{M})$ uniformly hyperbolic. This implies the statement of Propositions \ref{proposition_existence_almost-Anosov_structure} for $(\phi,M_\Gamma)$, simply by pulling-back under $H$ the atlas $\mathcal{D}_{\circ}$ and the Riemannian metric $\vert\cdot\vert_\circ$ obtained for $\mathring{M}$. 

Since the $\widehat{P}:\widehat{\Sigma}\to\widehat{\Sigma}$ is pseudo-Anosov, it has an associated pair $(\Fs,\mu_\mathrm{s})$ and $(\Fu,\mu_\mathrm{u})$ of transverse foliations equipped with transverse measures and a stretching factor $0<\lambda<1$. 
Since there are no singularities on the open surface $\Sigma_\Delta$, this pair of transverse foliations provides a translation atlas 
$$\mathcal{D}_\Delta=\{\varphi_i:U_i\to\R^2\}_{i\in I},\text{ where }\{U_i\}_{i\in I}\text{ is an open cover of }\Sigma_\Delta,$$
such that, on each coordinate neighborhood, the foliations $\Fs$ and $\Fu$ correspond to the foliations of the plane by horizontal and vertical lines, and the transverse measures $\mu_\mathrm{s}$ and $\mu_\mathrm{s}$ correspond to integrate the 1-forms $|dx|$ and $|dy|$ in $\R^2$, respectively. This translation atlas defines a smooth structure in $\Sigma_\Delta$ and a Riemannian metric $|\cdot|_\Delta^2={dx^2+dy^2}$. 
In the local coordinates of this atlas, the first return $P:\Sigma_\Delta\to\Sigma_\Delta$ takes the form of a homeomorphism  
$$\varphi_i\circ P\circ\varphi_j^{-1}:(x,y)\mapsto\pm
\begin{pmatrix}
\lambda & 0\\
0 & \lambda^{-1}
\end{pmatrix}
(x,y)+\tau_{ij}$$ 
between open sets in $\R^2$. It follows that $P$ is smooth and $DP:T\Sigma_\Delta\to T\Sigma_\Delta$ preserves a splitting $T\Sigma_\Delta=\Es\oplus\Eu$ given in local coordinates by $\Es=\R\times 0$ and $\Eu=0\times\R$. Moreover, with the metric $|\cdot|_\Delta$ this splitting is uniformly hyperbolic.

Consider the suspension flow 
$(\mathring{\phi}_t,\mathring{M})=\textrm{suspension}\left(P:\Sigma_\Delta\to\Sigma_\Delta\right)$.
With the structure of smooth Riemannian manifold on $\Sigma_\Delta$ defined above, the flow
$\phi$ is the flow generated by the smooth vector field $\partial/\partial t$ on the smooth (open) manifold
$$\mathring{M}=\left(\Sigma_\Delta\times\R\right)/{(z,t)\sim(P(z),t-1)}.$$
Let us denote by $\mathcal{D}_\circ$ the smooth atlas on $\mathring{M}$, that is obtained tensoring the charts in $\mathcal{D}_\Delta$ with the standard structure of $\R$. The $DP$-invariant splitting on the surface induces a splitting of the form
$T\mathring{M}=\Es\oplus\Eu\oplus\spa\{\partial/\partial t\}$
that is invariant by $D\mathring{\phi}_t$. In addition, for each point $(z,t)$ of $\Sigma_\Delta\times\R$, the expression 
$|\cdot|_\circ^2=\lambda^{-2t}dx^2+\lambda^{2t}dy^2+dt^2$
defines a Riemannian metric that pushes-down to the quotient manifold $\mathring{M}$ and induces a metric $|\cdot|_\circ$ that coincides with the metric on $\Sigma_\Delta$ along a fixed global transverse section $\Sigma_\Delta\hookrightarrow\mathring{M}$. The invariant splitting of $T\mathring{M}$ is uniformly hyperbolic with respect to this metric, and it follows from its definition that $|D\phi_t|_{\Es}|_{\circ}=\lambda^t$ and $|D\phi_t|_{\Eu}|_{\circ}=\lambda^{-t}$. 
\end{proof}

\begin{rem}
Observe that the Riemannian metric $\vert\cdot\vert_\Gamma$ constructed above does not extend onto the closed curves in $\Gamma$. To see this, fix some $\gamma_i\in\Gamma$ and let $x_{ij}\in\Delta$ be one of the points in $\widehat{\Sigma}$ obtained after blowing down the components of $\iota^{-1}(\gamma_i)$ of $\partial\Sigma$. Given $r>0$, consider the closed loop $\sigma_r$ in $\widehat{\Sigma}$ that is the boundary of the disk of radius $r$ centered at $x_{ij}$, in the (singular) Euclidean metric $|\cdot|_\Delta^2={ds^2+du^2}$ defined above. Then, observe that the length of $\sigma_r$ goes to zero when $r\to 0$. However, the inclusion $\iota(\sigma_r)$ gives a family of simple closed curves in $M$ converging uniformly to $\gamma_i$ when $r\to 0$, and this implies that $\gamma_i$ must have length zero with respect to any extension of $\vert\cdot\vert_\Gamma$ onto $\Gamma$. 
\end{rem}

The almost Anosov atlases constructed upon Birkhoff sections in Proposition \ref{proposition_existence_almost-Anosov_structure} have an additional property, namely, they are \emph{affine} atlases and the action of the flow is by \emph{affine transformations}. This can be seen directly from the fact that the suspension flow associated to $P:\Sigma_\Delta\to\Sigma_\Delta$ is the quotient of the 1-parameter family $\tilde{\phi}_t(p,s)=(p,s+t)$ acting on $\Sigma_\Delta\times\R$, under the affine discrete group generated by $(p,s)\mapsto(P(p),s-1)$. 

\begin{question}
\label{question_affine_almost_anosov_structures}
Suppose that, in the complement of a finite set $\Gamma$ of periodic orbits, there is an almost Anosov structure $\mathcal{D}_\Gamma$ that in addition is affine, such that the flow acts by affine transformations. Is it true that $\Gamma$ bounds a Birkhoff section?
\end{question}

\subsection{Normal form in a neighborhood of $\gamma\in\Gamma$}\label{subsection_normal_form}
We describe here a normal form for the vector field $X_\Gamma$ given in Proposition \ref{proposition_existence_almost-Anosov_structure}, in a neighborhood of each $\gamma\in\Gamma$. Fix some periodic orbit $\gamma\in\Gamma$ and consider some small tubular neighborhood $W$. We assume that the invariant local manifolds are orientable, so every small tubular neighborhood $W$ is partitioned in four quadrants $W_i$, $i=1,\dots,4$. The following proposition gives a normal form for the vector field $X_\Gamma$ on the punctured neighborhood $W\setminus\gamma$. The normal form is  constructed by gluing along the boundaries (in a non-trivial way) the four quadrants of a saddle type periodic orbit, generated by an affine vector field in $\R^2\times\R/\Z$. 

Recall that the first return to $\Sigma$ is pseudo-Anosov and we denote by $0<\lambda<1$ its stretching factor. Also recall that $\gamma$ is the image of a number $p=p(\gamma,\Sigma)$ of connected components of $\partial\Sigma$. This means that in a small neighborhood $W$ of $\gamma$ the surface $\Sigma\cap W$ splits as $p$ different local Birkhoff sections at $\gamma$, each of them with the same linking number $n=n(\gamma,\Sigma)$ and multiplicity $m=m(\gamma,\Sigma)$. 

\begin{proposition}\label{lema_normal_coordinates}
Let $(\phi,M)$ be transitive topologically Anosov flow on a closed orientable 3-manifold and let $\iota:(\Sigma,\partial\Sigma)\to(M,\Gamma)$ be a Birkhoff section. Consider the smooth atlas $\mathcal{D}_\Gamma$ and the smooth vector field $X_\Gamma$ on $M\backslash\Gamma$ induced by the Birkhoff section (cf. Proposition \ref{proposition_existence_almost-Anosov_structure}). Then, for every $\gamma\in\Gamma$ there exists a small tubular neighborhood $W$, divided in four quadrants $W_i$, $i=1,\dots,4$, and a systems of smooth charts 
\begin{equation}\label{cartas_Pi-prop_4.3}
\Pi_i:W_i\setminus\gamma\to\left(\D_i\setminus\{0\}\right)\times\R/\Z,\ i=1,\dots,4 
\end{equation}
where $\D_i$, $i=1,\dots,4$ are the four quadrant of the plane $\R^2$ obtained by splitting along the vertical and horizontal axis, ordered in counterclockwise fashion, satisfying that:
\begin{enumerate}
\item $D\Pi_*(X_\Gamma)=X_{(\lambda,p,n)}$, where $n=n(\gamma,\Sigma)$, $p=p(\gamma,\Sigma)$ and $X_{(\lambda,p,n)}:\R^2\times\R/\Z\to\R^3$ is the vector field 
\begin{equation}\label{normal_form_X_Gamma}
X_{(\lambda,p,n)}(x,y,z)=\left(\log(\lambda)x,-\log(\lambda)y,\frac{1}{np}\right).
\end{equation}

\item The charts $\Pi_i$ send each connected component of $\Sigma\cap (W_i\backslash\gamma)$ isometrically onto a surface of the form $\left(U\times 0\right)\cap\left(\D_i\backslash \{0\}\right)\times\left\{\frac{k}{np}\right\}$, where $U$ is an open neighborhood of $0\in\R^2$, $k\in\{0,\dots,np-1\}$. 

\item It is verified that:
\begin{align*}
& \Pi_1(W_1)\cap\Pi_2(W_2)=\{0\}\times(0,+\infty)\times\R/\Z\\
& \Pi_2(W_2)\cap\Pi_3(W_3)=(-\infty,0)\times\{0\}\times\R/\Z\\
& \Pi_3(W_3)\cap\Pi_4(W_4)=\{0\}\times(-\infty,0)\times\R/\Z\\
& \Pi_4(W_4)\cap\Pi_1(W_1)=(0,+\infty)\times\{0\}\times\R/\Z\\
\end{align*}
and
\begin{align}\label{combinatorics_charts_Pi}
& \Pi_2\circ\Pi_1^{-1}:(0,y,z)\mapsto(0,y,z)\nonumber\\
& \Pi_3\circ\Pi_2^{-1}:(x,0,z)\mapsto(x,0,z) \\
& \Pi_4\circ\Pi_3^{-1}:(0,y,z)\mapsto(0,y,z)\nonumber \\
& \Pi_1\circ\Pi_4^{-1}:(x,0,z)\mapsto\left(x,0,z+\frac{m}{n}\right)\nonumber.
\end{align}
\end{enumerate}
\end{proposition}

\begin{rem}\label{remark_cross-shaped}
The charts defined in Proposition \ref{lema_normal_coordinates} send the orbit segments of $\phi$ that lie inside each punctured quadrant $W_i\setminus\gamma$ onto the orbit segments of the flow generated by the vector field \eqref{normal_form_X_Gamma} inside the quadrant $\D_i\times\R/\Z$ \emph{preserving the time parameter}. Thus, we can reconstruct the vector field $X_\Gamma$ in $W\backslash\gamma$ by gluing the four pieces 
$$\left(X_{(\lambda,p,n)}\ ,(\D_i\backslash \{0\})\times\R/\Z \right),\ i=1,\dots,4$$
along their boundaries in the way specified in (\ref{combinatorics_charts_Pi}). Since the vector field $X_{(\lambda,p,n)}$ is invariant by vertical translations, the gluing map $\Pi_1\circ\Pi_4^{-1}$ preserves $X_{(\lambda,p,n)}$ and we get a well-defined vector field in the quotient manifold, which is homeomorphic to a solid torus with an essential closed curve on the interior removed.
\end{rem}

\begin{proof}[Proof of Proposition \ref{lema_normal_coordinates}]
We assume first that $p(\gamma,\Sigma)=1$. Choose a small tubular neighborhood $W$ of $\gamma$ such that $B=\Sigma\cap W$ is a connected local Birkhoff section, equipped with a first return map $P_B:U\backslash\gamma\to B\backslash\gamma$ (cf. Definition \ref{defn_local_Birkhoff_section} for notations). Denote by $\widehat{B}$ the disk obtained by blowing down $B$ along $\gamma$ and let $\widehat{P}_B:\widehat{U}\to\widehat{B}$ be the corresponding map. Since $\widehat{P}_B$ is a local homeomorphism with a fixed point $q\in\widehat{U}$ of multi-saddle type (with $2n(\Sigma,\gamma)$ prongs), by suspending we obtain a germ $(\widehat{\phi},\widehat{W})_{\widehat{\gamma}}$ of a flow $\widehat{\phi}$ around a multi-saddle periodic orbit $\widehat{\gamma}\subset\widehat{W}$, where $\widehat{W}\simeq\D^2\times\s$ (cf. Section \ref{section_preliminaries}). 
Up to shrinking the neighborhoods $W$ and $\widehat{W}$ if necessary, there is a local orbital equivalence 
$$\left(\phi,W\backslash\gamma\right)\to\left(\widehat{\phi},\widehat{W}\backslash\widehat{\gamma}\right)$$
between the respective complements of $\gamma$ and $\widehat{\gamma}$, that takes the surface $B\backslash\gamma$ onto $\widehat{B}\backslash\{q\}$. In particular, observe that the local invariant manifolds of $\widehat{\gamma}$ partition $\widehat{W}\backslash\widehat{\gamma}$ into four quadrants $\widehat{W}_i$, $i=1,\dots,4$. The smooth atlas and the metric in $W\backslash\gamma$ given in Proposition \ref{proposition_existence_almost-Anosov_structure} are defined by pulling back those on $\widehat{W}\backslash\widehat{\gamma}$. Thus, it suffices to check \ref{lema_normal_coordinates} using $\left(\widehat{\phi},\widehat{W}\backslash\widehat{\gamma}\right)$.

Following the description in Section \ref{section_Birk_sect_and_equiv}, the disk $\widehat{B}$ is divided into $4n$ quadrants, that we enumerate as $\widehat{B}_1,\dots,\widehat{B}_{4n}$ following the positive circular order. On each $\widehat{W}_i$, $i=1,\dots,4$ there are contained exactly $n$ quadrants of the form $\widehat{B}_{i+4j}$, $j=0,\dots,n-1$. All of these quadrants share the common vertex $q$. 
Consider the sets $\widehat{B}_{i+4j}^*=\widehat{B}_{i+4j}\setminus\{q\},$ that we call \emph{punctured} quadrants. These are a collection of pairwise disjoint surfaces, properly embedded in the manifold $\widehat{W}_i\backslash\widehat{\gamma}$, each one homeomorphic to a disk. By Proposition \ref{prop_combinatorial_relation_projections_first_return} then $\widehat{P}$ permutes cyclically these surfaces and $\widehat{P}^n$ preserves each of them. So each $\widehat{B}_{i+4j}^*$ is a transverse section for $\widehat{\phi}_t$ on $\widehat{W}_i\setminus\widehat{\gamma}$, with first return defined for points near $q$ and return time constant and equal to $n$. It follows that there is a smooth time-preserving conjugacy between $\left(\widehat{\phi}_t\ ,\widehat{W}_i\setminus\widehat{\gamma}\right)$ and the flow induced by $z\mapsto z+t$ in 
\begin{equation}\label{suspension_quadrante_pinchado}
\widehat{B}_{i+4j}^*\times\R/_{(z,t)\mapsto(\widehat{P}^n(z),t-|n|)}.
\end{equation}

In addition, since $\widehat{P}_B$ is the restriction of a pseudo-Anosov homeomorphism $\widehat{P}:\widehat{\Sigma}\to\widehat{\Sigma}$ with stretching factor $\lambda$ and $q$ is a $2n$-prong, there exists a map $\varphi_q:(\widehat{B},q)\to(\R^2,0)$ that is smooth outside $\{q\}$, obtained by integrating the pair of transverse measures $(ds,du)$ on $\widehat{B}$ along a path with one extremity on $q$, that satisfies:
\begin{enumerate}
\item $\varphi_q:\widehat{B}\setminus\{q\}\to\R^2\setminus\{0\}$ is a $n$-fold covering,

\item $A\circ\varphi_q(p)=\varphi_q\circ\widehat{P}(p)$, for every $p$ in a neighborhood of $q$, where 
$A=
\begin{pmatrix}
\lambda & 0\\
0 & \lambda^{-1}
\end{pmatrix}
$.

\item For every $i=1,\dots,4$,  the chart $\varphi_q$ sends $\widehat{B}_{i+4j}$ isometrically onto its image in $\D_i\subset\R^2$.
\end{enumerate}
In particular, the chart $\varphi_q$ provides a smooth conjugation quadrant by quadrant:
$$
\begin{tikzcd}
\widehat{B}_{i+4j}\arrow[r, "\widehat{P}^n"]\arrow[d, "\varphi_q"] & \widehat{B}_{i+4j}\arrow[d, "\varphi_q"] \\
\D_i\arrow[r, "A^n"]     & \D_i
\end{tikzcd}\ \text{for every}\ i=1,\dots,4.
$$
and therefore the flow in (\ref{suspension_quadrante_pinchado}) above is conjugated to the flow generated by $z\mapsto z+t$ in 
\begin{equation}
\D_i\times\R/_{(z,t)\mapsto(A^n(z),t-|n|)}.
\end{equation}
 
Since this latter flow is smoothly time-preserving equivalent to the germ of the vector field $X_{(\lambda,n)}(x,y,z)=(\log(\lambda)x,-\log(\lambda)y,1/n)$ on the quadrant $\D_i\times\R/\Z$, we obtain a family of charts
$$\Pi_i^j:\widehat{W}_i\backslash\widehat{\gamma}\to(\D_i\backslash\{0\})\times\R/\Z,\ j=0,\dots,n-1,$$
satisfying the desired properties. Each chart is induced from one of the maps $\widehat{P}^n:\widehat{B}^*_{i+4j}\to\widehat{B}^*_{i+4j}$.

We may assume that the coordinate charts $\Pi_i$ send each surface $\widehat{B}^*_{i+4j}$ inside a plane of the form $\D_i\times\{j/n\}$, for every $j=0,\dots,n-1$. Chose $\widehat{B}^*_1$ in the quadrant $\widehat{W}_1$ and $\Pi_1$ such that $\Pi_1:\widehat{B}^*_1\mapsto\D_1\times\{0\}$. Observe that the chart $\varphi_q$ extends over the union 
$$\widehat{B}^*_1\cup\widehat{B}^*_2\cup\widehat{B}^*_3\cup\widehat{B}^*_4$$
sending each $\widehat{B}^*_i$ onto $\D_i\times\{0\}$. Thus, we can coherently extend $\Pi_1$ to the adjacent quadrants $\widehat{W}_i\setminus\widehat{\gamma}$, $i=2,3,4$, in such a way that each $\Pi_i$ sends the surface $\widehat{B}^*_i$ to $\D_i\times\{0\}$. The quadrant $\widehat{B}^*_4$ contained in $\widehat{W}_4$ is adjacent to some quadrant $\widehat{B}^*_5$ in $\widehat{W}_1$. By Proposition \ref{prop_combinatorial_relation_projections_first_return} we see that $\widehat{B}^*_5$ coincides with $\D_1\times \{m/n\}$, provided $\Pi_1(\widehat{B}^*_1)=\D_1\times\{0\}$. Thus, the \emph{holonomy defect} of $\Pi_4\circ\Pi_1^{-1}$ on the common boundary $\widehat{W}_4\cap\widehat{W}_1$ is a vertical translation by $\frac{m}{n}$. 

This completes the proof of Proposition \ref{lema_normal_coordinates} assuming that $p(\gamma,\Sigma)=1$. If there are more than one boundary components of $\Sigma$ that cover $\gamma$, the argument is the same but the first return map to each $\widehat{B}^*_{i+4j}$ changes by $\widehat{P}^{p(\gamma,\Sigma)\cdot n}$, so we must modify the parameters of $X_{(n,\lambda)}$ in the appropriate way.
\end{proof}


\section{Construction of the hyperbolic models}\label{section_hyperbolic_models}
Consider a transitive topologically Anosov flow $(\phi,M)$. In this section we construct another flow $(\psi,N)$ which is smooth, Anosov, volume-preserving and orbitally equivalent to the first one. We start by sketching the general argument, and we develop the steps in the following subsections. 

\subsection{Proof of Theorem \textbf{A}}\label{subsection_proof_thm_A}
Choose a tame Birkhoff section $\iota:(\Sigma,\partial\Sigma)\to(M,\Gamma)$ with first return $P:\mathring{\Sigma}\to\mathring{\Sigma}$ for the flow $(\phi,M)$, as given in Theorem \ref{thm_Fried-Brunella} and Proposition \ref{prop_existence_tame_Birkhoff_sections}. Then, the constructions in propositions \ref{proposition_existence_almost-Anosov_structure} and \ref{lema_normal_coordinates} associates: 

\begin{enumerate} 

\item A smooth atlas $\mathcal{D}_\Gamma$ and a Riemannian metric $|\cdot|_\Gamma$ on the manifold $M\setminus\Gamma$ such that, up to reparametrization, the restriction of $\phi$ onto $M\backslash\Gamma$ is generated by a smooth vector field $X_\Gamma$ and preserves a uniformly hyperbolic splitting $\Es_\Gamma\oplus\Ec_\Gamma\oplus\Eu_\Gamma$. Moreover,
\begin{align*}
&|D{\phi}_t(v)|_\Gamma=\lambda^t|v|_\Gamma,\ \forall\ v\in\Es_\Gamma,\ t\geq 0,\\
&|D{\phi}_t(v)|_\Gamma=\lambda^{-t}|v|_\Gamma,\ \forall\ v\in\Eu_\Gamma\ t\leq 0,
\end{align*} 
where $0<\lambda<1$ is the stretching factor of the first return map to the Birkhoff section.

\item For every orbit $\gamma\in\Gamma$ there is a small tubular neighborhood $W$, divided in four quadrants $W_i$, $i=1,\dots,4$, and a system of smooth charts
\begin{equation}\label{normal_coordinates_section_5}
\Pi_i:W_i\setminus\gamma\to (\D_i\setminus\{0\})\times\R/\Z,\ i=1,\dots,4
\end{equation}
called \emph{normal coordinates}, verifying that:
\begin{enumerate}
\item $D\Pi_i:X_\Gamma\mapsto X_{(\lambda,n,p)}$, $D\Pi_i:\Es_\Gamma\mapsto\R\times 0\times 0$ and $D\Pi_i:\Eu_\Gamma\mapsto 0\times\R\times 0$, where 
\begin{align}\label{normal_form_section_5}
X_{(\lambda,n,p)}(x,y,z)=\left(\ \log(\lambda)x\ ,\ -\log(\lambda)y\ ,\ \frac{1}{np}\ \right)\ \ \text{on}\ \ \ \R^2\times\R/\Z ;
\end{align}

\item The charts $\Pi_i$ send the quadrants of the local Birkhoff section $\Sigma\cap W\backslash\gamma$ isometrically onto surfaces of the form $\left(\D_i\setminus\{0\}\right)\times \left\{\frac{k}{np}\right\}$ with $k=0,\dots,np-1$;

\item Along the corresponding domains of intersection between the four quadrants, we have
\begin{align*}
&\Pi_2\circ\Pi_1^{-1}:(x,y,z)\mapsto(x,y,z)\\
&\Pi_3\circ\Pi_2^{-1}:(x,y,z)\mapsto(x,y,z)\\
&\Pi_4\circ\Pi_3^{-1}:(x,y,z)\mapsto(x,y,z)\\
&\Pi_1\circ\Pi_4^{-1}:(x,y,z)\mapsto\left(x,y,z+\frac{m}{n}\right);
\end{align*}
\end{enumerate}
where $p=p(\gamma,\Sigma)$, $n=n(\gamma,\Sigma)$ and $m=m(\gamma,\Sigma)$ are the combinatorial parameters of the Birkhoff section at $\gamma$ (see Section \ref{section_Birk_sect_and_equiv} for definitions).

\end{enumerate}

\noindent
\textbf{Assumption:} 
For simplicity, from now on and for the rest of the section we assume that $\Gamma$ consists in only one periodic orbit $\gamma$.  We assume as well that the local invariant manifolds of $\gamma$ are orientable (cf. Remark \ref{remark_existence-orientable_loc_inv_mflds}). The general case can be derived from the present case, by applying the following construction on a neighborhood of each curve $\gamma\in\Gamma$.

From the one side, using the charts given above in \eqref{normal_coordinates_section_5} we define a family $R(r_1,r_2)\subset M$ of tubular neighborhoods of $\gamma$, depending on two parameters $0<r_2<r_1$. From the other side, we consider the affine vector field $X_{(\lambda,n,p)}$ given in equation \eqref{normal_form_section_5} above, defined for all points in $\R^2\times\R/\Z$, and we define a family $\mathbb{V}(r_1,r_2)\subset\R^2\times\R/\Z$ of tubular neighborhoods of the saddle type periodic orbit $\gamma_0=\{0\}\times\R/\Z$, depending on two parameters $0<r_2<r_1$. This is the content of Section \ref{subsection_cross-shaped_neighbourhoods}. 

Now, for some fixed parameters $0<r_2<r_1<1$, consider the smooth manifold $M_R(r_1,r_2)$ obtained from $M$ by removing the interior of $R(r_1,r_2)$, equipped with the restriction of the vector field $X_\Gamma$. The construction is the following:
\begin{itemize}
\item[\textbf{(1)}] Using the system of normal coordinates $\{\Pi_i:i=1,\dots,4\}$ we define a diffeomorphism 
$$\varphi:\partial M_R(r_1,r_2)\to\partial \mathbb{V}(r_1,r_2)$$
from the boundary of $M_R(r_1,r_2)$ to the boundary of $\mathbb{V}(r_1,r_2)$, that depends on the signature of the multiplicity $m(\gamma,\Sigma)$. Gluing along the boundaries with $\varphi$ produces a closed manifold 
$$N=N(r_1,r_2)\coloneqq M_R(r_1,r_2)\sqcup_\varphi\mathbb{V}(r_1,r_2)$$ 
endowed with a smooth atlas. The vector fields $X_\Gamma$ on $M_R(r_1,r_2)$ and $X_{(\lambda,n,p)}$ on $\mathbb{V}(r_1,r_2)$ match together along the boundary and induce a smooth vector field $Y$ on $N$. In addition, the flow associated to $Y$ preserves a smooth volume form. This is the content of Section \ref{subsection_construction_smooth_model}. 

\item[\textbf{(2)}] Let $(\psi,N)$ be the flow generated by the vector field $Y$. In Section \ref{subsection_smooth_model_is_Anosov} we show that, for sufficiently small values of the parameters $0<r_2<r_1<1$, this flow is Anosov. This is basically a consequence of the hyperbolicity of $X_\Gamma$ away from $\gamma$ and the hyperbolicity of $X_{(\lambda,n,p)}$ near $\gamma_0$, and the argument is carried out using the so called \emph{cone-field criterion}. Remark that, since $\psi$ is smooth Anosov and preserves a smooth volume form, then it is ergodic (cf. \cite{katok-hasselblatt}).  

\item[\textbf{(3)}] Finally, in Section \ref{subsection_smooth_model_is_equivalent}, we show that $(\psi,N)$ is orbitally equivalent to the original topologically Anosov flow $(\phi,M)$. In order to prove this, we show that both flows can be equipped with adequate Birkhoff sections, where we can check the criterion in Theorem \textbf{B}.  
\end{itemize}

Theorem \textbf{A} follows then from the statements (1),(2) and (3) above.

\subsection{Cross-shaped neighborhoods}\label{subsection_cross-shaped_neighbourhoods} 
We describe here a class of compact tubular neighborhoods of a saddle type periodic orbit, that we call \emph{cross-shaped neighborhood} and are  used in the course of the proof of Theorem \textbf{A}. We do it for both a regular periodic orbit of an affine vector field in euclidean space, and for the singular orbit $\gamma\in\Gamma$ using the charts in \eqref{normal_coordinates_section_5} above.

\subsubsection{Affine model on $\R^2\times\R/\Z$}\label{subsection_linear_local_model}
Given $0<\lambda<1$ and two integers $n,p\geq 1$, consider the vector field $X=X_{(\lambda,n,p)}$ on $\R^2\times\R/\Z$ defined by 
$X_{(\lambda,n,p)}(x,y,z)=\left(\log(\lambda)x,-\log(\lambda)y,\frac{1}{np}\right).$
The corresponding flow is given by $\phi_t^{X_{(\lambda,n,p)}}(x,y,z)=(\lambda^tx,\lambda^{-t}y,z+t/np)$, so it splits as a linear flow on $\R^2$ times a translation flow on $\R/\Z$. 
\begin{rem}\label{remark_properties_of_affine_X}
It follows that:
\begin{enumerate}
\item The plane $\R^2\times\{0\}$ is a global transverse section, the first return map is the linear transformation $(x,y)\mapsto(\lambda^{np}x,\lambda^{-np}y)$, and the returning time is constant and equal to $np$. 
\item The non-wandering set consists in one saddle type, hyperbolic, periodic orbit 
$\gamma_0$, that coincides point-wise with the set $\{0\}\times\R/\Z$. The stable and unstable manifolds are the cylinders
\begin{align*}
&\Ws(\gamma_0)=\R\times\{0\}\times\R/\Z\\
&\Wu(\gamma_0)=\{0\}\times\R\times\R/\Z
\end{align*}
\item The action of $D\phi_t^{X_{(\lambda,n,p)}}$ on $\R^3$ preserves a splitting $\Es\oplus\Ec\oplus\Eu$ into three line bundles
\begin{align*}
& \Es=\R\times \{0\}\times \{0\},\ 
 \Eu=\{0\}\times \R\times\{0\}\ \text{and}\  
 \Ec=\spa\{{X_{(\lambda,n,p)}}\},
\end{align*}
and for every $p\in\R^2\times\R/\Z$ we have that
\begin{align*}
& \Vert D\phi_t^{X_{(\lambda,n,p)}}(p)\cdot v\Vert=\lambda^t\cdot\Vert v\Vert,\ \forall\ v\in\Es\text{ and }t\geq 0, \\
& \Vert D\phi_t^{X_{(\lambda,n,p)}}(p)\cdot v\Vert=\lambda^{-t}\cdot\Vert v\Vert,\ \forall\ v\in\Eu\text{ and }t\leq 0,
\end{align*}
where $\Vert\cdot\Vert$ is the standard Euclidean norm on $\R^2\times\R/\Z$. 
\end{enumerate}
\end{rem}

\paragraph*{Cross-shaped neighborhood}
Start with the standard partition in four quadrants $\{\D_i:i=1,\dots,4\}$ of the plane $\R^2$ along the vertical and horizontal axes. Given two real numbers $0<r_2<r_1<1$ consider the region $Q_1=Q_1(r_1,r_2)\subset \D_1$ delimited by the segments:
\begin{alignat*}{2}
&(1)\ \ w_1^\mathrm{s}     &&=[0,r_1]\times\{0\}\\
&(2)\ \ w_1^\mathrm{u}     &&=\{0\}\times[0,r_1]\\
&(3)\ \ J_\mathrm{in}^1  &&=\{r_1\}\times[0,r_2]\\
&(4)\ \ J_\mathrm{out}^1 &&=[0,r_2]\times\{r_\}\\
&(5)\ \ l_1       &&=\text{ segment of the hyperbola }xy=r_1r_2\text{ that connects }(r_1,r_2)\text{ with }(r_2,r_1).
\end{alignat*}

\begin{figure}
\centering
\begin{subfigure}[b]{0.45\textwidth}
\centering
\includegraphics[scale=0.256]{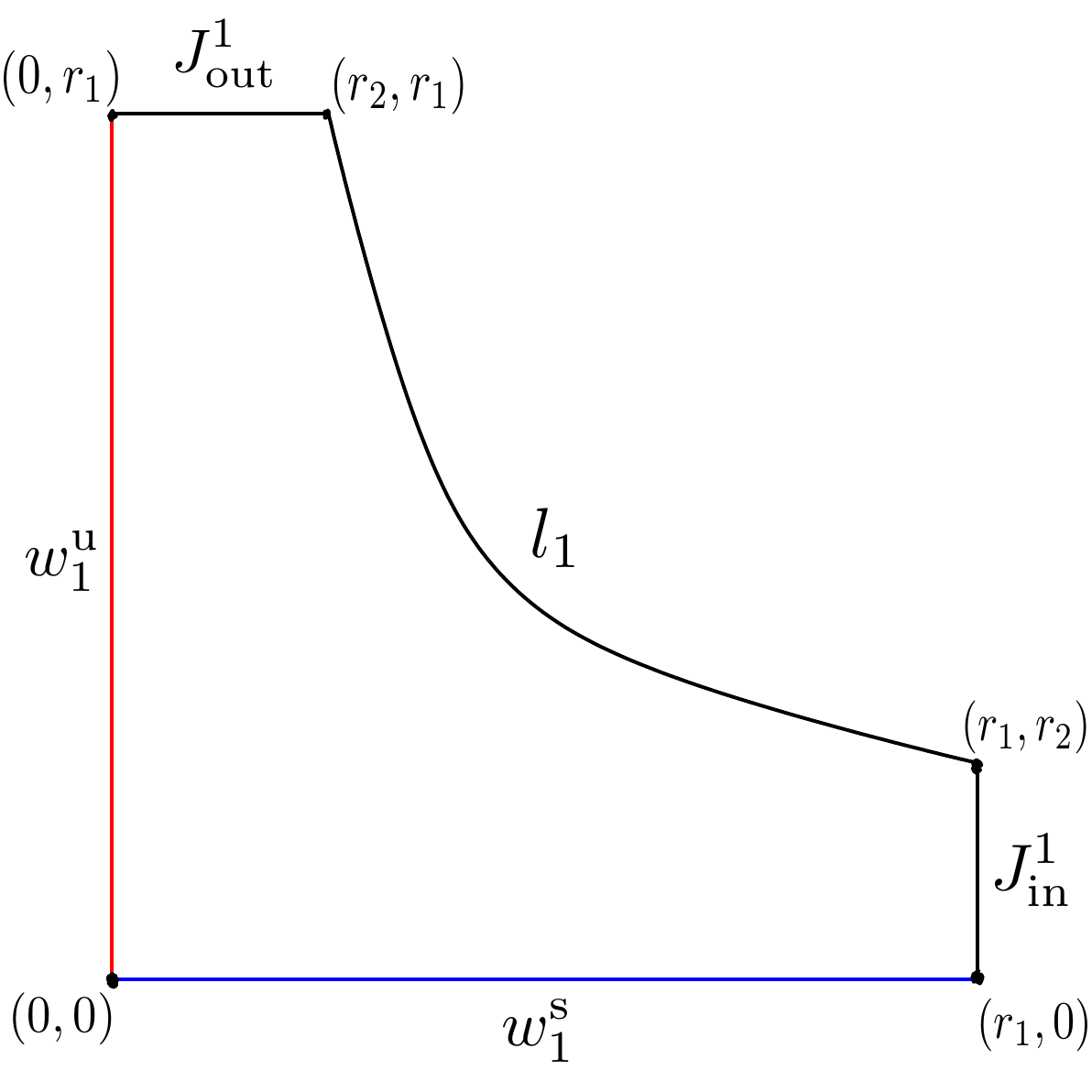}
\caption{Region $Q_1(r_1,r_2)$.}
\end{subfigure}
\hfill
\begin{subfigure}[b]{0.45\textwidth}
\centering
\includegraphics[scale=0.144]{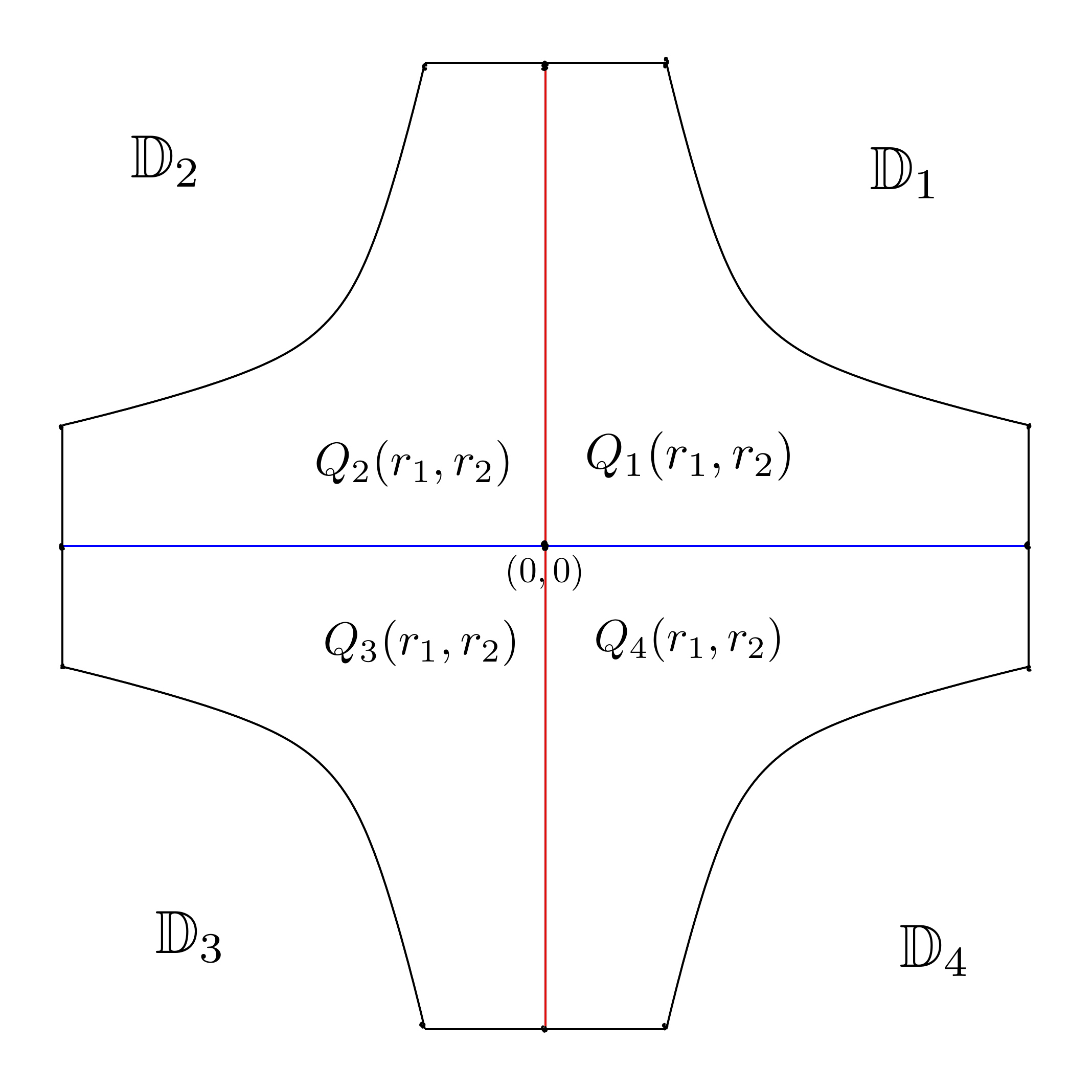}
\caption{Neighborhood $Q(r_1,r_2)$. }
\end{subfigure}
\caption{The neighborhood $Q(r_1,r_2)$}
\label{fig_region_Q}
\end{figure}

\noindent
By analogy we define the corresponding regions $Q_i\subset \D_i$ in each quadrant $i=2,3,4$. The union of these four regions determines a compact neighborhood $Q=Q(r_1,r_2)$ of $0\in\R^2$, as in Figure \ref{fig_region_Q}. 

\begin{figure}
\centering
\begin{subfigure}[b]{0.45\textwidth}
\centering
\includegraphics[scale=0.28]{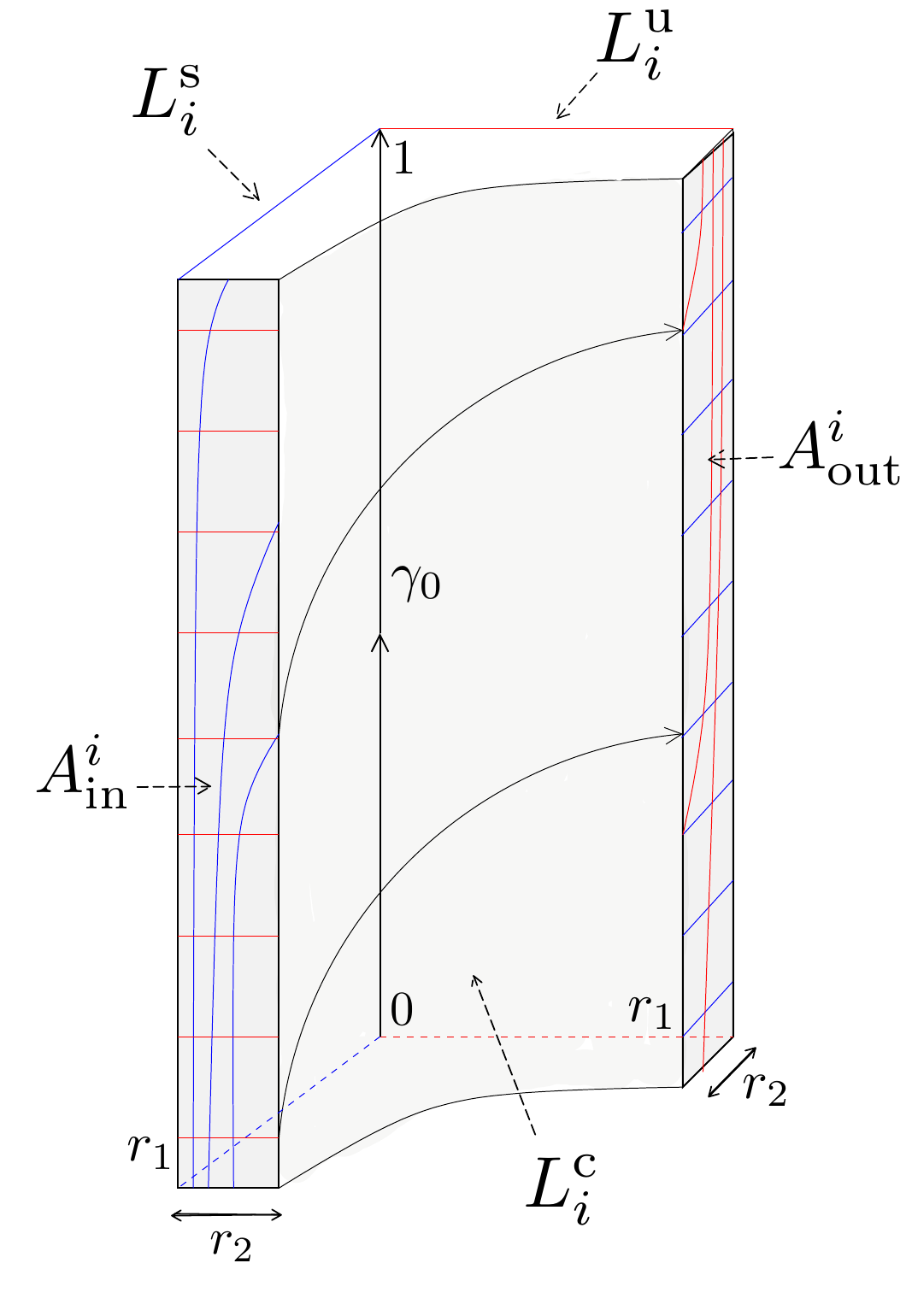}
\caption{Region $\mathbb{V}_i(r_1,r_2)$.}
\label{fig_region_Vi}
\end{subfigure}
\hfill
\begin{subfigure}[b]{0.45\textwidth}
\centering
\includegraphics[scale=0.24]{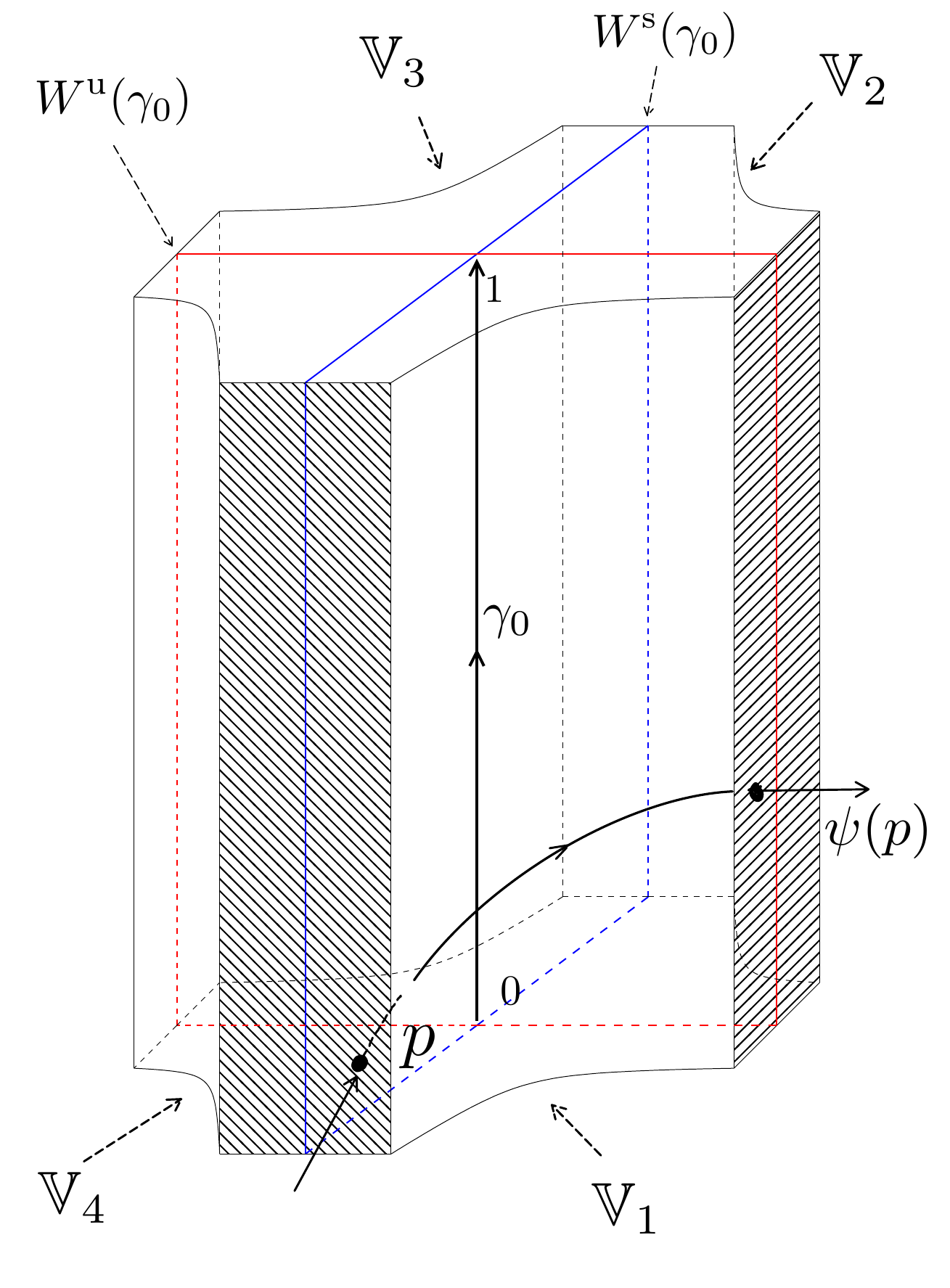}
\caption{Cross shaped neighborhood $\mathbb{V}(r_1,r_2)$. }
\label{fig_region_V}
\end{subfigure}
\caption{Cross-shaped neighborhood.}
\label{fig_cross-shaped_nbhd}
\end{figure}

Define the sets $\mathbb{V}_i=\mathbb{V}_i(r_1,r_2)\subset \D_i\times\R/\Z$ and $\mathbb{V}=\mathbb{V}(r_1,r_2)\subset\R^2\times\R/\Z$ by 
\begin{align}
& \mathbb{V}_i(r_1,r_2)\coloneqq Q_i(r_1,r_2)\times\R/\Z,\text{ for }i=1\dots,4\label{region_Vi}\\
& \mathbb{V}(r_1,r_2)\coloneqq Q(r_1,r_2)\times\R/\Z.\label{region_V}
\end{align}
The set $\mathbb{V}(r_1,r_2)$ is a compact, regular, tubular neighborhood of the periodic orbit $\gamma_0$. We call it a \emph{cross-shaped neighborhood}. It is decomposed as the union of the four regions $\mathbb{V}_i(r_1,r_2)$, $i=1,\dots,4$, glued along their stable/unstable boundaries, as in Figure \ref{fig_cross-shaped_nbhd}.

Each of $\mathbb{V}_i$ is homeomorphic to a solid torus and its boundary is composed of five annuli:
\begin{alignat*}{2}
&(1)\ L^\mathrm{s}_i     && =\ w^\mathrm{s}_i\times\R/\Z\\
&(2)\ L^\mathrm{u}_i     && =\ w^\mathrm{u}_i\times\R/\Z\\
&(3)\ A_\mathrm{in}^i  && =J_\mathrm{in}^i\times\R/\Z\\
&(4)\ A_\mathrm{out}^i && =J_\mathrm{out}^i\times\R/\Z\\
&(5)\ L^\mathrm{c}_i       && =l_i\times\R/\Z. 
\end{alignat*}
The annuli $L_i^\mathrm{s}$, $L_i^\mathrm{u}$ and $L_i^\mathrm{c}$ are tangent the vector field ${X_{(\lambda,n,p)}}$, being the first two contained in the stable and unstable manifolds of $\gamma_0$, respectively. The annuli $A_\mathrm{in}^i$ and $A_\mathrm{out}^i$ are the \emph{entrance} and \emph{exit} annuli, respectively. For $i=1$, these two annuli corresponds to the sets
\begin{align*}
A_\mathrm{in}^1&=\{(r_1,r,z):0\leq r\leq r_2,\ z\in\R/\Z\}\\
A_\mathrm{out}^1&=\{(r,r_1,z):0\leq r\leq r_2,\ z\in\R/\Z\}.
\end{align*} 
There is a diffeomorphism $\psi:A_\mathrm{in}^1\backslash\Ws(\gamma_0)\to A_\mathrm{out}^1\backslash\Wu(\gamma_0)$ of the form
$\psi:p\mapsto q=\phi^{X_{(\lambda,n,p)}}(\tau(p),p),$
which sends each entrance point $p$ onto the point $q$ determined by the intersection $A_\mathrm{out}^1\cap\mathcal{O}^+(p)$. The same holds for the other entrance-exit pair of annuli on the boundary. 

The following statement is elementary and will be used later. 

\begin{lem}
Given $0<r_2<r_1<1$, the map $\psi:\{r_1\}\times(0,r_2]\times\R/\Z\to(0,r_2]\times\{r_1\}\times\R/\Z$ defined above satisfies that: For every $p=(r_1,r,z)$, then $\psi(p)=(r,r_1,z+\tau(p)/np)$ and $\tau(p)=\tau(r)=\frac{\log(r/r_1)}{\log(\lambda)}$.
\end{lem}

\subsubsection{Singular orbit in the boundary of a Birkhoff section}\label{subsection_almost-anosov_model}
For each $\gamma\in\Gamma$, the local coordinates $\{\Pi_i:i=1,\dots,4\}$ given in \eqref{normal_coordinates_section_5} allow to construct a tubular neighborhood $R=R(r_1,r_2)$ of the orbit $\gamma$ which depends on two parameters $0<r_2<r_1<1$, in the following way:

\paragraph*{Cross-shaped neighborhood} 
For each $i=1,\dots,4$ consider the sets $\mathbb{V}_i(r_1,r_2)$ in $\R^2\times\R/\Z$ defined in \eqref{region_Vi} above. Let $0<r_2<r_1<1$ be such that each $\mathbb{V}_i(r_1,r_2)$ is contained in the image of $\Pi_i$. Define  
\begin{align}
& R_i(r_1,r_2)=\overline{\left\{p\in W_i\backslash\gamma:\Pi_i(p)\in\mathbb{V}_i(r_1,r_2)\right\}}.\label{region_R}\\ 
& R(r_1,r_2)=\bigcup_i R_i(r_1,r_2), i=1,\dots,4.\label{region_Ri}
\end{align}
We say that $R=R(r_1,r_2)$ is a \emph{cross-shaped neighborhood} of $\gamma$. 

\begin{figure}
\centering
\begin{subfigure}[b]{0.45\textwidth}
\centering
\includegraphics[scale=0.25]{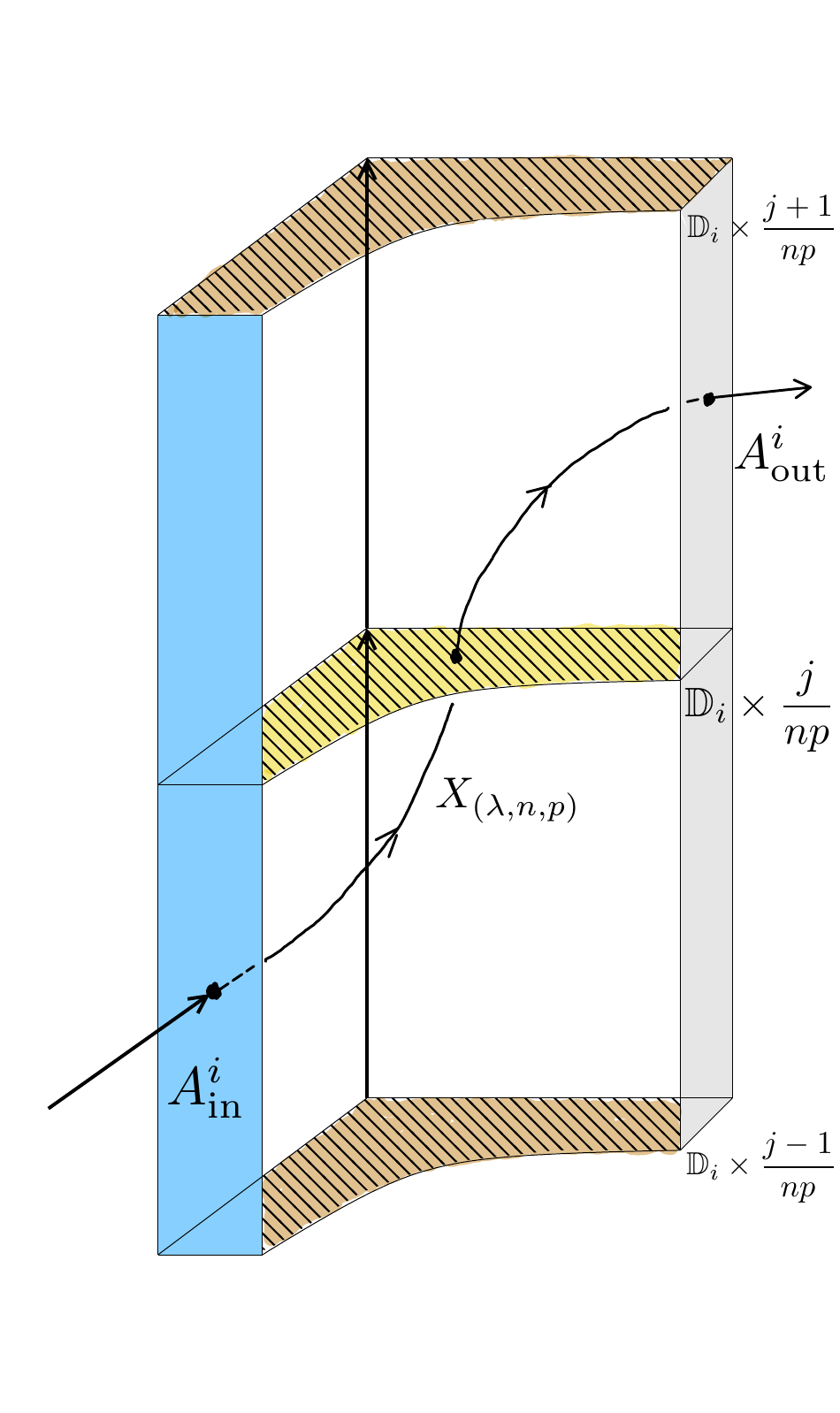}
\caption{$\mathbb{V}_i(r_1,r_2)$ and $X_{(\lambda,n,p)}$.}
\label{fig_region_Ni_Fried}
\end{subfigure}
\hfill
\begin{subfigure}[b]{0.45\textwidth}
\centering
\includegraphics[scale=0.199]{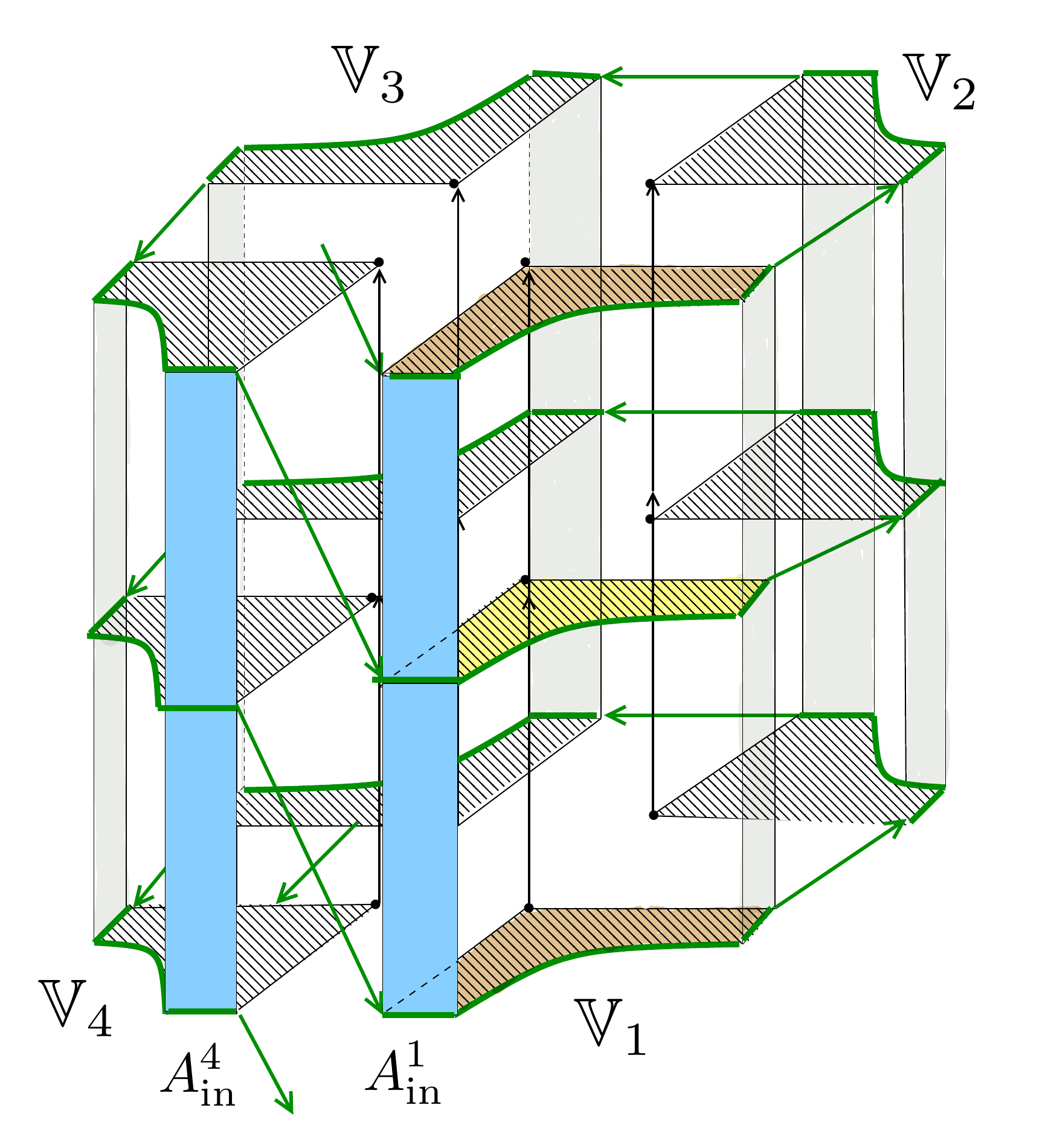}
\caption{$R(r_1,r_2)$}
\label{fig_region_N_Fried}
\end{subfigure}
\caption{The neighborhood $R(r_1,r_2)$, obtained by gluing the four pieces $\mathbb{V}_i(r_1,r_2)$.}
\label{fig_prong-shaped_nbhd}
\end{figure}

Following Remark \ref{remark_cross-shaped} in Section \ref{section_almost-Anosov}, we can describe the vector field $X_\Gamma$ on $R(r_1,r_2)\setminus\gamma$ in the following way: Consider the four pairs 
$$\left(\mathbb{V}_i(r_1,r_2),X_{(\lambda,n,p)}\right),\ i=1,\dots,4$$
consisting in the manifold $\mathbb{V}_i(r_1,r_2)$ equipped with the vector field $X_{(\lambda,n,p)}$ as in figure \eqref{fig_region_Ni_Fried}.
Then $\left(R(r_1,r_2)\setminus\gamma,X_\Gamma\right)$ is equivalent to the manifold obtained by gluing:
\begin{enumerate}
\item $\mathbb{V}_1$ with $\mathbb{V}_2$ along the boundary $L^\mathrm{u}_1\to L^\mathrm{u}_2$ with the map $(0,y,z)\mapsto(0,y,z)$,
\item $\mathbb{V}_2$ with $\mathbb{V}_3$ along the boundary $L^\mathrm{s}_2\to L^\mathrm{s}_3$ with the map $(x,0,z)\mapsto(x,0,z)$,
\item $\mathbb{V}_3$ with $\mathbb{V}_4$ along the boundary $L^\mathrm{u}_3\to L^\mathrm{u}_4$ with the map $(0,y,z)\mapsto(0,y,z)$,
\item $\mathbb{V}_4$ with $\mathbb{V}_1$ along the boundary $L^\mathrm{s}_4\to L^\mathrm{s}_1$ with the map $(0,y,z)\mapsto(0,y,z+\frac{m}{n})$.
\end{enumerate}
In figure (\ref{fig_region_N_Fried}) we see this for $m=-1$, $n=2$ and $p=1$. The associated quotient space is a solid torus. Observe that the curve $0\times\R/\Z$ is a $n$-fold covering of its image under the quotient projection. In the complement of this curve, the vector field $X_{(\lambda,n,p)}$ is invariant by $z$-translations, so we get a well-defined vector field on the quotient manifold. 

The boundary $\partial R$ is decomposed in eight smooth annuli, four of them tangent to the vector field, two where $X_\Gamma$ is transverse and points inward the neighborhood and two other where $X_\Gamma$ is transverse and points outward. Using the coordinates $\Pi_i$ we can identify the sets 
$$\Pi_i:\partial R\cap R_i\to\partial \mathbb{V}\cap\mathbb{V}_i=A_\mathrm{in}^i\cup L_i^\mathrm{c}\cup A_\mathrm{out}^i\subset\D_i\times\R/\Z.$$
In particular, the union $A_\mathrm{in}^1\cup A_\mathrm{in}^4$ is one of the eight annuli that forms $\partial R$, where the flow traverse inwardly. In coordinates \eqref{normal_coordinates_section_5}, these two annuli correspond to the sets
\begin{align*}
A_\mathrm{in}^1&=\{(r_1,r,z):0\leq r\leq r_2,\ z\in\R/\Z\}\\
A_\mathrm{in}^4&=\{(r_1,r,z):-r_2\leq r\leq 0,\ z\in\R/\Z\},
\end{align*}
glued along the boundary $\{r_1\}\times 0\times\R/\Z$ with a vertical translation by $m/n$. 

\subsection{Construction of the smooth model}\label{subsection_construction_smooth_model}
In this subsection we construct the smooth model $(\psi,N)$ associated to the topologically Anosov flow $(\phi,M)$. 

We recall the general assumption that $\Gamma$ consists in only one periodic orbit $\gamma$ and we set $n=n(\gamma,\Sigma)$, $m=m(\gamma,\Sigma)$ and $p=p(\gamma,\Sigma)$. Consider a smooth decreasing function $\rho:[0,1]\to[0,1]$ such that:
\begin{enumerate}
\item $\rho(t)=1$,  for $0\leq t\leq \frac{1}{3}$,
\item $\rho(t)=0$,  for $\frac{2}{3}\leq t\leq 1$,
\item $\rho'(t)<0$ and $\alpha(t)=|pm\log(\lambda)\rho'(t)|<\frac{1/2}{3t^2-t}$, for $\frac{1}{3}\leq t\leq \frac{2}{3}$,
\end{enumerate}

Observe that, for all choices of parameters $p,m,\log(\lambda)$, it is always possible to find $\rho$ satisfying item 3., since the function $t\mapsto\frac{1/2}{3t^2-t}$ has a pole of order $\sim 1/(t-1/3)$ at $t=1/3$. This condition is used in the proof of Lemma \ref{lema_conos_1} in Section \ref{subsection_smooth_model_is_Anosov}, in order to show that the flow is Anosov.

\begin{figure}[t]
\begin{center}
\includegraphics[width=\textwidth]{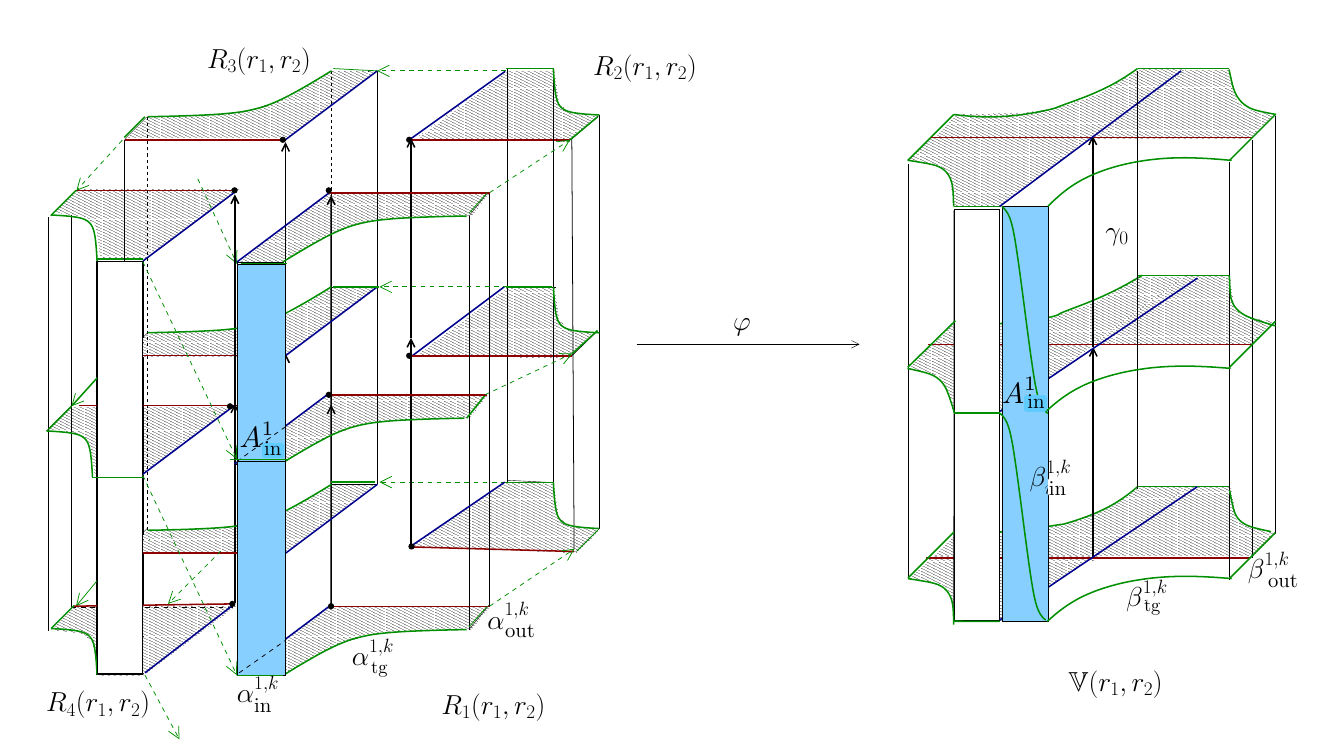}
\caption{Gluing map $\varphi:\partial M_R(r_1,r_2)\to\mathbb{V}(r_1,r_2)$ with $m=-1$, $n=2$ and $p=1$.}
\label{fig_glueing_map}
\end{center}
\end{figure}

\paragraph*{Gluing map}
Let $\mathbb{V}(r_1,r_2)$ and $R(r_1,r_2)$ be the tubular neighborhoods of $\gamma$ and $\gamma_0$ defined in \eqref{region_Vi}-\eqref{region_V} and \eqref{region_R}-\eqref{region_Ri}, respectively. Let $M_R(r_1,r_2)\coloneqq M\backslash\interior(R(r_1,r_2))$. Define $\varphi:\partial M_R(r_1,r_2)\to\partial\mathbb{V}(r_1,r_2)$ in the following way: For every $p\in\partial M_R(r_1,r_2)$ choose $1\leq i\leq 4$ such that $p$ belongs to the $i$-th quadrant and denote its coordinates by $(x,y,z)=\Pi_i(p)\in\R^2\times\R/\Z$. Then:
\begin{enumerate}
\item If $m=m(\gamma,\Sigma)<0$ then 
\begin{equation}\label{glueing_map_negative_linking_}
\varphi(p)=\left\{
\begin{array}{cll}
 \left(x\ ,\ y\ ,\ z\right) &;\text{ if }p\notin A_\mathrm{in}^1,\\
 \left(r_1\ ,\ y\ ,\ z+\frac{|m|}{n}\rho\left(\frac{y}{r_2}\right)\right) &;\text{ if }p\in A_\mathrm{in}^1\text{ and }(r_1,y,z)=\Pi_1(p).
\end{array}\right.
\end{equation}

\item If $m=m(\gamma,\Sigma)>0$ then
\begin{equation}\label{glueing_map_positive_linking_}
\varphi(p)=\left\{
\begin{array}{cll}
 \left(x\ ,\ y\ ,\ z\right) &;\text{ if }p\notin A_\mathrm{in}^4,\\
 \left(r_1\ ,\ y\ ,\ z+\frac{|m|}{n}\rho\left(-\frac{y}{r_2}\right)\right) &;\text{ if }p\in A_\mathrm{in}^4\text{ and }(r_1,y,z)=\Pi_4(p).
\end{array}\right.
\end{equation}
\end{enumerate}
In Figure \ref{fig_glueing_map} we represent the map $\varphi$ in normal coordinates for the case $m=-1$, $n=2$, $p=1$. It is a \emph{twist map} over the annulus $A_\mathrm{in}^1$ and the identity on its complement. 

\begin{proposition}\label{lemma_defn_psi}
For every couple of small parameters $0<r_2<r_1<1$, the map $\varphi$ given in \eqref{glueing_map_negative_linking_}-\eqref{glueing_map_positive_linking_} above is a well-defined diffeomorphism. Consider the manifold 
$N(r_1,r_2)=M_R(r_1,r_2)\sqcup_\varphi\mathbb{V}(r_1,r_2)$ obtained by gluing $M_R(r_1,r_2)$ and $\mathbb{V}(r_1,r_2)$ along their boundaries using $\varphi$. Then, there exists a smooth atlas $\mathcal{D}$ and a smooth vector field $Y$ on $N(r_1,r_2)$ such that:
\begin{enumerate}[$(i)$]
\item The inclusion $\iota:M_R(r_1,r_2)\hookrightarrow N(r_1,r_2)$ is a diffeomorphism onto its image and $D\iota$ maps the vector field $X_\Gamma$ to the vector field $Y$,
\item The inclusion $\iota:\mathbb{V}(r_1,r_2)\hookrightarrow N(r_1,r_2)$ is a diffeomorphism onto its image and $D\iota$ maps the vector field $X_{(\lambda,n,p)}$ to the vector field $Y$. 
\end{enumerate}
We denote by $\{\psi_t:N(r_1,r_2)\to N(r_1,r_2)\}_{t\in\R}$ the flow generated by $Y$.
\end{proposition}

\noindent
\textbf{Notation:}
For simplicity we set $M_0=M_R(r_1,r_2)$ and $M_1=\mathbb{V}(r_1,r_2)$, and for each $j=0,1$ we denote by $X_j$ the corresponding vector field. We denote by $N=M_0\sqcup_\varphi M_1$ and by $\iota_j:M_j\to N$ the natural inclusion maps. The two components are glued along the subset $M_0\cap M_1\simeq\partial M_0\simeq\partial M_1$ that is homeomorphic to a two dimensional torus.

\begin{proof}
In coordinates \eqref{normal_coordinates_section_5} the map $\varphi$ consists in apply a twist map, supported on the annuli $A_\mathrm{in}^1$ or $A_\mathrm{in}^4$ depending on the signature of the multiplicity. This assumption is important at Subsection \ref{subsection_smooth_model_is_equivalent} for the purpose of constructing a Birkhoff section. Along the present section, we stay in the case $m<0$ and $\textsl{supp}(\varphi)=A_\mathrm{in}^1$, being analogous the other one. 

The partition of the neighborhood $W$ in quadrants decomposes $\partial R$ in four regions, individually homeomorphic to the annulus. To show that $\varphi$ is well-defined we have to check its definition along the intersections of these four regions. It is well-defined over the intersections $A_\mathrm{out}^1\cap A_\mathrm{out}^2$, $A_\mathrm{in}^2\cap A_\mathrm{in}^3$ and $A_\mathrm{out}^3\cap A_\mathrm{out}^4$ because $\Pi_{i+1}\circ\Pi_i^{-1}=id$, for $i=1,2,3$. To check that it is well-defined on the intersection $A_\mathrm{in}^4\cap A_\mathrm{in}^1$ let $p$ be a point in this intersection with coordinates $\Pi_4(p)=(r_1,0,z)$ and $\Pi_1(p)=(r_1,0,z+m/n)$. Then applying \eqref{glueing_map_negative_linking_} we have 
\begin{align*}
&\text{for}\ p\in A_\mathrm{in}^4, &&\  \varphi(p)=id\circ\Pi_4(p)=(r_1,0,z),\\
&\text{for}\ p\in A_\mathrm{in}^1, &&\ \varphi(p)=\Pi_1(p)+(0,0,|m|/n)=(r_1,0,z).
\end{align*}
Thus $\varphi$ is a well-defined bijection between the boundary of $M_0$ and the boundary of $M_1$, and since the function $\rho$ above is constant in a neighborhood of $0$, we deduce that $\varphi$ is a diffeomorphism.

The quotient space $N=M_0\sqcup_\varphi M_1$ is then a $3$-manifold, and since each $M_i$ is smooth there is a smooth atlas defined over $N\backslash(M_0\cap M_1)$, induced by the inclusions $\iota_j:\interior(M_j)\hookrightarrow N$. To extend this smooth atlas along $M_0\cap M_1$ it suffices to provide a family
\begin{equation}\label{transformaciones_Phi}
\left\{\Phi_p:T_pM_0\to T_{\varphi(p)}M_1, p\in\partial M_0\right\}
\end{equation}
satisfying that
\begin{enumerate}[(i)]
\item $\Phi_p:T_pM_0\to T_{\varphi(p)}M_1$ is a linear isomorphism, $\forall\ p\in\partial M_0$,
\item $\Phi_p|_{T_p\partial M_0}=D\varphi$,
\item the map $p\mapsto\Phi_p$ varies smoothly.
\end{enumerate}
If in addition we want that the vector fields $X_j$, $j=1,2$ match together along the boundary, we need that
\begin{itemize}
\item[(iv)] $\Phi_p:X_0(p)\mapsto X_1(\varphi(p))$, $\forall\ p\in\partial M_0$.
\end{itemize}

Now, choose a point $p\in\partial M_0$ belonging to the $i$-th quadrant. If $p\notin A_\mathrm{in}^1$, since $\varphi\circ\Pi_i^{-1}=id$ in a neighborhood of $p$, $D\Pi_i:X_0\mapsto X_{(\lambda,n,p)}$ and $(M_1,X_1)=(\mathbb{V},X_{(\lambda,n,p)})$, we see that 
\begin{equation}\label{transf_Phi=id}
\Phi_p:T_pM_0\to T_{\varphi(p)}M_1\text{ given by }\Phi_p\circ D\Pi_i^{-1}=id_{\R^3} 
\end{equation}  
satisfies (i)-(iv) above. If $p\in A_\mathrm{in}^1$, since the vector field $X_0$ is transverse to the boundary of $M_0$ along the annulus $A^1_\mathrm{in}$ and points outward, and $X_1$ is transverse to the boundary of $M_1$ along $\varphi(A_\mathrm{in}^1)$ and points inward, then 
\begin{equation}\label{transf_Phi_soporte}
\Phi_p:T_pM_0\to T_{\varphi(p)}M_1\text{ given by }\Phi_p|_{T_p\partial M_0}=D\varphi(p)\text{ and }\Phi_p(X_0(p))=X_1(p) 
\end{equation}  
satisfies the desired properties (i)-(iv). Since $\varphi$ coincides with a vertical translation near the boundary components of $A_\mathrm{in}^1$ and $X_{(\lambda,n,p)}$ is invariant by these translations, then the transformation in \eqref{transf_Phi_soporte} satisfies $\Phi_p\circ D\Pi_1^{-1}=id$ for points $p$ near $\partial A_\mathrm{in}^1$, so it is compatible with \eqref{transf_Phi=id} above. 

We refer to Chapter 13 of \cite{brocker-janich_diff-top} for more details on how to glue manifolds and vector fields along the boundary.
\end{proof}

\subsubsection{The action of the flow $\{\psi_t:N\to N\}_{t\in\R}$}\label{subsection_action_flow_psi_t}
The action of the flow $\psi$ on a point $p$ in the manifold $N=M_0\sqcup_\varphi M_1$ can be described as an alternated composition of integral curves of $X_j$ on each component $M_j$ and the gluing map $\varphi$. 

\begin{lem}\label{lema_iterated_action_formula}
Given $p\in\interior(M_0)$ and $t\geq 0$ such that $\psi_t(p)\in \interior(M_0)$, there exist $0<t_1,\dots,t_{l+1}<t$ and $0<s_1,\dots,s_l<t$ satisfying that $t_k\geq T$ for $k=2,\dots,l$, $s_k\geq T$ for $k=1,\dots,l$, their sum is  $\sum_{k=0}^{l}t_k+\sum_{k=1}^{l}s_k=t$ and
\begin{align*}
\psi_t(p)=\phi_{t_{l+1}}^0\circ\left(\varphi^{-1}\circ\phi_{s_l}^1\circ\varphi\right)\circ                 \cdots\circ\left(\varphi^{-1}\circ\phi_{s_2}^1\circ\varphi\right)\circ\phi_{t_2}^0\circ \left(\varphi^{-1}\circ\phi_{s_1}^1\circ\varphi\right)\circ\phi_{t_1}^0(p).
\end{align*}
The set of the $s_i$'s may be empty ($l=0$) and in that case $t_1=t$. Analogous formulations hold for $\psi$-orbit segments starting and ending in any combination of components $M_j$, $j=0,1$.
\end{lem}

\begin{proof}
Let $p\in\interior(M_0)$ and $t>0$ such that $\psi_t(p)\in\interior(M_0)$. If the orbit segment $[p,\psi_t(p)]$ is disjoint from $\interior(M_1)$, then it is completely contained in $M_0$ and the statement follows just by setting $t_1=t$ (and an empty set of $s_i$'s). If the orbit segment $[p,\psi_t(p)]$ intersects $\interior(M_1)$, then the positive semi-orbit of $p$ will follow the $\phi^0$-trajectory in $M_0$ until the first time $t_1>0$, when it hits $\partial M_0$ transversally in a point $p_1=\phi_{t_1}^0(p)$. Then, $p_1$ is identified with $\varphi(p_1)\in\partial M_1$ and the orbit switches to the component $M_1$, concatenating with the $\phi^1$-trajectory starting on the point $q_1=\varphi(p_1)$. Since the $\phi^1$-orbits starting at $q_1$ intersects $\interior(M_0)$ in a future time, there exists a first time $s_1>0$ when it hits $\partial M_1$ transversally in a point $\phi_{s_1}^1(q_1)$, and switch to the component $M_0$ concatenating with the $\phi^0$-orbit starting at $p_2=\varphi^{-1}(\phi_{s_1}^1(q_1))$. 

Proceeding inductively, we can decompose the orbit segment $[p,\psi_t(p)]$ as an alternated composition of a finite number $l+1$ of $\phi^0$-orbit segments of length $t_i$ and $\phi^1$-orbit segments of length $s_i$, using the gluing map ($\varphi$ or $\varphi^{-1}$) to switch from one component to the other. We remark that there exists a constant $T>0$ such that, for both $j=0,1$, 
\begin{align*}
& \min\left\{t>0:\exists\ p\in\partial M_j\text{ such that }\phi_t^j(p)\in\partial M_j\text{ and }\phi_s^j(p)\in\interior(M_j),\ \forall\ 0< s< t\right\}\geq T.
\end{align*}
This implies that all $s_i$'s, and all $t_i$'s except $t_1$ and $t_{l+1}$, are greater that a constant $T>0$.
\end{proof}

We describe now the derivative action $D\psi_t:TN\to TN$ on the tangent bundle. Recall that on each component $(M_j,X_j)$, $j=0,1$, there is a splitting~$TM_j=\Es_j\oplus \Ec_j\oplus \Eu_j$ and a Riemannian metric $|\cdot|_j$ (cf. \eqref{normal_form_section_5} for $(M_0,X_0)=(M_\Gamma,X_\Gamma)$ and Remark \ref{remark_properties_of_affine_X} for $(M_1,X_1)=(\mathbb{V},X_{(\lambda,n,p)})$), that is invariant by $D\phi_t^j(p)$ (for every $p\in M_j$ and $t\in\R$ such that the $\phi^j$-orbit segment between $p$ and $\phi_t^j(p)$ is entirely contained in $M_j$) and uniformly hyperbolic with contraction rate $0<\lambda<1$, where $\Es_j$ is the contracting bundle, $\Eu_j$ the expanding bundle and $\Ec_j=\spa\{X_j\}$. 

\paragraph*{Framing on TN}
To study the action of $\{D\psi_t\}_{t\in\R}$ on $TN$ we define a (non-continuous) splitting by
\begin{align}\label{non-cont-framing_N}
& T_pN=H^\mathrm{s}(p)\oplus H^\mathrm{c}(p)\oplus H^\mathrm{u}(p)=
\left\{
\begin{array}{lll}
 \Es_1(p)\oplus\Ec_1(p)\oplus\Eu_1(p),\ \text{if}\ p\in \interior(M_1),\\
 \Es_0(p)\oplus\Ec_0(p)\oplus\Eu_0(p),\ \text{if}\ p\in N\setminus\interior(M_1),
\end{array}\right. \ \forall p\in N.
\end{align}

\begin{rem}
We can see how does $\Es_0\oplus\Ec_0\oplus\Eu_0$ matches $\Es_1\oplus\Ec_1\oplus\Eu_1$ along $M_0\cap M_1$ using the identifications $\Phi_p:T_pM_0\to T_{\varphi(p)}M_1$ given in \eqref{transformaciones_Phi}, \eqref{transf_Phi=id}, \eqref{transf_Phi_soporte}. For every $p\in \partial M_0$ the map $\Phi_p$ sends $X_0(p)\mapsto X_1(\varphi(p))$, so it follows that the bundle $H^\mathrm{c}=\spa\{Y\}$ is continuous. If $p\notin A_\mathrm{in}^1$ then by \eqref{normal_form_section_5}, \eqref{transf_Phi=id} and 3. on Remark \ref{remark_properties_of_affine_X}, we obtain that $\Phi_p$ sends $\Es_0(p)\mapsto \Es_1(\varphi(p))$ and $\Eu_0(p)\mapsto \Eu_1(\varphi(p))$, obtaining continuity of the decomposition \eqref{non-cont-framing_N} in $p$.  Nevertheless, for points $p\in A_\mathrm{in}^1$ the two decompositions $\Es_0(p)\oplus\Ec_0(p)\oplus\Eu_0(p)$ and $\Es_1(\varphi(p))\oplus\Ec_1(\varphi(p))\oplus\Eu_1(\varphi(p))$ do not agree via $\Phi_p$ in general, so the framing \eqref{non-cont-framing_N} above may be non-continuous in $p$.
\end{rem}

\begin{lem}\label{lema_accion_Dpsi}
Given $p\in N$ and $t>0$, there exist $0<t_1,\dots,t_{l+1}<t$ satisfying $t_k\geq T$ for $k=2,\dots,l$ and points $p_1,\dots,p_l\in A_\mathrm{in}^1$ (empty in case $l=0$) such that $D\psi_t:T_pN\to T_{\psi_t(p)}N$ is an iterated composition of the form:
\begin{align}\label{Dpsi_forma_gral}
& D\psi_t(p)=\Psi_{t_{l+1}}\circ\Phi_{p_l}\circ\Psi_{t_l}\cdots\circ\Phi_{p_2}\circ\Psi_{t_2}\circ\Phi_{p_1}\circ \Psi_{t_1}(p),
\end{align}
where the transformations $\Psi_\tau$ and $\Phi_{p}$ are the following:

\begin{enumerate}
\item $\Psi_\tau$ denotes the application $T_p N \to T_{\psi_\tau(p)} N$, defined for all the couples $(\tau,p)$ satisfying that $[p,\psi_\tau(p)]\cap A_\mathrm{in}^1=\emptyset$ by the expression:
\begin{equation}\label{mapa_Psi}
\Psi_t: aY(p)+be_\mathrm{s}(p)+ce_\mathrm{u}(p)\mapsto aY(\psi_\tau(p))+\lambda^t b e_\mathrm{s}(\psi_\tau(p))+\lambda^{-t} c e_\mathrm{u}(\psi_\tau(p)),
\end{equation}
where $\{Y,e_\mathrm{s},e_\mathrm{u}\}$ is a continuous frame along the curve $[p,\psi_\tau(p)]$ s.t. $e_\mathrm{s}(\psi_\theta(p))\in H^\mathrm{s}(\psi_\theta(p))$ and $e_\mathrm{u}(\psi_\theta(p))\in H^\mathrm{u}(\psi_\theta(p))$ are unitary vectors, for every $0\leq \theta\leq \tau$.

\item $\Phi_p:T_pM_0\to T_{\varphi(p)}M_1$ are the transformations defined in \eqref{transf_Phi_soporte} for $p\in A_\mathrm{in}^1$. Given $p\in A_\mathrm{in}^1$ with coordinates $\Pi_1(p)=(r_1,r,z)$, $0\leq r\leq r_2$, consider two basis 
$$\mathcal{B}_0=\{X_0(p),e_0^\mathrm{s}(p),e_0^\mathrm{u}(p)\}\ \text{and}\ \mathcal{B}_1=\{X_1(\varphi(p)),e_1^\mathrm{s}(\varphi(p)),e_1^\mathrm{u}(\varphi(p))\}$$
of $T_pM_0$ and $T_{\varphi(p)}M_1$ respectively, where the $e_j^\mathrm{s}$, $e_j^\mathrm{u}$ are unitary vectors contained the spaces $E_j^\mathrm{s}$, $E_j^\mathrm{u}$ and oriented as in Figure \ref{fig_su-coordinates} below. Then 
\begin{equation}\label{lema_cmbio_de_lado_en_coordenadas_hiperbolicas}
_{\mathcal{B}_1}\left(\Phi_p\right)_{\mathcal{B}_0}=
\begin{pmatrix}
1 & \frac{K(r,r_2)}{|\log(\lambda)|r_1} & \frac{K(r,r_2)}{|\log(\lambda)|r}\\
0 & -K(r,r_2)+1 & -K(r,r_2)\frac{r_1}{r}\\
0 & K(r,r_2)\frac{r}{r_1} & K(r,r_2)+1
\end{pmatrix}\, \text{;}\ K(r,r_2)=-\left\vert p m\log(\lambda)\rho'(r/r_2)\right\vert\frac{r}{r_2}.\ 
\end{equation}
\end{enumerate} 
\end{lem}

\begin{proof}
The expression in \eqref{Dpsi_forma_gral} follows Lemma \ref{lema_iterated_action_formula} by taking derivatives and using that, when switching from one component $M_j$ to the other, the tangent spaces are identified with the transformations \eqref{transformaciones_Phi}-\eqref{transf_Phi=id}-\eqref{transf_Phi_soporte}. For $k=1,\dots,l$ denote by $p_k=\psi_{t_k+s_{k-1}}(p_{k-1})$ and $q_k=\psi_{s_k}(p_k)$, so the orbit of $p$ switches from $M_0$ to $M_1$ on the points $p_k$, and from $M_1$ to $M_0$ on the points $q_k$. We obtain that
\begin{align*}
D\psi_t(p)=D\phi_{t_{l+1}}^0 \circ \Phi_{q_l} \circ D\phi_{s_l}^1 \circ \Phi_{p_l} \circ \cdots \circ \Phi_{q_2} \circ\phi_{s_2}^1 \circ \Phi_{p_2} \circ \phi_{t_2}^0 \circ \Phi_{q_1} \circ \phi_{s_1}^1 \circ \Phi_{p_1} \circ \phi_{t_1}^0(p).
\end{align*}
Since the points $q_k$ are not contained in the support $A_\mathrm{in}^1$ of the gluing map $\varphi$, then $\Phi_{q_k}=id$ by \eqref{transf_Phi=id} and we can reduce the previous expression by only considering the points $p_k$ (and corresponding times $t_k$) where the orbit segment $[p,\psi_t(p)]$ intersects $A_\mathrm{in}^1$. Since each transformation $D\phi_{\tau}^j$ acts on $TM_j$ in the form specified in \eqref{mapa_Psi}, we deduce \eqref{Dpsi_forma_gral}.

To show \eqref{lema_cmbio_de_lado_en_coordenadas_hiperbolicas}, consider $p\in A_\mathrm{in}^1$ with coordinates $\Pi_1(p)=(r_1,r,z)$, $0\leq r\leq r_2$. Consider the basis of $\R^3$ given by
$$\{e_1=(1,0,0),e_2=(0,1,0),e_3=(0,0,1)\}.$$
Then, by \eqref{normal_form_section_5} and 3. on \ref{remark_properties_of_affine_X} it follows that $\mathcal{B}_0=\{X_{(\lambda,n,p)}(p),e_1,e_2\}$ and $\mathcal{B}_1=\{X_{(\lambda,n,p)} (\varphi(p)),$ $e_1,e_2\}$ are basis satisfying the hypothesis above. Thus, the matrix \eqref{lema_cmbio_de_lado_en_coordenadas_hiperbolicas} is the expression in these basis of $\Phi_p\circ D\Pi_1^{-1}:\R^3\to\R^3$. Now, since $\Phi_p$ sends $X_{(\lambda,n,p)}(\Pi_1(p))\mapsto X_{(\lambda,n,p)}(\varphi(p))$ and coincides with $D\varphi$ on $TA_\mathrm{in}^1=0\times\R\times\R$, in the basis $\mathcal{C}=\{X_{(\lambda,n,p)},e_2,e_3\}$ we have that
\begin{equation}\label{lema_cambio_de_lado}
_\mathcal{C}\left(\Phi_p\right)_{\mathcal{C}}=
\begin{pmatrix}
1 & 0 & 0\\
0 & 1 & 0\\
0 & -\frac{\kappa(r/r_2)}{r_2} & 1
\end{pmatrix},\text{ where }\kappa(r/r_2)=\frac{|m|}{n}\left\vert\rho'(r/r_2)\right\vert.
\end{equation}

The expression \eqref{lema_cmbio_de_lado_en_coordenadas_hiperbolicas} is finally obtained by applying a change of basis between $\mathcal{C}$ and $\mathcal{B}_j$. Namely, if we write the change of basis matrices
\begin{equation*}
_{\mathcal{B}_1}\left(I\right)_{\mathcal{C}}=
\begin{pmatrix}
1 & 0 & \frac{1}{X_3(\varphi(p))}\\
0 & 0 & -\frac{X_1(\varphi(p))}{X_3(\varphi(p))}\\
0 & 1 & -\frac{X_2(\varphi(p))}{X_3(\varphi(p))}
\end{pmatrix}
\ \text{and}\ 
_\mathcal{C}\left(I\right)_{\mathcal{B}_0}=
\begin{pmatrix}
1 & \frac{1}{X_1(p)} & 0\\
0 & -\frac{X_2(p)}{X_1(p)} & 1\\
0 & -\frac{X_3(p)}{X_1(p)} & 0
\end{pmatrix},
\end{equation*}
where $X_1$, $X_2$, $X_3$ denotes the components of the vector field $X_{(\lambda,n,p)}$ in \eqref{normal_form_section_5}, 
we obtain that
\begin{equation*}
_{\mathcal{B}_1}\left(\Phi_p\right)_{\mathcal{B}_0}
=
\begin{pmatrix}
1 & 0 & \frac{1}{X_3}\\
0 & 0 & -\frac{X_1}{X_3}\\
0 & 1 & -\frac{X_2}{X_3}
\end{pmatrix}
\begin{pmatrix}
1 & 0 & 0\\
0 & 1 & 0\\
0 & -\frac{\kappa}{r_2} & 1
\end{pmatrix}
\begin{pmatrix}
1 & \frac{1}{X_1} & 0\\
0 & -\frac{X_2}{X_1} & 1\\
0 & -\frac{X_3}{X_1} & 0
\end{pmatrix}
=
\begin{pmatrix}
1 & \frac{\kappa}{r_2}\frac{X_2}{X_1X_3} & -\frac{\kappa}{r_2}\frac{1}{X_3}\\
0 & -\frac{\kappa}{r_2}\frac{X_2}{X_3}+1 & -\frac{\kappa}{r_2}\frac{X_1}{X_3}\\
0 & -\frac{\kappa}{r_2}\frac{X_2^2}{X_1X_3} & \frac{\kappa}{r_2}\frac{X_2}{X_3}+1
\end{pmatrix}.
\end{equation*}
Finally, recalling that $X_1(r_1,r,z)=\log(\lambda)r_1$, $X_2(r_1,r,z)=-\log(\lambda)r$, $X_3(r_1,r,z)=1/np$, that $\varphi$ leaves invariant the components of the vector field, and that $K(r,r_2)=-|n\log(\lambda)\kappa(r/r_2)|$, we obtain the expression \eqref{lema_cmbio_de_lado_en_coordenadas_hiperbolicas} from the statement of the lemma. 
\end{proof}

\subsubsection{Invariant measure}
\begin{proposition}\label{prop_existence_of_smooth_measure}
The flow $\{\psi_t:N\to N\}_{t\in\R}$ in Proposition \ref{lemma_defn_psi} preserves a smooth volume form.
\end{proposition}

\begin{proof}
To prove this we define a smooth volume form $\omega_i$ on each $M_i$, $i=0,1$, invariant by the action of $\phi^i$, and we prove that they match together along $M_0\cap M_1$ under the gluing map $\varphi$. 

For each $p\in M_0$ consider a positively oriented basis $\{X_0(p),e_\mathrm{s}^0(p),e_\mathrm{u}^0(p)\}$ of $T_pM_0$, according to the decomposition $\Es_0(p)\oplus\Ec_0(p)\oplus\Eu_0(p)$, and such that $|e_\mathrm{s}^0(p)|_0=|e_\mathrm{u}^0(p)|_0=1$. Define a volume form $\omega_0$ in $M_0$ by the expression $\omega_0(p)=dX_0(p)\wedge de_\mathrm{s}^0(p)\wedge de_\mathrm{u}^0(p)$. Since by construction the splitting $TM_0=\Es_0\oplus\Ec_0\oplus\Eu_0$ is smooth then $\omega_0$ is a smooth volume form on $M_0$, and from \eqref{mapa_Psi} it follows that:
\begin{equation*}
\begin{split}
\left(\left(\phi_t^0\right)^*\omega_0\right)(p) & =\left(\left(\phi_t^0\right)^* dX_0\right)(p)\wedge\left(\left(\phi_t^0\right)^* de_\mathrm{s}^0\right)(p)\wedge\left(\left(\phi_t^0\right)^* de_\mathrm{u}^0\right)(p) \\
& = dX_0(p)\wedge \left(\lambda^t\cdot de_\mathrm{s}^0\right)(p)\wedge \left(\lambda^{-t}\cdot de_\mathrm{u}^0\right)(p) =\omega_0(p) 
\end{split}
\end{equation*}
for every $p\in M_0$ and $t\geq 0$ such that $[p,\phi_t^0(p)]$ is contained in $M_0$.

For each point $p=(x,y,z)\in M_1=\mathbb{V}(r_1,r_2)$ consider the basis $\{X_1(p),e_\mathrm{s}^1(p),e_\mathrm{u}^1(p)\}$, where $e_\mathrm{s}^1(p)=(1,0,0)$ and  $e_\mathrm{u}^1(p)=(0,1,0)$. Define $\omega_1(p)=dX_1(p)\wedge de_\mathrm{s}^1(p)\wedge de_\mathrm{s}^2(p)$. Then, this is a smooth volume form and by Remark \ref{remark_properties_of_affine_X} we have
\begin{equation*}
\begin{split}
\left(\left(\phi_t^1\right)^*\omega_1\right)(p) & =\left(\left(\phi_t^1\right)^* dX_1\right)(p)\wedge\left(\left(\phi_t^1\right)^* de_\mathrm{s}^1\right)(p)\wedge\left(\left(\phi_t^1\right)^* de_\mathrm{u}^1\right)(p) \\
& = dX_1(p)\wedge \left(\lambda^t\cdot de_\mathrm{s}^1\right)(p)\wedge \left(\lambda^{-t}\cdot de_\mathrm{u}^1\right)(p) =\omega_1(p) 
\end{split}
\end{equation*}
for every $p\in M_1$ and $t\geq 0$ such that $[p,\phi_t^1(p)]$ is contained in $M_1$.

Finally, observe that by \eqref{lema_cmbio_de_lado_en_coordenadas_hiperbolicas}, for every $p\in\partial M_0$ the transformation $\Phi_p:T_pM_0\to T_{\varphi(p)}M_1$ has determinant $\det(\Phi_q)=1$ in the basis $\{X_i(p),e_\mathrm{s}^i(p),e_\mathrm{u}^i(p)\},\ i=1,2$. Thus
$$\left(\left(\Phi_p\right)^*\omega_1\right)(p)=\det(\Phi_q)\cdot\omega_0(p)=\omega_0(p),$$
for every $p\in \partial M_0$.
From the three expressions above and the decomposition of $D\psi_t$ in products of $\Psi_t$ and $\Phi_q$ given in Lemma \ref{lema_accion_Dpsi}, we conclude that the $3$-form $\omega$ on $N=M_0\sqcup_\varphi M_1$ defined by $\omega(p)=\omega_i(p)$ for $p\in M_i$ is a well-defined smooth volume form satisfying $\psi_t^*(\omega)=\omega$.
\end{proof}

\subsection{Hyperbolicity}\label{subsection_smooth_model_is_Anosov}

\begin{proposition}\label{prop_smooth_model_is_Anosov}
If $0<r_2<r_1<1$ are chosen sufficiently small, then the smooth flow $\{\psi_t:N(r_1,r_2)\to N(r_1,r_2)\}_{t\in\R}$ constructed in Proposition \ref{lemma_defn_psi} is Anosov. 
\end{proposition}

We prove Proposition \ref{prop_smooth_model_is_Anosov} by showing that the associated \emph{Linear Poincar\'e Flow} is uniformly hyperbolic. This technique consists in the following: 

\paragraph*{Linear Poincar\'e flow}
The \emph{normal bundle} associated to the vector field $Y$ generating the flow $\{\psi_t:N\to N\}_{t\in\R}$, is the smooth 2-dimensional vector bundle $H\to N$, obtained by fiberwise quotienting the tangent spaces $T_pN$ by the subspaces $H^\mathrm{c}(p)=\spa\{Y(p)\}$. Since for every $t\in\R$, $p\in N$ the derivative action of the flow on $TN$ satisfy $D\psi_t(Y(p))=Y(\psi_t(p))$, there is an induced $\R$-action on the normal bundle, that commutes with $\psi_t:N\to N$, and for every $t\in\R$ is defined by
\begin{equation*}
D\psi_t:H\to H\ \text{such that}\ D\psi_t(\bar{v})=\overline{D\psi_t(v)},\ \text{for every}\ \bar{v}=v+H^\mathrm{c},
\end{equation*}
where $\bar{v}$ denotes the class $v+H^\mathrm{c}$ in the quotient $H=TN/H^\mathrm{c}$. Note that we use the same symbol $D\psi_t$ to denote either the action on $TN$ or the action on $H$. This action is called the \emph{linear Poincar\'e flow.}

The linear Poincar\'e flow captures the relevant information of the derivative action of the flow on the tangent bundle, by modding out the direction where the action is trivial. In particular, as noticed by Doering in \cite{Doering} (see also \cite{hozoori_2024_regularity_WL} and \cite{Asaoka_projectively_Anosov}), a smooth flow is Anosov (according to Definition \ref{defn_Anosov_flow}) if and only its associated linear Poincar\'e flow preserves a splitting of the normal bundle into two transverse sub-bundles, where norms are uniformly contracted/expanded under forward iteration. We will use this approach and show the following proposition that, according to Proposition 1.1 of \cite{Doering}, is equivalent to Proposition \ref{prop_smooth_model_is_Anosov} above.\footnote{We thank S. Hozoori for indicating us this argument and its references.} 

\begin{proposition}\label{corolario_existe_splitting_dominado_planos}
If $0<r_2<r_1<1$ are chosen sufficiently small, then the normal bundle $H\to N$ associated to $\{\psi_t:N\to N\}_{t\in\R}$ splits as the direct sum $H=F^{\mathrm{s}}\oplus F^{\mathrm{u}}$ of two continuous line bundles $F^{\mathrm{s}}\to N$ and $F^{\mathrm{u}}\to N$, satisfying that:
\begin{enumerate}
\item $D\psi_t(p)\left(F^{\mathrm{u}}(p)\right)=F^{\mathrm{u}}(\psi_t(p))$ and $D\psi_t(p)\left(F^{\mathrm{s}}(p)\right)=F^{\mathrm{s}}(\psi_t(p))$, for every $p\in N$, $t\in\R$;
\item Given a fixed Riemannian metric $\Vert\cdot\Vert$ on the bundle $H$, there exists constants $L>0$ and $\mu>1$ such that 
\begin{align*}
&\Vert D\psi_{t}(p)\cdot \bar{v}\ \Vert\geq L\mu^{t}\Vert \bar{v}\Vert,\ \forall\ p\in N,\ \bar{v}\in F^{\mathrm{u}}(p),\ t\geq 0,\\
&\Vert D\psi_{-t}(p)\cdot \bar{v}\Vert\geq L\mu^{t}\Vert \bar{v}\Vert,\ \forall\ p\in N,\ \bar{v}\in F^{\mathrm{s}}(p),\ t\geq 0.
\end{align*}
\end{enumerate}
\end{proposition}

\begin{rem}
In the definition of the gluing map $\varphi$ in  \eqref{glueing_map_negative_linking_} and \eqref{glueing_map_positive_linking_} there is an auxiliary smooth function $\rho$ and two parameters $0<r_2<r_1<1$. The maximum size of these real parameters depend on the charts $\Pi_i$, but they can be chosen arbitrarily small. We assume from now on that $\rho$ is fixed and we will adjust the parameters $r_1,r_2$ to satisfy the proposition. 
\end{rem}

For proving Proposition \ref{corolario_existe_splitting_dominado_planos}, we use the \emph{cone field criterion}; see \cite{katok-hasselblatt} for a precise statement. Let us introduce the following terminology:

\paragraph*{Stable and unstable slopes}
Given a point $p\in N$ consider a \emph{positively oriented} basis of $T_pN$ of the form $\{Y(p), e_\mathrm{s}(p), e_\mathrm{u}(p)\}$, where $e_\mathrm{s}(p)$ and $e_\mathrm{u}(p)$ are unitary vectors contained in the spaces $H^\mathrm{s}(p)$ and $H^\mathrm{u}(p)$ of the decomposition \eqref{non-cont-framing_N}, respectively. Given a vector $v=aY(p)+be_\mathrm{s}(p)+ce_\mathrm{u}(p)$ in $T_pN$ we define its $\mathrm{u}$-\emph{slope} and $\mathrm{s}$-\emph{slope} respectively as
\begin{align*}
& \Delta_\mathrm{u}(v)=\frac{b}{c}\ \text{and}\ \Delta_\mathrm{s}(v)=\frac{c}{b}.
\end{align*}
There are two possibilities for choosing a positive basis as before, but the slope is unchanged by switching this choice, so the $\mathrm{u}$,$\mathrm{s}$-slopes are well-defined. Observe that $\mathrm{u}$- and $\mathrm{s}$-slopes of vectors in $TN$ are unchanged by adding a component collinear with $Y$, and hence $\Delta_\mathrm{u}(\bar{v})=\Delta_\mathrm{u}(v)$ and $\Delta_\mathrm{s}(\bar{v})=\Delta_\mathrm{s}(v)$ are well-defined for vectors $\bar{v}=v+H^\mathrm{c}$ in the normal bundle $H(p)=T_pN/H^\mathrm{c}(p)$. 

\paragraph*{Cone distributions in $H$}
Given two real numbers that $-\infty<\delta_0<\delta_1<+\infty$ we define two cone distributions on the normal bundle, as follows: 
\begin{align*}
& \Cu(p;\delta_0,\delta_1)=\left\{\bar{v}\in H(p)\ :\ \delta_0\leq \Delta_\mathrm{u} (v) \leq\delta_1\right\}\\
& \Cs(p;\delta_0,\delta_1)=\left\{\bar{v}\in H(p)\ :\ \delta_0\leq\Delta_\mathrm{s} (v) \leq\delta_1\right\}.
\end{align*}
At each point $p\in N$, these sets form a pair of $1$-dimensional cones contained in the normal space $H(p)=T_pN/H^\mathrm{c}(p)$, and complementary in case $\delta_0$, $\delta_1$ have modulus not too big. 

\paragraph*{Norm in $H$}
For every vector $v=aY(p)+be_\mathrm{s}(p)+ce_\mathrm{u}(p)$ in $T_pN$ we define its \emph{$\mathrm{su}$-norm} by the expression $\Vert v\Vert_\mathrm{su}=\sqrt{b^2+c^2}$. It clearly induces a norm on each normal space by setting $\Vert\bar{v}\Vert_\mathrm{su}=\Vert v\Vert_\mathrm{su}$, for every $\bar{v}\in H(p)=T_pN/H^\mathrm{c}(p)$. 

We will show that for some adequate slope values $\delta_0^\mathrm{u},\delta_1^\mathrm{u}$ and $\delta_0^\mathrm{s},\delta_1^\mathrm{s}$ there is a pair of $\mathrm{u}$ and $\mathrm{s}$-cones satisfying the cone field criterion under the action of $\{D\psi_t:H\to H\}_{t\in\R}$. We remark that the slope functions, the Riemannian metric and the cone distributions 
$p\mapsto\Cu(p;\delta_0^\mathrm{u},\delta_1^\mathrm{u})$ and $p\mapsto\Cs(p;\delta_0^\mathrm{s},\delta_1^\mathrm{s})$ are not continuous as functions of $p\in N$. This is due to the discontinuities of the splitting $TN=H^\mathrm{s}\oplus H^\mathrm{c}\oplus H^\mathrm{u}$ over the set $A_\mathrm{in}^1$ given at \eqref{non-cont-framing_N}. Nevertheless, this poses no obstructions for applying the criterion. 

\paragraph*{Choice of slopes}
Define 
\begin{align}\label{funciones_delta_i}
& \delta_0^\mathrm{u}=-\frac{3r_1}{r_2},\ \delta_1^\mathrm{u}=\frac{r_2}{3r_1}\ \text{ and }\ \delta_0^\mathrm{s}=-\frac{3r_1}{2r_2},\ \delta_1^\mathrm{u}=\frac{2r_2}{3r_1}.
\end{align} 

The cone field criterion consists in checking the following three statements: 

\begin{lem}\label{lema_conos_1}
Let $\delta_0^\mathrm{u},\delta_1^\mathrm{u},\delta_0^\mathrm{s},\delta_1^\mathrm{s}:N\to\R$ be the quantities defined in (\ref{funciones_delta_i}) above. If $0<r_2<r_1<1$ are sufficiently small, then there exists $T>0$ such that: For every $p\in N$, 
\begin{align*}
&D\psi_{t}\left(\Cu\left(p\ ;\ \delta_0^\mathrm{u},\delta_1^\mathrm{u}\ \right)\right)\subset
\Cu\left(\psi_{t}(p)\ ;\ \delta_0^\mathrm{u},\delta_1^\mathrm{u}\ \right),\ \forall\ t\geq T,\\
&D\psi_{-t}\left(\Cs\left(p\ ;\ \delta_0^\mathrm{s},\delta_1^\mathrm{s}\ \right)\right)\subset
\Cs(\psi_{-t}(p)\ ;\ \delta_0^\mathrm{s},\delta_1^\mathrm{s}\ )),\ \forall\ t\geq T.
\end{align*}
\end{lem}

\begin{lem}\label{lema_conos_2}
For the parameters $0<r_2<r_1<1$ given in the previous lemma, it is satisfied the following: There exists a constant $L_0>0$ such that, for every $p\in N$ and every $t\geq 0$ then
\begin{align*}
&\left\vert\ \Delta_\mathrm{u}(D\psi_{t}(p)\cdot v_2)-\Delta_\mathrm{u}(D\psi_{t}(p)\cdot v_1)\ \right\vert\leq 
\lambda^{2t}L_0,\ \forall\ \bar{v}_1,\bar{v}_2\in\Cu(p;\delta_0^\mathrm{u},\delta_1^\mathrm{u}),\\
&\left\vert\Delta_\mathrm{s}(D\psi_{-t}(p)\cdot v_2)-\Delta_\mathrm{s}(D\psi_{-t}(p)\cdot v_1)\right\vert\leq
\lambda^{2t}L_0,\ \forall\ \bar{v}_1,\bar{v}_2\in\Cs(p;\delta_0^\mathrm{s},\delta_1^\mathrm{s}).
\end{align*}
\end{lem}

\begin{lem}\label{lema_conos_3}
By shrinking the parameters $0<r_2<r_1<1$ of the previous lemmas if necessary, it is satisfied that: There exist constants $L>0$ and $\mu>1$ such that, for every $p\in N$ and every $t\geq 0$ then
\begin{align*}
&\left\Vert D\psi_t(p)\cdot v\ \right\Vert_{\mathrm{su}}\geq L\mu^{t}\Vert v\Vert_{\mathrm{su}},\ \forall\ \bar{v}\in\Cu(p;\delta_0^\mathrm{u},\delta_1^\mathrm{u}),\\
&\left\Vert D\psi_{-t}(p)\cdot v\right\Vert_{\mathrm{su}}\geq L\mu^{t}\Vert v\Vert_{\mathrm{su}},\ \forall\ \bar{v}\in\Cs(p;\delta_0^\mathrm{s},\delta_1^\mathrm{s}).
\end{align*}
\end{lem}

\begin{proof}[Proof of Proposition \ref{corolario_existe_splitting_dominado_planos}]
Define
\begin{align*}
& F^{\mathrm{u}}(p)=\bigcap_{k\geq 0}D\psi_{kT}\left(\ \Cu(\ \psi_{-kT}(p)\ ;\ \delta_0^\mathrm{u},\delta_1^\mathrm{u}\ )\ \right),\\
& F^{\mathrm{s}}(p)=\bigcap_{k\geq 0}D\psi_{-kT}\left(\ \Cs(\ \psi_{kT}(p)\ ; \delta_0^\mathrm{s},\delta_1^\mathrm{s}\ )\ \right).
\end{align*}
This is a pair of $1$-dimensional cones (decreasing intersection of cones by Lemma \ref{lema_conos_1}) contained in $H(p)$. The slope function $\Delta_\mathrm{u}$ allows to identify each cone $\Cu(p;\delta_0^\mathrm{u},\delta_1^\mathrm{u})$ with the closed interval $[\delta_0^\mathrm{u},\delta_1^\mathrm{u}]\subset\R$. By Lemmas \ref{lema_conos_1} and \ref{lema_conos_2}, $F^\mathrm{u}(p)$ corresponds with a nested intersection of compact segments whose diameter tends to zero, so it is a non-empty cone that in fact reduces to a 1-dimensional subspace. The same considerations apply for $F^\mathrm{s}(p)$. 

Now, by Lemma \ref{lema_conos_3} we deduce the expansion/contraction property stated at item 2. of Proposition \ref{corolario_existe_splitting_dominado_planos} for each $F^{\mathrm{u}}$ and $F^{\mathrm{s}}$. Since $\Vert D\psi_t(p)\cdot v\Vert_{\mathrm{su}}\to+\infty$ for vectors $\bar{v}\in F^{\mathrm{u}}(p)$ and $\Vert D\psi_{t}(p)\cdot w\Vert_{\mathrm{su}}\to 0$ for vectors $\bar{w}\in F^{\mathrm{s}}(p)$, for $t\to +\infty$, then $F^{\mathrm{u}}(p)\cap F^{\mathrm{s}}(p)=\{0\}$ and hence the normal bundle splits as the direct sum $H(p)=\FFs(p)\oplus \FFu(p)$.  

Finally, observe that if another pair of subspaces $E^{\mathrm{s}}(p)$ and $E^{\mathrm{u}}(p)$ in $H(p)$ satisfies the conditions of \ref{corolario_existe_splitting_dominado_planos}, then 
$E^{\mathrm{s}}(p)=F^{\mathrm{s}}(p)$ and $E^{\mathrm{u}}(p)=F^{\mathrm{u}}(p)$. In particular, we obtain the invariance 
$$D\psi_t:F^{\mathrm{u}}(p)\mapsto F^{\mathrm{u}}(\psi_t(p)),\ \ D\psi_t:F^{\mathrm{s}}(p)\mapsto F^{\mathrm{s}}(\psi_t(p)),\ \ \forall\ t\in\R.$$ 

The continuity of the bundles $p\mapsto F^{\mathrm{s}}(p)$ and $p\mapsto F^{\mathrm{u}}(p)$ is automatic: this property is true for every pair of invariant line  bundles satisfying items 1. and 2. of Proposition \ref{corolario_existe_splitting_dominado_planos}, see \cite[Chapter 6.4]{katok-hasselblatt}.
\end{proof}

We now proceed to the proof of Lemmas \ref{lema_conos_1}, \ref{lema_conos_2} and \ref{lema_conos_3}. 

\paragraph*{Action on the normal bundle}
There is a decomposition of the normal bundle of the form $H=\bar{H}^\mathrm{s}\oplus\bar{H}^\mathrm{u}$, induced form $TN=H^\mathrm{s}\oplus H^\mathrm{c}\oplus H^\mathrm{u}$. By Lemma \ref{lema_accion_Dpsi} the action of $D\psi_t$ on $H$ is an alternated composition of transformations of the form $\Psi_t$ and $\Phi_p$. For every $p\in N$ let $\{\bar{e}_\mathrm{s}(p),\bar{e}_\mathrm{u}(p)\}$ be a positive basis of $H$ induced from two unitary vectors in $H^\mathrm{s}$ and $H^\mathrm{u}$. We have:

\begin{figure}[t]
\begin{center}
\includegraphics[scale=0.46,angle=0]{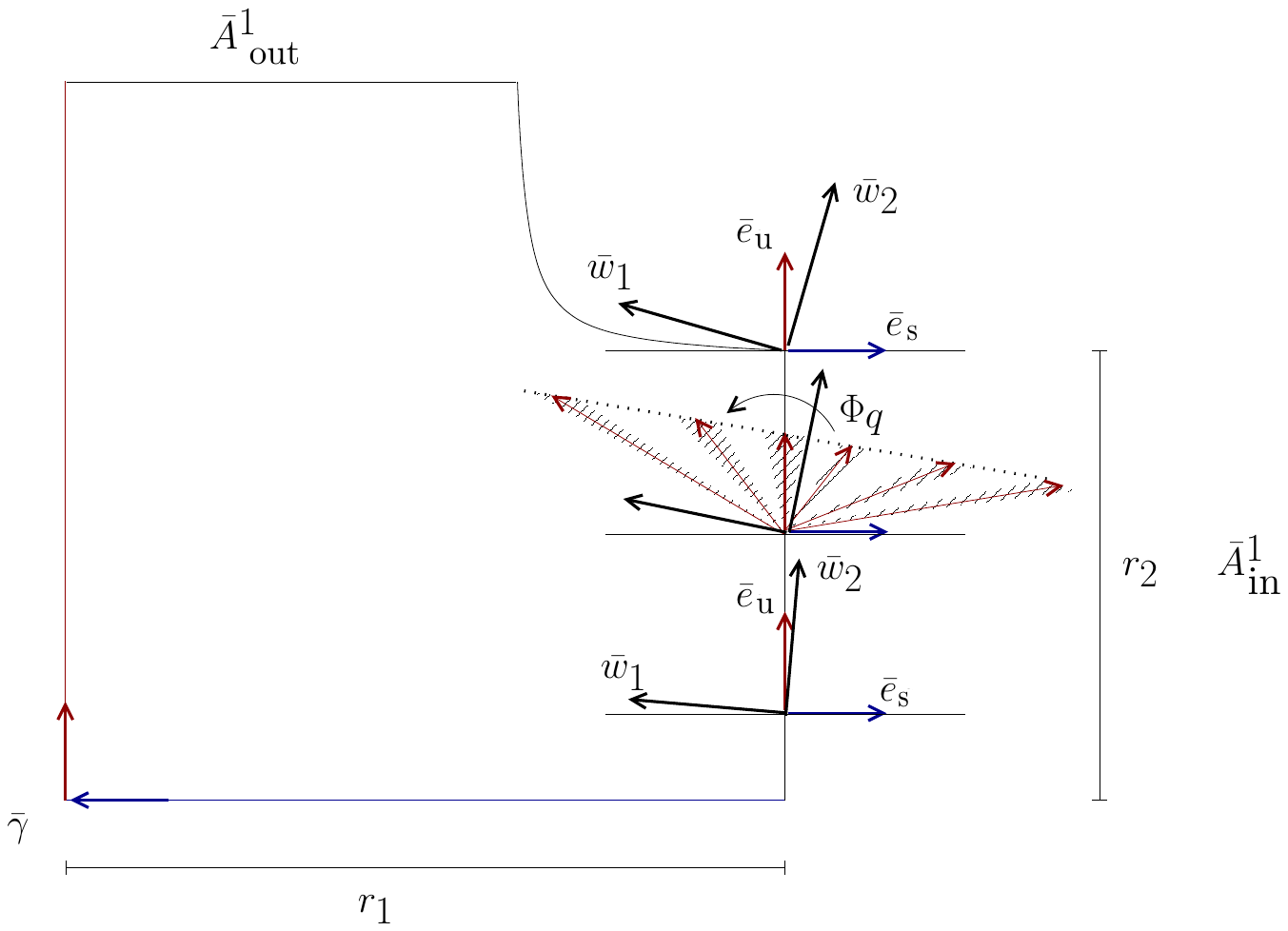}
\caption{The map $\Phi_p$ in $\mathrm{su}$-coordinates.}
\label{fig_su-coordinates}
\end{center}
\end{figure}

\begin{enumerate}
\item If $p\notin A_\mathrm{in}^1$ and $t>0$ satisfies that $[p,\psi_t(p)]\cap A_\mathrm{in}^1=\emptyset$, then the matrix associated to the action $D\psi_t(p)$ on the normal bundle, in the basis $\{\bar{e}_\mathrm{s},\bar{e}_\mathrm{u}\}$, is given by
\begin{equation}\label{Psi_matrix_su}
\left(\Psi_t\right)_{\mathrm{su}}=
\begin{pmatrix}
\lambda^t & 0\\
0 & \lambda^{-t}
\end{pmatrix}.
\end{equation}

\item Let $p\in A_\mathrm{in}^1$ with coordinates $\Pi_1(p)=(r_1,r,s)$. By \eqref{lema_cmbio_de_lado_en_coordenadas_hiperbolicas} in Lemma \ref{lema_accion_Dpsi} we see that that the action of $\Phi_p$ in the normal bundle, in the basis $\{\bar{e}_\mathrm{s},\bar{e}_\mathrm{u}\}$, is given by the matrix
\begin{equation}\label{su_action_Phi}
\left(\Phi_p\right)_{\mathrm{su}}=
\begin{pmatrix}
K(r,r_2)+1 & K(r,r_2)\frac{r_1}{r}\\
-K(r,r_2)\frac{r}{r_1} & -K(r,r_2)+1
\end{pmatrix},
\end{equation} 
where $K(r,r_2)=-(cte)\cdot|\rho'(r/r_2)|\frac{r}{r_2}$ and $(cte)=p|m||\log(\lambda)|$. 
Remark that $K$ is a non-positive function. Thus, we obtain a family of transformations parametrized over $0\leq r\leq r_2$. Observe that $\Phi_p$ is non-trivial just for points $p=(r_1,r,s)$ with $\frac{r_2}{3}\leq r\leq \frac{2r_2}{3}$ due to the definition of $\rho:[0,1]\to\R$. Since \eqref{su_action_Phi} has determinant equal to one and trace equal to two, it has a double eigenvalue equal to one. The vector $\bar{w}_1=-\bar{e}_\mathrm{s}+\frac{r}{r_1}\bar{e}_\mathrm{u}$ is an eigenvector for the action of this matrix on the normal bundle. Consider the vector $\bar{w}_2=\frac{r}{r_1}\bar{e}_\mathrm{s}+\bar{e}_\mathrm{u}$. In the orthogonal basis $\{\bar{w}_1,\bar{w}_2\}$ the transformation \eqref{su_action_Phi} takes the form 
\begin{equation}\label{jordan_form_Phi}
\left(\Phi_p\right)_{\{\bar{w}_1,\bar{w}_2\}}=
\begin{pmatrix}
1 & \eta(r)\frac{r_1}{r_2}\\
0 & 1
\end{pmatrix}
\end{equation}
for some continuous function $\eta:[0,r_2]\to[0,\infty)$ with support in $[\frac{r_2}{3},\frac{2r_2}{3}]$. We illustrate the action of the matrix (\ref{su_action_Phi}) in Figure \ref{fig_su-coordinates}.
\end{enumerate}

\begin{figure}[t]
     \begin{center}  
     \begin{subfigure}[b]{0.49\textwidth}
         \begin{center}
		 \includegraphics[width=0.9\textwidth,keepaspectratio]{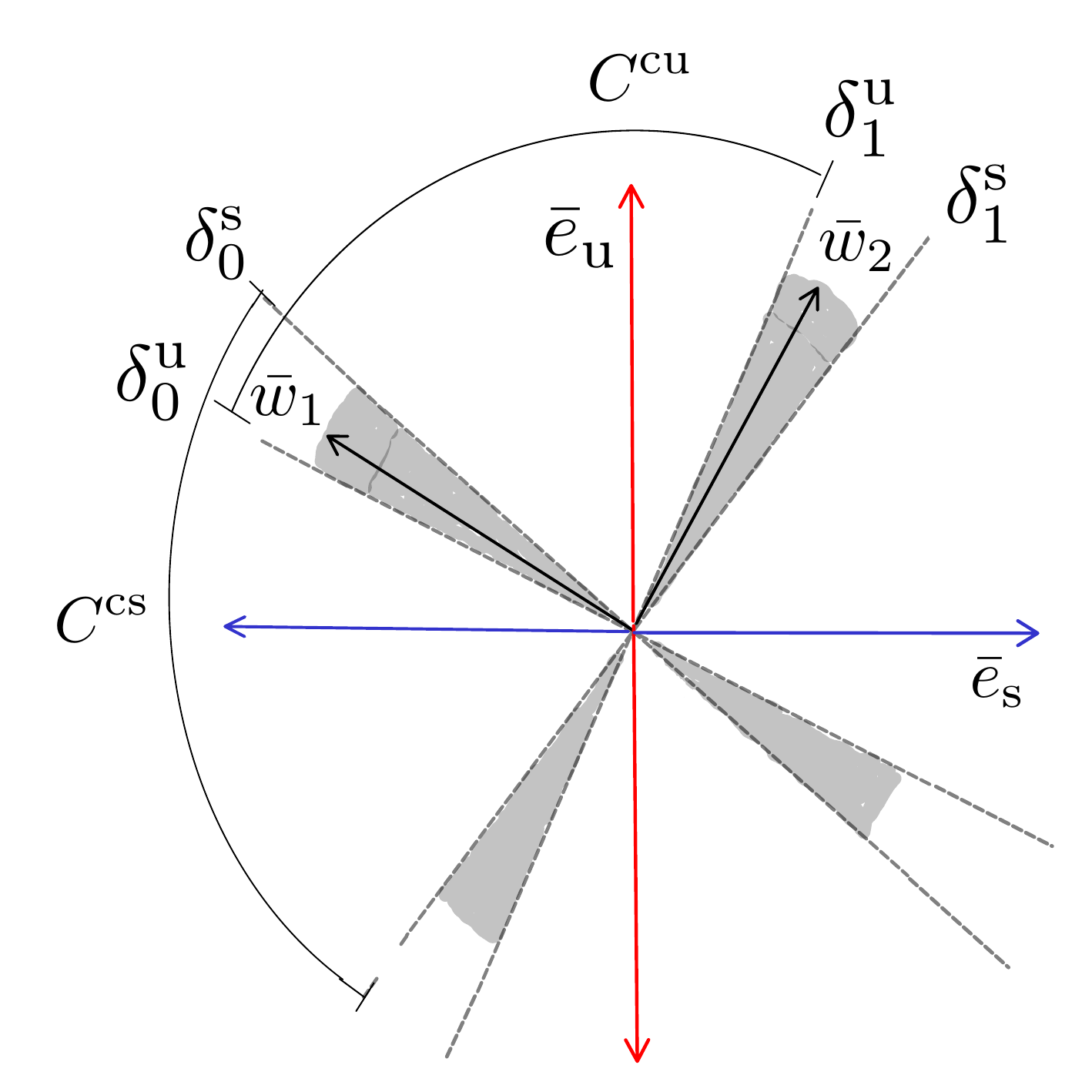}
         \vspace*{5mm}
         \caption{\ }
         \label{fig_slope_functions(2)}
		 \end{center}
     \end{subfigure}
     \hfill
     \hfill
     \begin{subfigure}[b]{0.49\textwidth}
         \begin{center}
	 	 \includegraphics[width=\textwidth,keepaspectratio]{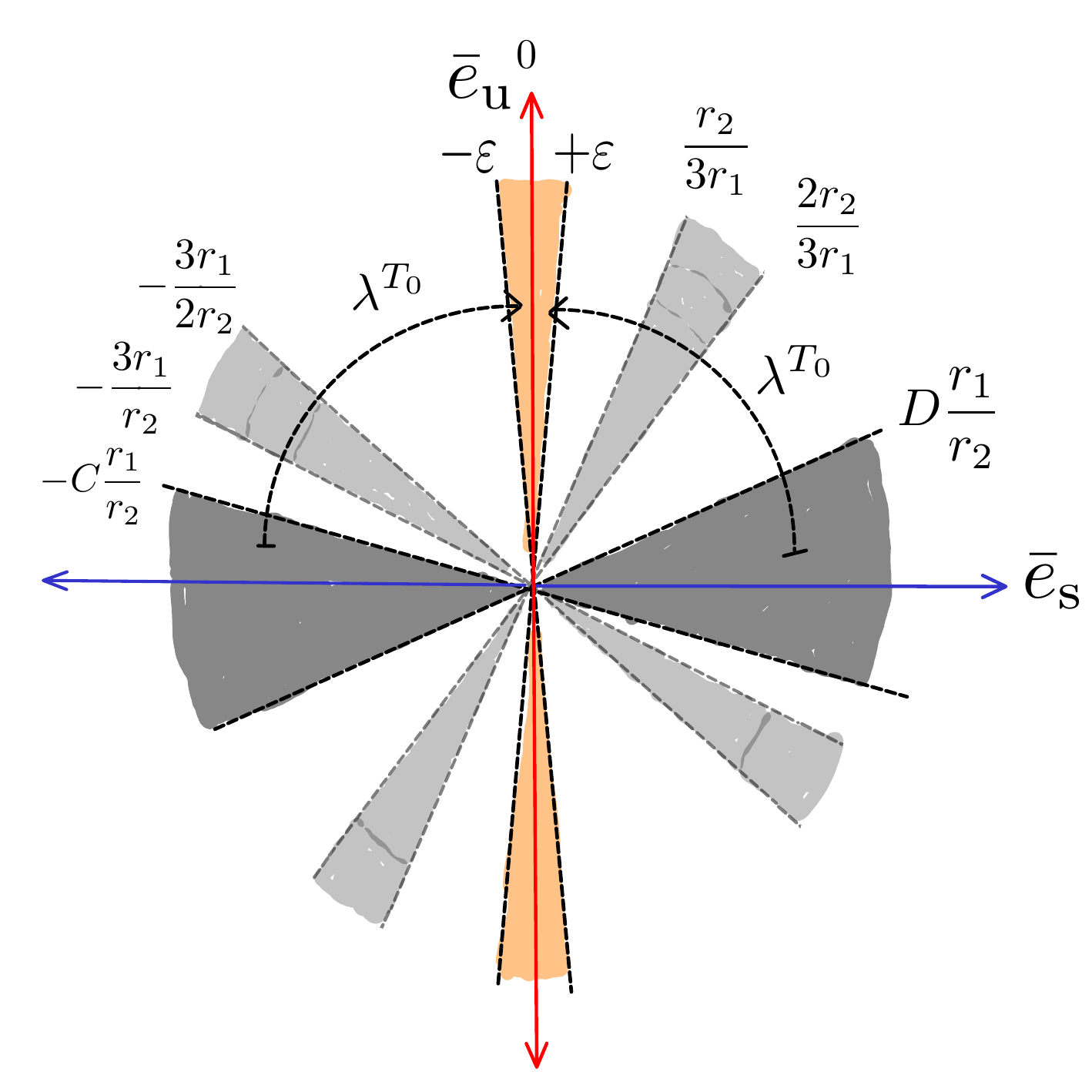}
         \caption{\ }
         \label{fig_slope_functions(1)}
         \end{center}
     \end{subfigure}
     \caption{Cones in the normal bundle}
     \label{fig_slope_functions}
     \end{center}
\end{figure}

\begin{proof}[Proof of Lemma \ref{lema_conos_1}]
We make the proof for the $\mathrm{u}$- case, being analogous the other one. We use the decomposition of $D\psi_t$ as product of matrices \eqref{Psi_matrix_su} and \eqref{su_action_Phi} acting on the normal bundle. To show that the combined action of these matrices preserves the $\mathrm{u}$-cone field, the key fact is that the all the entries of the matrix \eqref{su_action_Phi} are bounded, if we keep bounded the ratio $r_2/r_1$.

We study first how the transformations $\Phi_q$ given in \eqref{su_action_Phi} act on the cones $C^{\mathrm{u}}(q;\delta_0^\mathrm{u},\delta_1^\mathrm{u})$, for every point $q=(r_1,r,s)\in A_\mathrm{in}^1,\ 0\leq r\leq r_2$. Recall from \eqref{jordan_form_Phi} above that each non-trivial $\Phi_q$ has an eigenvector $\bar{w}_1(q)$ (double eigenvalue equal to $1$) whose $\mathrm{u}$-slope is bounded between $-3r_1/r_2$ and $-3r_1/2r_2$, as represented in figure \eqref{fig_slope_functions(2)}. 

From the one hand, there exists a positive constant $\varepsilon=\varepsilon(r_2/r_1)<\min\{3r_1/2r_2,r_2/3r_1\}$, just depending on the ratio $r_2/r_1$, such that:
\begin{equation}\label{conos-cotas_1-seccion-5}
-\frac{3r_1}{r_2}<\Delta_\mathrm{u}(\Phi_q(v))<\frac{r_2}{3r_1},\ \forall\ q\in A_\mathrm{in}^1\ \text{and}\ v\in T_qN\ \text{s.t.}\ -\varepsilon<\Delta_\mathrm{u}(v)<\varepsilon.
\end{equation}
To see this, by \eqref{su_action_Phi} we have that 
$$\Delta_\mathrm{u}(\Phi_q(\bar{e}_\mathrm{u}))=-\left(\frac{-K(r,r_2)}{1-K(r,r_2)}\right)\frac{r_1}{r}.$$
Now, since $K(r,r_2)=-(cte)\cdot|\rho'(r/r_2)|\cdot r/r_2$, it can be checked that the expression between parentheses on the right side of the previous equality is bounded in absolute value by a constant $0<K_0<1$ not depending on $r,r_1,r_2$. Since $K(r,r_2)$ is not trivial just for $r_2/3\leq r\leq 2r_2/3$, we obtain that $-3r_1/r_2<(-3r_1/r_2)\cdot K_0\leq \Delta_\mathrm{u}(\Phi_q(\bar{e}_\mathrm{u}))<0$, for every $0\leq r\leq r_2$. So, taking $\varepsilon<(3r_1/r_2)\cdot (1-K_0)$ it follows that $-3r_1/r_2<\Delta_\mathrm{u}(\Phi_q(v))<0$, for every $v\in T_qN$ satisfying $-\varepsilon<\Delta_\mathrm{u}(v)\leq 0$. It follows directly that $\Delta_\mathrm{u}(\Phi_q(v))<r_2/3r_1$, for every $v\in T_qN$ satisfying $0\leq \Delta_\mathrm{u}(v)\leq \varepsilon$. 

From the other hand, there exist constants $C>3$ and $D>1/3$, only depending on $p$, $m$, $\log(\lambda)$ and $\rho$, satisfying that
\begin{equation}\label{conos-cotas-seccion-5}
-C\frac{r_1}{r_2}<\Delta_\mathrm{u}(\Phi_q(v))<D\frac{r_2}{r_1},\ \forall\ q\in A_\mathrm{in}^1\ \text{and}\ v\in T_qN\ \text{s.t.}\ -\frac{3r_1}{r_2}<\Delta_\mathrm{u}(v)<\frac{r_2}{3r_1}.
\end{equation}
To see this, for every vector $v$ of $\mathrm{u}$-slope $-3r_1/r_2\leq \Delta_\mathrm{u}(v)\leq 0$, by (\ref{su_action_Phi}) we have: 
\begin{equation*}
-\frac{3r_1}{r_2}\cdot\left(\frac{K(r,r_2)(1-\frac{r_2}{3r})+1}{-K(r,r_2)(1-\frac{3r}{r_2})+1}\right) \leq \Delta_\mathrm{u}(\Phi_q(v))\leq 0.
\end{equation*}
Again since $K(r,r_2)=-(cte)\cdot|\rho'(r/r_2)|\cdot r/r_2$, it can be checked that the expression between parenthesis on the left side of the previous equation is bounded, by a bound which is independent from $r_2/r_1$. To see this, denote $t=r/r_2$ and let $\alpha(t)=(cte)\cdot|\rho'(t)|$. Recall from Section \ref{subsection_construction_smooth_model} the particular condition we have required on the function $\rho$, namely, that $\alpha(t)=|pm\log(\lambda)\rho'(t)|<\frac{1/2}{3t^2-t}$ for every $1/3\leq t\leq 2/3$. Using this condition, the expression between parenthesis above becomes
\begin{equation*}
\left(\frac{K(r,r_2)(1-\frac{r_2}{3r})+1}{-K(r,r_2)(1-\frac{3r}{r_2})+1}\right)=
\frac{\alpha(t)(1/3-t)+1}{\alpha(t)(t-3t^2)+1}\leq \frac{\alpha(t)(1/3-t)+1}{1/2}\leq \frac{(1/3)\Vert \alpha \Vert_\infty+1}{1/2},
\end{equation*}
from which we see that there exists a constant $C>0$ (necessarily $C>3$) such that:
\begin{equation*}
\forall\ q\in A_\mathrm{in}^1\text{ with }-\frac{3r_1}{r_2}\leq\Delta_\mathrm{u}(v)\leq 0\ \Rightarrow\  -C\frac{r_1}{r_2}\leq \Delta_\mathrm{u}(\Phi_q(v))\leq 0. 
\end{equation*}
An analogous reasoning applies for vectors satisfying $0\leq \Delta_\mathrm{u}(v)\leq 2r_2/3r_1$, giving the constant $D>1/3$.

We have now three cone distributions $\Cu\left(p;-\varepsilon,\varepsilon\right)\subset\Cu\left(p;\delta^\mathrm{u}_0,\delta^\mathrm{u}_1\right)\subset\Cu\left(p;-C\frac{r_1}{r_2},D\frac{r_2}{r_1}\right)$
defined on each $H(p)=T_pN/H^\mathrm{c}(p)$, as in figure \eqref{fig_slope_functions(1)}. Let $T_0=T_0(r_1,r_2)>0$ be such that 
\begin{equation}\label{cota_general}
-\varepsilon < -\lambda^{2T_0} C\frac{r_1}{r_2}<0<\lambda^{2T_0}D\frac{r_2}{r_1}<\varepsilon,
\end{equation}
so $T_0(r_1,r_2)$ is the time needed for a transformation $\Psi_t$ to send a cone of the form $\Cu(p;-Cr_1/r_2,$ $Dr_2/r_1)$ inside $\Cu(\psi_t(p);-\varepsilon,\varepsilon)$ and depends just on the ratio $r_1/r_2$.

Let $T_1=T_1(r_1,r_2)$ be defined by $$T_1(r_1,r_2)=\min\left\{t>0:\exists\ p\in A_\mathrm{in}^1 \text{ s.t. }p,\psi_t(p)\in\mathrm{supp}(\varphi)\subset A_\mathrm{in}^1\right\}.$$ This is the minimal returning time of points in the region $\{r_1\}\times[r_2/3,2r_2/3]\times\R/\Z\subset A_\mathrm{in}^1$ (where $\Phi_q$ is non-trivial) onto itself. We claim that:

\begin{cl}
If we shrink $r_1$ and $r_2$ keeping constant the ratio $r_2/r_1$, then $T_1(r_1,r_2)$ tends to infinity.
\end{cl}

Assuming this claim, define $T=2T_0(r_1,r_2)$ and shrink both $0<r_2<r_1$ (keeping constant $r_2/r_1$) in such a way that $T_1(r_1,r_2)>T_0(r_1,r_2)$. Let $p\in N$ be a point and let $t\geq T$. Then:
\begin{enumerate}

\item If $[p,\psi_t(p)]\cap A_\mathrm{in}^1=\emptyset$ then the action of $D\psi_t$ on $H$ is just the transformation 
$\Psi_t$. By \eqref{Psi_matrix_su} and \eqref{cota_general} we have:
\begin{align*}
&\Psi_{t}\left(\Cu\left(p;\delta^\mathrm{u}_0,\delta^\mathrm{u}_1\right)\right)=\Cu\left(\psi_{t}(p);-\lambda^{2t}\frac{3r_1}{r_2},\lambda^{2t}\frac{r_2}{3r_1}\right) \subset \Cu\left(\psi_t(p);\delta^\mathrm{u}_0,\delta^\mathrm{u}_1\right).
%
\end{align*}

\item If $[p,\psi_t(p)]\cap A_\mathrm{in}^1$ contains exactly one point $p_1=\psi_{t_1}(p)$, then $D\psi_t=\Psi_{t_2}\circ\Phi_{p_1}\circ\Psi_{t_1}$ with $t_1\geq T_0$ or $t_2\geq T_0$. Using \eqref{conos-cotas_1-seccion-5}, \eqref{conos-cotas-seccion-5} and \eqref{cota_general}, in the first case we have
\begin{align*}
&\Psi_{t_2}\circ\Phi_{p_1}\circ\Psi_{t_1}\left(\Cu\left(p;\delta^\mathrm{u}_0,\delta^\mathrm{u}_1\right)\right)
\subset\Psi_{t_2}\circ\Phi_{p_1}\left(\Cu\left(p_1;-\varepsilon,\varepsilon\right)\right)\\
&\subset\Psi_{t_2}\left(\Cu\left(\varphi(p_1);\delta^\mathrm{u}_0,\delta^\mathrm{u}_1\right)\right)\subset\Cu\left(\psi_{t}(p);\delta^\mathrm{u}_0,\delta^\mathrm{u}_1\right),
%
\end{align*}
while in the second
\begin{align*}
&\Psi_{t_2}\circ\Phi_{p_1}\circ\Psi_{t_1}\left(\Cu\left(p;\delta^\mathrm{u}_0,\delta^\mathrm{u}_1\right)\right)
\subset\Psi_{t_2}\circ\Phi_{p_1}\left(\Cu\left(p_1;\delta^\mathrm{u}_0,\delta^\mathrm{u}_1\right)\right)\\
&\subset\Psi_{t_2}\left(\Cu\left(\varphi(p_1);-C\frac{r_1}{r_2},D\frac{r_2}{r_1}\right)\right)\subset\Cu\left(\psi_{t}(p);-\varepsilon,\varepsilon\right).
%
\end{align*}

\item If $[p,\psi_t(p)]\cap A_\mathrm{in}^1$ consists in $l\geq 2$ points $p_1,\dots,p_l$, write $D\psi_t=\Psi_{t_{l+1}}\circ\Phi_{p_l}\circ\cdots\circ\Phi_{p_1}\circ\Psi_{t_1}$, where $t_k\geq T_1$ for $k=2,\dots,l$. Then by \eqref{conos-cotas_1-seccion-5} and \eqref{conos-cotas-seccion-5} we see that $\Psi_{t_{k+1}}\circ\Phi_{p_k}$ sends the cone $\Cu(p_k;\delta^\mathrm{u}_0,\delta^\mathrm{u}_1)$ inside $\Cu(p_{k+1};-\varepsilon,\varepsilon)$, for $k=1,\dots,l-1$, from where we can conclude the desired property by finite iteration. 
\end{enumerate}
In any case, $D\psi_t(p)$ sends the $\mathrm{u}$-cone in $p$ inside the corresponding $\mathrm{u}$-cone in $\psi_{t}(p)$, for $t\geq T$.

\begin{figure}[t]
\begin{center}
\includegraphics[scale=0.4]{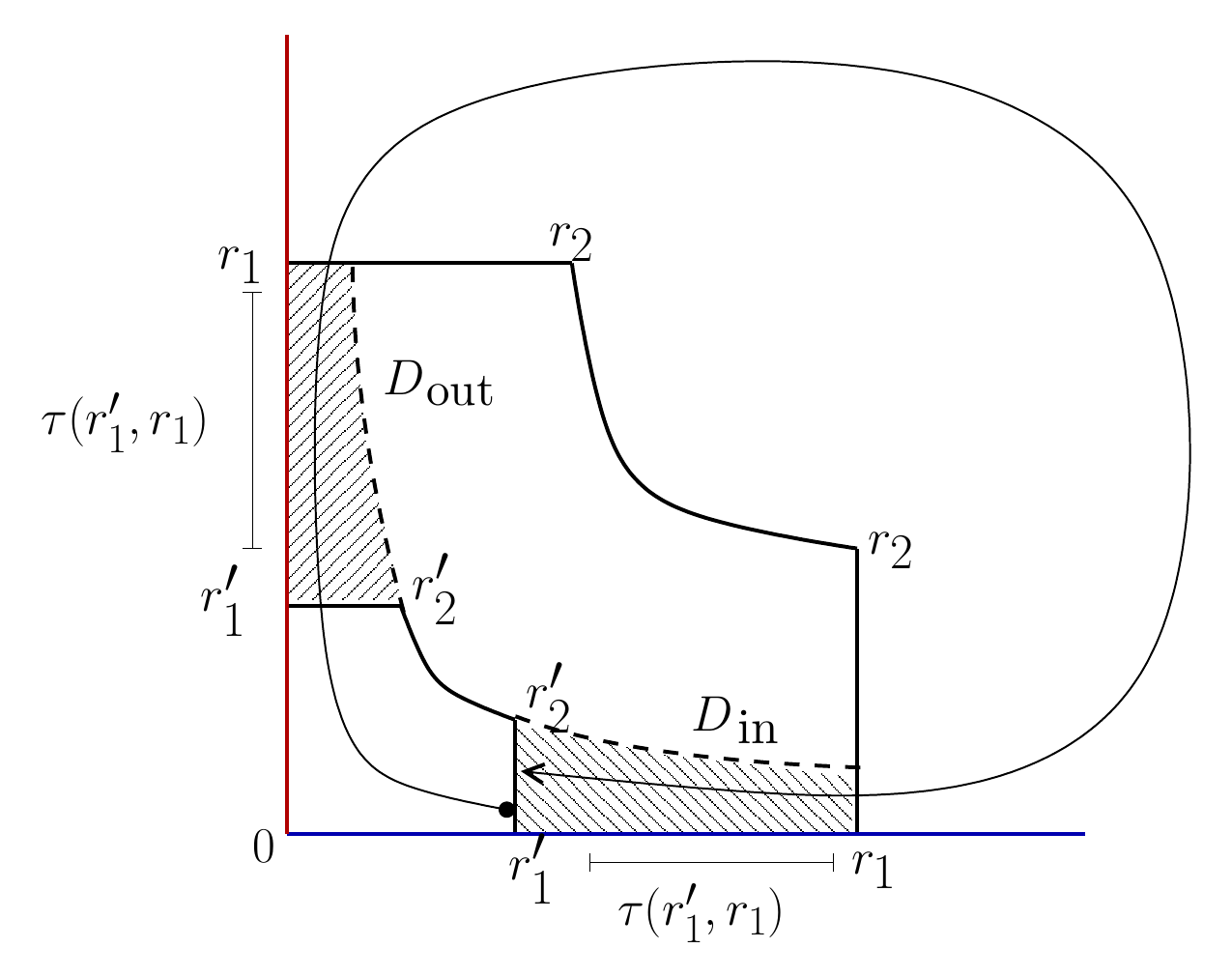}
\caption{The returning time to the annulus at position $r_1\times[0,r_2]\times\R/\Z$ increases when $r_2/r_1\to 0$.}
\label{fig_relation_r1r2}
\end{center}
\end{figure}

To prove the claim consider $r_2'<r_2$ and $r_1'<r_1$ satisfying that $r_2'/r_1'=r_2/r_1$. Then, the manifold $N(r_1',r_2')$ is obtained replacing the original cross-shaped region for a smaller $\mathbb{V}(r_1',r_2')$ in $\R^2\times\R/\Z$. See Figure \ref{fig_relation_r1r2}. Let $p=(r_1',r,s)$ be a point in the entrance annulus ${r_1'}\times[0,r_2']\times\R/\Z$. Observe that if the future orbit of $p$ comes back to the same annulus, it must traverse the regions $D_\mathrm{in}$, $D_\mathrm{out}$ in the present figure. This means that the returning time $T_1(r_1',r_2')$ is bounded from below by twice the time $\tau$ for traversing each of these regions. In addition, observe that the time for traversing each $D_i$ from the entrance boundary to the exit one is given by $\lambda^{\tau}r_1=r_1'$. Therefore, we obtain
$$T_1(r_1',r_2')\geq 2\tau(r_1',r_1)=\frac{2}{\log(\lambda)}\log\left(\frac{r_1'}{r_1}\right),$$ 
which proves the claim. \qedhere
\end{proof}

\begin{proof}[Proof of Lemma \ref{lema_conos_2}]
We use the constants $\varepsilon$, $C$, and $D$ defined in the proof of Lemma \ref{lema_conos_1}. Let $\bar{v}_1,\bar{v}_2\in\Cu(p;\delta_0^\mathrm{u},\delta_1^\mathrm{u})$. By (\ref{Psi_matrix_su}) the forward action of each transformation $\Psi_t$ has the effect of contract the difference of $\mathrm{u}$-slopes between $\bar{v}_1$ and $\bar{v}_2$ by a ratio exponential in $t$. More precisely,
\begin{equation}\label{contraction_slopes}
\left\vert\Delta_\mathrm{u}(\Psi_t(v_2))-\Delta_\mathrm{u}(\Psi_t(v_1))\right\vert
=\lambda^{2t}\left\vert\Delta_\mathrm{u}(v_2)-\Delta_\mathrm{u}(v_1)\right\vert,\ \text{for all}\ t\geq 0.
\end{equation}

For a point $q\in A_\mathrm{in}^1$ consider two vectors $\bar{v}_i=b_i\bar{e}_\mathrm{s}+c_i\bar{e}_\mathrm{u}$ in $H(q)$ and let $(r_1,r,s)$ be the coordinates of $q$. First, we claim that if $|\Delta_\mathrm{u}(v_2)-\Delta_\mathrm{u}(v_1)|<\varepsilon$ then 
\begin{equation}\label{contraction_slopes_2}
\left\vert\Delta_\mathrm{u}(\Phi_q(v_2))-\Delta_\mathrm{u}(\Phi_q(v_1))\right\vert\leq \left\vert\Delta_\mathrm{u}(v_2)-\Delta_\mathrm{u}(v_1)\right\vert.
\end{equation}
To see this, consider a vector $\bar{v}=b\bar{e}_\mathrm{s}+c\bar{e}_\mathrm{u}$ and denote its $\mathrm{u}$-slope by $\delta=b/c$. From the expression of $\Phi_q$ is $\mathrm{su}$-coordinates in \eqref{su_action_Phi} we have that 
\begin{equation*}
\Delta_\mathrm{u}(v)=\frac{(b/c)(K(r,r_2)+1)+K(r,r_2)\frac{r_1}{r}}{-(b/c)K(r,r_2)\frac{r}{r_1}+(1-K(r,r_2))}.
\end{equation*}
Denote by $\delta\mapsto g(\delta)$ the family of rational functions obtained by setting $\delta=b/c$. Then
$$g'(\delta)=\frac{1}{\left(-\delta K(r,r_2)\frac{r}{r_1}+(1-K(r,r_2))\right)^2}<1,\ \text{for all}\ -r_1/r\leq \delta\leq 0.$$
We conclude the claim \eqref{contraction_slopes_2}, since $\varepsilon<r_1/r$ for every $r$ where $\Phi_q\neq \mathrm{id}$. Second, by \eqref{conos-cotas-seccion-5} we can see that $g'$ has no poles for for slopes $-3r_1/r_2\leq \delta\leq 0$, so there is a constant $K_0>0$ such that
\begin{equation}\label{contraction_slopes_3}
\left\vert\Delta_\mathrm{u}(\Phi_q(v_2))-\Delta_\mathrm{u}(\Phi_q(v_1))\right\vert\leq K_0\left\vert\Delta_\mathrm{u}(v_2)-\Delta_\mathrm{u}(v_1)\right\vert,\ \forall\ \bar{v}_1,\bar{v}_2\in\Cu(q;\delta^\mathrm{u}_0,\delta^\mathrm{u}_1).
\end{equation}

Given $p\in N$ and $t\geq 0$, by Lemma \ref{lema_accion_Dpsi} there exist $t_1,\dots,t_{l+1}>0$ and $p_1,\dots,p_{l+1}\in A_\mathrm{in}^1$, such that $t=t_1+\cdots t_{l+1}$, $t_n\geq T_1(r_1,r_2)$ for $n=2,\dots,l$ (cf. proof of \ref{lema_conos_1}) and the action of $D\psi_t(p)$ on $H$ decomposes as a product of transformations $\Psi_{t_n}$ and $\Phi_{p_n}$. Fix two vectors $\bar{v}_1,\bar{v}_2\in\Cu(p,\delta^\mathrm{u}_0,\delta^\mathrm{u}_1)$. Using \eqref{contraction_slopes} and \eqref{contraction_slopes_3} we see that
$$\left\vert\Delta_\mathrm{u}\right(\Phi_{p_1}\circ\Psi_{t_1}(v_2)\left)-\Delta_\mathrm{u}\right(\Phi_{p_1}\circ\Psi_{t_1}(v_1)\left)\right\vert\leq \lambda^{t_1} K_0 \left\vert\Delta_\mathrm{u}(v_2)-\Delta_\mathrm{u}(v_1)\right\vert.$$
For each $i=1,2$ define $\bar{v}_i^1=\Phi_{p_1}\circ\Psi_{t_1}(\bar{v}_i)$ and for $n\geq 2$ define $\bar{v}_i^n=\Phi_{t_n}\circ\Psi_{p_{n}}(\bar{v}_i^{n-1})$. Since $\bar{v}_1^1,\bar{v}_2^1\in\Cu(\varphi(p_1;-Cr_1/r_2,Dr_2/r_1))$ and since $t_n\geq T_1$ for $n\geq 2$, it is verified that $\bar{v}_i^n\in\Cu(-\varepsilon,\varepsilon),\ \forall\ n\geq 2$ (cf. proof of Lemma \ref{lema_conos_1}). In particular, using \eqref{contraction_slopes_2} we have
$$\left\vert\Delta_\mathrm{u}(v_2^{n+1})-\Delta_\mathrm{u}(v_1^{n+1})\right\vert\leq \lambda^{2t_n}\left\vert\Delta_\mathrm{u}(v_2^n)-\Delta_\mathrm{u}(v_1^n)\right\vert,\ \forall\ n\geq 1,$$
from where we deduce
\begin{equation*}
\left\vert\Delta_\mathrm{u}(D\psi_t(p)\cdot v_2)-\Delta_\mathrm{u}(D\psi_t(p)\cdot v_1)\right\vert\leq \lambda^{2\left(t_1+\cdots t_{l+1}\right)}\cdot K_0\cdot 
\left\vert\Delta_\mathrm{u}(v_2)-\Delta_\mathrm{u}(v_1)\right\vert.
\end{equation*} 
Finally, since $|\Delta_\mathrm{u}(v_2)-\Delta_\mathrm{u}(v_1)|\leq \delta^\mathrm{u}_1-\delta^\mathrm{u}_0$, we can bound the right side of this inequality by $\lambda^{2t}L_0$, where $L_0=K_0(\delta^\mathrm{u}_1-\delta^\mathrm{u}_0)>0$. 

This completes the $\mathrm{u}$-case. The analogous reasoning applies for $\Ccs(p;\delta_0^\mathrm{s},\delta_1^\mathrm{s})$ and the backward action of the flow. \qedhere
\end{proof}

\begin{proof}[Proof of Lemma \ref{lema_conos_3}]
Let $p\in N$, $\bar{v}\in T_pN/H^\mathrm{c}(p)$ and $t\geq 0$. Let us call $\delta=\Delta_\mathrm{u}(v)=b/c$ and assume that $-\infty<\delta<\infty$. 

From the one hand, let $q\in A_\mathrm{in}^1$ and $\bar{v}\in T_qN/H^\mathrm{c}(q)$. Since $q\mapsto\Phi_q$ has compact support, by taking $R_0=\min\left\{\left(\Vert\Phi_{q}^{-1}\Vert_{\infty}\right)^{-1}:q\in A_\mathrm{in}^1\right\}>0$ we have
\begin{equation}\label{growing_vectors_Phi_ii}
\Vert\Phi_q(v)\Vert_\mathrm{su}>R_0 \Vert v\Vert_\mathrm{su},\ \text{for all}\ \bar{v}\in T_qN/H^\mathrm{c}(q).
\end{equation}
Moreover, consider $0<\varepsilon<3r_1/2r_2$ defined in the proof of Lemma \ref{lema_conos_1}. From the expression \eqref{jordan_form_Phi} it follows that
\begin{equation}\label{growing_vectors_Phi}
\Vert\Phi_q(v)\Vert_{\mathrm{su}}>\Vert v\Vert_{\mathrm{su}},\ \text{for all}\ \bar{v}\in\Cu(q;-\varepsilon,\varepsilon).
\end{equation}

From the other hand, since the transformations $\Psi_t$ on the normal bundle correspond to the hyperbolic matrices \eqref{Psi_matrix_su}, it follows that for every $\bar{v}\in TN/H$ there exists a quantity $0<Q(\delta)=(\delta^2+1)^{-1}\leq 1$ such that
\begin{equation}\label{growing_vectors_Psi}
\Vert\Psi_t(v)\Vert_{\mathrm{su}}\geq Q(\delta)\lambda^{-t}\Vert v\Vert_{\mathrm{su}},\ \forall\ t\geq 0.
\end{equation}
To see an expansion on the $\mathrm{su}$-norm by an application of $\Psi_t$ on $\bar{v}$, we need to wait some time depending on the $\mathrm{u}$-slope $\delta=\Delta_\mathrm{u}(v)$. It can be seen that, if $|\delta|\leq 1$ then $|\Psi_t(v)|_{\mathrm{su}}>|v|_{\mathrm{su}}$, $\forall\ t\geq 0$, and if $|\delta|>1$ then $|\Psi_t(v)|_{\mathrm{su}}>|v|_{\mathrm{su}}$ if and only if $\Delta_\mathrm{u}(\Psi_t(v))>1/|\delta|$. 

Consider $T_0=T_0(r_1,r_2)>0$ such that 
$$-\frac{1}{C}\frac{r_2}{r_1}<-\lambda^{2T_0(r_1,r_2)}C\frac{r_1}{r_2}<0<\lambda^{2T_0(r_1,r_2)}D\frac{r_2}{r_1}<\frac{1}{D}\frac{r_1}{r_2},$$
where $C,D>0$ are the constants defined along the proof of Lemma \ref{lema_conos_1}. Defined in this way, $T_0$ is the minimal time needed to see an expansion of $\Vert\Psi_t(v)\Vert_{\mathrm{su}}$ for vectors $\bar{v}$ in the cone $\Cu(p;-Cr_1/r_2,$ $Dr_2/r_1)$. 
For such a choice of $T_0$ we have that
\begin{equation*}
\Vert\Psi_{T_0}(v)\Vert_{\mathrm{su}}\geq Q_0\lambda^{-T_0}\Vert v\Vert_{\mathrm{su}}>\Vert v\Vert_{\mathrm{su}},\ \text{for all}\ -Cr_1/r_2 <\Delta_\mathrm{u}(v)< Dr_2/r_1,
\end{equation*}
where $Q_0=\min\left\{Q(-Cr_1/r_2),Q(Dr_2/r_1)\right\}$. 
We define $\mu=Q_0^{1/T_0}\lambda^{-1}$, which satisfies $1<\mu<1/\lambda$. 

Since $T_0(r_1,r_2)$ depends just on $r_1/r_2$, we can shrink both $r_1$ and $r_2$, keeping constant its ratio, in such a way that the minimal returning time $T_1=T_1(r_1,r_2)$ of points in the support of the surgery to itself is greater that $T_0(r_1,r_2)$, cf. proof of Lemma \ref{lema_conos_1}. 

Assuming this condition, let $p\in N$, $t\geq 0$, $\bar{v}\in\Cu(p;\delta_0^\mathrm{u},\delta_1^\mathrm{u})$ and consider the decomposition
$D\psi_t(p)\cdot \bar{v}=\Psi_{t_{l+1}}\cdot\Phi_{p_l}\cdots\Phi_{p_1}\cdot\Psi_{t_1}(\bar{v})$, where $l\geq 0$, $t_k\geq T_1\geq T_0$ for $k=2,\dots,l$ and $p_k\in A_\mathrm{in}^1$. Define $\bar{v}^1=\Phi_{p_1}\circ\Psi_{t_1}(\bar{v})$ and for $k=2,\dots,l$ define $\bar{v}^k=\Psi_{t_k}\circ\Phi_{p_{k}}(\bar{v}^{k-1})$. It follows that
\begin{enumerate}[(i)]
\item By \eqref{growing_vectors_Psi} and \eqref{growing_vectors_Phi_ii} then
$$\Vert v^1\Vert_{\mathrm{su}}\geq \left(\Vert\Phi_{p_1}^{-1}\Vert_{\infty}\right)^{-1} \Vert\Psi_{t_1}(v)\Vert_{\mathrm{su}}\geq R_0 Q_0\cdot\lambda^{-t_1}\Vert v\Vert_{\mathrm{su}}.$$

\item Since $\bar{v}^1\in\Cu(p_1;-Cr_1/r_2,Dr_2/r_1)$ and $t_k\geq T_0$ for $k=2,\dots,l$ then $\bar{v}^k\in\Cu(p_k;-\varepsilon;\varepsilon)$ (cf. proof of Lemma \ref{lema_conos_1}). By \eqref{growing_vectors_Phi} and \eqref{growing_vectors_Psi} then
$$\Vert v^k\Vert_{\mathrm{su}}\geq Q_0\lambda^{-t_k}\Vert v^{k-1}\Vert_{\mathrm{su}},\ \text{for}\  k=2,\dots,l.$$
\end{enumerate}
Therefore, putting together (i) and (ii) above in the decomposition $D\psi_t=\Psi_{t_{l+1}}\cdot\Phi_{p_l}\cdots\Phi_{p_1}\cdot\Psi_{t_1}$, and using the fact that $0\leq l\leq t/T_1$, we obtain 
$$\Vert D\psi_t(p)\cdot v\Vert_{\mathrm{su}}\geq R_0Q_0^{l+1}\lambda^{-(t_1+\cdots+t_{k+1})}\Vert v\Vert_{\mathrm{su}}\geq 
(R_0Q_0)Q_0^{t/T_1}\lambda^{-t}\Vert v\Vert_{\mathrm{su}}\geq L\cdot\mu^{t}\Vert v\Vert_{\mathrm{su}}.
$$
for some constant $0<L< R_0Q_0$. \qedhere
\end{proof}

\subsection{Equivalence with the original flow}\label{subsection_smooth_model_is_equivalent}
The general setting at the beginning of Section \ref{section_hyperbolic_models} is a transitive topologically Anosov flow $(\phi,M)$ equipped with a Birkhoff section $\iota:(\Sigma,\partial\Sigma)\to (M,\Gamma)$ which, for simplicity, we have assumed to have only one orbit $\gamma\in\Gamma$. In Section \ref{subsection_construction_smooth_model} we have constructed a smooth flow $(\psi,N)$ depending on the combinatorial parameters $p(\gamma,\Sigma)$, $n(\gamma,\Sigma)$, $m(\gamma,\Sigma)$ of the Birkhoff section and on two real parameters $0<r_2<r_1<1$. Then, in Section \ref{subsection_smooth_model_is_Anosov} we have showed that if $r_1,r_2$ are chosen small enough, the flow $(\psi,N)$ is Anosov. 

We show here that if $0<r_2<r_1<1$ are sufficiently small such that $(\psi,N)$ is Anosov, then $(\psi,N)$ is orbitally equivalent to the original flow $(\phi,M)$. We do it by checking the criterion provided in Theorem \textbf{B}.

\begin{proposition}\label{prop_subsection_smooth_model_is_equivalent}
The flow $(\psi,N)$ constructed in Proposition \ref{lemma_defn_psi} is orbitally equivalent to the original flow 
$(\phi,M)$, provided $0<r_2<r_1<1$ are sufficiently small.
\end{proposition}

To check the criterion in Theorem \textbf{B}, the flow $(\phi,M)$ is already endowed with a Birkhoff section and we have to show that $(\psi,N)$ admits a (homeomorphic) Birkhoff section, that carries the same first return action on fundamental group and the same combinatorial data on the boundary. 

It is convenient to recall the general construction of the manifold $N=M_0\sqcup_\varphi M_1$. Using the charts $\{\Pi_i:i=1,\dots,4\}$ from Proposition \ref{lema_normal_coordinates} we have defined a family of compact tubular neighborhoods $R(r_1,r_2)\subset M$ of the orbit $\gamma$. Then $N$ is constructed by removing from $M$ the interior of such a tubular neighborhood, obtaining a manifold $M_0$ with boundary $\partial M_0=\partial R(r_1,r_2)$, and then gluing along the boundary the solid torus $M_1=\mathbb{V}(r_1,r_2)$. The gluing map $\varphi:\partial M_0\to\partial M_1$ given in definition \eqref{glueing_map_negative_linking_}-\eqref{glueing_map_positive_linking_} has support contained in the annulus $A_\mathrm{in}^1$ inside the first quadrant of $R(r_1,r_2)$, or in $A_\mathrm{in}^4$ in the fourth quadrant, depending on the signature of $m=m(\gamma,\Sigma)$.

\begin{lem}\label{lema_birk_sect_for_psi_t}
There exists an immersion $\zeta:\Sigma\to N$ verifying that:
\begin{enumerate}[i.]
\item $\zeta:(\Sigma,\partial\Sigma)\to(N,\gamma_0)$ is a Birkhoff section that sends the boundary components of $\Sigma$ onto the periodic orbit $\gamma_0$ contained in $M_1$.

\item If we call $\Sigma'=\zeta(\Sigma)$, then $n(\gamma_0,\Sigma')=n(\gamma,\Sigma)$, $m(\gamma_0,\Sigma')=m(\gamma,\Sigma)$ and $p(\gamma_0,\Sigma')=p(\gamma,\Sigma)$.

\item $\zeta(\Sigma)\cap M_0\equiv \iota(\Sigma)\cap M_0$.
\end{enumerate}
\end{lem}

\begin{proof}
We prove the lemma for the case $m=m(\gamma,\Sigma)<0$, the other case being analogous. For simplicity, we assume as well that $p(\gamma,\Sigma)=1$.

Since $N=M_0\sqcup_\varphi M_1$ and $M=M_0\sqcup_{id} R(r_1,r_2)$, we consider the component $M_0$ as included in both $N$ and $M$. Define $\Sigma_0=\iota(\Sigma)\cap M_0$, which is the part of the original Birkhoff section that lies outside the neighborhood $R(r_1,r_2)$. To prove the lemma, we show that $\Sigma_0$ can be extended inside the manifold $M_1$, adding an helicoidal-like surface $S$ that connects $\partial\Sigma_0$ with $\gamma_0$, in such a way to obtain the desired Birkhoff section.

Let $\alpha=\partial\Sigma_0=\Sigma_0\pitchfork \partial M_0$. By construction of $R(r_1,r_2)$, this curve is a piecewise smooth, simple, closed curve in $\partial R(r_1,r_2)$. The partition into quadrants of $R(r_1,r_2)$ originates a partition of $\alpha$ into segments
$$\left\{\alpha_\mathrm{in}^{i+4j},\alpha_\mathrm{tg}^{i+4j},\alpha_\mathrm{out}^{i+4j};\ i=1,\dots,4,\ k=0,\dots,4n-1\right\}$$
of constant $\R/\Z$-coordinate, as in the left part of figure \eqref{fig_glueing_map}. We choose the supra-index $i+4j$ in such a way that each $\alpha^{i+4j}_{*}$ is contained in the quadrant $R_i(r_1,r_2)$, $i=1,\dots,4$, and we use the sub-index $in,\ tg,\ out$ to indicate if the segment belongs to the region of $\partial R(r_1,r_2)$ where the flow enters, is tangent or escapes the neighborhood $R(r_1,r_2)$, respectively. 

Let $\beta=\varphi(\alpha)$, which is a piecewise smooth, simple, closed curve contained in $\partial\mathbb{V}(r_1,r_2)$. In an analogous fashion, the curve $\beta$ can be decomposed in a concatenation of compact segments
$$\left\{\beta_\mathrm{in}^{i+4j},\beta_\mathrm{tg}^{i+4j},\beta_\mathrm{out}^{i+4j};\ i=1,\dots,4,\ k=0,\dots,4n-1\right\},$$
as in figures \eqref{fig_glueing_map} or \eqref{fig_glueing_map_Birk_section_3}. Since $\varphi$ restricts to a monotonous twist on the annulus $A_\mathrm{in}^1$ and to the identity on the complement, we find that each segment $\beta_\mathrm{in}^{1+4j}$ (in the first quadrant) admits a parametrization with monotonous $\R/\Z$-coordinate, while every other segment $\beta_{*}^{i+4j}=\varphi(\alpha_{*}^{i+4j})=\alpha_{*}^{i+4j}$ has constant $\R/\Z$-coordinate.

\begin{cl}
There exists an immersion $S\hookrightarrow\mathbb{V}(r_1,r_2)$, where $S$ is a compact annulus, satisfying that:
\begin{enumerate}

\item One boundary component of $S$ coincides with $\beta=S\cap\partial\mathbb{V}(r_1,r_2)$,

\item The immersion is actually a \emph{local} Birkhoff section at $\gamma_0$ for the flow $\psi$, with $n(\gamma_0,S)=n(\gamma,\Sigma)$ and $m(\gamma_0,S)=m(\gamma,\Sigma)$.
\end{enumerate}
\end{cl}

\begin{proof}[Proof of the claim]
Define
$S=\left\{(\theta x,\theta y,z)\in\R^2/\times\R/\Z\ :\ (x,y,z)\in\beta,\ 0\leq\theta\leq 1\right\}.$ 
Then $S$ is the surface obtained by joining with a straight segment each point $(x,y,z)\in\beta$ with the point $(0,0,z)\in\gamma_0$. It is homeomorphic to a compact annulus and clearly $\partial S=\beta\cup\gamma_0$. 

The set $S$ can be decomposed as the union of some smooth horizontal surfaces, each one isometric to the region $Q(r_1,r_2)$ defined in Section \ref{subsection_linear_local_model}, and some smooth non-horizontal bands, as we see in figure \eqref{fig_glueing_map_Birk_section_3}. Each segment $\beta_\mathrm{in}^{1+4j}$ belongs to the boundary of a band, while every other segment composing $\beta$ is contained in the boundary of a horizontal surface.  

\begin{figure}
\begin{center}
\includegraphics[scale=0.35,angle=0]{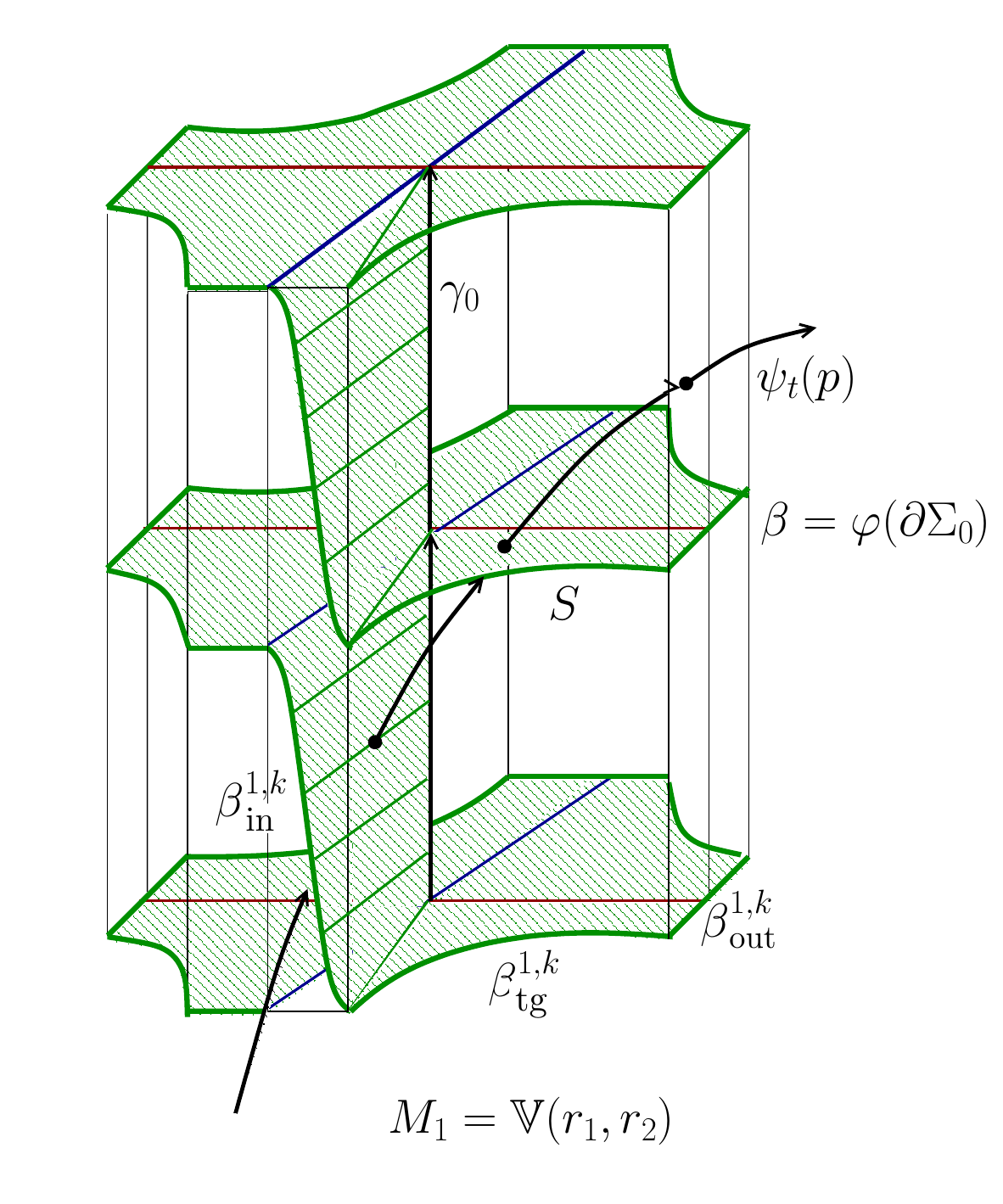}
\caption{Local Birkhoff $S$ section around $\gamma_0\subset\mathbb{V}(r_1,r_2)$.}
\label{fig_glueing_map_Birk_section_3}
\end{center}
\end{figure}

In the complement of $\gamma_0$ the surface $S$ is transverse to the vector field $X_1$. To see this, recall that $X_1$ is defined to be the vector field 
$(x,y,z)\mapsto\left(\log(\lambda)x,-\log(\lambda)y,{1}/{n}\right)$
in $\R^2\times\R/\Z.$ Since the third component is nowhere zero, it follows that it is transverse to $S$ along the horizontal parts. For the transversality along the bands, consider a point $p=(r_1,r,s)\in S\cap A_\mathrm{in}^1$. On $A_\mathrm{in}^1$ there is a parametrization $r\mapsto\beta(r)=\left(r_1,r,z_0+(|m|/n)\cdot\rho(r/r_2)\right)$, where $z_0$ is some constant. From the definition of $S$ we have that $T_pS$ is generated by the vectors 
$$\frac{\partial}{\partial r}(p)=(0,1,-\kappa(r))\ \ \text{and}\ \ \frac{\partial}{\partial\theta}(p)=(r_1,r,0),\text{ where }\kappa(r)=\frac{|m|}{n}|\rho'(r/r_2)|\geq 0.$$
It follows that
\begin{align*}
X_1(p)\wedge\frac{\partial}{\partial r}(p)\wedge\frac{\partial}{\partial \theta}(p)=
\begin{vmatrix}
\log(\lambda)r_1 & 0          & r_1\\
-\log(\lambda)r  & 1          & r  \\
1/n            & -\kappa(r) & 0    
\end{vmatrix}
=-\frac{r_1}{n}+2\kappa(r)\log(\lambda)r_1r<0,
\end{align*}
so we conclude that $X_1$ is transverse to $TS$ along the the arc $\beta_\mathrm{in}^{1+4j}$. Then, since the surface $S$ is the image of a $1$-parameter family of horizontal homotheties 
$(x,y,z)\mapsto(\theta x,\theta y,z),$
and since $X_1$ is invariant under these transformations, this implies the transversality of $X_1$ with $TS$ along all the interior of the band.\footnote{It is here that it is important to consider a separate definition \eqref{glueing_map_negative_linking_}-\eqref{glueing_map_positive_linking_} for $\varphi$, discriminating by the signature of the linking number. Observe that if we apply the same formula for $\varphi$ in the first quadrant with positive $n$, then the bands switch to a bad position and transversality cannot be guaranteed.} 

By construction the surface $S$ is orientable, and it is not hard to see that $X_1$ traverse each band or horizontal surface always in the same sense. Thus, the orbits of $\psi$ are everywhere (topologically) transverse to the interior of $S$. In a meridian/longitude basis $\{a,b\}$ of $H_1(M_1\backslash\gamma_0)$ the homological coordinates of $\beta$ are $[\beta]=n(\gamma,\Sigma)\cdot a +m(\gamma,\Sigma)\cdot b$, from where we deduce the linking number and multiplicity of $S$. Since $n(\gamma,\Sigma)\neq 0$, there exists a tubular neighborhood $O$ of $\gamma_0$ such that every point in $O$ intersects the surface $S$ in a uniformly bounded time. Therefore, $S$ is a local Birkhoff section at $\gamma_0$. This completes the claim. 
\end{proof}

To complete the proof of the lemma, consider the set $\Sigma'\coloneqq \Sigma_0\cup S$
that is contained in $N$. Then $\Sigma'$ is the image of a continuous immersion $\zeta:\Sigma\to N$, which is an embedding on $\mathring{\Sigma}$ and coincides with $\iota$ over $\Sigma_0$. From the construction it can be seen that its interior $\mathring{\Sigma}'=\Sigma'\backslash\partial\Sigma'$ is (topologically) transverse to the flow lines. To prove that this is actually a Birkhoff section, it remains to show that all the $\psi$-orbits intersect $\Sigma'$ in a uniformly bounded time. For this, consider two tubular neighborhoods $O_1\subset O_0\subset N$ of $\gamma_0$ satisfying that $M_0\cap M_1$ is contained in the interior of $O_0\backslash O_1$, and that there exists $T_1,T_0>0$ such that $[p,\psi_{T_1}(p)]\cap S\neq\emptyset$ for every $p\in O_1$ and $[p,\psi_{T_0}(p)]\cap \Sigma_0\neq\emptyset$ for every $p\in N\backslash O_0$. Then, since the neighborhoods $O_1\subset O_0$ are nested germs of a saddle type periodic orbit, we see that there exits some $T_2>0$ such that $[p,\psi_{T_2}(p)]$ is not contained in $O_0\backslash O_1$, for every $p\in O_0\backslash O_1$. Then, taking $T>\max\{T_0+T_2,T_1+T_2\}$, we deduce that $[p,\psi_{T}(p)]\cap \Sigma'\neq\emptyset$, for every $p\in N$. This completes the proof of the lemma. 
\end{proof}

For the next statement, it is convenient to consider $\Sigma$ as a fixed compact surface, for which we have two immersions $\iota:(\Sigma,\partial\Sigma)\to(M,\gamma)$ and $\zeta:(\Sigma,\partial\Sigma)\to(N,\gamma_0)$, that are Birkhoff sections for $(\phi,M)$ and $(\psi,N)$, respectively. Denote by $P:\mathring{\Sigma}\to\mathring{\Sigma}$ and $P':\mathring{\Sigma}\to\mathring{\Sigma}$ the corresponding first return maps.  

\begin{lem}\label{lema_birk_sect_for_psi_t-fund_group}
The homeomorphisms $P$ and $P'$ define the same action on $\pi_1(\mathring{\Sigma})$. 
\end{lem}

\begin{proof}
By construction, both flows $\{\phi_t:M\to M\}_{t\in\R}$ and $\{\psi_t:N\to N\}_{t\in\R}$ induce the same foliation by orbit segments in $M_0$. Moreover, if $p\in M_0$ and $t\geq 0$ satisfy that $[p,\psi_t(p)]\subset M_0$ or $[p,\phi_t(p)]\subset M_0$ then $\psi_s(p)=\phi_s(p)$, for $0\leq s\leq t$. Let $U\subset\Sigma_0$ be a collar neighborhood of $\partial\Sigma_0$ and consider the subsurface $\Sigma_U=\Sigma_0\backslash U$. If we take $U$ big enough, then it is satisfied that every point $p\in\Sigma_U$ has a first return $\psi_{\tau}(p)$ in the interior $\mathring{\Sigma_0}$ and all the orbit segment $[p,\psi_\tau(p)]$ is contained in $M_0$. Therefore, $P(p)=P'(p)$, for every $p\in\Sigma_U$.

Since there is a deformation retraction $\Sigma\to\Sigma_0\to\Sigma_U$, every homotopy class $[\alpha]$ in $\pi_1(\Sigma)$ has a representative $\alpha_0\subset\Sigma_U$, and from the previous paragraph it follows that $P(\alpha_0)=P'(\alpha_0)$. Thus, $P_*([\alpha])=P_*'([\alpha])$, for every $[\alpha]\in \pi_1(\Sigma)$.
\end{proof}

\begin{proof}[Proof of Proposition \ref{prop_subsection_smooth_model_is_equivalent}]
Choose $0<r_2<r_1$ small enough such that $\{\psi_t:N\to N\}_{t\in\R}$ is Anosov. By Lemma \ref{lema_birk_sect_for_psi_t} above, we have two topologically Anosov flows $(\phi,M)$ and $(\psi,N)$ equipped with Birkhoff sections $\iota:(\Sigma,\partial\Sigma)\to(M,\gamma)$ and $\zeta:(\Sigma,\partial\Sigma)\to(N,\gamma_0)$ respectively. By construction $n(\gamma,\Sigma)=n(\gamma_0,\Sigma')$, $m(\gamma,\Sigma)=m(\gamma_0,\Sigma')$ and $p(\gamma,\Sigma)=p(\gamma_0,\Sigma')$. By Lemma \ref{lema_birk_sect_for_psi_t-fund_group}, the identity map on $\Sigma$ conjugates the corresponding first return actions on fundamental groups, so the hypotheses of Theorem \textbf{B} hold. We conclude that $(\phi,M)$ and $(\psi,N)$ are orbitally equivalent flows.
\end{proof}


\backmatter
\bibliographystyle{smfplain}
\bibliography{bibliografia_thesis}
\end{document}

%% file: fig_1_2_1.pdf_t
\begin{picture}(0,0)%
\includegraphics{fig_1_2_1.pdf}%
\end{picture}%
\setlength{\unitlength}{4144sp}%
\begingroup\makeatletter\ifx\SetFigFont\undefined%
\gdef\SetFigFont#1#2#3#4#5{%
  \reset@font\fontsize{#1}{#2pt}%
  \fontfamily{#3}\fontseries{#4}\fontshape{#5}%
  \selectfont}%
\fi\endgroup%
\begin{picture}(14556,9426)(868,-9454)
\put(12466,-2221){\makebox(0,0)[lb]{\smash{{\SetFigFont{25}{30.0}{\rmdefault}{\mddefault}{\updefault}{\color[rgb]{0,0,0}$x$}%
}}}}
\put(9496,-5551){\makebox(0,0)[lb]{\smash{{\SetFigFont{29}{34.8}{\rmdefault}{\mddefault}{\updefault}{\color[rgb]{0,0,0}$U$}%
}}}}
\put(12151,-5821){\makebox(0,0)[lb]{\smash{{\SetFigFont{29}{34.8}{\rmdefault}{\mddefault}{\updefault}{\color[rgb]{0,0,0}$\Sigma$}%
}}}}
\put(1486,-9016){\makebox(0,0)[lb]{\smash{{\SetFigFont{29}{34.8}{\rmdefault}{\mddefault}{\updefault}{\color[rgb]{0,0,0}$\hat{\mathcal{U}}$}%
}}}}
\put(11116,-8971){\makebox(0,0)[lb]{\smash{{\SetFigFont{29}{34.8}{\rmdefault}{\mddefault}{\updefault}{\color[rgb]{0,0,0}$\mathcal{U}$}%
}}}}
\put(3016,-7441){\makebox(0,0)[lb]{\smash{{\SetFigFont{29}{34.8}{\rmdefault}{\mddefault}{\updefault}{\color[rgb]{0,0,0}$(u,0)$}%
}}}}
\put(2971,-2941){\makebox(0,0)[lb]{\smash{{\SetFigFont{29}{34.8}{\rmdefault}{\mddefault}{\updefault}{\color[rgb]{0,0,0}$(u,\tau(u))$}%
}}}}
\put(13501,-7306){\makebox(0,0)[lb]{\smash{{\SetFigFont{29}{34.8}{\rmdefault}{\mddefault}{\updefault}{\color[rgb]{0,0,0}}%
}}}}
\put(9946,-4696){\makebox(0,0)[lb]{\smash{{\SetFigFont{29}{34.8}{\rmdefault}{\mddefault}{\updefault}{\color[rgb]{0,0,0}$u$}%
}}}}
\put(11476,-5011){\makebox(0,0)[lb]{\smash{{\SetFigFont{29}{34.8}{\rmdefault}{\mddefault}{\updefault}{\color[rgb]{0,0,0}$P(u)$}%
}}}}
\put(5266,-4606){\makebox(0,0)[lb]{\smash{{\SetFigFont{29}{34.8}{\rmdefault}{\mddefault}{\updefault}{\color[rgb]{0,0,0}$\varphi$}%
}}}}
\end{picture}%

%% file: main_-_Apr_10.bbl
\providecommand{\bysame}{\leavevmode ---\ }
\providecommand{\og}{``}
\providecommand{\fg}{''}
\providecommand{\smfandname}{\&}
\providecommand{\smfedsname}{\'eds.}
\providecommand{\smfedname}{\'ed.}
\providecommand{\smfmastersthesisname}{M\'emoire}
\providecommand{\smfphdthesisname}{Th\`ese}
\begin{thebibliography}{10}

\bibitem{agol-tsang_Veering-triangulations_AGT}
{\scshape I.~Agol {\normalfont \smfandname} C.~C. Tsang} -- {\og Dynamics of
  veering triangulations: infinitesimal components of their flow graphs and
  applications\fg}, \emph{Algebr. Geom. Topol.} \textbf{24} (2024), no.~6,
  p.~3401--3453.

\bibitem{aoki-hiraide}
{\scshape N.~Aoki {\normalfont \smfandname} K.~Hiraide} -- \emph{Topological
  theory of dynamical systems}, North-Holland Mathematical Library, vol.~52,
  North-Holland Publishing Co., Amsterdam, 1994.

\bibitem{asaoka_invariant_volumes}
{\scshape M.~Asaoka} -- {\og On invariant volumes of codimension-one {A}nosov
  flows and the {V}erjovsky conjecture\fg}, \emph{Invent. Math.} \textbf{174}
  (2008), no.~2, p.~435--462.

\bibitem{Asaoka_projectively_Anosov}
\bysame , {\og Regular projectively {A}nosov flows on three-dimensional
  manifolds\fg}, \emph{Ann. Inst. Fourier (Grenoble)} \textbf{60} (2010),
  no.~5, p.~1649--1684.

\bibitem{asaoka2022oriented}
{\scshape M.~Asaoka, C.~Bonatti {\normalfont \smfandname} T.~Marty} -- {\og
  Oriented {B}irkhoff sections of {A}nosov flows\fg}, \emph{arxiv 2212.06483}
  (2022).

\bibitem{barbot_charact_flots_Anosov}
{\scshape T.~Barbot} -- {\og Caract\'erisation des flots d'{A}nosov en
  dimension 3 par leurs feuilletages faibles\fg}, \emph{Ergodic Theory Dynam.
  Systems} \textbf{15} (1995), no.~2, p.~247--270.

\bibitem{bonatti-guelman}
{\scshape C.~Bonatti {\normalfont \smfandname} N.~Guelman} -- {\og Axiom {A}
  diffeomorphisms derived from {A}nosov flows\fg}, \emph{J. Mod. Dyn.}
  \textbf{4} (2010), no.~1, p.~1--63.

\bibitem{bonatti-wilkinson}
{\scshape C.~Bonatti {\normalfont \smfandname} A.~Wilkinson} -- {\og Transitive
  partially hyperbolic diffeomorphisms on 3-manifolds\fg}, \emph{Topology}
  \textbf{44} (2005), no.~3, p.~475--508.

\bibitem{brocker-janich_diff-top}
{\scshape T.~Br\"{o}cker {\normalfont \smfandname} K.~J\"{a}nich} --
  \emph{Introduction to differential topology}, Cambridge University Press,
  Cambridge-New York, 1982.

\bibitem{brunella_thesis}
{\scshape M.~Brunella} -- {\og Expansive flows on three - manifolds.\fg},
  \smfphdthesisname, SISSA, Trieste, 1992.

\bibitem{brunella_top_equivalence_criterion}
\bysame , {\og On the topological equivalence between {A}nosov flows on
  three-manifolds\fg}, \emph{Comment. Math. Helv.} \textbf{67} (1992), no.~3,
  p.~459--470.

\bibitem{brunella_expansivos_seifert}
\bysame , {\og Expansive flows on {S}eifert manifolds and on torus bundles\fg},
  \emph{Bol. Soc. Brasil. Mat. (N.S.)} \textbf{24} (1993), no.~1, p.~89--104.

\bibitem{Doering}
{\scshape C.~I. Doering} -- {\og Persistently transitive vector fields on
  three-dimensional manifolds\fg}, \emph{Dynamical systems and bifurcation
  theory ({R}io de {J}aneiro, 1985)} \textbf{160} (1987), p.~59--89.

\bibitem{farrell-jones_Anosov_exotic}
{\scshape F.~T. Farrell {\normalfont \smfandname} L.~E. Jones} -- {\og Anosov
  diffeomorphisms constructed from {$\pi _{1}\,{\rm Diff}\,(S^{n})$}\fg},
  \emph{Topology} \textbf{17} (1978), no.~3, p.~273--282.

\bibitem{FLP}
{\scshape A.~Fathi, F.~Laudenbach {\normalfont \smfandname} V.~Po\'{e}naru} --
  \emph{Thurston's work on surfaces}, Mathematical Notes, vol.~48, Princeton
  University Press, Princeton, NJ, 2012.

\bibitem{fenley_Anosov}
{\scshape S.~R. Fenley} -- {\og Anosov flows in {$3$}-manifolds\fg}, \emph{Ann.
  of Math. (2)} \textbf{139} (1994), no.~1, p.~79--115.

\bibitem{fried_surgery}
{\scshape D.~Fried} -- {\og Transitive {A}nosov flows and pseudo-{A}nosov
  maps\fg}, \emph{Topology} \textbf{22} (1983), no.~3, p.~299--303.

\bibitem{gerber-katok}
{\scshape M.~Gerber {\normalfont \smfandname} A.~Katok} -- {\og Smooth models
  of {T}hurston's pseudo-{A}nosov maps\fg}, \emph{Ann. Sci. \'{E}cole Norm.
  Sup. (4)} \textbf{15} (1982), no.~1, p.~173--204.

\bibitem{ghys_Classification_Anosov_circle_bundles}
{\scshape E.~Ghys} -- {\og Flots d'{A}nosov sur les {$3$}-vari\'{e}t\'{e}s
  fibr\'{e}es en cercles\fg}, \emph{Ergodic Theory Dynam. Systems} \textbf{4}
  (1984), no.~1, p.~67--80.

\bibitem{goodman_surgery}
{\scshape S.~Goodman} -- {\og Dehn surgery on {A}nosov flows\fg}, in
  \emph{Geometric dynamics ({R}io de {J}aneiro, 1981)}, Lecture Notes in Math.,
  vol. 1007, Springer, Berlin, 1983, p.~300--307.

\bibitem{handel-thurston}
{\scshape M.~Handel {\normalfont \smfandname} W.~P. Thurston} -- {\og Anosov
  flows on new three manifolds\fg}, \emph{Invent. Math.} \textbf{59} (1980),
  no.~2, p.~95--103.

\bibitem{hiraide-Expansivos}
{\scshape K.~Hiraide} -- {\og Expansive homeomorphisms of compact surfaces are
  pseudo-{A}nosov\fg}, \emph{Osaka J. Math.} \textbf{27} (1990), no.~1,
  p.~117--162.

\bibitem{hozoori_2024_regularity_WL}
{\scshape S.~Hozoori} -- {\og Regularity and persistence in non-{W}einstein
  {L}iouville geometry via hyperbolic dynamics\fg}, \emph{arXiv 2409.15592}
  (2024).

\bibitem{iakovoglou_Thesis}
{\scshape I.~Iakovoglou} -- {\og Classifying anosov flows in dimension 3 by
  geometric types\fg}, \smfphdthesisname, Université de Bourgogne
  Franche-Comt\'e, 2023, Hal: tel-04402978.

\bibitem{inaba-matsumoto_expansive_flows}
{\scshape T.~Inaba {\normalfont \smfandname} S.~Matsumoto} -- {\og Nonsingular
  expansive flows on {$3$}-manifolds and foliations with circle prong
  singularities\fg}, \emph{Japan. J. Math. (N.S.)} \textbf{16} (1990), no.~2,
  p.~329--340.

\bibitem{katok-hasselblatt}
{\scshape A.~Katok {\normalfont \smfandname} B.~Hasselblatt} --
  \emph{Introduction to the modern theory of dynamical systems}, Encyclopedia
  of Mathematics and its Applications, vol.~54, Cambridge University Press,
  Cambridge, 1995.

\bibitem{lewowicz_exp}
{\scshape J.~Lewowicz} -- {\og Expansive homeomorphisms of surfaces\fg},
  \emph{Bol. Soc. Brasil. Mat. (N.S.)} \textbf{20} (1989), no.~1, p.~113--133.

\bibitem{lewowicz_anal}
{\scshape J.~Lewowicz {\normalfont \smfandname} E.~Lima~de S\'{a}} -- {\og
  Analytic models of pseudo-{A}nosov maps\fg}, \emph{Ergodic Theory Dynam.
  Systems} \textbf{6} (1986), no.~3, p.~385--392.

\bibitem{mosher_monograph}
{\scshape L.~Mosher} -- {\og Laminations and flows transverse to finite depth
  foliations\fg}, \emph{Unpublished manuscript} (1996).

\bibitem{Palis_saddle-connection-stability}
{\scshape J.~Palis} -- {\og A differentiable invariant of topological
  conjugacies and moduli of stability\fg}, in \emph{Dynamical systems, {V}ol.
  {III}---{W}arsaw}, Ast\'{e}risque, vol. No. 51, Soc. Math. France, Paris,
  1978, p.~335--346.

\bibitem{paternain_expansive_flows}
{\scshape M.~Paternain} -- {\og Expansive flows and the fundamental group\fg},
  \emph{Bol. Soc. Brasil. Mat. (N.S.)} \textbf{24} (1993), no.~2, p.~179--199.

\bibitem{paternain_geodesic_flows}
\bysame , {\og Expansive geodesic flows on surfaces\fg}, \emph{Ergodic Theory
  Dynam. Systems} \textbf{13} (1993), no.~1, p.~153--165.

\bibitem{plante_sol_mafld_classification}
{\scshape J.~F. Plante} -- {\og Anosov flows, transversely affine foliations,
  and a conjecture of {V}erjovsky\fg}, \emph{J. London Math. Soc. (2)}
  \textbf{23} (1981), no.~2, p.~359--362.

\bibitem{shannon_thesis}
{\scshape M.~Shannon} -- {\og Dehn surgeries and smooth structures on
  3-dimensional transitive anosov flows\fg}, \smfphdthesisname, Université de
  Bourgogne Franche-Comt\'e, 2020, Hal: tel-02951219.

\end{thebibliography}
